\documentclass[a4paper, reqno, 10pt]{amsart}
\usepackage{amsmath, amssymb, amscd, amsthm, amsfonts,tikz-cd,enumerate,leftindex,bm}
\usepackage{graphicx}
\usepackage{hyperref}
\usepackage{geometry}
\usepackage{dsfont}
\usepackage{enumitem}
\usepackage{xcolor}
\usepackage{stmaryrd}
\usepackage{comment}
\usepackage{imakeidx}
\graphicspath{ {images/} }
\usepackage{cleveref}

\usepackage[backend=bibtex,style=alphabetic]{biblatex}
\usetikzlibrary{decorations.markings}

\makeatletter
\tikzcdset{
  open/.code     = {\tikzcdset{hook, circled};},
  closed/.code   = {\tikzcdset{hook, slashed};},
  open'/.code    = {\tikzcdset{hook', circled};},
  closed'/.code  = {\tikzcdset{hook', slashed};},
  circled/.code  = {\tikzcdset{markwith = {\draw (0,0) circle (.375ex);}};},
  slashed/.code  = {\tikzcdset{markwith = {\draw[-] (-.4ex,-.4ex) -- (.4ex,.4ex);}};},
  markwith/.code ={
    \pgfutil@ifundefined%
    {tikz@library@decorations.markings@loaded}%
    {\pgfutil@packageerror{tikz-cd}{You need to say %
      \string\usetikzlibrary{decorations.markings} to use arrows with markings}{}}{}%
    \pgfkeysalso{/tikz/postaction = {
      /tikz/decorate,
      /tikz/decoration={markings, mark = at position 0.5 with {#1}}}
    }
  },
}
\makeatother
\newtheorem{theorem}{Theorem}[section]
\newtheorem*{theorem*}{Theorem}
\newtheorem{lemma}[theorem]{Lemma}

\newtheorem{prop}[theorem]{Proposition}
\newtheorem{corol}[theorem]{Corollary}
\newtheorem*{prop*}{Proposition}
\theoremstyle{definition}
\newtheorem{defn}[theorem]{Definition}
\newtheorem{ex}[theorem]{Example}
\newtheorem{rem}[theorem]{Remark}

\newcommand{\zz}{\mathbb{Z}}
\newcommand{\fp}{\mathbb{F}_p}

\newcommand{\h}{\mathrm{H}}
\newcommand{\dist}{\mathrm{Dist}}
\newcommand{\ve}{\mathrm{V}}
\newcommand{\w}{\mathrm{W}}

\newcommand{\wn}{\w_n}
\newcommand{\ox}{\mathcal{O}_{X}}
\newcommand{\oo}{\mathcal{O}}
\newcommand{\oxx}{\mathcal{O}_{\mathcal{X}}}
\newcommand{\spec}{\mathrm{Spec}}

\newcommand{\enk}{\mathrm{End}_k}

\newcommand{\fil}{\mathrm{Fil}}
\newcommand{\iend}{\mathcal{E}\mathit{nd}}
\newcommand{\mm}{\mathcal{M}}
\newcommand{\nn}{\mathcal{N}}
\newcommand{\ext}{\mathrm{Ext}}
\newcommand{\ind}{\mathrm{Ind}}
\newcommand{\F}{\Phi_A}
\newcommand{\G}{\mathbf{G}}

\numberwithin{equation}{section}

\title{
 {On Hodge--Witt cohomology 
 of Drinfeld's upper half space 
 over a finite field}
}

\author{
{Mattia Tiso}}
\makeindex
\begin{document}
\maketitle

\begin{abstract}
             In this dissertation we study the Hodge-Witt cohomology of the $d$-dimensional Drinfeld's upper half space $\mathcal{X} \subset \mathbb{P}_k^d$ over a finite field $k$. We consider the natural action of the $k$-rational points $G$ of the linear group $\mathrm{GL}_{d+1}$ on $\h^0(\mathcal{X},\wn\Omega_{\mathbb{P}_k^d}^i)$, making them natural $\wn(k)[G]$-modules. To study these representations, we introduce a theory of differential operators over the Witt vectors for smooth $k$-schemes $X$, through a quasi-coherent sheaf of $\wn(k)$-algebras  $\mathcal{D}_{\wn(X)}$. We apply this theory to equip suitable local cohomology groups arising from $\h^0(\mathcal{X},\wn\oo_{\mathbb{P}_k^d})$ with a $\Gamma(\mathbb{P}_k^d,\mathcal{D}_{\wn(\mathbb{P}_k^d)})$-module structure. Those local cohomology groups are naturally modules over some parabolic subgroup of $\mathrm{GL}_{d+1}(k)$, and we prove that they are finitely generated $\Gamma(\mathbb{P}_k^d,\mathcal{D}_{\wn(\mathbb{P}_k^d)})$-modules.
\end{abstract}

\tableofcontents


\section*{Introduction}

Let $k$ be a finite field of characteristic $p > 0$. The $d$-dimensional Drinfeld's half space $\mathcal{X}^{}$ over $k$ is the open affine $k$-subvariety of $\mathbb{P}_k^d$ obtained by taking the complement of all $k$-rational hyperplanes $\mathbb{P}(H)$ in $\mathbb{P}_k^d$, i.e.,
\begin{equation}\label{Drinf_half_space}
    \mathcal{X}^{} := \mathbb{P}_k^d \backslash \bigcup_{H \subsetneq k^{d+1}} \mathbb{P}(H).
\end{equation}
We omit the dimension $d$ from the notation, assuming it  is implicitly fixed.
The finite group $G$ of $k$-rational points of the general linear algebraic group $\mathrm{GL}_{d+1,k}$ acts naturally on $\mathcal{X}^{}$ by permuting the $k$-rational hyperplanes in the complement $\mathcal{Y}^{}:=\bigcup_{H \subsetneq k^{d+1}} \mathbb{P}(H)$. 
Set $\G_k=\mathrm{GL}_{d+1,k}$ and let $\mathcal{E}$ be a $\textbf{G}_k$-equivariant vector bundle on $\mathbb{P}_k^d$. By functoriality, the cohomology group $\h^0(\mathcal{X}, \mathcal{E})$ inherits a canonical structure of a $k[G]$-module. In the case where $\mathcal{E}$ is the module of $i$-th differential forms $\Omega^i$ on $\mathbb{P}_k^d$ for some $i \geq 0$, we may further consider, for any natural number  $n\geq 1$, the global sections Hodge-Witt cohomology $\h^0(\mathcal{X}, \wn\Omega^i)$, which is the main object of investigation in this thesis. It has, analogously, a natural structure of a $\wn(k)[G]$ module, where $\wn(k)$ is the ring of Witt vectors of level $n$ of $k$. This study generalizes  the work of a preprint by Orlik \cite{Orlik21} for $n>1$ in the corresponding cases. 
  The main result we achieve is given by the following proposition:
 
 \begin{theorem*}[\Cref{mainfinitegeneration}]
Assume that $\mathrm{char}(k)\neq 2$. 
 The $P_j$-module $\tilde{\h}^{d-j}_{\mathbb{P}_k^j}(\mathbb{P}_k^d,\wn(\mathcal{O}_{\mathbb{P}_k^d}))$ admits a submodule $N_{n,j}$ that is a finitely generated $P_j$-module over $\wn(k)$ 
and a $\wn(k)$-linear epimorphism of $\mathcal{D}_n$-modules
\begin{equation*}
    \rho_{n,j}: \mathcal{D}_n \otimes_{\wn(k)}N_{n,j} \twoheadrightarrow \tilde{\h}^{d-j}_{\mathbb{P}_k^j}(\mathbb{P}_k^d,\wn(\mathcal{O}_{\mathbb{P}_k^d})).
\end{equation*}
\end{theorem*}
 All objects, like the $\wn(k)$-algebra $\mathcal{D}_n$ and the group $P_j$ will be introduced later on. In the case $n=1$ the reader may compare the Proposition above with similar results which have been discussed by Kuschkowitz (\cite[Proposition 2.5.1.3]{Ku16}) and in a preprint by Orlik (\cite[Proposition 3.11]{Orlik21}). In the formulation of the Theorem above, a different cohomology group appears instead of the global section cohomology. This is a consequence of successive reductions of the initial problem.
  
 In general, computing the representations $\h^0(\mathcal{X}, \mathcal{E})$  reduces to study certain (subgroups of) local cohomology groups, namely $\tilde{\h}^{d-j}_{\mathbb{P}_k^j}(\mathbb{P}_k^d, \mathcal{E})$, via a spectral sequence argument appearing already in \cite{Orl01}. 
These subgroups have a structure of modules of the maximal parabolic group $P_j \subset G$, stabilizing $\mathbb{P}_k^j$, and they are not finitely generated $k$-vector spaces. Also, their algebraic nature is not completely clear. 
To gain more information, one may consider the natural structure of $\dist(\G_k)$-module on $\tilde{\h}^{d-j}_{\mathbb{P}_k^j}(\mathbb{P}_k^d, \mathcal{E})$. We ask whether it fulfills some finiteness condition. 

Over a $p$-adic field $K$ a similar problem arises by considering the rigid analytic Drinfeld upper half space and, respectively, the $K$-rational points of the general linear group $\G_K$ acting on it. The case of the canonical bundle has been originally studied by Schneider and Teitelbaum (\cite{ST02}) and generalized by Orlik (\cite{Orl08}) for equivariant vector bundles. In Orlik's approach there are, analogously, some local cohomology groups equipped with a structure of modules over the enveloping algebra $\mathcal{U}(\mathrm{Lie}(\G_K))$. In this setting, each of those local cohomology groups is a $K[P(K)]$-module for the corresponding maximal parabolic subgroup $P \subset \G_K(K)$ and it is generated over $\mathcal{U}(\mathrm{Lie}(\G_K))$ by a finitely generated $K[P(K)]$-submodule. 
Unfortunately, a similar property is not satisfied in characteristic $p$, as observed in \cite{Ku16}. 
The strategy we adopted to overcome this problem is to replace the distribution algebra with the ring of differential operators $\mathcal{D}_1=\mathcal{D}(\mathbb{P}_k^d)$ as suggested by Orlik in \cite{Orlik21}.

The Hodge-Witt cohomology of $\mathcal{X}$ can be seen as a generalization of the cohomology of the $k$-modules $\h^0(\mathcal{X}, \Omega^i_{\mathbb{P}_k^d})$ by means of the $\wn(k)$-modules $\h^0(\mathcal{X}, \wn\Omega^i_{\mathbb{P}_k^d})$, where for any $n \geq 1$, $\wn\Omega^i_{\mathbb{P}_k^d}$ denotes the (Hodge-)Witt module of differentials, appearing in the $i$-th degree of the De Rham-Witt complex of $\mathbb{P}_k^d$. By functoriality, there is an action of $G$ on the Hodge-Witt cohomology groups, inducing a $\wn(k)[G]$-module structure. We prove that, similarly to the cases described above, the problem of computing the cohomology reduces to computing $\wn(k)$-submodules of certain local cohomology modules. This is a geometric phenomenon that only depends  on $\mathcal{X}$, giving rise to naturally non-finitely generated $\wn(k)[P_j]$-modules $\tilde{\h}^{d-j}_{\mathbb{P}_k^{j}}(\mathbb{P}_k^d,\wn\Omega^i_{\mathbb{P}_k^d})$. More precisely, we prove it in the following proposition:
\begin{prop*}[\Cref{prop_red_problem}]
Let $\mathcal{F}=\wn\Omega_{\mathbb{P}^d}^i$ for some $i \geq 0$. There is a spectral sequence $E_1^{-r,s}$ converging to $\h^{s-r}_{\mathcal{Y}}(\mathbb{P}^d,\mathcal{F})$ and degenerating at the $E_2$-page, such that:
 \begin{equation}
 E_2^{0,j}=\h^j(\mathbb{P}^d,\mathcal{F}) \quad j \geq 2,
 \end{equation}
 and the terms $E_2^{-j+1,j}$ for $j \geq 1$ appear as an extension of certain $\wn(k)[G]$-modules:
 \begin{equation}
 0 \rightarrow E_{2,\sim}^{-j+1,j} \rightarrow E_{2}^{-j+1,j} \rightarrow E_{2,\mathrm{w.s.}}^{-j+1,j} \rightarrow 0,
 \end{equation}
 where, the following equality hold:
 \begin{equation}
 E_{2,\sim}^{-j+1,j}= \mathrm{Ind}^G_{P_{(d+1-j,j)}}(\tilde{\h}_{\mathbb{P}^{d-j}}^j(\mathbb{P}^d,\mathcal{F}) \otimes_{\wn(k)} {\leftindex_n {\mathrm{St}}}_j^{\vee})
 \end{equation}
 \begin{equation}
 E_{2,\mathrm{w.s.}}^{-j+1,j}=\h_{}^j(\mathbb{P}^d,\mathcal{F})\otimes_{\wn(k)} ({\leftindex_n {\upsilon}}^{G}_{P_{(d+1-j,1^{j})}})^{\vee},
 \end{equation}
 for any $1\leq j \leq d$, and finally
 \begin{equation}
 E_2^{0,1}=E_1^{0,1}=\mathrm{Ind}^G_{P_{(d,1)}}\h^1_{\mathbb{P}^{d-1}}(\mathbb{P}^d,\mathcal{F}).
 \end{equation}
 \end{prop*}
Here, ${\leftindex_n {\upsilon}}^{G}_{P_{(d+1-j,1^{j})}}$ is the generalized Steinberg representation of $G$ associated to $P_{(d+1-j,1\dots,1)}$ ($1$ appears $j$ times) over $\wn(k)$ and ${\leftindex_n {\mathrm{St}}}_j$ is the standard Steinberg representation of $G$ over $\wn(k)$ (see \Cref{sec_generalized_steinb}).

On the algebraic side, unlike the case $n = 1$, we do not have (a priori) a natural action of the distribution algebra or the ring of differential operators at our disposal. Therefore, we develop, just for the sake of application, a suitable theory of differential operators over the Witt vectors, introducing a $\wn(k)$-algebra $\mathcal{D}_{\wn(X)}(X)$ for smooth $k$-varieties $X$.

 Here, we must mention that a more general theory was going to be introduced in a recent work of Dodd \cite{dodd2024} that appeared while this thesis was being written. In particular, the author of this thesis independently addressed the problem and provided analogous definitions. However, the idea of proving the relations in \Cref{relliftderivation} was inspired by the analogous one in \cite{dodd2024}. Then, we construct a sheaf $\mathcal{D}_{\wn(X)}$ similar to the one in \cite[Definition 2.33]{dodd2024}.
 
  Even if the techniques are different, both agree on the main idea of defining Witt differential operators as a restriction of differential operators with additional properties (i.e. for which \eqref{HS-formula} holds, classically called Hasse-Schimdt derivations) on smooth lifts.  Although similar, our construction is given locally, considering local parameters of smooth algebras, while the one of Dodd is more intrinsic. Moreover, we make consistently use of a map $\tilde{w}_n$, as a replacement of a map $\tilde{F}^n$, introduced and studied originally in the work of Illusie and Raynaud \cite{Illusie1983} and successively extended by  Berthelot et al. in \cite{BHR12}. The strategy in \cite{dodd2024} is analogous, thus there are similarities in computations, but the author does not mention any relation with $\tilde{F}^n$  and its properties, as we do, for example, in the proof of existence of Witt differential operators (cf. \Cref{liftder}). Also, it is worth mentioning that, in contrast to Dodd, we explicitly do not construct any canonical Witt differential operator (cf. \cite[Definition 2.8]{dodd2024}).

Roughly speaking, the main feature of this theory is that any Hasse-Schmidt $k$-linear differential operator over a smooth $k$-algebra $A$ admits a compatible lift to some $\wn(k)$-linear differential operator of $\wn(A)$. This is more precisely described in \Cref{liftder}. Also, those lifts satisfy some compatibilities with Verschiebung, Frobenius and Restriction maps as proved in  \Cref{relliftderivation}. Moreover, thanks to the property \eqref{derversch}, we apply this theory to describe the $\wn(k)[P_j]$-modules above for the case of the Witt vectors cohomology (i.e., for $i = 0$, where $\wn\Omega^i_{\mathbb{P}^d} = \wn\oo_{\mathbb{P}^d}$), and prove that the group $\tilde{\h}^{d-j}_{\mathbb{P}^{j}}(\mathbb{P}^d,\wn\oo_{\mathbb{P}^d})$ has a structure both of a $\wn(k)[P_j]$-module and of a $\mathcal{D}_{\wn(\mathbb{P}_k^d)}(\mathbb{P}_k^d)$-module generated by a finite $\wn(k)[P_j]$-submodule, that is precisely the meaning of the first proposition above.
\\ [2ex]
We will explain the structure of this paper in more detail. In \Cref{sec_diff_op}, we recall definitions and basics properties of Grothendieck's differential operators. We focus on the properties of $\mathcal{D}(A)$ when $A$ is a $\fp$-algebra and give some examples. Moreover, we introduce the notion of \textit{crystalline Weyl algebra}, to be thought as an integral version of the Weyl algebra, giving an explicit description of the module of differential operators in characteristic $p$. Furthermore, we deduce a relation between differential operators of a smooth $\w(k)$-scheme $X$ and its smooth nilpotent thickenings $X_n$ over $\wn(k)$ (cf. \Cref{prop_diff_thickening}).
\\[2ex]
In \Cref{sec_de Rham_Witt_cx}, we recall the construction of the de Rham-Witt complex for a $k$-scheme $X$ ($k$ perfect). We additionally consider $X$ equipped with an action of a finite group $G$ and we define and discuss the notion of $G$-equivariant $\wn\ox$-modules. In particular, by the universal property of de Rham-Witt complex, we deduce that any Hodge-Witt module $\wn\Omega_X^i$ ($i \geq 0$) is $G$-equivariant. Furthermore, following \cite{BHR12}, we explain the classical computation of the de Rham-Witt complex for the affine space of dimension $d$ and how to compute the Hodge-Witt cohomology of the projective space (equipped with the natural action of $G=\mathrm{GL}_{d+1,k}(k)$) of dimension $d$,  (cf. \cref{sec_hodge_witt_proj}). For completeness, we also introduce the map $\tilde{F}^n$, in the way defined in \cite{Illusie1983}, and we provide a self-contained elementary proof of  the known \Cref{tildeF^n_iso_prop}, classically deduced as a consequence of a more involved theory that here we do not investigate. 
 We further give an introduction of a less known concept of \textit{Witt line bundles} (cf. \cref{sec_witt_line}), some particular case of locally free $\wn\ox$-modules of rank one (following Tanaka in \cite{TH22}). Then, we prove some functoriality properties and extend the construction of the map $\tilde{F}^n$ for Witt line bundles as well. Then, as an example, we furnish a computation of the cohomology of Witt line bundles of the projective space of dimension $d$ as a $G$-module (cf. \cref{sec_coho_witt_line_proj}).
\\[2ex]
In \Cref{secwittdiffop}, we introduce a theory of Witt differential operators. If $A$ is a smooth $k$-algebra of dimension $d$, it admits compatible lifts to some smooth $\wn(k)$-algebras $A_n$. The ring of differential operators $\mathcal{D}(A)$ as $k$-algebra, is locally generated by  operators $\partial_i^{[r]}$ (for $i=1,\dots,d$, $r \geq 0$) satisfying certain relations (called Hasse-Schmidt derivations) which we treat in an appendix (c.f.  \eqref{rel_partial_H-S}). By a lifting property of smooth morphisms (cf. \Cref{cor_smooth_lift_der_over_wn}), there exist compatible $\wn(k)$-differential operators $\partial_{i,n}^{[r]} \in \mathcal{D}(A_n)$ lifting $\partial_i^{[r]}$, such that the analogous property \eqref{rel_partial_H-S} holds. Furthermore, there is a ring monomorphism $\tilde{w}_{n-1} \colon \wn(A) \rightarrow A_{n}$. Then, a Witt differential operator is given by restriction of such $\partial_{i,n}^{[r]}$ to $\wn(A)$ via $\tilde{w}_{n-1}$ (cf. \Cref{liftder}). Moreover, for any $\partial_{i,n}^{[r]} \in \mathcal{D}(A_n)$ lifting $\partial_{i}^{[r]} \in \mathcal{D}(A)$, the restriction does not depend on the chosen lifting, but only on that differential operator in  characteristic $p$ (cf. \Cref{uniquenesslift}). Furthermore, we prove some additional properties of Witt differential operators in order to relate them with Frobenius, Verschiebung and restriction maps on  Witt vectors of level $n$ (cf. \Cref{relliftderivation}). Finally, we define a sheaf of Witt differential operators $\mathcal{D}_{\wn(X)}$ for a smooth $k$-scheme $X$, such that for $n=1$, it agrees with Grothendieck's sheaf of differential operators (cf. \Cref{lem_presh_of_wittdiff_is_sheaf}).
\\[2ex]
In \Cref{chap_drinf_hodge_witt}, we consider the $\wn(k)[G]$-module cohomology $\h^0(\mathcal{X}, \mathcal{F})$ where $\mathcal{F}$ may be a Witt line bundle on $\mathbb{P}_k^d$ or one of the Hodge-Witt module $\wn\Omega^i_{\mathbb{P}_k^d}$. We explain how to reduce the computation of the latter by studying certain submodules of the local cohomology $\h^{d-j}_{\mathbb{P}_k^j}(\mathbb{P}^d, \mathcal{F})$ for any $j=0,\dots, d$ (cf. \Cref{prop_red_problem}). To prove the aforementioned Proposition, we adapt a Orlik's resolution of the constant sheaf $\zz_{\mathcal{Y}}$ (introduced in \cite{Orl01}) for the case of the Witt scheme $\wn(X)$, via the natural universal homeomorphism $X \rightarrow \wn(X)$ (induced in the affine case by the canonical projection $\wn(A) \rightarrow A$). We see that the dual (as a $\wn(k)$-module) generalized Steinberg (free finitely generated ) representations over $\wn(k)$ appear in the $G$-module structure of $\h^0(\mathcal{X},\mathcal{F})$. It follows that the only unknown module structures are given by $\tilde{\h}^{d-j}_{\mathbb{P}_k^j}(\mathbb{P}^d, \mathcal{F})$.
\\[2ex]
The latter will be investigated as an application in the last section. In order to proceed, we need a link between local cohomology groups and $\mathcal{D}$-modules (\Cref{sec_crys_bb}). In characteristic $0$, this goes back to Beilinson and Bernstein \cite{beiber81}. An action of a  reductive algebraic group $\G$ on a smooth scheme $X$, induces a natural $\dist(\G)$-module structure on the global section of $\ox$. Moreover, any element of $\dist(\G)$ acts as a differential operator on $\ox(X)$, inducing a natural map of associative algebras $\phi^{\ox} \colon \dist(\G)\rightarrow \Gamma(X,\mathcal{D}_X)$. More generally, this reasoning works for arbitrary quasi-coherent $\G_k$-equivariant sheaves in any characteristic (cf. \Cref{lie_alg_map_induced_by_lineariz} and \Cref{prop_bb_map_constr} hold in any characteristic, and indeed it coincides with the construction in \cite[11, (11.1.6)]{HTT08} in characteristic $0$, by identifying $\mathcal{U}(\mathfrak{g})=\dist(\G)$). The main difference between characteristic $0$ and $p$ is that $\phi^{\ox}$ is not in general  surjective (cf. \cite{smith86} for a counterexample in the case $\G=\mathrm{SL}_{2,k}$), while this is true in characteristic $0$ (cf. \cite[Theorem 11.2.2 (i)]{HTT08}). Using the theory of Witt differential operators, we can define a suitable (Teichm\"uller) lift in $\mathcal{D}_{\wn(X)}$ of differential operators in $\mathcal{D}_{X}$ (cf. \Cref{prop_teich_lift_diff}). 
\\[2ex]
 As an application, we consider the natural $\mathcal{D}_n:=\Gamma(\mathbb{P}_k^d,\mathcal{D}_{\wn(\mathbb{P}_k^d)})$-module structure (given on the left by evaluation of differential operators on functions) on the  cohomology subgroups  $\tilde{\h}^{d-j}_{\mathbb{P}_k^j}(\mathbb{P}_k^d,\wn\oo_{\mathbb{P}_k^d})$ together with the natural action of the finite parabolic subgroup $P_j:=P_{(j+1,d-j)}$ of $ G$. We ask the following: Does exist a $P_j$-equivariant $\wn(k)$-submodule $N_{n,j}$, such that it generates  $\tilde{\h}^{d-j}_{\mathbb{P}_k^j}(\mathbb{P}_k^d,\wn\oo_{\mathbb{P}_k^d})$ as $\mathcal{D}_n$-module? We answer positively  in \Cref{mainfinitegeneration}. The proof is constructive and only requires some elementary properties of Witt vectors. Then, the statement reduces to the characteristic $p$ case ($n=1$), by properties of Witt differential operators.


\section*{Acknowledgements}
This paper is the author's Ph.D. Thesis with slight corrections. I want to thank my advisor Sascha Orlik, to give me the opportunity to study this challenging and interesting topic. Moreover, I want to thank my co-advisor Kay R\"ulling for his constant help and support in learning and introducing me to several questions regarding the theory of Witt vectors. Also, I am grateful to my colleagues and friends Fei Ren, Andreas Bode, Dennis Peters, Daan Van Sonsbeek for enlightening conversations. 

The author conducted this project as part of the research
training group \textit{GRK 2240: Algebro-Geometric Methods in Algebra, Arithmetic and Topology}
which is funded by the Deutsche Forschungsgemeinschaft.

 \section*{Notation}
     In the course of this exposition, $p$ is a prime number and $k$, when not specified, will denote a finite field in characteristic $p>0$.\\
     The set of natural numbers $\mathbb{N}$ contains $0$. \\
     For any integers $a,b$ with $a < b$, we denote the range of integers between $a$ and $b$ by $[a,b]:=\{i \in \zz \mid a \leq i \leq b \}$.
     For any natural number $d \geq 1$, bold symbols $\mathbf{i}=(i_1,\dots,i_d),\mathbf{j}=(j_1,\dots,j_d),$ etc., are vectors of the abelian group $\mathbb{Z}^d$ or they might denote just $d$-uples. It will be clear from the context. We write $\mathbf{i} \leq \mathbf{j}$, if and only if $i_l \leq j_l$ for any $1 \leq l \leq d$. 
     
     The bold symbols $\G_k,\textbf{B}_k,\textbf{T}_k, \textbf{P}_k$ will denote algebraic groups over $k$. Sometimes the index $k$ is omitted, when it is clear. When the index $\zz$ appears, it means that the algebraic groups are defined over $\zz$. All rings are assumed to be commutative with unit, unless we are talking of the differential operators algebra, the Weyl algebra, the crystalline Weyl algebra, or the enveloping algebra of some Lie algebra, which are generally non commutative.
     
      The symbols, $G,B,T,P$ denotes the $k$-rational points of the respective algebraic groups. Also, when it is not specified, $\G$ is $\mathrm{GL}_{d+1}$ and generally always  a reductive group; $\textbf{T},\textbf{B}$ are fixed maximal torus and Borel subgroups of $\G$ and $\textbf{P}_I$ is the parabolic subgroup  associated to a subset $I \subset \Delta$, where $\Delta$ is the set of simple roots for the root system $\Phi(\textbf{T},\textbf{B})$ of $\G$. In the case of $\G=\mathrm{GL}_{d+1}$,
      $\Phi=\{\alpha_{ij}=\epsilon_i-\epsilon_j \mid 0 \leq i \neq j \leq d\}$ and 
      $\Delta=\{\alpha_i:=\alpha_{i,i+1} \mid i=0,\dots,d-1 \}$ where
      $\epsilon_i \in X(\textbf{T})=\mathrm{Hom}(\mathbf{T}, \mathbb{G}_m)$ is the character sending $\mathbf{T}(A) \ni (t_0,\dots,t_d) \mapsto t_i \in A$, for any $k$-algebra $A$. Also, the set of positive roots $\Phi^+ \subset \Phi$ consists of the elements $\alpha_{ij}$ with $i<j$, while $\Phi^{-}$ is its complement.

       The gothic symbols $\mathfrak{g}$, $\mathcal{U}(\mathfrak{g})$ denotes  respectively the Lie algebra and the enveloping algebra of $\G$, while $\dist(\G)$ is its distribution algebra.
       
        The letters $X$, $Y$ denote $k$-schemes of finite type. The index $(-)_A$ for any $k$-algebra $A$ denotes the base change along the structure morphism $k \rightarrow A$.
        
        For any $k$-schemes $X$, $Y$ and $Z$ we adopt the following conventions:
       we denote by $\text{pr}_1 \colon Y \times_k X \rightarrow X$ the canonical projection given by $(y,x) \mapsto x$, for $y \in Y, x \in X$, $\text{pr}_{12}: Z \times_k Y \times_k X \rightarrow Y \times_k X$ the projection $(z,y,x) \mapsto (y,x)$ for $z \in Z, y  \in Y$, $x \in X$, and $\text{pr}_2=\text{pr}_1\circ \text{pr}_{12}$.

         The $p$-typical ring of Witt vectors of length $n \geq 1$ will be denoted by $\wn(A)$ and for any $i$, $\wn\Omega_{X}^i$ is the $\wn\ox$-module appearing at the $i$-th degree of the De Rham-Witt complex of $X$. 
         
          We also consider the action of $\mathrm{GL}_{d+1,k}$ on $\mathbb{P}_k^d$ given on points by
    \begin{equation}
 \begin{array}{cc}
       \mathrm{GL}_{d+1,k}(A) \times \mathbb{P}_k^d(A) \longrightarrow \mathbb{P}_k^d(A)\\
      \qquad \qquad (g,[x_0:\dots:x_d]) \longmapsto [x_0:\dots :x_d]g^{-1}
     \end{array}
 \end{equation}
 where $g \in GL_{d+1,k}(A)$ and $[x_0 : \dots x_d] \in  \mathbb{P}_k^d(A)$ for any $k$-algebra $A$.

 Let $A$ be a $\fp$-algebra and $X$ a scheme over $\spec(A)$ (simply said over $A$). Let $\mathrm{Frob}^{\#} \colon A \rightarrow A, \quad x \mapsto x^p$ be the Frobenius morphism of $A$ and $\mathrm{Frob}$ be the induced map of spectra. Let $X^{(p)}:=X\times_{\mathrm{Frob}, \spec(A)} \spec(A)$ and consider the pullback square 
 \begin{equation}
     \begin{tikzcd}
X \arrow[rd, "F_{X/A}", dashed] \arrow[rrd, "F_X", bend left] &                                        &             \\
                                                                  & X^{(p)} \arrow[d] \arrow[r, "W"] & X \arrow[d] \\
                                                                  & \spec(A) \arrow[r, "\mathrm{Frob}"]    & \spec(A)   
\end{tikzcd}
 \end{equation}
where $W$ is the map of schemes induced by the pullback construction.
  The \textit{relative Frobenius} $F_{X/A} \colon X \rightarrow X^{(p)}$ on $X$ respect to $A$ is the map of $A$-schemes given by the universal property of pullback diagrams.\\
  Lastly, the \textit{absolute Frobenius} is given by $F_X=W \circ F_{X/A}$. 
  
  Let $\bar{k}$ be an algebraic closure of $k$ and $X$ be a $k$-scheme with base change $X\times_k \bar{k}=:{X}_{\bar{k}}$. The \textit{geometric Frobenius} is the $\bar{k}$-scheme morphism $F_X \times id_{\bar{k}} \colon {X}_{\bar{k}} \rightarrow {X}_{\bar{k}}$.
  If $X=\mathrm{GL}_{d,k} \subset \mathbb{A}_k^{d^2}$ , then the \textit{standard geometric Frobenius} of ${X}_{\bar{k}}$ is the restriction of the geometric Frobenius of $\mathbb{A}_{\bar{k}}^{d^2}$ to the open subscheme ${X}_{\bar{k}}$.

  For a topological space $X$ and an abelian sheaf $\mathcal{F}$ on $X$, as usual we denote by $\h^i(X,\mathcal{F})$ be the $i$-th derived functor of the global sections $\Gamma(X,-)$ with value in $\mathcal{F}$. If $U \subset X$ is an open subset with closed complement $Z=X \backslash U$, then $\h^i_Z(X,\mathcal{F})$ denotes the $i$-th derived functor of the global section $\Gamma_{Z}(X,-)$ with support in $Z$ and value in $\mathcal{F}$. For such a triple $(Z,U;X)$, we refer to the induced long exact sequence on group cohomology
  of
  \begin{equation}
      0 \rightarrow \Gamma_Z(X,\mathcal{F}) \rightarrow \Gamma(X,\mathcal{F}) \rightarrow \Gamma(U,\mathcal{F}) \rightarrow 0
  \end{equation}
    as the associated long exact sequence of the couple $(Z,U;X)$:
    \begin{equation}
       \dotsc \rightarrow \h^i_Z(X,\mathcal{F}) \rightarrow \h^i(X,\mathcal{F}) \rightarrow \h^i(U,\mathcal{F}) \rightarrow \h^{i+1}_Z(X,\mathcal{F}) \rightarrow \dotsc
  \end{equation}
 We denote the subgroups $\ker(\h^i_Z(X,\mathcal{F}) \rightarrow \h^i(X,\mathcal{F}))$ by $\tilde{\h}_Z^{i}(X,\mathcal{F})$.

 The for any abelian sheaf $\mathcal{F}$ on $X$, and for any integer $i$, the sheaf associated to the presheaf given by assigning for any open $U \subset X$
 \begin{equation*}
      U \mapsto \h^{i}(U,\mathcal{F})
 \end{equation*}
 will be denoted by $\mathcal{H}^i(\mathcal{F})$. Analogously, for a closed subset $Z \subset X$, the presheaf given by 
 \begin{equation*}
     U \mapsto \h_Z^{i}(U,\mathcal{F})
 \end{equation*}
 will be denoted by $\mathcal{H}_Z^i(\mathcal{F})$.

 For a complex of abelian sheaves $(\mathcal{F}^{\bullet},d)$, we denote by $$h^i(\mathcal{F}^{\bullet}):=\frac{\ker(d \colon \mathcal{F}^i \rightarrow \mathcal{F}^{i+1})}{d\mathcal{F}^{i-1}}$$ the $i$-th sheaf cohomology for any $i \in \zz$.

\section{Grothendieck's differential operators}\label{sec_diff_op}
Let $k$ be a commutative ring. In this section, we recall the definition of differential operators in the sense of \cite[Ch. IV, Sec. 16.8]{EGA4} and discuss some properties found in \cite{smith86} for the characteristic $p$ case.
\subsection{Basic definitions and properties} 
Let $A$ be a commutative, unitary  $k$-algebra. Then, define the $A$-algebra $\mathcal{D}(A)$ given by  $\mathcal{D}(A)=\bigcup_{n=0}^{\infty} \mathcal{D}_n(A)$, where 
\begin{equation*}
    \mathcal{D}_m(A):=\{\theta \in \enk(A) \mid [a_0,[a_1,\dots,[a_m,\theta]\dots]]=0 \quad\forall a_0,\dots, a_m \in A\}.
\end{equation*}
Here $(\enk(A),+,\circ)$ is the algebra of $k$-linear endomorphisms of $A$  and $A \subset \enk(A)$ is identified with the left (or right) multiplication morphism; the bracket $[-,-]: \enk(A) \times \enk(A) \rightarrow \enk(A)$ is the map sending $(\theta,\eta)$ to $\theta \circ \eta-\eta \circ \theta$. We recall that the filtration 
$\mathcal{D}_m(A) \subset  \mathcal{D}_{m+1}(A)$ makes $\mathcal{D}(A)$ a filtered $k$-algebra. \\
For any affine scheme (of finite type) $X=\spec(A)$, set $\mathcal{D}(X):=\mathcal{D}(A)$. Then, the  notion of differential operators sheafifies (for example by \cite[Theorem 3.2.5]{Will98}).   
\begin{defn}
For  a $k$-scheme $X$, $\mathcal{D}_X$ is the  unique quasi-coherent  $\ox$-module given by $\Gamma(U,\mathcal{D}_X)=\mathcal{D}(U)$ for any Zariski open affine $U \subset X$.
\end{defn}
 Moreover, $\mathcal{D}_X$ is also equipped with a filtration, by setting
 \begin{equation*}
     \fil_m\mathcal{D}_X(U):=\{\theta \in \mathcal{E}\mathit{nd}(\ox)(U) \mid [\theta,a] \in \fil_{m-1}\mathcal{D}_X(U), \forall a \in \ox(U)\}.
 \end{equation*}
 It is straightforward to see that $\fil_m\mathcal{D}_X(X)=\mathcal{D}_m(X)$ for $X$ affine. Moreover, we have a decomposition $\fil_1\mathcal{D}_X=\ox\oplus \mathcal{T}_X$ where $\mathcal{T}_X$ is the \textit{tangent sheaf},  given explicitly by
 \begin{equation*}
    \mathcal{T}_X(U)=\{\theta \in \mathcal{E}\mathit{nd}(\ox)(U) \mid \theta(ab)=a\theta(b)+b\theta(a) \quad \forall a,b \in \ox(U)\}
\end{equation*}
for  any open affine $U$ of $X$.

 Assume now, that  $X$ is equipped with an action of a linear algebraic group $\G$ over $k$. Then, the structure sheaf $\ox$ as well as $\Omega_{X}^i$  for any $i \geq 1$ are $\G$-equivariant in the sense of \cite[1, 3, Def. 1.3]{Mumf94} or  \cite[I, 0.2]{BerLun94}. If $G$ is an abstract group, with multiplication map $m \colon G \times G \rightarrow G$, the notion of linearization for quasi-coherent modules on $X$ can be analogously formulated.
 Consider $(G,m)$ as the constant $k$-group $\coprod_{g \in G} \spec{(k)}$. Assume that $G$ acts on $X$ via an action $\sigma \colon G \times_k X \rightarrow X$. Denote by $\text{pr}_1 \colon G \times_k X \rightarrow X$ the canonical projection, $\text{pr}_{12}: G \times_k G \times_k X \rightarrow G \times_k X$ the projection $(g_1,g_2,x) \mapsto (g_2,x)$ for any $g_1,g_2 \in G$, $x \in X$, and $\text{pr}_2=\text{pr}_1\circ \text{pr}_{12}$.  
 \begin{defn}
 A  quasi-coherent $\oo_{X}$-module $\mathcal{F}$ is said to be $G$-\textit{equivariant} (or $G$-\textit{linearizable}) if there exists an isomorphism (called $G$-\textit{linearization}) of $\oo_{G \times_k X}$-modules
 \begin{equation}
     \phi: \sigma^*\mathcal{F} \xrightarrow{\sim} \text{pr}_1^* \mathcal{F}
 \end{equation}
 such that the following 
     \begin{equation}\label{diag_linearization}
\begin{tikzcd}
(1_G \times \sigma)^*\text{pr}_1^*\mathcal{F} \arrow[r, "\text{pr}_{12}^*\phi"]        &  \text{pr}_2^*\mathcal{F}                                                  \\
(1_G \times \sigma)^*\sigma^*\mathcal{F} \arrow[u, "(1_ G \times \sigma)^*\phi"] \arrow[r, "="] & (m \times 1_ X)^*\sigma^*\mathcal{F} \arrow[u, "(m \times \text{id}_ X)^*\phi"']
\end{tikzcd}
\end{equation}
is a commutative diagram of $\oo_{G \times_k G \times_k X }$-modules.
 \end{defn}
 A $G$-linearization induces a canonical $G$-action on the global sections of $\mathcal{F}$, thus on each cohomology group with coefficients in $\mathcal{F}$. Indeed,
 since $G$ is the constant $k$-group associated to an abstract group, the definition above is equivalent to say that the collection of $k$-scheme isomorphisms $(\phi_g \colon \sigma_g^*\mathcal{F} \rightarrow \mathcal{F})_{g \in G}$, where  $\sigma_g \colon \{g\} \times X \xrightarrow{i_g} G \times_k X \xrightarrow{\sigma} X$ and 
 $\phi_g =i_g^*\phi$,
 satisfies the property for which the following diagram
 \begin{equation}\label{dia_linearization_points}
            \begin{tikzcd}
\sigma_{g_1g_2}^*\mathcal{F} \arrow[rd, "\sigma_{g_2}^*\phi_{g_1}"] \arrow[rr, "\phi_{g_1g_2}"] &                                                    & \mathcal{F} \\
  & \sigma_{g_2}^*\mathcal{F} \arrow[ru, "\phi_{g_2}"] &            
\end{tikzcd}
      \end{equation}
commutes for any $g_1,g_2 \in G$.
Then, for any $g \in G$, the adjoint morphism $\mathcal{F} \rightarrow \sigma_{g*}\sigma_g^{*}\mathcal{F}$ induces a morphism 
\begin{equation}
    \h^i(X,\mathcal{F}) \rightarrow \h^i(X,\sigma_{g*}\sigma_g^{*}\mathcal{F})=\h^i(G \times X,\sigma_g^{*}\mathcal{F} )\xrightarrow[\sim]{\phi_g} \h^i(X,\mathcal{F}).
\end{equation}
In this way $\h^i(X,\mathcal{F})$ has a structure of a $G$-module.
 \begin{lemma}\label{lemma_D_X_is_G_equi}
 If $G$ acts on a $k$-variety $X$, then $\fil_m\mathcal{D}_X$ and $\mathcal{D}_X$ are  $G$-equivariant quasi-coherent $\ox$-modules for any $m\geq 1$.
 \end{lemma}
 \begin{proof}
 It is sufficient to prove the statement for $\fil_m\mathcal{D}_X$, $m \geq 1$. As abuse of notation, for any $g \in G$, the isomorphism
 $\sigma_g(g,-) \colon X \rightarrow X$ will be simply denoted by $\sigma_g$. We will prove that for any $g \in G$, the morphism $\sigma_g $ induces a $\ox$-module isomorphism
 \begin{equation}\label{isofasci}
 \phi_g : \fil_m\mathcal{D}_X \xrightarrow{\sim}\sigma_{g*}\fil_m\mathcal{D}_X.
  \end{equation}
 It suffices to prove it for any open affine $U \subset X$. Let $\sigma_g^{\#}  \colon \ox \rightarrow (\sigma_{g})_*\ox$ be the canonical map induced by $\sigma_g$. For each $\eta \in \fil_m\mathcal{D}_X(U)$ let
 $\phi_g(\eta)$ be defined by 
 \begin{equation} \label{isofasciformula}
 \phi_g(\eta)(f):=\sigma_{g^{}}^{\#}(\eta(\sigma_{g^{-1}}^{\#}(f))) \in \ox(g^{-1}.U) \quad \text{for any }
  f \in \ox(g^{-1}.U).\end{equation}\footnote{We adopt the convention that $g.$ means to apply $\sigma_g(g,-)\colon X \rightarrow X$}
The $\ox$-module homomorphism $\phi_g$ is an isomorphism, with inverse given by $\phi_{g^{-1}}$. \\
 We prove that $\phi_g(\eta) \in \sigma_{g*}\fil_m\mathcal{D}_X(U) $. We proceed by induction on $m$.
 Recall that
 \begin{equation*}
     \fil_m\mathcal{D}_X(U)=\{\eta \in \iend(\ox)(U) \mid [\eta, F] \in \fil_{m-1}\mathcal{D}_X(U), \quad \forall F \in \ox(U)\}
 \end{equation*}
 and 
 \begin{equation*}
     \sigma_{g*}\fil_{m}\mathcal{D}_X(U)=\{\theta \in \iend(\ox)(g^{-1}.U) \mid [\theta,f] \in \fil_{m-1}\mathcal{D}_X(g^{-1}U) \quad \forall f\in \ox(g^{-1}.U)\}.
  \end{equation*}
  For $m=1$, we have that $\fil_1\mathcal{D}_X(U)=\ox(U) \oplus \mathcal{T}_X(U).$ We just need to prove that $\phi_g(\eta)$ is a derivation if $\eta$ is a derivation.
  Indeed, for $f,h \in \ox(g^{-1}.U)$, we get
  \begin{align*}
     \phi_g(\eta)(fh)= \sigma_g^{\#}(\eta(\sigma_{g^{-1}}^{\#}(fh))=\sigma_g^{\#}(\eta(\sigma_{g^{-1}}^{\#}(f)\sigma_{g^{-1}}^{\#}(h)))\\
     =\sigma_g^{\#}(\eta(\sigma_{g^{-1}}^{\#}(f))\sigma_{g^{-1}}^{\#}(h)) + \sigma_g^{\#}(\eta(\sigma_{g^{-1}}^{\#}(h))\sigma_{g^{-1}}^{\#}(f))\\
     =\phi_g(\eta)(f)\cdot h+f\cdot \phi_g(\eta)(h).		
 \end{align*}
        For a generic $\eta \in \fil_m\mathcal{D}_X(U)$, we have
  \begin{align*}
     [\phi_g(\eta),f](h)= \phi_g(\eta)(fh)-f\phi_g(\eta)(h) \\=	
     \phi_g(\eta)(f)\cdot h+f\cdot \phi_g(\eta)(h)-f\phi_g(\eta)(h) \\
     =\phi_g(\eta)(f)\cdot h\\
     =\sigma_{g}^{\#}([\eta,\sigma_{g^{-1}}^{\#}(f)](h)).
  \end{align*}
  In particular, this shows that 
  \begin{equation*}
  [\phi_g(\eta),f]=\phi_g([\eta,\sigma_{g^{-1}}^{\#}(f)]) \in \fil_{m-1}\mathcal{D}_X(g^{-1}.U)
  \end{equation*}
  by inductive hypothesis.\\
  We are left to verify the cocycle condition \eqref{dia_linearization_points}. It follows by the straightforward equality $\phi_{g_1g_2}(\eta)=\phi_{g_2}(\phi_{g_1}(\eta))$ for any $g_1,g_2 \in G$ and $\eta \in  \fil_m\mathcal{D}_X$.
 \end{proof}
 Now suppose that the $k$-scheme $X$ is equipped with an action ${}^k{\sigma} \colon \G \times_k X \rightarrow X$ of a reductive group $\G$ over $k$. Then, the analogous of \Cref{lemma_D_X_is_G_equi} holds true.
 Let $G=\G(k)$. Taking the base change with an algebraic closure $\bar{k}$ induces an action ${}^{\bar{k}}{\sigma} \colon \G_{\bar{k}} \times_{\bar{k}} X_{\bar{k}} \rightarrow X_{\bar{k}}$. We can identify $\bar{k}$-rational points with closed points of the scheme, $|\G_{\bar{k}}|=\G_{\bar{k}}(\bar{k})$. Moreover, suppose there is a $G$-stable open subscheme $U_{\bar{k}} \subset X_{\bar{k}}$. Thus, the restriction ${}^{\bar{k}}{\sigma}$ on $\G(k)$ induces an action $\bar{\sigma} \colon G \times_{\bar{k}} U_{\bar{k}} \rightarrow U_{\bar{k}}$. If there exists a $k$-rational structure $U$ for $U_{\bar{k}}$, then $\bar{\sigma}$ induces an action $\sigma \colon G \times_{k} U \rightarrow U$, via the natural map $U_{\bar{k}}=U \times_{k} \bar{k} \rightarrow U$.
 A $\G_{\bar{k}}$-linearization of $\mathcal{D}_{X_{\bar{k}}}$ can be defined similarly to the Formula \eqref{isofasciformula} by replacing the respective actions with ${}^{\bar{k}}\sigma_g$. Then it induces, by restriction, the $G$-linearization for $\mathcal{D}_U=\mathcal{D}_{X \mid U}$ defined in the lemma above.
 \begin{ex}
     When $X= \mathbb{P}_k^d$, $U=\mathcal{X}$ and $\G=\mathrm{GL}_{d+1,k}$ the $\G$-linearization on $\mathcal{D}_{X}$ induces a $G$-linearization of $\mathcal{D}_{\mathcal{X}}$ such that $G$ acts on $\mathcal{D}_{X}(\mathcal{X})$ via the isomorphism given by the Formula \eqref{isofasciformula}.
 \end{ex}
 \subsection{The Weyl algebra and crystalline Weyl algebra}
 Let $z_1,\dots,z_m$ be a set of variables together with symbols $\partial_{z_1},\dots,\partial_{z_m}$. 
 \begin{defn}
 A $m$-th \textit{Weyl algebra over} $k$ is a $k$-algebra isomorphic to $$k\{z_1,\dots,z_m,\partial_{z_1},\dots,\partial_{z_m}\}/J=: k[z_i \mid 1 \leq i \leq m]\langle \partial_{z_i} \mid 1 \leq i \leq m \rangle:=A_m(k)$$ for some $m$, where 
 $k\{z_1, \dots, z_m, \partial_{z_1},\dots,\partial_{z_m}\}$ is the free algebra generated by the symbols $z_i, \partial_{z_i}$ and $J$ is the ideal $$J=(z_iz_j-z_jz_i, z_i\partial_{z_j}-\partial
_{z_j}z_i +\delta_{ij}, \partial_{z_i}\partial_{z_j}-\partial_{z_j}\partial_{z_i},\text{ } 1 \leq i,j \leq m)$$
 \end{defn}
 Suppose that $k=\zz$ . Then consider the Weyl algebra over $\mathrm{Frac}(k)=\mathbb{Q}$,  $$A_m(\mathbb{Q})=\mathbb{Q}[z_i \mid 1 \leq i \leq m]\langle \partial_{z_i} \mid 1 \leq i \leq m \rangle.$$ Define the \textit{crystalline symbols} being the elements of $A_m(\mathbb{Q})$ given recursively by the integral relations $$\partial_i^{[r]}\partial_i^{[s]}=\binom{r+s}{r}\partial_i^{[r+s]}, \quad \forall r,s \in \mathbb{N}, $$ where 
 for $r=0$, $\partial_i^{[0]}:=1$, and for $r=1$, $\partial_i^{[1]}:=\partial_{z_i}$ for any $i$. Let $I \subset \{1,\dots,m\}$, and set $\partial_{I}^{[r_I]}:=\partial_{i_1}^{[r_{i_1}]}\cdots \partial_{i_{t}}^{[r_{i_t}]}$, also $z_I^{s_I}:=z_{i_1}^{s_1} \cdots z_{i_t}^{s_t}$ for $r_I,s_I \in \mathbb{N}^{t}$, $i_1 <i_2 < \dots < i_t \in I$.
 Define $S_m(\zz)$ be the $\zz$-submodule of $A_m(\mathbb{Q})$ generated by $z_I^{s_I}\partial_{I}^{[r_I]}$ for any $I$ and $s_I,r_I \in \mathbb{N}^{|I|}$. Note also, that the generators of $S_m(\zz)$ are linear independent over $\zz$. Moreover, over $\mathbb{Q}$ we have $\partial_i^{[r]}=\frac{1}{r!}\partial_{z_i}^{r}$. \\
 For any (commutative) ring $k$, the $k$-module $S_m(k)$ is the base change $S_m(\zz)\otimes k$. 
 Note that if $k \rightarrow k'$ is a ring map, then $S_m(k) \otimes_k k'=S_m(k')$ as $k'$-module.
 \begin{lemma}\label{tensorcrystallweyl}
 The $k$-module $S_m(k)$ is a $k$-algebra, and $S_m(k)\otimes_k k'\simeq S_m(k')$ as $k'$-algebras. Moreover, there is a natural isomorphism of $k$-algebras $S_1(k) \otimes_k \dots \otimes_k S_1(k) \xrightarrow{\sim} S_m(k)$ where the tensor product is taken $m$ times.
 \end{lemma}
 \begin{proof}
 By base change, it is sufficient to prove the statement for $k=\mathbb{\zz}$.\\
More precisely, we need to prove that for any $I$, $J \subset \{1,\dots,m\}$ and $r_I,s_I \in \mathbb{N}^{|I|}$ , $r'_J,s'_J \in \mathbb{N}^{|J|}$, we have that $z_I^{s_I}\partial_{I}^{[r_I]}z_J^{s'_J}\partial_{J}^{[r'_J]} \in S_m(k)$.  First of all, from the equality  $\partial_i^{[r_i]}\partial_j^{[r'_j]}=\partial_j^{[r'_j]}\partial_i^{[r_i]}$ and $z_iz_j=z_jz_i$, for any $i \neq j$, we can rearrange the product $z_I^{s_I}\partial_{I}^{[r_I]}z_J^{s'_J}\partial_{J}^{[r'_J]}$. 
  More precisely, if $i_1 < i_2 < \dots i_t \in I$,   $j_1 < j_2 < \dots j_{t'} \in J$ and $l_1< l_2 <\dots < l_{t''}$ are the elements of $I \cup J$, (with $t'' \leq t+t'$)  we have that 
 \begin{equation}\label{proofsm}
     z_I^{s_I}\partial_{I}^{[r_I]}z_J^{s'_J}\partial_{J}^{[r'_J]}=\prod_{i=1}^{t''}z_{l_i}^{s_{l_i}}\partial_{l_i}^{[r_{l_i}]}z_{l_i}^{s'_{l_i}}\partial_{l_i}^{[r'_{l_i}]}
 \end{equation}
 where we set $r_{l_i}$ (resp. $r',s,s'$) equal $0$ if such elements do not appear in $r_I$ (resp. $r'_J,s_I,s'_J$).  It follows that we can reduce the statement to $m=1$.
 Set $z_I=z$, $r_I=r, s_I=s, r'_J=r', s'_j=s'$. \\
 \textit{Claim:} For any $r,s,r',s' \in \zz_{\geq 0}$ with $s' \geq 1$ the relation
 \begin{equation}\label{claim_Lemma_Sm}
 z^s\partial_z^{[r]}z^{s'}\partial_z^{[r']}=z^s\partial_z^{[r-1]}z^{s'-1}\partial_z^{[r']}+z^{s+1}\partial_z^{[r]}z^{s'-1}\partial_z^{[r']} \end{equation}
 holds. 
 Indeed, for $r=1$, by the equality $\partial_z z-z\partial_z=1$ follows  
 $\partial_z z^{s'}=(1+z\partial_z)z^{s'-1}$, thus  
 $z^s\partial_z^{}z^{s'}\partial_z^{[r']}=z^sz^{s'-1}\partial_z^{[r']}+z^{s+1}\partial_z^{}z^{s'-1}\partial_z^{[r']}$. Moreover, we see that 
 \begin{align*}
     z^s\partial_z^{[r]}z^{s'}\partial_z^{[r']}=(1/{r!})z^s\partial_z^{r}z^{s'}\partial_z^{[r']}=
     (1/r)z^s\partial_z^{[r-1]}z^{s'-1}\partial_z^{[r']}+(1/r)z^s\partial_z^{[r-1]}z\partial_z z^{s'-1}\partial_z^{[r']}
 \end{align*}
 Thus, by applying an induction on $r$ with the term $z^s\partial_z^{[r-1]}z\partial_z$,  we get also 
 \begin{align*}
      (1/r)z^s\partial_z^{[r-1]}z\partial_z z^{s'-1}\partial_z^{[r']}= 
      (1-1/r)z^s\partial_z^{[r-1]}z^{s'-1}\partial_z^{[r']}+ z^{s+1}\partial_z^{[r]}z^{s'-1}\partial_z^{[r']}.
 \end{align*}
 The two latter equalities imply the claim. \\
 Set $g(r,s,r',s'):=z^s\partial_z^{[r]}z^{s'}\partial_z^{[r']}$. Then by definition, $g(r,s,r',0) \in S_1(k)$ for any non negative integers $r,s,r'$.   Now by induction on $s'\geq 1$,  any  $g(r,s,r',s')$ is a finite sum of elements $g(r_1,s_1,r_1',0)$ for some integers $r_1,s_1,r_1'$ by the relation \eqref{claim_Lemma_Sm} above. In particular,  $z^{s}\partial_z^{[r]}z^{s'}\partial_z^{[r']} \in S_1(k)$. \\
 For the last assertion, we notice that the natural  isomorphism of free polynomial algebras over $k$
 \begin{equation}
 k\{z_1,\partial_1\}\otimes_k \cdots \otimes_k k\{z_m,\partial_m\} \xrightarrow{\sim} k\{z_1,\cdots,z_m,\partial_1,\cdots,\partial_m\}
 \end{equation}
 induces an isomorphism
  \begin{equation*} 
 S_1(k) \otimes_k \dots \otimes_k S_1(k) \xrightarrow{\sim} S_m(k). \qedhere
 \end{equation*}
 \end{proof}
 \begin{defn} \footnote{We did not find any reference in the literature for calling such an object, so we took the freedom to give a name.}
 A $m$-th \textit{crystalline Weyl algebra} over a ring $k$ is a $k$-algebra isomorphic to $S_m(k)$ for some $m \in \mathbb{N}$.
  \end{defn} 
  
 \begin{lemma}
    Let $k$ be an integral domain. If $K=\mathrm{Frac}(k)$ is a field of characteristic zero, then  $S_m(k) \otimes _k K\simeq A_m(K)$
 \end{lemma}
 \begin{proof}
    It is  clear, because $\partial_i^{r}=r!\partial_i^{[r]}$ and $r!$ is invertible in $K$.
 \end{proof}
 \begin{ex} \label{example_diff_poly}
(i) Let  $\mathcal{D}(\mathbb{Q}[z_1,\dots,z_m])$ be the $\mathbb{Q}$-algebra of differential operators of $\mathbb{A}_{\mathbb{Q}}^m$. The $\mathbb{Q}$-linear derivations $\frac{\partial}{\partial z_i}=:\partial_i$ satisfy the relations $$z_i\partial_{j}-\partial
_{j}z_i +\delta_{ij}=0, \quad \partial_{i}\partial_{j}=\partial_{j}\partial_{i},$$ therefore we can identify
  $A_m(\mathbb{Q})=\mathcal{D}(\mathbb{Q}[z_1,\dots,z_m])$. Under this identification, a crystalline symbol is  a differential operator sending $\zz[z_1,\dots,z_m]$ to $\zz[z_1,\dots,z_m]$. Now, $\mathcal{D}(\zz[z_1,\dots,z_m])$ is the $\zz$-subalgebra of $\mathcal{D}(\mathbb{Q}[z_1,\dots,z_m])$ generated by those differential operators sending $\zz[z_1,\dots,z_m]$ to $\zz[z_1,\dots,z_m]$. Moreover,
 any   $\partial \in \mathcal{D}(\mathbb{Z}[z_1,\dots,z_m]) \subset \mathcal{D}(\mathbb{Q}[z_1,\dots,z_m]) $ may be uniquely written as 
 $$ \partial= \sum_{I}a_I\partial_{I}^{[r_I]}, \quad a_I \in \mathbb{Q}.$$
Since every crystalline symbol sends $\zz[z_1,\dots,z_m]$ to $\zz[z_1,\dots,z_m]$, it follows  $a_I \in \zz$. This shows that 
 $S_m(\zz)=\mathcal{D}(\zz[z_1,\dots,z_m])$.
 \\
 (ii) For any ring $k$, $S_m(k)=\mathcal{D}(\zz[z_1,\dots,z_m])\otimes k \subset \mathcal{D}(k[z_1,\dots,z_m])$: Indeed, it suffices to prove that any $\partial_i^{[r]} \in S_m(k)$ belongs to $\mathrm{Fil}_r\mathcal{D}(k[z_1,\dots,z_m])$. We can proceed by induction on $r\geq 0$, since $\partial_i^{[0]} \in k[z_1,\dots,z_m]=\mathrm{Fil}_0\mathcal{D}(k[z_1,\dots,z_m])$. For any $a,f \in k[z_1,\dots,z_m]$ we have 
 \begin{equation}
 [\partial_i^{[r]},a](f)=\partial_i^{[r]}(af)-a\partial_i^{[r]}(f)=\sum_{s=0}^{r-1}\partial_i^{[r-s]}(a)\partial_i^{[s]}(f)
 \end{equation}
 where the last equality follows by Formula \eqref{HS-formula}. Then, the claim follows since by inductive hypothesis, the operator $\partial_i^{[r-s]}(a)\partial_i^{[s]} \in \mathrm{Fil}_{r-1}\mathcal{D}(k[z_1,\dots,z_m]) $ for any $s \leq r-1$.
 \end{ex}
 
 When $k$ is a field of characteristic $p$, we have the following characterization:
 \begin{prop}[{cf. \cite[Theorem 2.7]{smith87}}]\label{Smith_theo_2.7}
 Let $k$ be a field of characteristic $p>0$. Set $A=k[z_1,\dots,z_m]$ and $A_n=k[z_1^{p^n},\dots, z_m^{p^n}]$. Then,
 \begin{equation}
     \mathcal{D}(A)=\bigcup_{n=0}^{\infty} \mathrm{End}_{A_n}(A).
 \end{equation}
 \end{prop}
 \begin{ex}\label{example_endomorph_gen_diff}
With the same notation of  \Cref{Smith_theo_2.7}, a basis for the $A_n$-module $\mathrm{End}_{A_n}(A)$ is given by the following maps : for any $\mathbf{s}=(s_1,\dots,s_m)$ such that $0 \leq s_r < p^n$ for any $1 \leq r \leq m$, let
 \begin{equation}
     \theta^{(n)}_{\mathbf{i}\mathbf{j}}(\mathbf{z}^{\mathbf{s}}):=\left\{\begin{array}{ll}
       \mathbf{z}^{\mathbf{i}}   & \text{if } \mathbf{s}=\mathbf{j}  \\
        0  & \text{otherwise} 
     \end{array}\right.
      \end{equation}
      where $\mathbf{i}=(i_1,\dots,i_m),\mathbf{j}=(j_1,\dots,j_m) \in \mathbb{N}^{m}$ are such that $0 \leq i_r,j_r < p^n$ and $\mathbf{z}^{\mathbf{i}}:=z_1^{i_1}\cdots z_m^{i_m}$. \\
Then, the maps $\theta^{(n)}_{\mathbf{i}\mathbf{j}}$ extends uniquely to $A_n$-linearly independent endomorphisms of $A$ and it is straightforward to check that $\theta^{(n)}_{\mathbf{i}\mathbf{j}}=\mathbf{z}^{\mathbf{i}}\partial_{\mathbf{z}}^{[\underline{p^n-1}]}\mathbf{z}^{\underline{p^n-1}-\mathbf{j}}$,
       where $\underline{p^n-1}=(p^n-1,\dots,p^n-1) \in \mathbb{N}^m$ and the sum of $m$-uples is taken component-wise. Therefore, $\theta^{(n)}_{\mathbf{i}\mathbf{j}} \in S_m(k) \subset \mathcal{D}(A)$ for any $n\geq 0$ and $\mathbf{i},\mathbf{j} \in \mathbb{N}^m$. In particular, by  \Cref{Smith_theo_2.7}, it follows that $\mathcal{D}(A)=S_m(k)$.  
 \end{ex}
 
 \begin{lemma}\label{smith2.7}
  Let $k$ be a field of characteristic $p>0$.
     If  $X=\spec(B)$ is a smooth affine $k$-scheme of dimension $m$, then  $\mathcal{D}(X)$ is locally isomorphic to a crystalline Weyl algebra. More precisely,  there exists an open covering $\mathcal{U}=\{U_i\}_i$ of $X$, where $U_i$ is a $k$-scheme étale over $ \mathbb{A}_k^m $, such that there is an isomorphism of $k$-algebras  $\mathcal{D}(U_i)\simeq S_m(k) \otimes_{\mathcal{O}_{\mathbb{A}_k^m}}\mathcal{O}_{U_i}$ for any $i$.
 \end{lemma}
 \begin{proof}
     Since $X=\spec(B)$ is smooth of dimension $m$, there exists an affine open cover $\mathcal{U}=\{U_i=\spec(B_i)\}$ together with étale maps $k[z_1,\dots,z_m] \rightarrow B_i$ . Then we have an isomorphism $\mathcal{D}(U_i)\simeq \mathcal{D}(k[z_1,\dots,z_m]) \otimes_{k[z_1,\dots,z_m]} B_i=S_m(k)\otimes_{k[z_1,\dots,z_m]}B_i$. The last equality is given by  \Cref{Smith_theo_2.7}.
 \end{proof}
 \begin{prop}\label{prop_diff_thickening}
 Assume that $k$ is a field of characteristic $p$. Let $X$ be a smooth scheme over $\w(k)$ and $X_n=X \times_{\spec(\w(k))} \wn(k)$. Then, there is a canonical isomorphism of $\wn(k)$-algebras
 \begin{equation}
 \mathcal{D}(X_n) \xleftarrow{\simeq} \mathcal{D}(X)  \otimes_{\w(k)} \wn(k).
\end{equation}  
 \end{prop}
 \begin{proof}
 It suffices to prove the statement locally, thus we can reduce to consider an affine $\w(k)$-scheme $X=\spec(B)$. Then, $X_n=\spec(B_n)$ where $B_n=B \otimes_{\w(k)} \wn(k)$, hence $B \rightarrow B_n$ is surjective. By smoothness, the modules $\Omega^1_{B/\w(k)}$ and $\Omega^1_{B_n/\wn(k)}$ are free respectively over $B$ and $B_n$ of the same rank. In particular, we have the canonical isomorphism of $B_n$-modules (by \cite[Tag 00RS, Lemma 10.131.7]{stacks-project} the following is an epimorphism and by smoothness, the $B_n$-modules have the same rank):
 \begin{equation}
 \Omega^1_{B/\w(k)} \otimes_B B_n \xrightarrow{\simeq} \Omega^1_{B_n/\wn(k)}
\end{equation} 
Moreover, since $B$ is smooth over $\w(k)$ we have for any $m\geq 1$ the following exact sequence of free $B$-modules (cf. \cite[(4.2.2)]{LT95}) :
\begin{equation}
0 \rightarrow \mathrm{Sym}_B^m(\Omega^1_{B/\w(k)})\rightarrow P_{B/\w(k)}^m\rightarrow P_{B/\w(k)}^{m-1}\rightarrow 0
\end{equation}
where $P_{B/\w(k)}^m$ denotes the $B$-module of principal part of order $m$ as defined in \cite[Definition 16.3.1]{EGA4} (where $P_{B/\w(k)}^0=B$) and $\mathrm{Sym}_B(-)$ is the  symmetric algebra functor. By induction on $m \geq 0$, and by the natural isomorphism 
\begin{equation}\mathrm{Sym}_B^m(\Omega^1_{B/\w(k)}) \otimes_B B_n \xrightarrow{\sim}\mathrm{Sym}_{B_n}^m(\Omega^1_{B_n/\wn(k)})
\end{equation}
 we get also the isomorphism of $B_n$-modules between  principal parts
\begin{equation}
P_{B/\w(k)}^m \otimes_B B_n \xrightarrow{\simeq}P^m_{B_n/\wn(k)}.
\end{equation} 
Thus the statement follows by the chain of identifications:
\begin{align*}
\mathrm{Hom}_{B_n}(P^m_{B_n/\wn(k)}, B_n)=\mathcal{D}_m(B_n) &\xleftarrow{\simeq} \mathrm{Hom}_{B_n}(P^m_{B/\w(k)} \otimes_B B_n, B_n) \\
 &=\mathrm{Hom}_{B}(P^m_{B/\w(k)}, B \otimes_B B_n) \\
 &=\mathrm{Hom}_{B}(P^m_{B/\w(k)}, B )\otimes_B B_n= \mathcal{D}_m(B) \otimes_B B_n
 \\
 &=\mathcal{D}_m(B) \otimes_{\w(k)} \wn(k).
\end{align*}
where from the second to third line we use that $P^m_{B/\w(k)}$ is free over $B$.
 \end{proof}
 \begin{rem}
     We notice that analogously to the  \Cref{example_diff_poly} (i), since $\w(k)$ is torsion free, then $S_m(\w(k))=\mathcal{D}(\w(k)[z_1,\dots,z_m])$. The proof of  \Cref{prop_diff_thickening} does not depend on  \Cref{Smith_theo_2.7}, therefore we may use it to deduce the equality $S_m(k)=\mathcal{D}(k[z_1,\dots,z_m])$. Moreover, by using the maps $\theta_{\mathbf{i}\mathbf{j}}^{(n)}$ constructed in  \Cref{example_endomorph_gen_diff}, we see that $S_m(k)=\bigcup_{n\geq 0}\mathrm{End}_{A_n}(A)$. Putting all together yields a proof of  \Cref{Smith_theo_2.7}. 
 \end{rem}

\section{Some classes of $\wn\ox$-modules}\label{sec_de Rham_Witt_cx}
In this chapter we introduce the main geometric objects of this thesis. We recall the main properties of the De Rham-Witt complex as defined in the absolute setting by Deligne and Illusie \cite{Illusie79} and we introduce a less standard concept of Witt line bundles, as recently studied in \cite{TH22} by the name of \textit{Teichm\"uller lift of line bundles}. These objects will be used as coefficients for  cohomology of schemes equipped with an action of a finite group $G$. Their natural structure of $\wn(k)[G]$-module will be further investigated in the particular case of the Drinfeld's upper half space in later chapters.
\subsection{The (Bloch-Deligne-Illusie) de Rham-Witt complex}
In this section we recall the notion of a de Rham-Witt complex of a $\fp$-scheme $X$, following the classical paper \cite{Illusie79}. We focus on some equivariant aspects that arise by assuming that $X$ is equipped with an action of a finite group $G$. In particular, the concept of $G$-linearization extends for $\wn\ox$-modules (see \Cref{appendix_witt_vectors} for the basic definitions of the ring of Witt vectors and its sheafification) and in the case of the de Rham-Witt complex, a natural $G$-linearization arises from that one on the de Rham complex.

\subsection{Action of $G$ on $\wn\ox$}
In the following we introduce the concept of $\ve$-pro-complex. In the next section we will see that the category of $\ve$-pro-complexes admits an initial object given by the de Rham-Witt complex. As warm up, we describe a $G$-linearization on $\wn\ox$ as a lift of the canonical $G$-linearization of $\ox$ in the category of $\ve$-pro-complexes. The same reasoning will be applied to the de Rham-Witt complex.
In this chapter $k$ will denote a perfect $\fp$-algebra.

\begin{defn}
Let $B$ be a ring and $A$ be a commutative $B$-algebra. Then, $(M^{\bullet},d)$ is said to be \textit{a differential graded $A$-algebra over $B$} if the following conditions hold:
\begin{itemize}
\item[i)] $(M^{\bullet},d)$ is a complex of abelian groups;
\item[ii)] For any $i \in \zz$, $M^i$ is a $A$-module and $d \colon M^i \rightarrow M^{i+1}$ is a $B$-linear map;
\item[iii)] For any $x \in M^i, y \in M^j$ , $i,j \in \zz$, the following relations hold:
\begin{equation}
xy=(-1)^{ij}yx, \quad d(xy)=(dx)y+(-1)^ixdy, \quad x^2=0 \text{ if the degree of } x \text{ is odd.}
\end{equation}
\end{itemize}
\end{defn}
Notice that the definition above without ii) corresponds to that of differential graded

 ($\zz$-)algebra (over $\zz$).

A $\ve$-pro complex is a projective system of (sheaves of) differential graded algebras (over $\zz$) satisfying certain relations.

\begin{defn}
Let $X$ be a $k$-scheme.
A $\ve$-pro-complex on $X$ consists of the following data:
\begin{itemize}
\item[a)] A projective system of sheaves of differential $\zz$-graded algebras over $\zz$ on $X$, $\{R \colon \mathcal{M}^{\bullet}_n \rightarrow \mathcal{M}^{\bullet}_{n-1}\}_{n \in \zz}$,
\item[b)] A collection of morphisms of sheaves of graded abelian groups on $X$, called \textit{Verschiebung maps} $\{\ve \colon \mathcal{M}^{\bullet}_n \rightarrow \mathcal{M}_{n+1}^{\bullet}\}_{n \in \zz}$,
\end{itemize}
satisfying the following conditions:
\begin{enumerate}
     \item For any $n,r \in \zz$, $\mathcal{M}_n^r$ is a quasi-coherent $\wn\ox$-module such that $\mathcal{M}_n^r=0$ if $n \leq 0$ or $r<0$; $\mathcal{M}_1^0$  is a sheaf of $k$-algebras  and  $\mathcal{M}_n^0=\wn(\mathcal{M}_1^0)$. Moreover, the maps $R \colon \w_{n+1}(\mathcal{M}_1^0) \rightarrow \w_{n}(\mathcal{M}_1^0)$ and $\ve \colon \w_{n}(\mathcal{M}_1^0)\rightarrow \w_{n+1}(\mathcal{M}_1^0)$ agree respectively with the canonical restriction and Verschiebung map of Witt vectors,
     \item The Verschiebung maps are compatible with the restrictions, i.e.  $R \ve=\ve R$ holds,
    \item $\ve(xdy)=\ve(x)d\ve(y)$ for $x,y \in \mathcal{M}_n^0$,
     \item $\ve(y)d[x]_{n+1}=\ve([x]_n^{p-1}y)d\ve([x]_n)$ for all $x \in \mathcal{M}_1^0$ and $y \in \mathcal{M}_n^{0}$.
\end{enumerate} 
A morphism of $\ve$-pro-complexes is a collection of graded differential algebra morphisms $f_n^{\bullet}: \mathcal{M}_n^{\bullet}\rightarrow\mathcal{N}_n^{\bullet}$ such that they are compatible with $R$ and $\ve$, and $f_n^0=\wn(f_1^0)$. 
\end{defn}
If we consider $\wn\ox$ itself as a trivial  graded differential algebra complex with degree (with respect to the index $r$ in the definition above) concentrated in $0$, with $d=0$ and $R$ and $\ve$ as usual, it is a $\ve$-pro-complex. As abuse of notation, we choose to omit the subscript of the Teichm\"uller representative, when it is clear from the context. \\
From now on, let us consider a finite group $G$ as constant $k$-group scheme and suppose $G$ acts on $X$ via $\sigma\colon G \times_k X \rightarrow X$.


Denote by $\text{pr}_1 \colon G \times_k X \rightarrow X$ the canonical projection, $\text{pr}_{12}: G \times_k G \times_k X \rightarrow G \times_k X$ the projection $(g_1,g_2,x) \mapsto (g_2,x)$ for $g_1,g_2 \in G$, $x \in X$, and $\text{pr}_2=\text{pr}_1\circ \text{pr}_{12}$.

For any morphism of schemes over $k$, $f \colon X \rightarrow Y$, and a $\wn\oo_Y$-module $\mathcal{F}$, define the pullback $f^*\mathcal{F}$ as the $\wn\ox$-module 
\begin{equation}
    f^*\mathcal{F}=f^{-1}\mathcal{F} \otimes_{f^{-1}\wn\oo_Y} \wn\ox.
\end{equation}

\begin{defn} 
Let $\mathcal{F}$ be a quasi-coherent $\wn\ox$-module. A $G$-\textit{linearization} on $\mathcal{F}$ is an isomorphism 
\begin{equation}
    \phi \colon \sigma^*\mathcal{F} \rightarrow \mathrm{pr}_1^*\mathcal{F}
\end{equation}
satisfying the following cocycle condition
\begin{equation}\label{cocycle_cond_wn}
    (m \times id_{X})^*\phi=\mathrm{pr}_{12}^*\phi \circ (\mathrm{id}_G \times \sigma)^*\phi.
\end{equation}
The $\wn\ox$-module $\mathcal{F}$ is said to be $G$-\textit{equivariant} if it has a $G$-linearization.
\end{defn}

\begin{defn}\label{lift}
Let $\mathcal{M}_n^{\bullet}$, $\mathcal{N}_n^{\bullet}$ be $\ve$-pro complexes on $X$ and let $f_1^{\bullet}: \mathcal{M}^{\bullet}_1\rightarrow\mathcal{N}^{\bullet}_1$ be a collection of $\ox$-module morphisms. A morphism of $\ve$-pro complexes $\tilde{f}_n^{\bullet}: \mm^{\bullet}_n \rightarrow \nn^{\bullet}_n$ is called a \textbf{lift of} $f_1^{\bullet}$ if $R^{n-1}\circ \tilde{f}_n^r=f_1^r \circ R^{n-1}$ at each $n,r \in \mathbb{Z}$.
\end{defn}
\begin{rem}
This definition says in particular that $\tilde{f}_1^{r}=f_1^{r}$ and for $r=0$ we have $\tilde{f}^0_n=\wn(f_1^0)$ is uniquely determined. 
\end{rem}
So we see that in the case of ($\wn\ox,0,R$) the notion of lifts is intrinsic in the description of $\wn$. Thus, any automorphism of $\ox$ lifts to a unique $\ve$-pro-complex automorphism of $\wn\ox$. In particular, $\sigma_g$ induces a $G$-linearization, given by $(\tilde{\phi}_g)_n=\wn(\phi_g)$.
\begin{rem}
The action given by $\tilde{\phi}_g$ is the "natural" one if we consider $\wn\ox \simeq \ox^n$ as a set.\footnote{we cannot ask more structure in order to get such an isomorphism, for example as $\ox$-modules.}
More precisely, we are saying that the rule $\tilde{\phi}_g(f_0,\dots,f_{n-1})=(\phi_g(f_0),\dots,\phi_g(f_{n-1}))$ defines a ring morphism, that simply is the functoriality of $\wn(-)$.
\end{rem}

Rephrasing last remark yields the following:
\begin{lemma} The map
$\tilde{\phi}_g \in \mathrm{End}(\wn\ox)$ is the unique morphism of $\ve$-pro complexes lifting $\phi_g \in \mathrm{End}(\ox)$.
\end{lemma}
We will see that the analogous lifting property (cf.  \Cref{unpropdrW}) holds for the de Rham-Witt complex.
\section{Definition and properties of the de Rham-Witt complex}
  Here we are going to recall definition and properties of the Bloch-Deligne-Illusie de Rham-Witt complex, following mainly \cite{Illusie79} and \cite{Langer2004DERC}.
  
\begin{prop} \label{derhamcompl_is_initial}
   For any commutative $k$-algebra $A$, there exists an object $(\wn\Omega_A^{\bullet}, R, \ve)$ in the category of $\ve$-pro-complexes on $\spec(A)$ such that $\wn\Omega_A^{\bullet}$ is a differential graded $\wn(A)$-algebra (dga) over  $\wn(k)$ and  for any other $\ve$-pro-complex $(M_n^{\bullet},R',\ve')$ such that $M_n^{\bullet}$ is a $\wn(A)$-dga over $\wn(k)$ there exists a unique map of $\ve$-pro-complexes
   \begin{equation}
       \wn\Omega_A^{\bullet} \rightarrow M_n^{\bullet}.
   \end{equation}
\end{prop}
\begin{proof}
    In the case $k=\fp$ this is \cite[I.1, Theorem 1.3]{Illusie79}. Since $k$ is perfect, the statement follows by loc. cit. I.1, Proposition 1.9.2.
\end{proof}
\begin{rem}[Universal property of the de Rham-Witt complex]\label{unpropdrW}
By \Cref{derhamcompl_is_initial}, in particular it follows that the functor $C \mapsto \wn\Omega_{C}^{\bullet}$ from $k$-algebras  to $\ve$-pro-complexes on $\spec(k)$ is fully faithful, i.e. for any $k$-algebras $A,B$ there is a natural bijection $$\mathrm{Hom}_{\ve-pro-c.(\spec(k))}(\wn\Omega_{A}^{\bullet},\wn\Omega_B^{\bullet}) \xrightarrow{\sim} \mathrm{Hom}_{k-alg}(A,B), \quad f \mapsto f_1^0.$$
\end{rem}
We recall the construction of $\wn\Omega_{A}^{\bullet}$.
Let us denote the de Rham complex (relative to $\zz$)  of any commutative ring $B$ by $\Omega_{B}^{\bullet}$. We proceed inductively on $n \geq 1$. Define $\w_1\Omega_A^{\bullet}:=\Omega_{A}^{\bullet}$ and suppose that we know $\w_i\Omega_A^{\bullet}$ for any $1 \leq i \leq n$ such that $\w_i\Omega_A^0=\w_i(A)$. The restriction $ R \colon \w_{n+1}(A) \rightarrow \wn(A)$ induces a dga morphism  $\Omega_{\w_{n+1}(A)}^{\bullet} \rightarrow \Omega_{\w_{n}(A)}^{\bullet}$. For any $n$ one can define a collection $\{I_{n}^{\bullet}\}_{n \in \mathbb{N}}$ of dga ideals for $\Omega_{\wn(A)}^{\bullet}$ and additive maps (Verschiebung) $ \ve \colon \wn\Omega_{A}^{\bullet} \rightarrow \Omega_{\w_{n+1}(A)}^{\bullet}/I_{n+1}^{\bullet}$ such that the following hold:
\begin{itemize}
 \item[1)] $\w_{n+1}\Omega_A^{\bullet}=\Omega_{\w_{n+1}(A)}^{\bullet}/I_{n+1}^{\bullet}$ and $\w_{n+1}\Omega_A^{0}=\Omega_{\w_{n+1}(A)}^{0}=\w_{n+1}(A)$;
    \item[2)] $\ve \colon \wn(A)=\wn\Omega_A^0 \rightarrow \w_{n+1}\Omega_A^0=\w_{n+1}(A)$ is the Verschiebung of Witt vectors and $ \ve(y)d[x]-\ve([x]^{p-1}y)d\ve([x]) \in I_{n+1}^1 \text{ for any } x \in A, y \in \wn(A)$ ;
       \item[3)] $\ve(adx_1\dots dx_i)=\ve(a)d\ve(x_1)\dots d\ve(x_i)$
for any $a, x_1,\dots,x_i \in \w_{n}(A)$, and  $adx_1\dots dx_i=\pi_{n}(a \otimes dx_1 \otimes \dots \otimes dx_i) \in \pi_{n}(\Omega^i_{\w_{n}(A)})$
where   $\pi_{n}:\Omega_{\w_{n}(A)}^{\bullet} \rightarrow \w_{n}\Omega_A^{\bullet} $ is the surjective canonical map;
\item[4)] $R(I^{\bullet}_{n+1})\subset I^{\bullet}_n$, thus it induces a map $R \colon \w_{n+1}\Omega_A^{\bullet} \rightarrow \wn\Omega_A^{\bullet}$. 
\end{itemize}
Moreover, $\ve$ is a map of complexes respect to $n$ and $R$: 
indeed condition 2) and 3) determine $\ve$ uniquely, then it follows from 4) and because $R$ is a dga morphism. 

The de Rham-Witt complex is equipped with a Frobenius operator lifting the Frobenius on the Witt vectors $F \colon \w_{n+1}(A) \rightarrow \wn(A)$. More precisely, there exists a unique morphism of projective system (respect to $n$ and $R$) of dga's 
\begin{equation}
F \colon \w_{n+1}\Omega_A^{\bullet} \rightarrow \wn\Omega_A^{\bullet}
\end{equation} 
such that 
\begin{center}
5)  $ F d[x]=[x]^{p-1}d[x]$;  \quad 
6)  $ Fd\ve=d \colon \wn(A) \rightarrow \wn\Omega^1_A$.
\end{center}
This is the content of the loc. cit. I , Theorem 2.17. Moreover, the following relations between $F,\ve,d, R$ hold (cf. loc.cit. I, Proposition 2.18):
\begin{align*}
&\text{7) } {F\ve}={\ve F}=p \colon \mathrm{W}_n \Omega_{\mathrm{A}}^i  \rightarrow \mathrm{W}_n \Omega_{\mathrm{A}}^i ;
 & \text{10) }  {F} d \mathrm{V}=d\colon \mathrm{W}_n \Omega_{\mathrm{A}}^i \rightarrow \mathrm{W}_n \Omega_{\mathrm{A}}^{i+1}; & 
\\
&\text{8) } d{F}=p {F} d: \quad \mathrm{W}_n \Omega_{\mathrm{A}}^i \rightarrow \mathrm{W}_{n-1} \Omega_{\mathrm{A}}^{i+1} ; 
& \text{11) }  x \mathrm{V} y=\mathrm{V}(y{F}(x) ), \quad x \in \mathrm{W}_n \Omega_{\mathrm{A}}^i, \quad y \in \mathrm{W}_{n-1} \Omega_{\mathrm{A}}^j. 
\\
&\text{9) }  \mathrm{V} d=p d \mathrm{V}\colon \mathrm{W}_n \Omega_{\mathrm{A}}^i \rightarrow \mathrm{W}_{n+1} \Omega_{\mathrm{A}}^i ; &
\end{align*}
For any $k$-scheme $X$ and any open affine $\spec(A)=U \subset X$ there is a unique quasi-coherent sheaf of $\wn\ox$-modules, namely  $\wn\Omega^i_{X}$,  such that 
\begin{equation}
  \wn\Omega^i_{X} \colon  U \longmapsto \Gamma(U,\wn\Omega^i_X):=\wn\Omega_A^i.
\end{equation}
We call  $\wn\Omega^i_{X}$ the \textit{$i$-th Witt differential module}, or \textit{$i$-th Hodge-Witt module}.

\begin{rem}
i) Since the Hodge-Witt modules are quasi-coherent, in particular, they behave well under localisation maps: more precisely, if $A \rightarrow B$ is a localisation map, then the natural map $\wn(B) \otimes_{\wn(A)} \wn\Omega_A^i \rightarrow \wn\Omega_B^i$ is an isomorphism. \\
ii) For any morphism of $k$-schemes $f \colon X \rightarrow Y$,$f^{-1}\wn\Omega_Y^{\bullet}$ is a $\ve$-pro-complex. Thus the natural map $f^{-1}\mathcal{O}_Y \rightarrow \ox$ induces a morphism of $\ve$-pro-complexes $f^{-1}\wn\Omega_Y^{\bullet} \rightarrow \wn\Omega_X^{\bullet}$.
In particular, for a point $i_x \colon x \hookrightarrow X$, from the equality $i_{x}^{-1}\ox=\mathcal{O}_{X,x}$ follows that there is a natural \footnote{The notation $\Omega^{\bullet}_{X,x}$ stands for $\Omega^{\bullet}_A$ where $A=\oo_{X,x}$.} isomorphism $i_x^{-1}\wn\Omega_{X}^{\bullet}=(\wn\Omega_{X}^{\bullet})_x \xrightarrow{\sim} \wn\Omega_{X,x}^{\bullet}$.
\end{rem} 
Moreover, Hodge-Witt modules behave well under étale morphisms of $k$-schemes. More precisely, the following holds (cf. loc. cit. I, Proposition 1.14):
\begin{prop}\label{prop_étale_derhamwitt}
Let $X \xrightarrow{f} Y$ be an étale morphism between $k$-schemes $X,Y$. Then, the natural map of $\wn\ox$-modules
\begin{equation}
    f^{*}\wn\Omega_{Y}^i \rightarrow \wn\Omega_X^i
\end{equation}
is an isomorphism.
\end{prop}

The action of $G$ on $X$ induces a $G$-linearization on the Hodge-Witt modules of $X$:

\begin{lemma}
If $G$ is a finite group acting on a $k$-scheme $X$, then for any $i \geq 0$,  $\wn\Omega^i_{X}$ is canonically a $G$-equivariant quasi-coherent $\wn\ox$-module.
\end{lemma}
\begin{proof}
Since $G$ is finite, $\sigma$ is a local isomorphism (because $\spec(k) \subset G$ is open). In particular, $\sigma$ is étale; the same is true for $\mathrm{pr}_1$ (for this  case we can just notice that $G \times_k X$ as a scheme over $X$ is a finite disjoint union of copies of $X$) thus by the \Cref{prop_étale_derhamwitt}  there are isomorphisms 
\begin{equation}
    \sigma^*\wn\Omega_X^i \xrightarrow{\sim} \wn\Omega^i_{G\times_k X}; \quad  \mathrm{pr}_1^*\wn\Omega_X^i \xleftarrow{\sim} \wn\Omega^i_{G\times_k X}
\end{equation}
Composing the two maps above we get an isomorphism $\phi\colon \sigma^*\wn\Omega_X^i \rightarrow  \mathrm{pr}_1^*\wn\Omega_X^i$. To verify the cocycle condition 
\begin{equation}
    (m \times id_{X})^*\phi=\mathrm{pr}_{12}^*\phi \circ (\mathrm{id}_G \times \sigma)^*\phi
\end{equation}
we first notice that it is well defined by the relations $\sigma \circ \mathrm{pr}_{12}= \mathrm{pr}_1 \circ (\mathrm{id}_G \times \sigma) $; $\mathrm{pr}_1\circ (m \times \mathrm{id}_X)=\mathrm{pr}_1 \circ \mathrm{pr}_{12}=\mathrm{pr}_2$ and (by definition of action) $\sigma \circ (m \times \mathrm{id}_X)=\sigma \circ (\mathrm{id}_G \times \sigma)$. Furthermore, consider $\phi$ as map of $\ve$-pro-complexes  varying $n$ and $i$. Denote by $\phi_{1}^0 \colon \sigma^*\ox \rightarrow  \mathrm{pr}_1^*\ox $ the map $\phi$ when $n=1$ and $i=0$. This is the natural $G$-linearization of the structure sheaf $\ox$, in particular the cocycle condition holds for $\phi_{1}^0$. By  \Cref{unpropdrW} the maps  $(m \times id_{X})^*\phi,\mathrm{pr}_{12}^*\phi,(\mathrm{id}_G \times \sigma)^*\phi $ are the unique morphisms of $\ve$-pro-complexes induced respectively by $(m \times id_{X})^*\phi_{1}^0,\mathrm{pr}_{12}^*\phi_{1}^0,(\mathrm{id}_G \times \sigma)^*\phi_{1}^0$, therefore the cocycle condition for 
$\phi$ follows again by universal property of the de Rham-Witt complex.
\footnote{Alternative argument: $\sigma$ induces an action $G \times \wn(X) \rightarrow \wn(X)$. The classical K\"ahler differentials $\Omega^i_{\wn(X)}$ are in this way $G$-equivariant and, by naturality, the projection $\pi_n$ is compatible with such linearizations. The cocycle condition follows from this compatibility.}
\end{proof}

Now suppose that $X$ is a smooth $k$-scheme. 
\begin{prop}[{\cite[I, Prop. 3.7 (a)]{Illusie79}}]
If $X/k$ is smooth of  dimension $N$ then $\w_{n+1}\Omega_X^i=0$ for $i>N$.
\end{prop}
 For every $n \geq 0$, $\w_{n+1}\Omega^{\bullet}_{X}$ is equipped with the following \textit{canonical filtration} of dga's:
\begin{equation}\label{canonicfilwitt}
\mathrm{Fil}^m\w_{n+1}\Omega^{\bullet}_{X}= \left \{\begin{array}{lc}\w_{n+1}\Omega^{\bullet}_{X} &\text{ if } m \leq 0 \\
\ker(\w_{n+1}\Omega^{\bullet}_{X} \xrightarrow{R^{n+1-m}}  \w_m\Omega^{\bullet}_{X}) &\text{ if }1 \leq m < n+1 \\
0 &\text{ if } m \geq n+1 
\end{array} \right.
\end{equation}
Denote by $F_X \colon X \rightarrow X$ the absolute Frobenius.
Then $F_{X*}\Omega_X^i$ is the sheaf of abelian groups $\Omega_X^i$ with a structure of  $\ox$-module induced by $F_X^{\sharp} \colon \ox \rightarrow F_{X *}\ox$. We recall the definition of the (inverse) Cartier operator.
\begin{prop}[{\cite[Theorem 7.2]{katz1970nilpotent}}]\label{prop_cartier_iso}
There is a unique $\ox$-module map 
\begin{equation}
C_X^{-1} \colon \Omega_X^i \rightarrow h^{i}(F_{X*}\Omega_X^{\bullet}):=\frac{\ker(F_{X*}\Omega_X^{i} \xrightarrow{d} F_{X*}\Omega_X^{i+1})}{d(F_{X*}\Omega_{X}^{i-1})}
\end{equation}
such that  $C_X^{-1}(dx)=[x^{p-1}dx] \in h^1$ for any local section $x \in \ox$, $C_X^{-1}(\eta \omega)=C_X^{-1}(\eta)C_X^{-1}(\omega)$ for $\eta \in \Omega_X^i, \omega \in \Omega_X^j$; and $C_X^{-1}=F_X^{\sharp}$ for $i=0$. Moreover, $C_X^{-1}$ is an isomorphism.
\end{prop}

We need $C_X^{-1}$ in order to define some abelian subsheaves 
 $$
B_n \Omega_{X}^i \subset Z_n \Omega_{X}^i \subset \Omega_{X}^i
$$
by letting
$$
\begin{gathered}
B_0 \Omega_{X }^i=0, \quad Z_0 \Omega_{X}^i=\Omega_{X}^i \\
B_1 \Omega_{X}^i=B \Omega_{X }^i=d \Omega_{X}^{i-1}, \quad Z_1 \Omega_{X}^i=Z \Omega_{X }^i=\operatorname{Ker}(d: \Omega_{X }^i \rightarrow \Omega_{X }^{i+1})
\end{gathered}
$$
and inductively on $n$, $B_{n+1} \Omega_{X}^i$, respectively $Z_{n+1} \Omega_{X}^i$ are the unique subsheaves such that 
\begin{equation}
C_X^{-1}(B_{n} \Omega_{X}^i) \xrightarrow[]{\simeq} B_{n+1} \Omega_{X}^i/B \Omega_{X}^i; \quad
C_X^{-1}(Z_{n} \Omega_{X}^i)\xrightarrow[]{\simeq} Z_{n+1} \Omega_{X}^i/B \Omega_{X}^i
\end{equation}
 We have the following relation between the graded module associated to the canonical filtration above and the sheaves $B_{n+1} \Omega_{X}^i,Z_{n+1} \Omega_{X}^i$:
\begin{prop}[{\cite[I, Corollary 3.9]{Illusie79}}]\label{prop_exact_seq_gr}
Let $X/k$ be smooth and for any $n,i \geq 0$, let  $\mathrm{gr}^n\w_{n+1}\Omega_X^i$ be the $n$-the graded piece of the filtration \eqref{canonicfilwitt}.  
Then,
\begin{itemize}
\item[a)] $\mathrm{gr}^n\w_{n+1}\Omega_X^i=\mathrm{Fil}^n\w_{n+1}\Omega_X^i=\ve^n\Omega_X^i+ d\ve^n\Omega_X^i$ where we identify $\Omega_X^i=\w_1\Omega_X^i$;
\item[b)] If $\mathrm{gr}^n\w_{n+1}\Omega_X^i$ is equipped with the structure of a $\ox$-module induced by $$F \colon \ox=\w_{n+1}\ox/\ve\w_{n}\ox \rightarrow \w_{n+1}\ox/p\w_{n+1}\ox, \quad (\text{Notice: } p(\mathrm{gr}^n\w_{n+1}\Omega_X^i)=0)$$ there is an exact sequence of $\ox$-modules:
\begin{equation}\label{exactsequhodgewitt}
0 \rightarrow F_{X*}^{n+1}\frac{\Omega_X^i}{B_n\Omega_X^i}\xrightarrow{\ve^n} \mathrm{gr}^n\w_{n+1}\Omega_X^i \xrightarrow{\beta_n} F_{X*}^{n+1}\frac{\Omega_X^{i-1}}{Z_n\Omega_X^{i-1}} \rightarrow 0
\end{equation}
where $\beta_n $ is the map sending an element of the form $\ve(x)+d\ve(y)$ to the class of $y$. Furthermore,
$F_{X*}^{n+1}\frac{\Omega_X^{i}}{B_n\Omega_X^{i-1}}$ and $F_{X*}^{n+1}\frac{\Omega_X^{i-1}}{Z_n\Omega_X^{i-1}}$ are locally free $\ox$-modules.
\end{itemize}
\end{prop}
\begin{rem}
 Assume that $G$ acts on $X$. Since $\sigma$ and $\mathrm{pr}_1$ are étale, in particular flat, morphisms, it follows that $\sigma^*,\mathrm{pr}_1^*$ are exact functors $\mathrm{QCoh}(G \times_k X)\rightarrow \mathrm{QCoh}(X)$. Since $\Omega^i_X$ is canonically a $G$-equivariant $\ox$-module for any $i$, and $d$ is a $G$-equivariant morphism, the same holds for $Z\Omega_X^i, B\Omega_X^i$. Furthermore, the Cartier operator is functorial on maps of $k$-schemes, hence it is $G$-equivariant. Thus, also $Z_{n+1}\Omega_X^i, B_{n+1}\Omega_X^i$ are such. We conclude the short exact sequence \eqref{exactsequhodgewitt} becomes $G$-equivariant.
 \end{rem}
\subsection{Description of $\wn\Omega^{\bullet}_{k[x_1,\dots,x_d]}$}
We are going to describe the De Rham-Witt complex of the $d$-dimensional affine space $\mathbb{A}^d_k$. \\
A \textit{weight function} is a map of sets $r\colon [1,d] \rightarrow \zz[1/p]_{\geq 0}$\footnote{$[1,d]=\{i \in \zz : 1 \leq i \leq d\}$ and $\zz[1/p]_{\geq 0}=\{ap^{n}: a,n \in \zz, a \geq 0\}$.}. The support $\mathrm{supp}(r) \subset [1,d]$ is the subset of elements $j$ such that $r(j):=r_j \neq 0$. We say that $r$ is nonzero if its support is nonempty. We fix a total order on
 $\mathrm{supp}(r)$ such that \footnote{In particular, for $a \neq b$ such that $v_p(r_a)=v_p(r_b)$ one can fix any order ($a \prec b$ or $b \prec a$) but it has to be the same after multiplying by $p^m$.} 
 \begin{equation}\label{ordersupp}
\forall \text{ } a,b \in \mathrm{supp}(r), \text{ } a \prec b \Leftrightarrow \left\{ 
\begin{array}{ll}
v_p(r_a) \leq v_p(r_b)  \text{ and }& \\
\text{the ordering of  } \mathrm{supp}(r) \text{ and } \mathrm{supp}(p^mr) \text{ agree } \forall m \in \zz. &
\end{array}\right.
 \end{equation}
 For any ordered subset $I\subset  \mathrm{supp}(r)$, we define the weight $r_I$ as the restriction of $r$ to $I$ and $0$ on the complement. Moreover, let
 \begin{equation} 
 t(I):=t(r_I):=\left\{ \begin{array}{lc}
 -\mathrm{min}_{a \in I}\{v_p(r_a)\} & \text{ if } r\neq 0 \\
 0  & \text{ if } r= 0
 	\end{array}\right., \qquad u(I):=u(r_I):=\mathrm{max}\{0,t(I)\}.
 	 \end{equation} We call $r$ \textit{integral} if and only if $t(r)\leq 0$. An integral weight $r$ is called \textit{primitive} if and only if $t(r)=0$. For any weight $r$, $p^{t(r)}r$ is primitive and $p^{u(r)}r$ is integral. Let $(I_0,\dots,I_i)$ be  a $(i+1)$-tuple of pairwise disjoint subsets of $\mathrm{supp}(r)$ satisfying the following properties: 
 \begin{itemize}
 \item[i)] $I_0 \sqcup I_1 \cdots \sqcup I_i=\mathrm{supp}(r)$;
 \item[ii)] $I_j \neq \emptyset$ if $i \geq j\geq 1$;
 \item[iii)] For any $j=0,\dots,i-1$, any element in $I_j$ is smaller then any element in $I_{j+1}$ with respect to the total order \eqref{ordersupp}.
 \item[iv)] For any $0 \leq j \leq i$, if $a,b \in I_j$ with $a\prec b$, then for every $ c \in \mathrm{supp}(r)$ with $a \prec c \prec b$, we have $c \in I_j$.
 \end{itemize}
 Denote by $\mathcal{P}^{(i)}_r$ be the set of such partitions. 
 \begin{rem}
 The conditions i) and ii) yield $|\mathrm{supp}(r)| \geq i$. Suppose 
 $|\mathrm{supp}(r)|=l \geq i$ and write $\mathrm{supp}(r)=\{a_1 \prec a_2 \prec \dots \prec a_l\}$. The property iv) implies that any $I_j$ with $|I_j|=c_j \neq 0$  is of the form $I_j=\{a_{s_j},a_{s_j+1},\dots a_{s_j+c_j-1} \}$. By ii), $c_j \neq 0$ and by iii), $s_j+c_j-1 < s_{j+1}$ for any $j=1,\dots,i-1$. Also, from i) it follows $s_j+c_j=s_{j+1}$. If $c_0=0$, then $a_1 \in I_1$ and the $i$-tuple $(c_1,\dots,c_i)$ is a partition of $l$ of length $i$ made of positive integers. Moreover, $c_1$ determines $I_1=\{a_1,\dots,a_{c_1}\}$ and so the set $I_{2}$ is uniquely determined by its size $|I_{2}|=c_{2}$. Inductively, the set $I_{j}$ is determined for any $j$. In the case $c_0 \neq 0$, $(c_0,\dots,c_i)$ is a partition of $l$ of length $i+1$ made of positive integers. Analogously, $c_0$ determines $I_0=\{a_1,\dots,a_{c_0}\}$, therefore any set $I_j$ with $j \geq 1$ is (inductively) uniquely determined by its size $c_j$. Hence, there is a bijection between the set $ \mathcal{P}^{(i)}_r$ and (ordered) partitions of $|\mathrm{supp}(r)|$ of length $i$ and $i+1$ made by positive integers, such that a $(i+1)$-tuple $(I_0,\dots,I_i) \in \mathcal{P}^{(i)}_r $ with $I_0=\emptyset$ (resp. $I_0 \neq \emptyset$) corresponds to a partition of length $i$ (resp. $i+1$). 
 \end{rem}
 
 Let $T_j:=[x_j] \in \wn(k[x_1,\dots,x_d])$ for any $1 \leq j \leq d$ and for any integral weight $r$, $T^r:=T_1^{r_1}\cdots T_d^{r_d}$. Let us define the following elements of $\wn(k[x_1,\dots,x_d])$ and $\wn\Omega^{1}_{k[x_1,\dots,x_d]}$:
 \begin{equation}
 e^0_n(1,r):=\ve^{u(r)}(T^{p^{u(r)}r}), \quad e_n^1(1,r)=\left\{\begin{array}{ll} d\ve^{u(r)}(T^{p^{u(r)}r}) & \text{if }r \text{ is not integral}\\
     F^{-t(r)}dT^{p^{t(r)}r}    & \text{if }r \text{ is integral.}
 \end{array}\right.
 \end{equation}

Let $P \in \mathcal{P}^{(i)}_r$. We combine those elements to get the following elements in $\wn\Omega^{i}_{k[x_1,\dots,x_d]}$: 
 \begin{equation}
 e_n(1,r,P)=\left\{\begin{array}{ll} e^0_n(1,r_{I_0})e^1_n(1,r_{I_1})\cdots e^1_n(1,r_{I_i}) & \text{if } I_0 \neq \emptyset\\
       e^1_n(1,r_{I_1})\cdots e^1_n(1,r_{I_i})  & \text{if } I_0=\emptyset.
 \end{array}\right.
 \end{equation}
 The elements $e_n(1,r,P)$ satisfy the following relations with the operators $F,\ve,d$ of the de Rham-Witt complex:
 \begin{align*}
 Fe_n(1,r,P)= &\left \{ \begin{array}{lcl}
 pe_n(1,p^{}r,P) &\text{ if }& I_0\neq \emptyset \text{ and } r \text{ not integral} \\
 e_n(1,p^{}r,P) &\text{ if }& I_0=\emptyset \text{ or } r \text{  integral}
\end{array}\right. \\
\ve e_n(1,r,P)= &\left \{ \begin{array}{lcl}
 e_n(1,p^{-1}r,P) &\text{ if }& I_0\neq \emptyset \text{ or } p^{-1}r \text{ integral} \\
  pe_n(1,p^{-1}r,P) &\text{ if }& I_0=\emptyset \text{ and } p^{-1}r \text{ not  integral}
\end{array}\right. \\
de_n(1,r,P)= &\left \{ \begin{array}{lcl}
0 & \text{ if }& I_0=\emptyset \\
e_n(1,r,(\emptyset,P)) &\text{ if }& I_0\neq \emptyset \text{ and } r \text{ not integral} \\
p^{-t(r)}e_n(1,r,(\emptyset,P)) &\text{ if }& I_0\neq \emptyset \text{ and } r \text{ integral}.
 \end{array}\right.
 \end{align*}
Then, the following holds:
\begin{prop}[{\cite[Proposition 2.17]{Langer2004DERC}}]
Every $\omega \in \wn\Omega^i_{k[x_1,\dots,x_d]}$ can be uniquely written as a finite sum of the form
\begin{equation}
\omega=\sum_{r,P \in \mathcal{P}_r^{(i)}}\eta_{r,P} \cdot e_n(1,r,P), \quad \eta_{r,P} \in \wn(k)
\end{equation}
where the sum runs over all weights $r$ such that $|\mathrm{supp}(r)| \geq i$ with $p^{n-1}r$ integral and all  partitions $P \in \mathcal{P}_r^{(i)}$. 
\end{prop}
\subsection{An isomorphism after Illusie-Raynaud}
Classically the de Rham-Witt complex is introduced to study crystalline cohomology. If $X$ is a smooth $k$-scheme admitting a smooth lift $X'$ over $\wn(k)$, there is a relation between Hodge-Witt modules and sheaf cohomology of the de Rham complex of $X'/\wn(k)$. This is discussed in \cite{Illusie1983}. The authors prove that for any $n \geq 1$, there are higher  Cartier isomorphisms $$C^{-n} \colon \wn\Omega_X^i \xrightarrow{\sim} h^i(\wn\Omega_X^{\bullet})$$ induced by the Frobenius map 
$F^{n}\colon \w_{2n}\Omega_X^i \rightarrow \wn\Omega_X^i$, such that for $n=1$, $C^{-n}$ agrees with the classical inverse Cartier operator. Furthermore, by comparison of crystalline cohomology and de Rham-Witt cohomology, there is a canonical $\wn(k)$-linear isomorphism $$h^i(\Omega_{X'/\wn}^{\bullet}) \xrightarrow{\sim} h^i(\wn\Omega_X^{\bullet}).$$
Taking its inverse and composing with $C^{-n}$, we get an isomorphism
\begin{equation}
\tilde{F}^{n} \colon \wn\Omega_X^i \rightarrow h^i(\Omega_{X'/\wn}^{\bullet}).
\end{equation}
When $i=0$, $\tilde{F}^n$ can be described explicitly (\cite[p.142, line 8]{Illusie1983}). Let $(x_1,\dots,x_n) \in \wn\ox$ and choose, respectively, some lifts $\tilde{x}_1,\dots, \tilde{x}_n \in \mathcal{O}_{X'}$. Then,
\begin{equation}
\tilde{F}^n \colon (x_1,\dots,x_n) \mapsto \sum_{i=0}^{n-1}p^i\tilde{x}_{i+1}^{p^{n-i}}.
\end{equation}
Let $\Phi \colon \wn(k) \rightarrow \wn(k)$ be the Witt vector Frobenius.
Here we don't introduce the higher Cartier operators, neither the crystalline comparison with de Rham-Witt cohomology. However, we will prove that the latter map above is a well defined $\wn(k)$-$\Phi^n$-semilinear isomorphism by elementary methods. We will only require the existence and injectivity of the classical Cartier operator in order to keep the proof self-contained as much as possible. We will point out that surjectivity of \eqref{tildeF^n} depends on $k$ being perfect (the motivation for such a notation relies on the more general statement of \cite[Prop. 8.4]{BHR12}). Notice that $p^n=0$ in $\mathcal{O}_{X'}$ and the topological spaces underlying $X$ and $X'$ are the same.
\begin{prop}[c.f. {\cite[Prop. 8.4]{BHR12}} and {\cite[Ch. III, sec. (1.5)]{Illusie1983}}]\label{tildeF^n_iso_prop}
Let $X$ be a smooth $k$-scheme together with a $\wn(k)$-smooth scheme $X'$ lifting $X$. Then, the map
\begin{equation}\label{tildeF^n}
\tilde{F}^n  \colon \wn\ox \rightarrow \oo_{X'}, \quad (x_1,\dots,x_n) \mapsto \sum_{i=0}^{n-1}p^i\tilde{x}_{i+1}^{p^{n-i}}
\end{equation}
is a  well defined $\Phi^n$-semilinear injective morphism of sheaves of $\wn(k)$-algebras, inducing an isomorphism onto $\ker(\mathcal{O}_{X'} \xrightarrow{d} \Omega^1_{X'})$.
\end{prop}
\begin{proof}
The map \eqref{tildeF^n} is well defined: if $\tilde{x}_{i+1}$ and $\tilde{\tilde{x}}_{i+1}$ are different lifts of $x_{i+1}$, i.e.
$\tilde{x}_{i+1} \equiv \tilde{\tilde{x}}_{i+1} \pmod p $, then 
$p^i\tilde{x}_{i+1}^{p^{n-i}} \equiv p^i\tilde{\tilde{x}}_{i+1}^{p^{n-i}} \pmod{p^n}$. Moreover, clearly $d \circ \tilde{F}^n \subset p^n\Omega^1_{X'}=0$, thus $\tilde{F}^n\wn\ox \subset \ker(\mathcal{O}_{X'} \xrightarrow{d} \Omega^1_{X'})$.

\eqref{tildeF^n} is a $\Phi^n$-semi-linear morphism of $\wn(k)$-modules:
We observe that as map of sets, it factorizes the $n$-th ghost map
\begin{equation}
\w_{n+1}\mathcal{O}_{X'} \xrightarrow{w_n} \mathcal{O}_{X'}, \quad (y_1,\dots,y_{n+1}) \mapsto \sum_{i=0}^n p^iy_{i+1}^{p^{n-i}} 
\end{equation} 
through the restriction morphism $\w_{n+1}(\mathcal{O}_{X'}) \rightarrow \w_{n+1}(\ox) \rightarrow \wn(\ox) $. Since the latter and $w_n$ are ring homomorphisms, the same is true for \eqref{tildeF^n}. If $a=(a_1,\dots,a_n) \in \wn(k)$, then $([a_1],\dots,[a_n],0)(\tilde{x}_1,\dots,\tilde{x}_n,0) \in \w_{n+1}\mathcal{O}_{X'}$ is a lift of $a\cdot (x_1,\dots,x_n)$, where $[a_i] \in \w_{n}(k)$ is the Teichm\"uller representative of $a_i \in k$. Therefore,
\begin{equation}
\tilde{F}^n(a\cdot (x_1,\dots,x_n))=w_n(([a_1],\dots,[a_n],0))w_n((\tilde{x_1},\dots,\tilde{x}_n,0))=\Phi^n(a)\tilde{F}^n((x_1,\dots,x_n)).
\end{equation} 

\eqref{tildeF^n} is injective: If $\tilde{F}^n((x_1,\dots,x_n))=0$, then reducing the expression $\pmod p$ we get $\tilde{x_1}^{p^n}\equiv 0 \pmod p$. Since $X$ is smooth (so reduced), $x_1=0$. Thus, we get an expression 
$p\tilde{x}_2^{p^{n-1}}+ \cdots +p^{n-1}\tilde{x}_{n}\equiv 0 \pmod{p^n}$. Multiplication by $p \colon \w_{n-1}(k) \rightarrow \wn(k)$ is injective, therefore by smoothness of $X'$, it implies
$$ \tilde{x}_2^{p^{n-1}}+ \cdots +p^{n-2}\tilde{x}_{n}\equiv 0 \pmod{p^{n-1}}.$$
Again reducing modulo $p$, we get $x_2=0$. Repeating the argument, we obtain that all $x_{i+1}$'s are $0$.

\eqref{tildeF^n} is surjective onto $\ker(d)$: let $y \in \mathcal{O}_{X'}$ such that $dy=0$. Suppose $n=2$. Then reducing $\pmod p$, we get $\bar{y} \in Z\Omega^0_X$. Since $Z\Omega^0_X$  is generated by the elements $x^p$ with $x \in \ox$, there exists  $\tilde{a}_1,\tilde{a}_2 \in \mathcal{O}_{X'}$ such that \begin{equation} \label{eq:y=a_1^p+pa_2}
y=\tilde{a}_1^p+p\tilde{a}_2.\end{equation} \footnote{Here we are using that $k$ is perfect: write $\bar{y}=\sum_{i} \alpha_ix_i^p=(\sum_{i}\beta_ix_i)^p$ for some $\alpha_i \in k$, $x_i \in \ox$ and $\beta_i=(\alpha_i)^{1/p} \in k$.} It follows by applying $d$ and dividing by $p$ (since again multiplication by $p$ is injective) to the relation \eqref{eq:y=a_1^p+pa_2}, that $a_1^{p-1}da_1+da_2=0 \in \ox$ (where $a_i=\tilde{a}_i \pmod p$ for $i=1,2$). Now the claim is that  $a_1^{p-1}da_1$ can not be a boundary \footnote{Heuristic motivation: we cannot integrate $x^{p-1}$ in characteristic $p$.}.
More precisely the following claim holds: 

\textit{Claim}: Let $a_1,\dots,a_{n-1} \in \ox$, then $ \sum_{i=1}^{n-1}a_i^{p^{n-i}-1}da_i \in B\Omega_X^1$ if and only if $da_i=0$ for any $i=1,\dots,n-1$. 

If we assume that, then $da_1=da_2=0$ and repeating the argument, we write $a_i=\bar{y}_i^p$ for some $\bar{y}_i \in \ox$ ($i=1,2$), thus  we get an expression (for $i=1,2$) $\tilde{a}_i=y_i^p+p\tilde{b}_i$ for some $y_i,\tilde{b}_i \in \oo_{X'}$ such that $y_i \pmod p=\bar{y}_i$. Replacing in \eqref{eq:y=a_1^p+pa_2}, we get an expression of the form $y=y_1^{p^2}+py_2^p$ for $y_1,y_2 \in \mathcal{O}_{X'} $. Then we can proceed by induction on $n\geq 2$. 

Reducing the expression $dy=0$ $\pmod{p^{n-1}}$, by the inductive hypothesis, there exist $y_1,\dots,y_{n-1} \in \mathcal{O}_{X'}/p^{n-1}$ such that 
$$y \equiv \sum_{i=0}^{n-2}p^iy_{i+1}^{p^{n-1-i}} \pmod{p^{n-1}}. $$
 Lifting every $y_{i+1}$ to some $\tilde{y}_{i+1} \in \mathcal{O}_{X'}$, it implies that there exists a $y_n \in \mathcal{O}_{X'}$ such that 
\begin{equation}
y - \sum_{i=0}^{n-2}p^i\tilde{y}_{i+1}^{p^{n-1-i}}=p^{n-1}y_n.
\end{equation}
Since, $dy=0$, applying $d$ to the equality above, we get that \begin{equation}\label{formula_dy=0}
p^{n-1}(\sum_{i=0}^{n-2}\tilde{y}_{i+1}^{p^{n-1-i}-1}d\tilde{y}_{i+1})=p^{n-1}d(-{y}_n).
\end{equation}
Therefore, 
$$\sum_{i=0}^{n-2}\bar{y}_{i+1}^{p^{n-1-i}-1}d\bar{y}_{i+1} \in B\Omega_X^1,$$
where $\bar{y}_{i+1}$ is the mod $p$ reduction of $y_{i+1}$.   By the claim it follows that $d\bar{y}_1=\dots=d\bar{y}_{n-1}=0$. Thus, by  \eqref{formula_dy=0}, $d\bar{y}_{n}=0 \in \Omega_X^1$ from which the statement follows.

\textit{Proof of the claim:} Of course if $da_i=0$ for any $i$, one implication is trivially verified.
For the other one, we observe that the statement is local, thus it suffices to be verified in the case of $X=\spec(A)$  for a smooth $k$-algebra $A$. 

We notice that the claim for $n=2$ is equivalent to say that the Cartier operator $C^{-1}:=C_{A}^{-1} \colon \Omega_{A^{}}^1 \rightarrow h^1(\Omega_A^{\bullet}) $ is injective, thus it holds by \Cref{prop_cartier_iso}. Moreover, $C^{-1}$ is $F$-semilinear over $A$, where $ F \colon A^{} \rightarrow A, a^{} \mapsto a^p$ is the absolute Frobenius. It follows, that 
\begin{equation}
C^{-1}(a^{p^r-1}da^{})=[(a^{p^{r}-1})^p]C^{-1}(da^{})= [a^{p^{r+1}-1}da]
\end{equation}
for any $a \in A$, $r \geq 0$. Thus, we have the following equality $\pmod {B\Omega_A^1}$:
\begin{equation}
0 \equiv \sum_{i=1}^{n-1}a_i^{p^{n-i}-1}da_i \equiv C^{-1}(da_{n-1}+\sum_{i=1}^{n-2}a_i^{p^{n-i-1}-1}da_i) \pmod{B\Omega_A^1}.
\end{equation}
Since $C^{-1}$ is injective, then follows that 
\begin{equation}
\sum_{i=1}^{n-2}a_i^{p^{n-i-1}-1}da_i^{} \in B\Omega^1_{A^{}}
\end{equation}
and by induction, it follows $da_1=\dots=da_{n-2}=0$. By applying the claim for $n=2$ again, it turns out that $da_{n-1}=0$. \qedhere

\end{proof}
\begin{rem}
    When $k$ is not a perfect $\fp$-algebra, the proof of \Cref{tildeF^n_iso_prop} shows that the map \eqref{tildeF^n} is still an injective morphism with image contained in $\ker(\mathcal{O}_{X'} \xrightarrow{d} \Omega^1_{X'})$ 
\end{rem}
\subsection{Hodge-Witt cohomology of $\mathbb{P}^d_k$ }\label{sec_hodge_witt_proj}
Let $d\geq 1$ be an integer and $k$ a perfect field of characteristic $p$. We are going to compute the Hodge-Witt cohomology of the projective space $P:=\mathbb{P}^d_k$.

The same proof of \Cref{prop_cohomo_of_proj_wn} can be found in \cite[Theorem 6.4]{BHR12} (we only avoid the derived category language used by Berthelot et al.). More classically the result has to be attributed to Gros in \cite[Corollaire 4.2.15]{Gro85} depending on relevant results of Ekedahl (cf. \cite[Corollary 1.1.3]{Eke85}). However, we found the first approach easier and more straightforward than the second one. The cost is that it requires the concept for a smooth and proper variety of being \textit{ordinary}, which is not used in the classical approach. Here, we will address these considerations, as well as a proof of \Cref{prop_cohomo_of_proj_wn}.

The following result for the Hodge cohomology is well known:
\begin{prop}\label{prop_cohomo_of_proj}
Let $i,j \geq0 $ be integers. Then, there is a natural isomorphism of $k$-modules
\begin{equation}
    \h^i(P,\Omega_P^j)=\left\{ \begin{array}{cc}
         k & \text{if } 0 \leq i=j \leq d \\
         0 & \text{otherwise} 
    \end{array}\right.
\end{equation}
\end{prop}
We will see that a similar computation for the Hodge-Witt cohomology holds true:
\begin{prop}\label{prop_cohomo_of_proj_wn}
    Let $i,j \geq0 $ be integers. Then, there is a natural isomorphism of $\wn(k)$-modules
\begin{equation}
    \h^i(P,\wn\Omega_P^j)=\left\{ \begin{array}{cc}
         \wn(k) & \text{if } 0 \leq i=j \leq d \\
         0 & \text{otherwise} 
    \end{array}\right.
\end{equation}
\end{prop}
What we actually prove is that  \Cref{prop_cohomo_of_proj_wn} descends from  \Cref{prop_cohomo_of_proj}.
The key point is a geometric property of the projective space $P$.
\begin{defn}
   A  smooth and proper scheme $X$ over $k$ is called \textit{ordinary} if
   \begin{equation}
       \h^i(X,B\Omega^j_{X/k})=0,\quad (B\Omega^j_{X/k}=d\Omega^{j-1}_{X/k})
   \end{equation} for any $i,j \geq 0$. 
 \end{defn}
 When the crystalline cohomology of $X$ over $\w=\w(k)$ is torsion free, the condition of being ordinary can be formulated in terms of the Hodge and Newton polygons associated to $\h^{j}_{crys}(X/\w)$. For the definition of Newton and Hodge polygons we refer the reader to \cite{Kat79} .
 \begin{prop}[{\cite[Proposition 7.3]{BK86}}]\label{prop_equiv_ordinary}
     Assume that $\h^i_{crys}(X/\w)$ is torsion free for any $i \geq 0$. Then, $X$ is ordinary if and only if the Hodge polygon associated to the numbers $h^{j}:=\mathrm{dim}_{k}\h^{i-j}(X,\Omega_{X/k}^j)$ coincides with the Newton polygon given by the Frobenius action on $\h^i_{crys}(X/\w)$ for any $i \geq 0$.
 \end{prop}
\begin{lemma}
    The projective space $P$ over $k$ is an ordinary scheme.
\end{lemma}
\begin{proof}
Let $\tilde{P}:=\mathbb{P}^d_{\w(k)}$.
The comparison of crystalline cohomology and de Rham cohomology of the projective space over $\w(k)$ ensures that $\h^i_{crys}(P/\w) \simeq \h_{dR}^i(\tilde{P}/\w)$ is torsion free. Thus we can apply  \Cref{prop_equiv_ordinary}.

The cohomology group $\h_{dR}^i(\tilde{P}/\w)$ vanishes if $i$ is odd such that $0 < i < 2d$, or $i>2d$. It follows that $h_{dR}^{i}=1$ for $i$ even and $h_{dR}^{j}=0$ in all other cases, where
 $h^i_{dR}$ denotes the rank of the cohomology $\h_{dR}^i(\tilde{P}/\w)$.  

The computation of Hodge cohomology groups gives
\begin{equation}
\h^{i}(\tilde{P},\Omega^j_{\tilde{P}/\w})= \left\{\begin{array}{cc}
      \h^{i/2}(\tilde{P},\Omega^{i/2}_{\tilde{P}/\w})  & \text{if } 0 \leq i \leq 2d \text{ is even}  \\
      0  &  \text{otherwise}
    \end{array}\right.
\end{equation}
and they are all torsion free modules. Thus, the degeneration of the Hodge spectral sequence (c.f. \cite[Corollaire 2.5]{DI87}) implies that $1=h^i_{dR}=h^{i/2}$. Thus, the Hodge polygon has slope $i/2$ and length $1$ for $0 \leq i \leq 2d$ even, and is trivial for other $i$'s. 
We notice that the absolute Frobenius on ${P}$ acts on the $i$-th crystalline cohomology as the multiplication by $p^{\lambda}=p^{i/2}$ or as the $0$ map.  Hence, the Newton polygon has numbers $\lambda=i/2$ and $\mathrm{mult}(\lambda)=1$, therefore agrees with the Hodge polygon.
\end{proof}
\begin{lemma}[{\cite[Lemma 6.2]{BHR12}}] \label{lemma_gr=}
Let us consider the canonical filtration on the de Rham-Witt complex $\w_{n+1}\Omega_P^{\bullet}$ for any $n\geq 1$. Then, for any $i,j \geq 0$ there are canonical isomorphisms
\begin{equation}\label{iso_gr=}
    \h^i(P,\Omega^j_P) \xrightarrow{\simeq} \h^i(P,\mathrm{gr}^n\w_{n+1}\Omega_P^j)
\end{equation}
induced by the map $\ve^n$ in the exact sequence \eqref{exactsequhodgewitt}. Let $F$ be the absolute Frobenius of $P$. When the cohomology on the left side is equipped with the $k$-module structure induced by the $\oo_P$-module structure of $F_*^{n+1}\Omega_P^j$ and on the right side with the one induced by the $\oo_P$-module structure given in  \Cref{prop_exact_seq_gr}, the isomorphisms \eqref{iso_gr=} are $k$-linear.
\end{lemma}
\begin{proof}
    Consider the following diagram of exact sequences of $\oo_P$-modules, where the vertical one is given by the short exact sequence \eqref{exactsequhodgewitt} :
    \begin{equation}
        \begin{tikzcd}
            &                                             &                                        & 0                                                                  &   \\
0 \arrow[r] & F^{n+1}_*Z_n\Omega_P^{j-1} \arrow[r, "a_n"] & F^{n+1}_*\Omega_P^{j-1} \arrow[r]      & F^{n+1}_*(\Omega_P^{j-1}/Z_n\Omega_P^{j-1}) \arrow[u] \arrow[r]      & 0 \\
            &                                             &                                        & \mathrm{gr}^n\w_{n+1}\Omega_P^j \arrow[u]                          &   \\
0 \arrow[r] & F^{n+1}_*B_n\Omega_P^{j} \arrow[r]          & F^{n+1}_*\Omega_P^{j} \arrow[r, "b_n"] & F^{n+1}_*(\Omega_P^{j}/B_n\Omega_P^{j}) \arrow[u, "\ve^n"] \arrow[r] & 0 \\
            &                                             &                                        & 0 \arrow[u]                                                        &  
\end{tikzcd}
    \end{equation}
    We prove by induction on $n$ that the maps $a_n$ and $b_n$ induce isomorphisms on the $i$-th cohomology for any $i \geq 0$. Indeed, for $n=1$  this is true because $P$ is ordinary.

    Generally, we have the following isomorphisms:
     \begin{equation}
        \h^i(P,B_n\Omega_P^j) \xrightarrow[\simeq]{C_P^{-1}}\h^i(P,B_{n+1}\Omega_P^j/B\Omega_P^j) \xleftarrow[\simeq]{} \h^i(P,B_{n+1}\Omega_P^j)
    \end{equation}
    where $C_P^{-1}$ is the inverse Cartier isomorphism and the other map is  induced by the natural projection. Since $P$ is ordinary, then it is an isomorphism. Thus, the claim for $b_n$ follows.

    For $a_n$ we can analogously consider the following commutative diagram, given by the functoriality of inverse Cartier isomorphism:
    \begin{equation}
      \begin{tikzcd}
{\h^i(P,Z_n\Omega_P^j)} \arrow[r, "C_P^{-1}"] \arrow[d, "\simeq"'] & {\h^i(P,Z_{n+1}\Omega_P^j/B\Omega_P^j)} \arrow[d] & {\h^i(P,Z_{n+1}\Omega_P^j)} \arrow[l, "\simeq"'] \arrow[d] \\
{\h^i(P,\Omega_P^j)} \arrow[r, "C_P^{-1}"]                         & {\h^i(P,Z\Omega_P^j/B\Omega_P^j)}                 & {\h^i(P,Z\Omega_P^j)} \arrow[l, "\simeq"]                           
\end{tikzcd}
    \end{equation}
    As before, the top and bottom line are isomorphisms. The left vertical map is isomorphism by inductive hypothesis, thus the right outer vertical map is isomorphism too, yielding the claim for $a_n$.
\end{proof}
\begin{proof}[Proof of  \Cref{prop_cohomo_of_proj_wn}]
By the exact sequence of $\w_{n+1}\oo_P$-modules \footnote{Here $\mathrm{gr}^n\w_{n+1}\Omega_P^j$ is considered naturally just as $\w_{n+1}\oo_P$-module, unlikely to the $\oo_p$-module structure given in previous Lemma. Analogously, $\w_n\Omega_P^j$ is a  $\w_{n+1}\oo_P$-module via the natural restriction $\w_{n+1}\oo_P \rightarrow \wn\oo_P$.} 
\begin{equation}
    0\rightarrow \mathrm{gr}^n\w_{n+1}\Omega_P^j \rightarrow \w_{n+1}\Omega_P^j\rightarrow \w_{n}\Omega_P^j \rightarrow 0,
\end{equation}
 The \Cref{lemma_gr=} together with  \Cref{prop_cohomo_of_proj}, imply that $\h^i(P,\wn\Omega_P^j)=0$ for $i \neq j$. Furthermore, when $i=j$, by same considerations, we get the following short exact sequence of $\w_{n+1}(k)$-modules on cohomology for any $i$:
\begin{equation}
    0\rightarrow \h^i(P,\mathrm{gr}^n\w_{n+1}\Omega_P^i) \rightarrow \h^i(P,\w_{n+1}\Omega_P^i) \rightarrow \h^i(P,\w_{n}\Omega_P^i) \rightarrow 0
\end{equation}
For any $n\geq 1$, define the map of abelian sheaves $$d\log_n \colon \oo_P^{\times} \rightarrow \wn\Omega_P^1,\quad  s \mapsto \frac{d[s]}{[s]}.$$
Passing to the first cohomology, the map $d\log_n$ induces a map of groups 
\begin{equation}
    cl_{P,n} \colon \zz=\mathrm{Pic}(P) \rightarrow \h^1(P,\wn\Omega_P^1), 
\end{equation}
such that $F(cl_{P,n+1}(1))=R(cl_{P,n+1}(1))=cl_{P,n}(1)$.
Let $cl_{P,n}(1)^i \in \h^i(P,\wn\Omega_P^i)$ be the $i$-th cup product (induced by the structure of dga of $\wn\Omega_P^{\bullet}$) of $cl_{P,n}(1)$. The multiplication by $cl_{P,n}(1)^i$ defines a $\wn(k)$-linear map $\wn(k) \rightarrow \h^i(P,\wn\Omega_P^i)$.

Notice that for $n=1$, $cl_{P,1}(1)=:cl_P(1)$ defines the first Chern class generating the cohomology $\h^1(P,\Omega_P^1)$. Similarly,  $cl_P(1)^i$ generates $\h^i(P,\Omega_P^i)$ for $i \geq 0$ (where $cl_P(1)^0=1$) (c.f. \cite[Tag 0FMG, Lemma 50.11.3]{stacks-project}).

We have the following commutative diagram
\begin{equation}
    \begin{tikzcd}
0 \arrow[r] & {\h^i(P,\mathrm{gr}^n\w_{n+1}\Omega_P^i)} \arrow[r]     & {\h^i(P,\w_{n+1}\Omega_P^i) } \arrow[r]      & {\h^i(P,\w_{n}\Omega_P^i)} \arrow[r] & 0 \\
0 \arrow[r] & F_*^nk \arrow[r, "\ve^n"] \arrow[u, "\ve^n \circ cl_P(1)^i"] & \w_{n+1}(k) \arrow[r] \arrow[u, "cl_{P,n+1}(1)^i"] & \wn(k) \arrow[r] \arrow[u, "cl_{P,n}(1)^i"]          & 0
\end{tikzcd}
\end{equation}
Here $F_*^nk$ denotes the group $k$ with the $\w_{n+1}(k)$-module structure induced by the iterated Frobenius on $\w_{n+1}(k)$, $F^n \colon \w_{n+1}(k) \rightarrow k$. In this way, the sequences on top and bottom lines are exact and the respective diagram is commutative, in the category of $\w_{n+1}(k)$-modules. 
By induction on $n$, the outer maps are isomorphisms, thus it is the map in the middle too.
\end{proof}
\begin{rem}
    The action of the finite group $G=\mathrm{GL}_{d+1}(k)$ on $P$ is by functoriality compatible with the $d\log_n$ map, thus, the identifications in  \Cref{prop_cohomo_of_proj_wn} are $G$-equivariant.  
\end{rem}
\section{Witt line bundles}\label{sec_witt_line}
   We introduce the notion of Witt line bundles in the sense of \cite[Section 3.2]{TH22}. There, the author defines the notion of \textit{ Witt divisorial sheaves}, that are sheaves  associated to a $\mathbb{R}$-divisor $D$ on  integral, normal, Noetherian $\fp$-schemes, and that one of \textit{ Teichm\"uller lift} for line bundles on general $\fp$-schemes.  Our terminology of \textit{Witt line bundles} coincides with the latter one. Moreover, when $X$ is an integral, normal, Noetherian $\fp$-scheme and $D$ is a Cartier divisor on $X$, both definitions agree.
   Let $X$ be an integral $k$-scheme. Assume there is  an integral $\w(k)$-scheme $\tilde{X}$ such that $X_n:=\tilde{X} \times_{\spec(\w(k))} \spec(\wn(k))$ is flat over $\wn(k)$ and $X_1=X$. We denote the total quotient ring of $X$ with $\mathcal{K}_{X}$ (cf. \cite[Ch. II, 6, \textit{Cartier Divisor}, Definition 1]{Har77} ). Notice that  since $X$ is integral, $\mathcal{K}_{X}$ agrees with the function field of $X$. Furthermore, for any scheme $Y$, let $\mathrm{Bun}_n(Y)$ denote the category of vector bundles of rank $n$ over $Y$. 
   
    Let $\mathcal{L}$ be a line bundle on $X$ given by the collection $\{(U_i)_{i \in I}, (f_{ij})_{(i,j) \in I \times I})\}$. By construction, $U_i$ and $f_{ij}$'s have the following properties:
   \begin{itemize}
   \item $(U_i)_{i \in I}$ is an affine Zariski open cover of $X$, 
  \item $f_{ij} \in \Gamma(U_{ij},\mathcal{O}_{X})^{\times}$ where $U_{ij}:=U_i \cap U_j$,
   \item $(f_{ij}f_{jk}f_{ki})_{\mid U_{ijk}}=1$ for any triple $(i,j,k) \in I \times I \times I$, where $U_{ijk}:=U_i\cap U_j\cap U_j$.
      \end{itemize}
To a given line bundle $ \mathcal{L}=\{(U_i),f_{ij}\}$ on $X$ we can associate the line bundle $\wn\mathcal{L}:=\{(U_i),[f_{ij}]\}$ on $\wn(X)$, defining a map

\begin{equation}\label{teich}
\begin{array}{cc}
\quad \mathrm{Bun}_{1}(X) \xrightarrow{\wn} \mathrm{Bun}_{1}(\wn(X)) \\
\mathcal{L} \longmapsto \wn\mathcal{L}.
\end{array}  
\end{equation}
 
 It is induced by the Teichm\"uller map $\mathcal{O}_{X}^{\times} \rightarrow \wn(\mathcal{O}_{X})^{\times}$. In particular, it gives
     a group homomorphism on the respective Picard groups:
   \begin{equation}
  \h^{1}(X,\ox^{\times})\simeq \mathrm{Pic}(X) \longrightarrow \mathrm{Pic}(\wn(X)) \simeq \h^{1}(X,\wn(\ox)^{\times}), \quad \mathrm{class}(\mathcal{L})\mapsto \mathrm{class}(\wn\mathcal{L}).
   \end{equation}
\begin{defn}
 The line bundle on $\wn(X)$ associated to  $\mathcal{L}$  via \eqref{teich} is the \textit{Witt line bundle}  $\wn\mathcal{L}$.
\end{defn} 
More concretely, if $\mathcal{L}=\mathcal{O}_{X}(D)$ is the line bundle associated to a Cartier divisor $D=\{(U_i,f_i)\}$, with $f_i \in \Gamma(U_i,\mathcal{K}_{X}^{\times})$ then $$\wn\mathcal{O}_{X}(D)_{\mid U_i}=\left[\frac{1}{f_i}\right]\wn\mathcal{O}_{U_i} \subset \wn(\mathcal{K}_{X}).$$

\begin{lemma}
For any integral $k$-scheme $X$, the association
\begin{equation*} 
\begin{array}{cc}
\quad\mathrm{Bun}_{1}(X) \xrightarrow{\wn} \mathrm{Bun}_{1}(\wn(X)) \\
 \mathcal{L} \longmapsto \wn\mathcal{L}
\end{array}
\end{equation*}
is functorial, i.e. for any morhpism of line bundles $E \colon \mathcal{L} \rightarrow \mathcal{M}$, there exists a natural morphism of line bundles $\wn(E)\colon \wn\mathcal{L} \rightarrow \wn\mathcal{M} $ such that $\wn(id_{\mathcal{L}})=id_{\wn\mathcal{L}}$ and for any line bundle morphism $S \colon \mathcal{M} \rightarrow \mathcal{N}$, we have $\wn(S \circ E)= \wn(S) \circ \wn(E)$.
\end{lemma}
\begin{proof}
We need to define how $\wn$ operates at level of morphisms. Let
$\mathcal{L}$ and $\mathcal{M}$ be line bundles over $X$, and let $[.]: \mathcal{K}_{X}^{\times} \rightarrow \wn(\mathcal{K}_{X})^{\times}$ be the Teichm\"uller map.
Let $E : \mathcal{L} \rightarrow \mathcal{M}$ be a morphism of line bundles, and $\mathcal{U}=(U_i)_{i \in I}$ be an affine open covering of $X$, such that $\mathcal{L}_{\mid U_i} = f_i\mathcal{O}_{U_i}\simeq \oo_{U_i}$ and
 $\mathcal{M}_{\mid U_i} = g_i\mathcal{O}_{U_i}\simeq \oo_{U_i}$, for some $f_i \in \Gamma(U_i,\mathcal{K}_{X}^{\times})$, $g_i \in \Gamma(U_i,\mathcal{K}_{X}^{\times})$. Moreover, the transition maps are multiplication by some invertible sections, namely $\oo_{ U_{ij}} \xrightarrow{\cdot f_{ij}} \oo_{ U_{ij}}$ and
 $\oo_{ U_{ij}} \xrightarrow{\cdot g_{ij}} \oo_{ U_{ij}}$, such that $f_{ij},g_{ij} \in \Gamma(U_{ij},\ox^{\times})$, $f_if_{ij}=f_j$ and $g_ig_{ij}=g_j$.
For any open $U_i$, define  $$\wn(E)_{U_i}: \wn\mathcal{O}_{U_i}\rightarrow \wn\mathcal{O}_{U_i}, $$ being the unique map of $\wn\mathcal{O}_{U_i}$-modules such that $\wn(E)_{U_i}(1)=[E_{\mid U_i}(1)]$ \footnote{As abuse of notation , here $1 \in f_i\mathcal{O}_{U_i} \simeq \oo_{U_i}$ denotes the unique element corresponding to $1 \in \oo_{U_i}$.}.
 Since $E$ is a morphism of line bundles, we have the following commutative diagram:
 \begin{equation}\label{diag_lemma_func_witt_bundle}
 \begin{tikzcd}
\mathcal{O}_{U_{ij}} \arrow[d, "\cdot f_{ij}"'] \arrow[r, "{E_{U_i}}_{\mid U_{ij}}"] & \mathcal{O}_{U_{ij}} \arrow[d, "\cdot g_{ij}"] \\
\mathcal{O}_{U_{ij}} \arrow[r, "{E_{U_j}}_{\mid U_{ij}}"]                      & \mathcal{O}_{U_{ij}}                    
\end{tikzcd}
 \end{equation}
 It implies that the following diagram commutes:
 \begin{equation}
 \begin{tikzcd}
\wn\mathcal{O}_{U_{ij}} \arrow[d, "\cdot {[f_{ij}]}"'] \arrow[r, "{\wn(E)_{U_i}}_{\mid U_{ij}}"] & \wn\mathcal{O}_{U_{ij}} \arrow[d, "\cdot {[g_{ij}]}"] \\
\wn\mathcal{O}_{U_{ij}} \arrow[r, "{\wn(E)_{U_j}}_{\mid U_{ij}}"]                      & \wn\mathcal{O}_{U_{ij}}.                    
\end{tikzcd}
 \end{equation}
Indeed, we have the following equalities:
\begin{multline*}
([g_{ij}]^{-1} \circ {\wn(E)_{U_j}}_{\mid U_{ij}} \circ [f_{ij}])(1)=[g_{ij}]^{-1}[f_{ij}]{\wn(E)_{U_j}}_{\mid U_{ij}}(1)= \\ 
 [g_{ij}]^{-1}[f_{ij}][E_{U_j}(1)]=[g_{ij}^{-1}f_{ij}E_{U_j}(1)]=[E_{U_i}(1)]=\wn(E)_{U_i}([1])
\end{multline*}
where the second-last equality follows by \eqref{diag_lemma_func_witt_bundle}.
It follows that $(\wn(E)_{U_i})_{i \in I}$ glue together, giving rise to a map of $\wn\mathcal{O}_{X}$-modules 
 $$ \wn(E): \wn\mathcal{L} \rightarrow \wn\mathcal{M}.$$ 
 It is also functorial. Indeed, by construction $\wn(id_{\mathcal{L}})=id_{\wn\mathcal{L}}$.
  Let 
 $$ \mathcal{L} \xrightarrow{E} \mathcal{M} \xrightarrow{S} \mathcal{N}$$
 be morphisms of line bundles on $X$, together with affine open cover  $\mathcal{U}=(U_i)_{i \in I}$, such that $\mathcal{L}_{\mid U_i}$, $\mathcal{M}_{\mid U_i}$ and $\mathcal{N}_{\mid U_i}$ are isomorphic to the trivial $\mathcal{O}_{U_i}$-module. Write $E_{U_i}(1)=s_i$ for some $s_i \in \mathcal{O}_{U_i}$. By construction, we have
 \begin{multline*}
 \wn(S)_{U_i}\circ \wn(E)_{U_i}(1)=\wn(S)_{U_i}([E_{\mid U_i}(1)])=[s_i]\wn(S)([1])=[s_iS_{U_i}(1)]=\\
 [S_{U_i}(s_i)]=[S_{U_i}(E_{U_i}(1))]=\wn(S \circ E)_{U_i}(1)
 \end{multline*}
 By linearity, it implies $\wn(S) \circ \wn(E)=\wn(S \circ E)$, thus the functoriality is verified.
\end{proof}

\begin{lemma}[c.f. {\cite[Prop. 3.7]{TH22}}]\label{wittcartier}
Let $D$ be a Cartier divisor on $X$. Let $\ve, R, F$ be respectively the Verschiebung, restriction and Frobenius (with $F \colon \w_{n+1}(\mathcal{K}_X) \rightarrow \w_{n}(\mathcal{K}_X)$) map of $\wn(\mathcal{K}_X)$ and $F_X$ the absolute Frobenius of $X$. Then, the following holds: 
	\begin{enumerate}
	\item[a)] $F(\w_{n+1}\mathcal{O}_{X}(D))\subset\w_{n}\mathcal{O}_{X}(pD)$ 
	\item[b)] $\ve(\w_{n}\mathcal{O}_{X}(pD))\subset\w_{n+1}\mathcal{O}_{X}(D)$ 
	\item[c)] $R(\w_{n+1}\mathcal{O}_{X}(D))\subset\wn\mathcal{O}_{X}(D)$.
\end{enumerate}
In particular, we have the short exact sequence of $\wn\ox$-modules
\begin{equation}\label{ses_witt_line_bundle}
    0 \rightarrow F_{X*}\w_{n}\mathcal{O}_{X}(pD) \xrightarrow{\ve} \w_{n+1}\mathcal{O}_{X}(D)
    \xrightarrow{R^n} \mathcal{O}_{X}(D) \rightarrow 0.
\end{equation}
Moreover, there exists a unique morphism of $\wn\mathcal{O}_{X}$-modules
    \begin{equation}\label{liftbundle}
        \tilde{F}^n : \wn\mathcal{O}_{X}(D) \rightarrow \mathcal{O}_{X_n}(p^nD)
    \end{equation}
    where $\mathcal{O}_{X_n}$ is considered as a $\wn\ox$-module via the the map \eqref{tildeF^n} and
    such that the composition $$\w_{n+1}\mathcal{O}_{X_n}(D) \xrightarrow{R} \wn\mathcal{O}_{X_n}(D) \rightarrow\w_{n}\mathcal{O}_{X}(D) \xrightarrow{\tilde{F}^n}\mathcal{O}_{X_n}(p^nD)$$ coincides with $F^n$.
\end{lemma} 
\begin{proof}
The statement is local, hence we can assume that $X=U_1$ for an affine $U_1=\spec(A)$, where $A$ is a $k$-algebra, and $D=(U_1,f)$ for $f \in \mathrm{Frac}(A)\backslash \{0\}$. Notice that $F([f])=[f]^p $ and
$\ve(a[f]^p)=[f]\ve(a)$ for $a \in \w_{n}\mathcal{O}_{U_1}$. So a),b) easily follow, and c) holds by definition. Moreover, $R$ is surjective because it is induced by the surjective map $\w_{n+1}(\mathcal{K}_X) \rightarrow \w_{n}(\mathcal{K}_X)$. To prove the \eqref{ses_witt_line_bundle}, we only need to verify that at level of stalks $\ker R^n \subset \mathrm{Im} \ve$. For, let $x \in X$ a point,  \eqref{ses_witt_line_bundle} has the form
\begin{equation}
     0 \rightarrow [1/f^p]\w_{n}\mathcal{O}_{X,x} \xrightarrow{\ve} [1/f]\w_{n+1}\mathcal{O}_{X,x}
    \xrightarrow{R^n} [1/f]\mathcal{O}_{X,x} \rightarrow 0.
\end{equation}
Let $a=(a_1,\dots,a_{n+1}) \in \w_{n+1}\mathcal{O}_{X,x} $ such that $R^n([1/f]a)=0$. Equivalently, 
$a_1/f=0$. Then, 
\begin{multline*}[1/f](a_1,\dots,a_{n+1})=(0,a_2/f^p,\dots,a_{n+1}/f^{p^{n+1}})= \\
\ve((a_2/f^p,\dots,a_{n+1}/f^{p^{n+1}}))=\ve([1/f^p](a_2,\dots,a_{n+1})).
\end{multline*}
For a generic section $f/g \in \mathcal{K}_{X}^{}$, the condition $\tilde{F}^n([g^{-1}]):=\tilde{g}^{-p^{n}}$ extends the map \eqref{tildeF^n} to a map of sheaves of rings $\wn(\mathcal{K}_{X})^{} \rightarrow\mathcal{K}_{X_n}^{}$. We prove that its restriction to $\wn\oo_X(D) \subset \wn(\mathcal{K}_X)^{}$ gives the searched map. Indeed,
for \eqref{liftbundle} we can suppose that $X$ has a model $\tilde{X}$ over $\w(k)$, thus we can assume that $U_1$ is given by base change of a Zariski open $V$ of $\tilde{X}$. We denote by $U_n$ the base change of $V$ with $\spec(\wn(k))$. Take a lift $\tilde{f} \in \mathcal{O}_{U_n}^{}$ of $f \in \mathcal{O}_{U_1}^{}$. Then the map defined in \eqref{tildeF^n} is such that $\tilde{F}^n([f])=\tilde{f}^{p^n}$ and does not depend on the chosen lifting. Since  $\wn\mathcal{O}_{U_1}(D)=[1/f]\wn\mathcal{O}_{U_1}$, $\mathcal{O}_{U_n}(p^nD)=\tilde{f}^{-p^n}\mathcal{O}_{U_n}$ and  $\tilde{F}^n$ is homomorphism of sheaves of rings, it extends to a such $\wn\mathcal{O}_{X}$-module morphism in \eqref{liftbundle}.  
\end{proof}

\begin{lemma}
Let $X$ be a Noetherian, integral, separated $k$-scheme. Let $\mathcal{U}$ be an affine open cover of $X$, and let $\mathcal{L}$ be a line bundle on $X$. Then, there is a canonical isomorphism of groups
\begin{equation}\label{wittcocechco}
   \check{\h}^i(\mathcal{U},\wn\mathcal{L}) \xrightarrow{\simeq} \h^i(X,\wn\mathcal{L}),
\end{equation}
where $ \check{\h}^{i}$ is the $i$-th \v{C}ech cohomology.
\end{lemma}
\begin{proof}
If $D$ is the Cartier divisor associated to $\mathcal{L}=\mathcal{O}_{X}(D)$, by Lemma \eqref{wittcartier} there is a short exact sequence of abelian sheaves: 
\begin{equation}\label{s.e.s_wnL}
0\rightarrow \w_{n-1}\mathcal{O}_{X}(pD)\xrightarrow{\ve} \w_{n}\mathcal{O}_{X}(D) \rightarrow \mathcal{O}_{X}(D) \rightarrow 0.
\end{equation}
For $n=1$, the statement is a particular case of the analogous result for quasi-coherent cohomology. Using the long exact sequence on cohomology \footnote{Since $\wn\mathcal{L}$ is a quasi-coherent module on the separated scheme $\wn(X)$, a long exact sequence for \v{C}ech cohomology associated to the short exact sequence \eqref{s.e.s_wnL} exists as constructed in the proof of \cite[III, Theorem 4.5]{Har77}.} and the 4-Lemma, the isomorphism \eqref{wittcocechco} follows by induction on $n$.
\end{proof}
\begin{lemma}\label{Witt_line_commute_with_pullback}
Let $f \colon X \rightarrow Y$ be a morphism of integral $k$-schemes and $E \colon \mathcal{L} \rightarrow \mathcal{M}$ be a morphism of line bundles over $Y$. Then, for any $n \geq 1$, we have that $f^*\wn(\mathcal{L})=\wn(f^*\mathcal{L})$ (respectively for $\mathcal{M}$) and
\begin{equation}
f^*\wn(E)=\wn(f^*E) \colon f^*\wn\mathcal{L} \rightarrow f^*\wn\mathcal{M}.
\end{equation}
\end{lemma}
\begin{proof}
The statement is local, thus we assume $X=\spec(B)$, $Y=\spec(A)$ where $A,B$ are two integral $k$-algebras. Then $\mathcal{L}, \mathcal{M}$ are $A$-modules of rank one, respectively  $L\simeq s_LA$, $M \simeq s_MA$, for some $s_L,s_M \in \mathrm{Frac}(A)$. Then 
$f^*\mathcal{L}=B \otimes_{A,f} L \simeq h_LB$ and
$f^*\mathcal{M}=B \otimes_{A,f} M \simeq h_MB$ where $h_L=f(s_L)=1\otimes s_L$, $h_M=f(s_M)=1\otimes s_M$.
Then, by definition, $\wn(f^*\mathcal{L})=[h_L]\wn(B)$ (analogously for $\mathcal{M}
$). Thus, $$f^*\wn\mathcal{L}=\wn(B) \otimes_{\wn(A),\wn(f)} [s_L]\wn(A) \simeq [f(s_L)]\wn(B)=[h_L]\wn(B)=\wn(f^*\mathcal{L})$$
and analogously for $\mathcal{M}$.\\
 Moreover, by linearity $\wn(f^*E)=f^*\wn(E)$ if and only if $\wn(f^*E)([h_L])=f^*\wn(E)([h_L])$. By construction we have
\begin{align*}
&\wn(f^*E)([h_L])=[f^*E(h_L)]=1 \otimes [E(s_L)]  \\
&f^*\wn(E)([h_L])=1 \otimes \wn(E)([h_L])=1 \otimes [E(s_L)].
\end{align*}
Thus the equality is verified.
\end{proof}

Assume that $X$ is equipped with an action $\sigma$ of a group $G$ and  the line bundle $\mathcal{L}$ is equipped with a $G$-linearization in the sense of  \Cref{sec_diff_op}.
Explicitly, an isomorphism 
\begin{equation}
\phi \colon \sigma^*\mathcal{L} \xrightarrow{\sim} \mathrm{pr}_1^*\mathcal{L}
\end{equation}
of line bundles over ${G \times X}$ is given such that 
\eqref{diag_linearization} holds.
By  \Cref{Witt_line_commute_with_pullback} the functor $\wn$ commutes with pullback, thus  the isomorphism $\phi$ induces an isomorphism $\wn(\phi)$ satisfying the cocycle condition \eqref{cocycle_cond_wn} for $\wn\mathcal{L}$.
 In other words,
\begin{lemma}
If $\mathcal{L}$ has a $G$-linearization, then it lifts to a $G$-linearization of $\wn\mathcal{L}$ for any $n\geq1$, such that the canonical projection $\wn\mathcal{L}\xrightarrow{R^{n-1}}\mathcal{L}$ is $G$-equivariant.
 $\square$
\end{lemma}

\subsection{Cohomology of Witt line bundles on $\mathbb{P}_k^d$}\label{sec_coho_witt_line_proj}
 Let $d \geq 1$ be an integer. We are going to compute the   cohomology of Witt line bundles for $\mathbb{P}^d_k$.  Denote by $G$  the group of $k$-rational points of $\mathrm{GL}_{d+1,k}$ acting on $\mathbb{P}_k^d$.

 Recall that we have the following $G$-equivariant isomorphisms of $k$ -modules for any integers $i,r$:
\begin{equation}
\h^i\left(\mathbb{P}_k^d, \mathcal{O}_{\mathbb{P}_k^d}(r)\right)=\left\{\begin{array}{ccc}
\left(k\left[z_0, \ldots, z_d\right]\right)_r & \text { if } & i=0 \\
0 & \text { if } & i \neq 0, d \\
\mathrm{Hom}_k(\left(k\left[z_0, \ldots, z_d\right]\right)_{-d-1-r},k) & \text { if } & i=d
\end{array}\right.
\end{equation}
where the index $(-)_r$ denotes the $r$-th homogeneous degree part of the respective graded module. Let $\oo:=\oo_{\mathbb{P}_k^d}$.
 \begin{lemma}\label{positivelinebundle}
 Let $b \geq 0$ be a non negative integer and fix an integer $n \geq 1$. Then, for any $d \geq 1$
 \begin{equation*}
 \h^i(\mathbb{P}^d_k,\wn\mathcal{O}(b))=0 \qquad \forall i > 0.
 \end{equation*}
 For $b >0$ holds that $\h^0(\mathbb{P}^d_k,\wn\mathcal{O}(b))\neq 0$ and $\h^0(\mathbb{P}^d_k,\wn\mathcal{O})=\wn(k)$.
 \end{lemma}
\begin{proof}
The non vanishing assertion for the global section is clear (since it holds for $n=1$). To prove the vanishing of the cohomology groups $\h^i$ for $i>0$, consider the short exact sequences of abelian sheaves 
 $$ 0 \rightarrow \w_{n-1}\mathcal{O}(bp)\xrightarrow{\ve} \wn\mathcal{O}(b) \rightarrow \mathcal{O}(b)\rightarrow 0.$$ Since
  $\h^i(\mathbb{P}^d,\mathcal{O}(b))=0$ whenever $i >0$, using the corresponding long exact sequence and by induction on 
  $n$, we have $\h^i(\mathbb{P}^d,\w_{n}\mathcal{O}(b))=0$ for any $i >0$ . The last equality follows since 
$\h^0(\mathbb{P}^d_k,\wn\mathcal{O})=\wn(\h^0(\mathbb{P}^d_k,\mathcal{O}))=\wn(k)$.
\end{proof}

 \begin{lemma}\label{negativelinebundle}
 Let $a < 0$ be a negative integer and fix an integer $n \geq 1$. Then, for any $d \geq 1$,
 \begin{equation}
 \h^i(\mathbb{P}^d_k,\wn\mathcal{O}(a))=0 \qquad \forall i \not = d.
 \end{equation}
 Moreover, if  $d> -p^{n-1}a-1$, then
 \begin{equation}
 \h^d(\mathbb{P}^d_k,\wn\mathcal{O}(a))=0.
 \end{equation}
 When $d \leq -p^{n-1}a-1$, the cohomology groups $$\ve^i\h^d(\mathbb{P}^d_k,\w_{n-i}\mathcal{O}(p^ia)), \qquad i=0,\dots,n-1$$ form a non trivial descending filtration of $\wn(k)$-sub-modules $\mathrm{Fil}^{\bullet}$ of $\h^d(\mathbb{P}^d_k,\wn\mathcal{O}(a))$, such that $$\mathrm{gr}^i\mathrm{Fil}^{\bullet} \simeq \h^0(\mathbb{P}_k^d,\mathcal{O}(-p^{i}a-d-1))^{\vee}$$
 with the $k$-vector structure induced by $F^i$, where $F$ is the absolute Frobenius of $\mathbb{P}^d_k.$ 
 \end{lemma}
 \begin{proof}
 
  We have short exact sequence of $\wn\oo_{\mathbb{P}^d_k}$-modules
 $$ 0 \rightarrow F_*\w_{n-1}\mathcal{O}(ap)\xrightarrow{\ve} \wn\mathcal{O}(a) \rightarrow \mathcal{O}(a)\rightarrow 0.$$ Using the corresponding long exact sequence, we see, by induction on $n$, that  $$\h^i(\mathbb{P}_k^d,\wn\mathcal{O}(a))=0$$ for $i \not = d$. 
We have to determine $\h^d(\mathbb{P}_k^d,\wn\mathcal{O}(a))$. For any negative $a$ the long exact sequence has the form 
\begin{multline*}
0=\h^{d-1}(\mathbb{P}_k^d,\mathcal{O}(a))\rightarrow \h^d(\mathbb{P}_k^d,\w_{n-1}\mathcal{O}(pa)) \\ \rightarrow \h^d(\mathbb{P}_k^d,\w_{n}\mathcal{O}(a)) \rightarrow \h^d(\mathbb{P}_k^d,\mathcal{O}(a)) \\ \rightarrow\h^{d+1}(\mathbb{P}_k^d,\w_{n-1}\mathcal{O}(pa))\rightarrow \dots,
\end{multline*}
where the last term appearing above is $0$. This means that we have a short exact sequence 
\begin{equation*}
0\rightarrow F_*\h^d(\mathbb{P}_k^d,\w_{n-1}\mathcal{O}(pa))  \rightarrow \h^d(\mathbb{P}_k^d,\w_{n}\mathcal{O}(a)) \rightarrow \h^0(\mathbb{P}_k^d,\mathcal{O}(-a-d-1))^{\vee} \rightarrow 0,
\end{equation*}
since $\h^d(\mathbb{P}_k^d,\mathcal{O}(a))\simeq \h^0(\mathbb{P}_k^d,\mathcal{O}(-a-d-1))^{\vee}$, by Serre duality.
For any $a$ and $n$ fixed, define $V_i^{(n)}(a):=\h^d(\mathbb{P}_k^d,\w_{n-i}\mathcal{O}(p^ia))$ for $i=0,\dots,n-1$ and $V_n^{(n)}(a):=0$. Then, we have a descending chain of $\wn(k)$-modules 
\begin{equation*}
\mathrm{Fil}^{\bullet}:\quad \h^d(\mathbb{P}_k^d,\w_{n}\mathcal{O}(a))=V_0^{(n)}(a) \supset F_*V_1^{(n)}(a) \supset \dots \supset F^{n-1}_*V_{n-1}^{(n)}(a) \simeq F^{n-1}_*\h^d(\mathbb{P}_k^d,\mathcal{O}(p^{n-1}a))  
\end{equation*}
 such that $\mathrm{gr}^i\mathrm{Fil}^{\bullet}$ are the $k$-vector spaces $F^i_*\big(V_i^{(n)}(a)/F_*V_{i+1}^{(n)}(a)\big) \simeq F^i_*\h^0(\mathbb{P}_k^d,\mathcal{O}(-p^ia-d-1))^{\vee}$ for $i=0,\dots,n-1$. This quotient is trivial when $p^ia > -d-1$. If we let 
 $$i_0(a):=\mathrm{max}\{i \mid p^ia > -d-1\},$$ the chain above has the first $i_0(a)$ terms all isomorphic and for $n>i>i_{0}(a)$ the chain is formed by non trivial proper submodules, so it is a non trivial filtration. When $i_0(a) \geq n-1$, that precisely happens when $p^{n-1}a > -d-1$, then the chain is stationary and in this case 
 $$\h^d(\mathbb{P}_k^d,\w_{n}\mathcal{O}(a))\simeq F^{n-1}_*\h^0(\mathbb{P}_k^d,\mathcal{O}(-p^{n-1}a-d-1))^{\vee}=0. \eqno\qedhere$$ 
 \end{proof}
 
\begin{rem}\label{remarkequivariant}
To compute the cohomology of Witt line bundles on the projective space, we essentially used: short exact sequences of $\wn\mathcal{O}$-modules that are $G$-equivariant, the group homomorphism \eqref{teich} and the Serre duality for projective space, which are $G$-equivariant, because are both functorial morphisms. This means that all the modules involved above are also $G$-modules, and the respective maps are morphism of $G$-modules.
\end{rem}
We can summarize the computation above by the following Proposition:
\begin{prop}
Let $\mathcal{L}$ be a line bundle on $\mathbb{P}^d_k$.  We have the following isomorphisms of $\wn(k)[G]$-modules
\begin{equation}
\h^i(\mathbb{P}^d_k,\wn\mathcal{L})=
\begin{cases}
0 & \text{if }   i \not=0,d \\ 
0 & \text{if }  i=d \text{, } \mathcal{L}=\mathcal{O}(a) \text{, } a<0 \text{, } d> -p^{n-1}a-1\\
\wn(k) & \text{if }  i=0\text{, } \mathcal{L}=\mathcal{O} \\
0 & \text{if }  i=0 \text{, } \mathcal{L}=\mathcal{O}(a) \text{, } a<0 \\
\h^0(\mathbb{P}^d_k,\wn\mathcal{O}(b)) \neq 0 & \text{if }  i=0 \text{, }\ \mathcal{L}=\mathcal{O}(b) \text{, } b>0\\
\not = 0 & \text{otherwise,}\\
\end{cases}
\end{equation}
where $G$ acts on $\wn(k)$ trivially.\\ Moreover, when $\mathcal{L}=\mathcal{O}(a)$ with $a<0$ and  $d \leq -p^{n-1}a-1$, then $$0 \not =F_*^{n-1}\h^0(\mathbb{P}_k^d,\mathcal{O}(-p^{n-1}a-d-1))^{\vee} \subset \h^d(\mathbb{P}^d_k,\wn\mathcal{L})$$
is a non trivial proper $\wn(k)[G]$-submodule.
\end{prop}
\begin{proof}
Since $\mathrm{Pic}(\mathbb{P}^d_k)=\mathbb{Z} \cdot \mathcal{O}(1)$, any line bundle $\mathcal{L}$ admits a $G$-linearization, because $\mathcal{O}(1)$ does. Here, we consider the natural linearization on $\mathcal{O}(1)$ induced by the natural action of $\mathrm{GL}_{d+1,k}$ on $\mathbb{P}_k^d$. So the respective cohomology groups are $G$-modules. By the  \Cref{negativelinebundle}, \Cref{positivelinebundle} and \Cref{remarkequivariant}, all the equalities follow.
\end{proof}
\section{Witt differential operators}\label{secwittdiffop}
Our next goal is to define a generalization of $\mathcal{D}_X$, considering the sheaf $\wn\ox$ in place of $\ox$, in such a way it can be seen as a lifted version of the sheaf of differentials over $\mathrm{W}_n$. \\
Let $k$ be a perfect field of characteristic $p>0$, and $d \geq 1$ an integer.
Let $A=k[t_1,\dots,t_d]$, $A_{n+1}=\w_{n+1}(k)[t_1,\dots,t_d]$ and $\Phi: \w_{n+1}(k) \rightarrow \w_{n+1}(k)$ the Frobenius morphism induced by the one on $k$. Consider the ring homomorphism defined by
the $n$-ghost map  $$w_n \colon \w_{n+1}(A_{n+1})\rightarrow A_{n+1}.$$ 
\begin{lemma}\label{Finjective}
Let $k$ be a reduced $\fp$-algebra. Let $B$ be a reduced $k$-algebra and assume for any $n \geq 1$, a projective system of flat $\wn(k)$-algebras $B_n$ is given such that $B_n/p^{n-1} \simeq B_{n-1}$ ($B_1:=B$). If $F$ is a lift over $B_n$ of the Frobenius $\sigma \colon B \rightarrow B$, then $F$ is injective.
\end{lemma}
\begin{proof}
Since $B$ is reduced, $\sigma$ is injective. Suppose there exists  a nonzero $x \in B_n$ such that $F(x)=0$. Since the reduction modulo $p$ of $F$ is $\sigma$, then $x \equiv 0 \pmod{p}$, so there exists a nonzero $x_1 \in B_n$ such that $x=px_1$ and so $pF(x_1)=0$. By flatness of $B_n$ over $\wn(k)$, this means that $F(x_1)=0 \in B_{n-1}$ (because $B_{n-1}=B_n \otimes_{\wn(k)} \w_{n-1}(k) \xrightarrow{p} B_n \otimes_{\wn(k)} \wn(k)=B_n$ is injective), thus reducing again modulo $p$, there is a nonzero $x_2 \in B_{n-1}$ such that $x_1=px_2$. Repeating the same argument, there should exist a nonzero $ x_{n} \in B$ such that $x_{n-1}=px_{n}=0 \in B_1=B$, then $x=px_1=p^2x_2= \dots=p^{n-1}x_{n-1}=0 \in B_n$; a contradiction.
\end{proof}

\begin{prop}
    There exists a unique $\Phi^{n}$-semilinear morphism over $\w_{n+1}(k)$,
    \begin{equation}
        \tilde{w}_n \colon \w_{n+1}(A) \rightarrow A_{n+1}
    \end{equation}
    factorizing $w_n \colon \w_{n+1}(A_{n+1})\rightarrow A_{n+1} $ via the restriction $\w_{n+1}(A_{n+1}) \rightarrow \w_{n+1}(A)$.

    Moreover, the following hold:
    \begin{itemize}
        \item[1)] For any lift $F \colon A_{n+1} \rightarrow A_{n+1}$ of the absolute Frobenius $\sigma \colon A \rightarrow A$, the relation
        \begin{equation}
            \tilde{F}^{n+1}=F \circ \tilde{w}_n
        \end{equation}
      (where $\tilde{F}^{n}$ is defined as in \eqref{tildeF^n} for any $n \geq 1$)  holds and does not depend on the choice of $F$. \item[2)] The map $\tilde{w}_n$ is injective and \textit{maps injectively} $\ve^i\w_{n+1-i}(A)$ into $p^iA_{n+1}$, i.e. it induces an isomorphism between $\ve^i\w_{n+1}(A)$ and $\tilde{w}_{n}\w_{n+1}(A) \cap p^iA_{n+1}$.
    \end{itemize}
\end{prop}
\begin{proof}
    Let $(f_1,\dots,f_{n+1}) \in \w_{n+1}(A)$ and choose some lifts $\tilde{f}_1,\dots,\tilde{f}_{n+1}$  of the respective $f_i$'s over $A_{n+1}$. Define
    \begin{equation}
    \tilde{w}_n \colon \w_{n+1}(A) \rightarrow A_{n+1}, \quad (f_1,\dots,f_{n+1}) \mapsto \tilde{f_1}^{p^{n}} +p\tilde{f_2}^{p^{n-1}} + \dots +p^n\tilde{f}_{n+1}. 
\end{equation}
The map $\tilde{w}_n$ is well defined: if $\tilde{f}_{i+1}$ and $\tilde{\tilde{f}}_{i+1}$ are different lifts of $f_{i+1}$, i.e.
$\tilde{f}_{i+1} \equiv \tilde{\tilde{f}}_{i+1} \pmod p $, then 
$p^i\tilde{f}_{i+1}^{p^{n-i}} \equiv p^i\tilde{\tilde{f}}_{i+1}^{p^{n-i}} \pmod{p^{n+1}}$.

It is clearly the unique map factorizing the ghost map $w_n$ and thus $\tilde{w}_n$ it is a $\Phi^{n}$-semilinear ring homomorphism. (compare with the proof of \Cref{tildeF^n_iso_prop}). 

By construction, for any Frobenius lift $F$, we have $\tilde{F}^{n+1}=F \circ \tilde{w}_n$. Further,
if $F_1,F_2$ are such lifts, then $F_1(\tilde{f_i})$ and $F_2(\tilde{f}_i)$  are both lifts of the same $\sigma(f_i)=f_i^p$. But $\tilde{w}_n$ does not depend on the choice of lifts, so 
\begin{multline*}
F_1(\tilde{w}_n(f_1,\dots,f_{n+1}))=F_1(\tilde{f_1})^{p^{n}} +pF_1(\tilde{f_2})^{p^{n}} + \dots +p^nF_1(\tilde{f}_{n+1})= \\
\tilde{w}_n(\sigma(f_1),\dots,\sigma(f_{n+1}))=F_2(\tilde{w}_n(f_1,\dots,f_{n+1})).
\end{multline*}
Moreover, since $F$ and $\tilde{F}^{n+1}$ are injective (cf. \Cref{Finjective},\Cref{tildeF^n_iso_prop}), the relation in 1) yields the injectivity of $\tilde{w}_n$. 

We are left to prove the isomorphism between $\ve^i\w_{n+1-i}(A)$ and $\tilde{w}_{n}\w_{n+1}(A) \cap p^iA_{n+1}$: Indeed, if $f=(f_1,\dots,f_{n+1}) \in \w_{n+1}(A)$ is such that $\tilde{w}_n(f)=\tilde{f}_1^{p^n}+p\tilde{f}_2^{p^{n-1}}+\dots+p^{n}\tilde{f}_{n+1}=p^iy$ for some $y \in A_{n+1}$ and some lift $\tilde{f}_s$ of $f_s$, then reducing inductively modulo $p^s$, we get $\tilde{f}_s \equiv 0 \pmod{p^{}}$ for any $1\leq s\leq i$. 
\end{proof}

Set
$X_1=X=\spec(A)$, 
$X_{n+1}=\spec(A_{n+1})$ and we denote the ring of differential operators by
    $\mathcal{D}(X_{n+1})$. By  \Cref{sec_diff_op}, recall that $\mathcal{D}(X_{n+1})=\mathcal{D}(\mathbb{A}^d_{\w(k)}) \otimes \w_{n+1}(k)= \mathcal{D}(\w(k)[t_1,\dots,t_d]) \otimes_{\w(k)} \w_{n+1}(k)$
where $\mathcal{D}(\mathbb{A}^d_{\w(k)})$ is the $\w(k)$-algebra generated by the differential operators described on the fraction field by $\partial_i^{[r]}:=\frac{1}{r!}\partial_{t_i}^r$ for any natural number $r\geq 0$.
The natural projection $\w(k) \xrightarrow{R_{n+1}} \w_{n+1}(k)$ induces a natural map $\mathcal{D}(\mathbb{A}^d_{\w(k)})\xrightarrow{R_{n+1}} \mathcal{D}(X_{n+1})$.
The image of $\partial_i^{[r]}$ under $R_{n+1}$ is $\partial_{i,n+1}^{[r]}$.
\begin{prop}\label{liftder}
\begin{itemize}
\item[(i)] Let $B$ be a perfect $\fp$-algebra, and set  $B_{n+1}:=\w_{n+1}(B)$. Let $C=B[t]$, $C_{n+1}=B_{n+1}[t]$. For any  $r \geq 0$, let $\partial_{n+1}^{[r]} \in \mathcal{D}(\mathbb{A}^1_{B_{n+1}}) $ be the $r$-th differential operator with respect to the variable $t$. Then, $\partial_{n+1}^{[r]} \circ \tilde{w}_n$ factorizes through $\tilde{w}_n$, i.e. there is a unique $\w_{n+1}(B)$-linear $r$-th order differential operator, still denoted by $\partial_{n+1}^{[r]} \colon \w_{n+1}(C)\rightarrow \w_{n+1}(C)$
such that the following diagram
\begin{equation}
\begin{tikzcd}
C_{n+1} \arrow[r, "{\partial_{n+1}^{[r]}}"]                & C_{n+1}                     \\
\w_{n+1}(C) \arrow[u, "\tilde{w}_n", hook] \arrow[r, "{\exists ! \partial_{n+1}^{[r]}}", dotted] & \w_{n+1}(C) \arrow[u, "\tilde{w}_n",  hook]
\end{tikzcd}
\end{equation}
commutes. 
\item[(ii)] For any $r \geq 0$, and $1 \leq i \leq d$, any  $\partial_{i,n+1}^{[r]} : A_{n+1} \rightarrow 
A_{n+1} $ factorizes through $\tilde{w}_n$, i.e. there is a unique $\w_{n+1}(k)$-linear $r$-th order differential operator, still denoted by $\partial_{i,n+1}^{[r]} \colon \w_{n+1}(A)\rightarrow \w_{n+1}(A)$
such that the following diagram
\begin{equation}
\begin{tikzcd}
A_{n+1} \arrow[r, "{\partial_{i,n+1}^{[r]}}"]                & A_{n+1}                     \\
\w_{n+1}(A) \arrow[u, "\tilde{w}_n", hook] \arrow[r, "{\exists ! \partial_{i,n+1}^{[r]}}", dotted] & \w_{n+1}(A) \arrow[u, "\tilde{w}_n",  hook]
\end{tikzcd}
\end{equation}
commutes.
\end{itemize}
\end{prop}
The proof of \Cref{liftder} will be given after introducing some notations and some preliminary computations. 
\begin{defn}
    Let $\partial \in \mathcal{D}(X_{n+1})$ such that there exists a $\w_{n+1}(k)$-linear differential operator $\partial_{\mid \w_{n+1}(A)} \colon \w_{n+1}(A) \rightarrow \w_{n+1}(A)$ making the diagram
    \begin{equation}
\begin{tikzcd}
A_{n+1} \arrow[r, "{\partial}"]                & A_{n+1}                     \\
\w_{n+1}(A) \arrow[u, "\tilde{w}_n", hook] \arrow[r, "{ \partial_{\mid \w_{n+1}(A)}}"] & \w_{n+1}(A) \arrow[u, "\tilde{w}_n",  hook]
\end{tikzcd}
\end{equation}
 commute.
 
 We call $\partial_{\mid \w_{n+1}(A)}$  \textit{the restriction to} $\w_{n+1}(A)$ of the respective differential operator $\partial \in \mathcal{D}(X_{n+1})$.
\end{defn}
\begin{rem}
The reason to use the injection given by $\tilde{w}_n$ instead of $\tilde{F}^{n+1}$ is simply that $\partial_{i,n+1} \circ \tilde{F}^{n+1}=0$ for any $\partial_{i,n+1} \in \mathcal{D}(X_{n+1})$.
\end{rem}

To prove  \Cref{liftder} we need some computations. Let us denote by  $v_p(.) \colon \zz \setminus \{0\} \rightarrow \zz$  the $p$-adic valuation and set $v_p(0)=+\infty$.
We will use frequently the following elementary (Legendre's) Formula (cf. \cite[Theorem 2.6.1]{Mol12}):
For any natural number $n$,
\begin{equation}
    v_p(n!)=\sum_{i=1}^{\infty} \left\lfloor \frac{n}{p^i} \right\rfloor.
\end{equation}
\begin{lemma}\label{binomval}
Let $z,w\not=0$ be  natural numbers such that $z \leq w   $. Then,
\begin{equation}
    v_p \left(\binom{w}{z}\right) \geq v_p(w)-v_p(z).
\end{equation}
\end{lemma}
\begin{proof}
The \Cref{binomval} is trivially true when $v_p(w) \leq v_p(z)$. So we can assume $v_p(w)> v_p(z)$. 
Any natural number $w$ can be written uniquely in the form $w=jp^{v_p(w)}$ where $j$ is a natural number such that $(j,p)=1$. Set $v_p(w)=:n$.  By assumption, $n>v_p(z)=:m$. Let $s=z/p^m$, then $(s,p)=1$, thus
\begin{equation}
\left\lfloor\frac{jp^{n-m}-s}{p^l}\right\rfloor=jp^{n-m-l}-1-\left\lfloor \frac{s}{p^l}\right\rfloor, \quad \forall n-m \geq l>0.
\end{equation}
For $l>n-m$ the inequality 
\begin{equation}
\left\lfloor\frac{jp^{n-m}-s}{p^l}\right\rfloor \leq \lfloor jp^{n-m-l} \rfloor - \left\lfloor \frac{s}{p^l}\right\rfloor
\end{equation}
holds.
We can write $w-z=p^m(jp^{n-m}-s)$. Then, by Legendre's Formula
\begin{multline}
v_p((w-z)!)=\frac{p^m-1}{p-1}(jp^{n-m}-s)+\sum_{l=1}^{\infty}\left\lfloor\frac{jp^{n-m}-s}{p^l}\right\rfloor=\\
\frac{p^m-1}{p-1}(jp^{n-m}-s)+\sum_{l=1}^{n-m}\left\lfloor\frac{jp^{n-m}-s}{p^l}\right\rfloor+\sum_{l>n-m}\left\lfloor\frac{jp^{n-m}-s}{p^l}\right\rfloor.
\end{multline}
By the equality and inequality above, it follows that
\begin{multline}
v_p((w-z)!) \leq \frac{p^m-1}{p-1}(jp^{n-m}-s)+\sum_{l=1}^{\infty}\left\lfloor\frac{jp^{n-m}}{p^l}\right\rfloor-\sum_{l=1}^{\infty}\left\lfloor\frac{s}{p^l}\right\rfloor -(n-m)\\=
\frac{p^m-1}{p-1}(jp^{n-m}-s)+v_p(j!)+j\frac{p^{n-m}-1}{p-1}-v_p(s!)-(n-m).
\end{multline} 
It implies that 
\begin{multline*}
v_p \left(\binom{w}{z}\right) \geq j\frac{p^n-1}{p-1}+v_p(j!)-v_p(z!)-jp^{n-m}\frac{p^m-1}{p-1}\\ + \frac{p^m-1}{p-1}s-v_p(j!)-j\frac{p^{n-m}-1}{p-1}+v_p(s!)+(n-m)
\\=v_p(s!)-v_p((sp^m)!)+s\frac{p^m-1}{p-1}+(n-m)=n-m. \qedhere
\end{multline*}
 
 \end{proof}
 \begin{proof}[ Proof of  \Cref{liftder}]
  In both cases (ii) and (i), the uniqueness of the map and the fact that it is a $\w_{n+1}(k)$-linear, respectively $\w_{n+1}(B)$-linear  $r$-th order differential operator is clear by construction (It is equal to $\tilde{w}_n^{-1} \circ \partial_{i,n+1}^{[r]} \circ \tilde{w}_n$, where $\tilde{w}_n^{-1}$ is the inverse of $\tilde{w}_n $ on the image. The linearity follows since $\tilde{w}_n$ is $\Phi^n$-semilinear, resp. $\Phi^n_{B}$-semilinear where $\Phi_B \colon \w_{n+1}(B) \rightarrow \w_{n+1}(B)$ is the Frobenius morphism induced by that one on $B$). We only need to prove that $\partial_{i,n+1}^{[r]}$ preserves $\w_{n+1}(C)$, resp. $\w_{n+1}(A)$.
  
We first notice that we can reduce to the case (i). Indeed, we have that $A=B'[t_i]$, and $A_{n+1}=B'_{n+1}[t_i]$ where here $B'$ and  $B'_{n+1}$ are respectively the polynomial algebra over $k$ and over $\w_{n+1}(k)$ in the variables $t_j$ for every $j \neq i$. 
 The differential operator $\partial_{i}^{[r]}$ is  $B'$-linear, of the form $\partial^{[r]}_{i,\zz} \otimes B'$ where $\partial^{[r]}_{i,\zz} \in \mathcal{D}(\zz[t_i])$.  Thus it is uniquely determined by the integral image of $t_i=:t$. Therefore, $\partial_{i,n+1}^{[r]}$ is a  $B'_{n+1}$-linear differential operator compatible with $\partial_{i}^{[r]}$. Let $B:=B'_{perf}$ be the perfect closure of $B'$ (cf. \cite[Section 4]{BG18}). 
 Then,  $\partial_{i,\zz}^{[r]} \otimes_{\zz} \w_{n+1}(B)=:\partial_{i,n+1}^{[r],perf} \colon \w_{n+1}(B)[t] \rightarrow \w_{n+1}(B)[t]$ defines a $\w_{n+1}(B)$-linear differential operator.
 We claim there is a unique $\w_{n+1}(B')$-linear differential operator $\partial_{i,n+1}^{[r]} \colon \w_{n+1}(B'[t]) \rightarrow  \w_{n+1}(B'[t])$ such that we have the following commutative diagram
 \begin{equation}
\begin{tikzcd}
B_{n+1}[t] \arrow[r, "{\partial_{i,n+1}^{[r],perf}}"]                & B_{n+1}[t]                     \\
\w_{n+1}(B_{}[t]) \arrow[u, "\tilde{w}_n", hook] \arrow[r, "{ \partial_{i,n+1}^{[r],perf}}"] & \w_{n+1}(B_{}[t]) \arrow[u, "\tilde{w}_n",  hook] \\
\w_{n+1}(B'[t]) \arrow[u, hook] \arrow[r,  "{\exists ! \partial_{i,n+1}^{[r]} }",dotted] & \w_{n+1}(B'[t]) \arrow[u,hook] 
\end{tikzcd}
\end{equation}
  Indeed, by (i) applied to $\mathbb{A}^1_{B_{}}$, the top square above exists. Denote by $\Phi_{B_{}} \colon \w_{n+1}(B_{}) \rightarrow \w_{n+1}(B_{})$ the Frobenius morphism induced by that one on $B_{}$. Any $f \in \w_{n+1}(B'[t])$ is a linear combination of elements $\ve^{l}(b_{}[t]^s)$ with $b_{} \in \w_{n+1-l}(B')$. However, $$\tilde{w}_n(\ve^{l}(b_{}[t]^s))=\Phi_{B_{}}^{n-l}(b_{})p^lt^{sp^{n-l}}=\Phi_{B'_{}}^{n-l}(b_{})\tilde{w}_n(\ve^{l}([t]^s)).$$ Thus, 
  \begin{multline*}
      (\tilde{w}_n^{-1} \circ \partial_{i,n+1}^{[r]}\circ \tilde{w}_n^{})(\ve^{l}(b_{}[t]^s))=\tilde{w}_n^{-1}\big(\Phi_{B'_{}}^{n-l}(b_{})\partial_{i,n+1}^{[r]}(\tilde{w}_n\ve^{l}([t]^s)) \big)=
      \\
      =\Phi_B^{-l}(b)(\tilde{w}_n^{-1} \circ \partial_{i,n+1}^{[r]}\circ \tilde{w}_n^{})(\ve^{l}([t]^s)).
  \end{multline*}
  We need to determine the image of the elements of the form $\ve^{l}([t]^s) \in \w_{n+1}(B[t])$ under $\partial_{i,n+1}^{[r]}$.  By construction $\partial_{i,n+1}^{[r],perf}(\ve^{l}([t]^s)) \in \w_{n+1}(B'[t])$. Indeed, since
  $$\partial_{i,n+1}^{[r],perf}(\tilde{w}_n(\ve^{l}([t]^s)))=\partial_{i,n+1}^{[r],perf}(p^lt^{sp^{n-l}})=p^l\binom{sp^{n-l}}{r}t^{sp^{n-l}-r}$$
  if $r \leq sp^{n-l}$, and $0$ otherwise,
   by \Cref{binomval}, we can write
      \begin{equation}
          p^l\binom{sp^{n-l}}{r}=  \left\{\begin{array}{lr}
           p^{n-v_p(r)}u,     & \text{if } v_p(r) \leq n-l \\
             p^lw,  & \text{if } v_p(r) > n-l, 
          \end{array}\right.
      \end{equation}
   for some $u,w \in \zz_{(p)}^{}$. It follows that
    \begin{equation}
       p^l\binom{sp^{n-l}}{r}t^{sp^{n-l}-r}=\left\{\begin{array}{lr}
          p^{n-v_p(r)}ut^{p^{v_p(r)}u'} ,& \text{if } v_p(r) \leq n-l  \\
        p^{l}wt^{p^{n-l}w'},   & \text{if } v_p(r) > n-l
      \end{array}\right.
        \end{equation}
        for some $u',w' \in \zz_{(p)}$. Hence,
        \begin{equation}
        \tilde{w}_n^{-1}\left(p^l\binom{sp^{n-l}}{r}t^{sp^{n-l}-r} \right)=  \left\{\begin{array}{lr}
                u\ve^{n-v_p(r)}([t]^{u'}), & \text{if } v_p(r) \leq n-l  \\
                w\ve^{l}([t]^{w'}), & \text{if } v_p(r) > n-l 
            \end{array} 
            \right.
        \end{equation}
    is an element of $\w_{n+1}(B'[t])$.    
  Since, $\Phi_B^{-l}(b)\ve^l([t]^{w'})=\ve^l(b[t]^{w'})$ lies in $\w_{n+1}(B'[t])$ and for $n-v_p(r)-l \geq 0$, also $\Phi_B^{-l}(b)\ve^{n-v_p(r)}([t]^{u'})=\ve^{n-v_p(r)}(\Phi_{B}^{n-v_p(r)-l}(b)[t]^{u'})$ lies in $\w_{n+1}(B'[t])$, then $\partial_{i,n+1}^{[r]}(\ve^l(b[t]^s)) \in \w_{n+1}(B'[t])$  and they determine the searched map. 
  
   Let us prove (i).\\
   Let $F \colon C_{n+1} \rightarrow C_{n+1}$ be a lift of the absolute Frobenius of $C$. Since $C$ is a reduced algebra over the reduced $\fp$-algebra $B$, $F$ is injective by  \Cref{Finjective}. Moreover, since $B$ is perfect,
   by  \Cref{tildeF^n_iso_prop}, $\tilde{F}^{n+1}$ is an isomorphism between $\w_{n+1}(C)$ and $\ker\big(C_{n+1} \xrightarrow{d} \Omega^1_{C_{n+1}/{B_{n+1}}}\big)$.  To prove  \Cref{liftder} we need to check that for any $f \in\w_{n+1}(C)$, then 
     \begin{equation}\label{rel_prop_liftder}
         \partial_{n+1}^{[r]}(\tilde{w}_{n}(f)) \in \tilde{w}_n(\w_{n+1}(C)).
     \end{equation}
     It suffices to verify that 
     \begin{equation}\label{eq:dF(partial)}
         dF(\partial_{n+1}^{[r]}(\tilde{w}_{n}(f)))=0 \quad \forall r \geq 0.
     \end{equation}
     Indeed, 
      by \Cref{tildeF^n_iso_prop}, \eqref{eq:dF(partial)} implies that $F(\partial_{n+1}^{[r]}(\tilde{w}_{n}(f))) \in \tilde{F}^{n+1}\w_{n+1}(C)$.  Since $\tilde{F}^{n+1}=F \circ \tilde{w}_n$,  by injectivity of $F$, the \eqref{rel_prop_liftder} follows.
    
     By linearity, it suffices to consider $f$ being of the form $\ve^i(a_{}[t]^j)$, and $a_{} \in \w_{n+1-i}(B)$ with $i<n+1$. The map $\tilde{w}_n \colon \w_{n+1}(C)=\w_{n+1}(B[t]) \rightarrow \w_{n+1}(B)[t]=C_{n+1}$ maps $[t] \mapsto t^{p^n}$. With those assumptions, we have, by definition, that
     \begin{equation}
       \tilde{w}_{n}(\ve^i(a_{}[t]^j))=\Phi_B^{n-i}(a_{})p^it^{jp^{n-i}}.  
     \end{equation}
     Then,
     \begin{equation*}
        F\Big(\partial^{[r]}\big( \Phi_B^{n-i}(a_{})p^it^{jp^{n-i}}\big)\Big)=
        \left\{
        \begin{array}{lc}
        p^i\Phi_B^{n+1-i}(a_{})\binom{ jp^{n-i}}{r}t^{jp^{n-i+1}-pr}, & \text{ if } r \leq jp^{n-i}\\
        0, & \text{ otherwise. }
        \end{array}\right.
     \end{equation*}
     and,
     \begin{multline*}
         d\left(p^i\Phi_B^{n+1-i}(a_{})\binom{jp^{n-i}}{r}t^{jp^{n+1-i}-pr}\right)= \\
                 = \left\{
        \begin{array}{lc}
        p^i\Phi_B^{n+1-i}(a_{})\binom{ jp^{n-i}}{r}(jp^{n+1-i}-pr)t^{jp^{n+1-i}-pr-1}dt, & \text{ if } r \leq jp^{n-i} \\
        0, & \text{ otherwise. }
        \end{array}\right.
     \end{multline*}
     Applying the result of  \Cref{binomval}, we get that 
     \begin{align*}
         v_p\left(p^i\binom{ jp^{n-i}}{r}(jp^{n+1-i}-pr)\right) &\geq i+n-i-v_p(r)+1+v_p(jp^{n-i}-r) \\
         &\geq n+1-v_p(r)+v_p(r) = n+1,
     \end{align*}
     if $v_p(r) \leq n-i$, while
     \begin{equation*}
         v_p\left(p^i\binom{ jp^{n-i}}{r}(jp^{n+1-i}-pr)\right)\geq i+1+v_p(jp^{n-i}-r) 
         \geq i+1+n-i = n+1,
     \end{equation*}
     if $v_p(r)> n-i$,
     implying the statement.
 \end{proof}
 \subsection{Properties of Witt differential operators}
 Let $X_{n+1}=\mathbb{A}^d_{\w_{n+1}(k)}=\spec(A_{n+1})$ as before. The composition of differential operators $\partial_{i,n+1}^{[r]} \circ \partial_{j,n+1}^{[s]} \in \mathcal{D}(X_{n+1})$ satisfies the following relations:
 \begin{align}
     \partial_{i,n+1}^{[r]} \circ \partial_{i,n+1}^{[s]}=\binom{r+s}{r}\partial_{i,n+1}^{[r+s]} \\
      \partial_{i,n+1}^{[r]} \circ \partial_{j,n+1}^{[s]}= \partial_{j,n+1}^{[s]} \circ \partial_{i,n+1}^{[r]} \quad \text{if } i \neq j.
 \end{align}
 Since the restriction to $\w_{n+1}(A)$, namely $\partial_{i,n+1 \mid \w_{n+1}(A)}^{[r]}$ agrees with $\tilde{w}_n^{-1}\circ \partial_{i,n+1}^{[r]} \circ \tilde{w}_n$, the relations above hold for the restriction too.
 \begin{rem}
 The composition could be zero depending on the $p$-adic valuation of $\binom{r+s}{r}$. For example, in the case of the $A=k[t]$ and $r=p, s=(p-1)p$ , by \Cref{binomval} follows that $(\partial_{t,n+1}^{[r]}\circ \partial_{t,n+1}^{[s]})^{\circ n+1}=0.$ 
  \end{rem}
  For any vector $\mathbf{r}=(r_{1},\dots,r_d) \in \mathbb{N}^{d}$, denote by $$\partial_{n+1}^{[\mathbf{r}]}:=\prod_{j=1}^d\partial_{j,n+1}^{[r_j]}=\partial_{1,n+1}^{[r_1]} \circ \dots \circ \partial_{d,n+1}^{[r_d]}.$$ Let $J:=\mathrm{supp}(\mathbf{r}) \subset \{1,\dots,d\}$ be the subset of indexes $j$ where $v_p(r_j) \neq 0$. Then, set $$v_p(\mathbf{r}):=\min_{j \in J}\{v_{p}(r_j)\}.$$

 \begin{lemma}\label{claimhighdiff}
 Suppose that $\delta \in \mathcal{D}(X_{n+1})$ is some differential operator  of order $q$. For any $x \in A_{n+1}$ the following relation holds: For any $m \geq 1$,
 \begin{equation}\label{eq:claimhighdiff}
 \delta(x^m) \equiv 0 \pmod{p^{v_p(m)-v_p(q)}},
 \end{equation}
 where for $m$ such that $0\leq v_p(m) \leq v_p(q)$, the \eqref{eq:claimhighdiff} means $\delta(x^m) \in A_{n+1}$.
 \end{lemma}
 \begin{proof}
 Let us write $m=jp^{v_p(m)}$ uniquely such that $(j,p)=1$. Then, by \cite[Proposition 9]{Nakai70}
 \begin{multline}
 \delta(x^m)=\binom{p^{v_p(m)}}{q}(x^{j})^{p^{v_p(m)}-q}\delta(x^{jq})+
 \\
 \sum_{s=1}^{q-1}(-1)^s\binom{p^{v_p(m)}}{q-s}\binom{p^{v_p(m)}-q+s-1}{s}(x^{j})^{p^{v_p(m)}-q+s}\delta(x^{j(q-s)}).
 \end{multline}
 Notice that when $q-s=p^{v_p(m)}$, the respective term in the sum is $0$. We can then assume that $q-s<p^{v_p(m)}$, thus $v_p(q-s)<v_p(m).$ We proceed by induction on $v_p(m)$. There is nothing to prove when $v_p(m)=0$. Assume that the \eqref{eq:claimhighdiff} is true for every $m'$ such that $v_p(m')<v:=v_p(m)$. Since $v_p(q-s)<v$, by induction $v_p(\delta(x^{j(q-s)}))\geq v_p(q-s)-v_p(q)$. Moreover, by \Cref{binomval}, $$v_p\left(\binom{p^{v_p(m)}}{q-s}\right)\geq v_p(m)-v_p(q-s).$$ Thanks to the sum above it follows that $$v_p(\delta(x^m))\geq \mathrm{min}_{s=1,\dots,q-1}\{v_p(m)-v_p(q),v_p(m)-v_p(q-s)+v_p(q-s)-v_p(q)\}=v_p(m)-v_p(q).$$ 
 Thus the statement is verified. 
 \end{proof}
 Any $\partial_{j,1}^{[r]} \in \mathcal{D}(X_1)$ is a differential operator satisfying \eqref{general_H-S_formula}. In particular, since $X_{n+1}$ is smooth over $\w_{n+1}(k)$, there is at least a differential operator $\tilde{\partial}$ on $\mathcal{D}(X_{n+1})$ satisfying the relation \eqref{general_H-S_formula} and lifting $\partial_{j,1}^{[r]}$ (cf. \Cref{lift_derivation_over_smooth} ii) ).  The restriction of $\partial$ to $\w_{n+1}(A)$  may depend a priori, on the chosen lift. In the next Proposition, we will see that this case indeed does not occur. 

 \begin{defn} Let $a^{(p^{n-i})}_{(c_1,\dots,c_m)}$ be the coefficient $\pmod{p^{n+1}}$ of the monomial $z_1^{c_1}\cdots z_{m}^{c_m}$ in the expansion of $(\sum_{j=1}^m z_j)^{p^{n-i}} \in \zz[z_1,\dots,z_m]$.
 \end{defn}
 \begin{lemma}\label{lem:before_indep_of_lift}
 Let $0\leq j \leq n$, and $m \geq 1$ be integers. Let $a_1,\dots, a_m, b_1,\dots,b_m \in \zz[x_1,\dots,x_d]$ such that $a_i \equiv b_i \pmod p$ for any $i=1,\dots,m$. Moreover, let $c_1,\dots, c_m \in \zz$ such that $c_1+\dots +c_m=p^{n-j}$.
 
 Then, $$p^ja^{(p^{n-j})}_{(c_1,\dots,c_m)}a_1^{c_1}\cdots a_m^{c_m} \equiv p^ja^{(p^{n-j})}_{(c_1,\dots,c_m)}b_1^{c_1}\cdots b_m^{c_m} \pmod{p^{n+1}}. $$
 \end{lemma}
 \begin{proof}
 Let $\epsilon= \mathrm{min}_{i=1,\dots,m}\{v_p(c_i)\}$. Without loss of generality, we can assume that $v_p(c_1)=\epsilon$. Observe \footnote{As abuse of notation $a^{(p^{n-j})}_{(c_1,\dots,c_m)}$ denotes the integer number, rather than its reduction modulo $p^{n+1}$.} that $\binom{p^{n-j}}{c_1} \mid a^{(p^{n-j})}_{(c_1,\dots,c_m)}$. By \Cref{binomval}, $v_p(a^{(p^{n-j})}_{(c_1,\dots,c_m)}) \geq v_p(\binom{p^{n-j}}{c_1}) \geq n-j-\epsilon$.
 
 Since $a_i \equiv b_i \pmod p$, it yields $a_i^{c_i} \equiv b_i^{c_i} \pmod {p^{\epsilon+1}}$, i.e. for any $i$, there is a $h_i \in \zz[x_1,\dots,x_d] $ such that $a_i^{c_i}=b_i^{c_i}+p^{\epsilon+1}h_i$. Hence, for some $H \in \zz[x_1,\dots,x_d]$, we have 
 \begin{align*}
 p^ja^{(p^{n-j})}_{(c_1,\dots,c_m)}a_1^{c_1}\cdots a_m^{c_m}=& \\
 p^ja^{(p^{n-j})}_{(c_1,\dots,c_m)} (b_1^{c_1}+p^{\epsilon+1}h_1) \cdots (b_m^{c_m}+p^{\epsilon+1}h_m)=& \\
 p^ja^{(p^{n-j})}_{(c_1,\dots,c_m)}b_1^{c_1}\cdots b_m^{c_m}
+p^ja^{(p^{n-j})}_{(c_1,\dots,c_m)}p^{\epsilon+1}H.
\end{align*} 
Moreover, $v_p(p^ja^{(p^{n-j})}_{(c_1,\dots,c_m)}p^{\epsilon+1})\geq j+n-j-\epsilon+\epsilon+1=n+1$, from which the assertion follows. 
 \end{proof}
 \begin{prop}\label{uniquenessliftHS}
      If $\partial_{n+1}^{[\mathbf{r}]} \colon A_{n+1} \rightarrow A_{n+1}$ is any $\w_{n+1}(k)$-linear differential operator lifting $\partial_{1}^{[\mathbf{r}]}$, satisfying the relation \eqref{general_H-S_formula} and such that the restriction to $\w_{n+1}(A)$ exists, then 
 ${\partial_{n+1}^{[\mathbf{r}]}}_{\mid \w_{n+1}(A)}$ does not depend on the choice of lift $\partial_{n+1}^{[\mathbf{r}]}$.
 \end{prop}
\begin{proof}

 
 Since it can be verified on each $j$-th coordinate separately, we can reduce to the case  $d=1$. 
 
      Let $0 \leq i  \leq n$ and $f=\ve^i([f_i]) \in \w_{n+1}(A)$, $\tilde{w}_n(f)=p^i\tilde{f}_i^{p^{n-i}}$ for some lift $\tilde{f}_i \in A_{n+1}$. Write $\partial_{n+1}^{[r]}:=\partial_{j,n+1}^{[r]}$. Then, by the  \eqref{general_H-S_formula}, it follows that we can write
 \begin{equation}\label{eq:hasse_diff_op_computation}
     \partial_{n+1}^{[r]}(p^i\tilde{f}_i^{p^{n-i}})= p^i\sum_{\substack{t_1 < \dots < t_m \\ c_1t_1+\dots + c_mt_m=r \\ c_1+ \dots +c_m=p^{n-i} }} a^{(p^{n-i})}_{(c_1,\dots,c_m)}\partial_{n+1}^{[t_1]}(\tilde{f})^{c_1}\dots \partial_{n+1}^{[t_{m}]}(\tilde{f})^{c_m}
 \end{equation}

  Because of $\partial_{n+1}^{[t_{j}]}(\tilde{f}) \equiv \partial_{1}^{[t_{j}]}({f}) \pmod{p}$, then
  for any lift 
 $\widetilde{\partial_{1}^{[t_{j}]}({f})} \in A_{n+1}$ of  $\partial_{1}^{[t_{j}]}({f})$, by \Cref{lem:before_indep_of_lift}
 $$ p^i a^{(p^{n-i})}_{(c_1,\dots,c_m)}\partial_{n+1}^{[t_1]}(\tilde{f})^{c_1}\dots \partial_{n+1}^{[t_{m}]}(\tilde{f})^{c_m} \equiv p^ia^{(p^{n-i})}_{(c_1,\dots,c_m)}\widetilde{\partial_{1}^{[t_1]}({f})}^{c_1}\dots \widetilde{\partial_{1}^{[t_{m}]}({f})}^{c_m} \pmod{p^{n+1}}.$$
 In particular, the \eqref{eq:hasse_diff_op_computation}, does not depend on choices of such lifts. 
\end{proof}
 
 \begin{corol}\label{uniquenesslift}
 \begin{itemize}
     \item[a)] 
 Any $\partial \in \mathcal{D}(X_{1})$ admits some lifting in $\mathcal{D}(X_{n+1})$ of the form \begin{equation}\label{sumofhs}
 \tilde{\partial}=\sum_{\mathbf{r}}c_{\mathbf{r}} \partial_{n+1}^{[\mathbf{r}]},
 \end{equation} where $c_{\mathbf{r}} \in A_{n+1}$ and $\partial_{n+1}^{[\mathbf{r}]}$ is obtained as product of lifts of $\partial_{j,1}^{[r_j]}$ for any $1 \leq j \leq d$.
 
 \item[b)] Furthermore, if the restriction $\tilde{\partial}_{\mid \w_{n+1}(A)}$ to $\w_{n+1}(A)$ exists, it does not depend on the lifts of $\partial_{j,1}^{[r_j]}$, for any $j$.
 
 \item[c)] Let $\tilde{\partial}_1, \tilde{\partial}_2 \in \mathcal{D}(X_{n+1})$ be liftings of $\partial$. Then, if the restriction of
 $\tilde{\partial}_1-\tilde{\partial}_2=\sum_{\mathbf{r}} a_{\mathbf{r}}\partial_{n+1}^{[\mathbf{r}]}$  to $\w_{n+1}(A)$ exists, it is $0$ if and only if $a_{\mathbf{r}} \in p^{v_{p}(\mathbf{r})+1}A_{n+1}$ for every $\mathbf{r}$.
 \end{itemize}
 \end{corol}
 
 \begin{proof}
 a) The existence of such liftings follows by  \Cref{lift_derivation_over_smooth} ii), because we can write $\partial=\sum_{\mathbf{r}} b_{\mathbf{r}} \partial_{1}^{[\mathbf{r}]}$ for some $b_{\mathbf{r}} \in A$, finitely many not zero. 

b) It follows by \Cref{uniquenessliftHS}.
 

 c) Let $f=(f_1,\dots,f_{n+1}) \in \w_{n+1}(A)$, then by \Cref{claimhighdiff} it follows that for any differential operator $\delta$ of order $q$,  $p^{i}\delta(\tilde{f}_{i+1}^{p^{n-i}}) \in p^{n-v_p(q)}A_{n+1}$.
 In particular, $\partial_{j,n+1}^{[r_j]}(\tilde{w}_n\w_{n+1}(A)) \subset p^{n-v_p(r_j)}A_{n+1}$. Thus, if $a_{\mathbf{r}} \in p^{v_p(\mathbf{r})+1}A_{n+1}$, since $$\prod_{j=1}^d \partial_{j,n+1}^{[r_j]}(\tilde{w}_n(\w_{n+1}(A))) \subset p^{n-v_p(r_j)}A_{n+1}$$ for any $j$, in particular $$a_{\mathbf{r}}\partial^{[\mathbf{r}]}_{n+1 \mid\w_{n+1}(A)}=0.$$ On the other hand, suppose to have a differential operator $\sum_{\mathbf{r}} a_{\mathbf{r}}\partial_{n+1}^{[\mathbf{r}]} \in \mathcal{D}(X_{n+1})$ that is $0$ restricted  to $\w_{n+1}(A)$. Assume that the set $S:=\{\mathbf{r} \in \mathbb{N}^{d} \mid a_{\mathbf{r}} \not \in p^{v_p(\mathbf{r})+1}A_{n+1} \}$ is not empty. Then, we have that
 \begin{equation*}
     \sum_{\mathbf{r} \in S} a_{\mathbf{r}}\partial^{[\mathbf{r}]}_{n+1 \mid \w_{n+1}(A)}=0.
 \end{equation*}
 Fix such a $\mathbf{r}=(r_1,\dots,r_d) \in S$. 
 In particular, $a_{\mathbf{r}} \neq 0$. 
 Let $f_{\mathbf{r}}=\ve^{n-v_p(\mathbf{r})}(\prod_{j=1}^d[t_j]^{r_jp^{-v_p(\mathbf{r})}})$. Then, for any other $\mathbf{s}=(s_1,\dots,s_d) \in \mathbb{N}^d$, we have
 \begin{equation*}
     \partial_{n+1}^{[\mathbf{s}]}(\tilde{w}_n(f_{\mathbf{r}}))=p^{n-v_p(\mathbf{r})}\binom{r_1}{s_1}\cdots \binom{r_d}{s_d}\prod_{j=1}^dt_j^{r_j-s_j}
 \end{equation*}
 where the equality above is meant to be $0$ if $s_j > r_j$ for some $j$. Thus, we have the following equality
 \begin{equation*}
     0=\sum_{\mathbf{s} \in S} a_{\mathbf{s}}\partial_{n+1}^{[\mathbf{s}]}(\tilde{w}_n(f_{\mathbf{r}}))= \sum_{\mathbf{s} \leq \mathbf{r}}p^{n-v_p(\mathbf{r})}\binom{r_1}{s_1}\cdots \binom{r_d}{s_d}\prod_{j=1}^dt_j^{r_j-s_j}a_{\mathbf{s}} \in A_{n+1}.
 \end{equation*}
 All the monomials $\prod_{j=1}^dt_j^{r_j-s_j}$ are different varying $\mathbf{s} \in S$. Thus, in particular for $\mathbf{s}=\mathbf{r}$, we have that 
 $$p^{n-v_p(\mathbf{r})}a_{\mathbf{r}}=0 \in A_{n+1}$$
 with $a_{\mathbf{r}} \neq 0$, implying that $a_{\mathbf{r}} \in p^{v_p(\mathbf{r})+1}A_{n+1}$; this is a contradiction to  $\mathbf{r} \in S$.
  \end{proof}
  \begin{rem}
     Note that we can always lift $\partial \in \mathcal{D}(A)$ to a $ \tilde{\partial} \in \mathcal{D}(\w_{n+1}(A))$ by lifting the coefficients of $\partial=\sum_{\mathbf{r}}b_{\mathbf{r}}\partial_1^{[\mathbf{r}]}$ to some $b^{(n+1)}_{\mathbf{r}} \in \w_{n+1}(A)$ (e.g. by taking the Teichm\"uller lift). 
 \end{rem}
   Another consequence of \Cref{claimhighdiff} is the following:
  \begin{lemma}\label{lem:image_res_diff_op}
  For any differential operator $\delta \in \mathcal{D}(X_{n+1})$ of order $q \geq 1$, with $v_p(q)\leq n$ and admitting a restriction to $\w_{n+1}(A)$, then 
  \begin{equation}\label{diff_order_q_versch}
  \delta_{\mid \w_{n+1}(A)}(\w_{n+1}(A)) \subset \ve^{n-v_p(q)}\w_{n+1}(A).
\end{equation}
\begin{proof}
Indeed, if $f=(f_1,\dots,f_{n+1}) \in \w_{n+1}(A)$ by \Cref{claimhighdiff} and  \Cref{liftder} we see that 
\begin{equation*} 
\tilde{w}_n(\delta(f))=\delta(\tilde{w}_n(f)) \in p^{n-v_p(q)}A_{n+1} \cap \tilde{w}_n\w_{n+1}(A).
\end{equation*}
Since $\tilde{w}_n$ maps injectively $\ve^{n-v_p(q)}\w_{n+1}(A)$ into $p^{n-v_p(q)}A_{n+1}$, then \begin{equation*}\delta(f) \in \ve^{n-v_p(q)}\w_{n+1}(A). \qedhere
\end{equation*} 
\end{proof}
  \end{lemma}
   In the following, we are going to prove some relations involving $\partial_{n+1}^{[r]}$ and $R,\F,\ve$, \\ where $\F \colon \w_{n+1}(A) \rightarrow \w_{n+1}(A)$ denotes the Witt vectors Frobenius induced by $\sigma$ (the Frobenius of $A$), and $R \colon \w_{n+1}(A) \rightarrow \wn(A)$ is the natural projection.
  

  \begin{lemma}\label{weakFpartial}
  Let $r \geq 0$ be an integer, and $1 \leq j \leq d$. Let $\partial_{j,1}^{[r]} \in \mathcal{D}(X)$. Then, for any $f \in A$
  \begin{equation}
   \partial_{j,1}^{[r/p]}(f)^p= \partial_{j,1}^{[r]}({f}^p), 
  \end{equation}
   where $\partial_{j,1}^{[r/p]}$ is meant to be $0$ if $p \nmid r$. Consequently, for any $\mathbf{r} \in \mathbb{N}^d$,
   \begin{equation}
      \partial_{1}^{[\mathbf{r}/p]}(f)^p= \partial_{1}^{[\mathbf{r}]}({f}^p).
   \end{equation}
   \end{lemma}

  \begin{proof} We fix such a $j$ and with abuse of notation, we simply write $\partial_{j,1}^{[r]}=:\partial_{1}^{[r]}$.  By Formula \eqref{general_H-S_formula}, we have
  \begin{equation*}
  \partial_1^{[r]}({f}^p)= \sum_{\substack{j_1+\dots + j_{p}=r \\ 0 \leq j_1,\dots ,j_{p} \leq r}}\partial_{1}^{[j_1]}({f})\dots \partial_{1}^{[j_{p}]}({f})=
 \partial_1^{[r/p]}(f)^p+ \sum_{\substack{j_1+\dots + j_{p}=r \\ 0 \leq j_1,\dots ,j_{p} \text{ not all equal } }}\partial_{1}^{[j_1]}({f})\dots \partial_{1}^{[j_{p}]}({f}).
  \end{equation*}
  Any partition of $r$ of length $p$ corresponds to a subset of $$\{(j_1,\dots,j_p)\mid j_1 +\dots +j_p=r, \text{ } j_1,\dots, j_r \geq 0\}$$ whose elements are the $p$-uples given by permuting elements of the partition. The corresponding subset to a partition not containing all equal elements  $j_1,\dots,j_p$ has cardinality being a multiple of $p$. Indeed, more precisely it is $\frac{p!}{c_1!\dots c_m!}$, where $m$ is the cardinality of the set $\{j_1,\dots,j_p\}=\{t_1,\dots,t_m\}$ (where $t_1,\dots,t_m$ are pairwise distinct)  and $c_l$ is the number of times for which $t_l$ appears in $(j_1,\dots,j_p)$. Since by assumptions, $c_l<p$ for any $l=1,\dots,m$, then $p \mid \frac{p!}{c_1!\dots c_m!}$. Thus, the big sum in the right-most hand side above is $0$.
  \end{proof}
  \begin{rem}
  If $\tilde{f} \in A_{n+1}$ is any lift of $f \in A $,  \Cref{weakFpartial} rephrases by the equivalence $$ \partial_{j,n+1}^{[r/p]}(\tilde{f})^p \equiv  \partial_{j,n+1}^{[r]}(\tilde{f}^p) \pmod{p}.$$
  \end{rem}
  \begin{lemma}
     For any $x,y \in \w_{n+1}(A)$ the following relations hold:
    \begin{align}
    & R\tilde{w}_n(x)= (F\circ \tilde{w}_{n-1} \circ R)(x) \pmod{p^nA_{n+1}} \label{formulaw_nw_n-1} \\
        &
      \tilde{w}_n(x) \equiv \tilde{w}_n(y) \pmod{p^n} \text{ } \Leftrightarrow x \equiv y \pmod{\ve^n\w_{n+1}(A)}.\label{equivwnmodp}
  \end{align}
  \end{lemma}
  \begin{proof}  
     Let $x=(x_1,\dots,x_{n+1}) \in \w_{n+1}(A)$, then $$R\tilde{w}_n(x)=\tilde{x}_1^{p^{n}}+\dots+p^{n-1}\tilde{x}_n^p=F(\tilde{w}_{n-1}(x_1,\dots,x_n))=(F \circ \tilde{w}_{n-1} \circ R)(x).$$
          For $x,y \in \w_{n+1}(A)$ , we have
  $\tilde{w}_n(x) \equiv \tilde{w}_n(y) \pmod{p^n}  \Leftrightarrow  R\tilde{w}_n(x)=R\tilde{w}_n(y)$, thus by the \eqref{formulaw_nw_n-1}, and injectivity of $\tilde{F}^{n}$, we get  
  \begin{equation*}
  R\tilde{w}_n(x)=R\tilde{w}_n(y) \Leftrightarrow \tilde{F}^{n}(Rx)=\tilde{F}^{n}(Ry) \Leftrightarrow Rx=Ry \Leftrightarrow x \equiv y \pmod{\ve^n\w_{n+1}(A)}. \qedhere
  \end{equation*}
  \end{proof}

  	\begin{prop}\label{relliftderivation}
  Let $r \geq 0$ be an integer and $1 \leq j \leq d$.
   Then the following relations hold:
  \begin{align}
    \partial_{j,n}^{[r/p]}\circ R = R \circ \partial_{j,n+1}^{[r]};\label{derrestr} \\
  \partial_{j,n+1}^{[r]}\circ \F=\F \circ \partial_{j,n+1}^{[r/p]}; \label{derfrob}\\
    \partial_{j,n+1}^{[r]}\circ \ve=\ve \circ \partial_{j,n}^{[r]}; \label{derversch}\\
    \partial_{j,n+1}^{[r]} (\ve^i\w_{n+1}(A)) \subset \ve^i\w_{n+1}(A) & \quad \text{for all } 0 \leq i < n+1;\label{derfiltr}
  \end{align}
  where $\partial_{j,n+1}^{[r/p]}$ is meant to be $0$ if $p \nmid r$ and $R=R_{n+1,n}$. Consequently, the analogous relations \text{\eqref{derrestr},\eqref{derfrob},\eqref{derversch},\eqref{derfiltr}} hold for $\partial_{n+1}^{[\mathbf{r}]}$ for any $\mathbf{r} \in \mathbb{N}^d$.
    \end{prop} 
    
    \begin{proof}
  As in the proof of \Cref{weakFpartial}, we fix a $j$ and omit its notation from the corresponding differential operator. The last relation follows by the commutative diagram in  \Cref{liftder} together with the fact that $\tilde{w}_n$ maps injectively $\ve^i\w_{n+1}(A)$ to $p^iA_{n+1}$.
    Firstly, assume $p \nmid r$. Then, $v_p(r)=0$ and by \Cref{lem:image_res_diff_op}, it follows that
    $\partial_{n+1}^{[r]}(\w_{n+1}(A)) \subset \ve^n\w_{n+1}(A) $, so that $R \circ \partial_{n+1}^{[r]}=0 \text{ (}\in \wn(A))$. Moreover, by the relation \eqref{eq:claimhighdiff}, 
    $\partial_{}^{[r]}(p^i\tilde{f}^{p^{n+1-i}}) \equiv 0 \pmod{p^{n+1}}$, for any $\tilde{f} \in A_{n+1}$, thus
    $\partial_{n+1}^{[r]} \circ \F =0.$
    
    Now, we are going to prove that $R \circ \partial_{n+1}^{[p]}=\partial_{n} \circ R$. 
    
     By linearity, it suffices to test the relation on the elements of the form $\ve^i([f]) \in \w_{n+1}(A)$. Then,
    \begin{equation}
    \tilde{w}_{n}\partial_{n+1}^{[p]}(\ve^i([f]))=p^i\partial_{n+1}^{[p]}(\tilde{f}^{p^{n-i}}).
    \end{equation}
    By  
    \Cref{weakFpartial}, the following hold \footnote{Notice that $\frac{1}{p^{n-i-1}}\binom{p^{n-i}}{p}=\frac{1}{(p-1)!}\prod_{j=1}^{p-1}(p^{n-i}-j) \equiv 1 \pmod p$.}:
    \begin{align}
    p^i\binom{p^{n-i}}{p} \equiv p^{n-1} \pmod{p^{n}};\\
    \partial_{n+1}^{[p]}(\tilde{f}^p) \equiv \partial_{n+1}(\tilde{f})^p
 \pmod p.
     \end{align}
    By \cite[Proposition 9]{Nakai70} and the congruences above,
    \begin{multline}
    p^i\partial_{n+1}^{[p]}(\tilde{f}^{p^{n-i}})\equiv_{p^n} p^i\binom{p^{n-i}}{p} \tilde{f}^{p^{n-i}-p}\partial_{n+1}^{[p]}(\tilde{f}^p) \equiv_{p^n}\\
    p^{n-1}\tilde{f}^{p^{n-i}-p}\partial_{n+1}(\tilde{f})^p \equiv_{p^n} \tilde{w}_{n}(\ve^{n-1}([f^{p^{n-i-1}-1}\partial_1(f)]))
    \end{multline}
    By the \eqref{equivwnmodp}, it follows
    \begin{equation}
    (R \circ \partial_{n+1}^{[p]})(\ve^i([f]))=\ve^{n-1}([f^{p^{n-i-1}-1}\partial_1(f)]).
    \end{equation}
    
    On the other hand,
    \begin{equation}
    \tilde{w}_{n-1}\partial_{n}^{}(R\ve^i([f]))=p^{n-1}\tilde{f}^{p^{n-i-1}-1}\partial_n(\tilde{f})=
    \tilde{w}_{n-1}(\ve^{n-1}([f^{p^{n-i-1}-1}\partial_1(f)]))
    \end{equation}
    so that 
    \begin{equation}
    R\partial_{n+1}^{[p]}(\ve^i([f]))= \ve^{n-1}([f^{p^{n-i-1}-1}\partial_1(f)])
        =\partial_n(R\ve^i([f])).
    \end{equation}
   
    Suppose now $r >p$. We can write
    \begin{equation}
     \tilde{w}_n\partial_{n+1}^{[r]}(\ve^i([f]))=\partial_{n+1}^{[r]}(p^i\tilde{f}^{p^{n-i}})=
   p^i\sum_{\substack{j_1+\dots + j_{p^{n-i}}=r \\ 0 \leq j_1,\dots ,j_{p^{n-i}}\leq r}}\partial_{n+1}^{[j_1]}(\tilde{f})\dots \partial_{n+1}^{[j_{p^{n-i}}]}(\tilde{f}).
   \end{equation}
   Assume $i<n$, otherwise, for $i=n$ the relation $R(\partial_{n+1}^{[r]}(\ve^n([f])))=0$ is trivial. The trick here is that to compute the sum above $\pmod{p^n}$, it suffices to assume  every $j_l \leq r/p$, for any $1 \leq l \leq p^{n-i}$. Indeed, if $\{t_1,\dots,t_m \}$ denotes the set $\{j_1,\dots,j_{p^{n-i}} \}$ where the  $t_l$'s are pairwise distinct, we have the relation (since $\partial_{n+1}^{[j_l]}$'s commute)
   \begin{equation}\label{expansionder}
   p^i\sum_{\substack{j_1+\dots + j_{p^{n-i}}=r \\ 0 \leq j_1,\dots ,j_{p^{n-i}}\leq r}}\partial_{n+1}^{[j_1]}(\tilde{f})\dots \partial_{n+1}^{[j_{p^{n-i}}]}(\tilde{f})=
       p^i\sum_{\substack{t_1 < \dots < t_m \\ c_1t_1+\dots + c_mt_m=r \\ c_1+ \dots +c_m=p^{n-i} }} a^{(p^{n-i})}_{(c_1,\dots,c_m)}\partial_{n+1}^{[t_1]}(\tilde{f})^{c_1}\dots \partial_{n+1}^{[t_{m}]}(\tilde{f})^{c_m}.
   \end{equation}
   For any such $m$-uple $(t_1<\dots< t_m)$, the coefficient $a^{(p^{n-i})}_{(c_1,\dots,c_m)} \in \zz$ is the same appearing in the corresponding monomial of the expansion of $(\partial_{n+1}^{[t_1]}(\tilde{f})+\dots+\partial_{n+1}^{[t_m]}(\tilde{f}))^{p^{n-i}}$. Indeed, in both case, it is computed as the number of ways to partition the set $\{1,\dots, p^{n-1}\}$ in $m$ subsets $C_1, \dots, C_m$ with $|C_i|=c_i$.
   Notice that we have the elementary equality $\pmod{p^n}$:
   \begin{equation} \label{elemtaryp^n}
   p^i(\partial_{n+1}^{[t_1]}(\tilde{f})+\dots+\partial_{n+1}^{[t_m]}(\tilde{f}))^{p^{n-i}}\equiv p^i({\partial_{n+1}^{[t_1]}(\tilde{f})}^p+\dots+{\partial_{n+1}^{[t_m]}(\tilde{f})}^p)^{p^{n-i-1}} \pmod{p^n}.
   \end{equation}
   In particular, the relation \eqref{elemtaryp^n} implies that
   $p^ia^{(p^{n-i})}_{(c_1,\dots,c_m)}\equiv 0 \pmod{p^n}$ if there exists some index $j$ with $1 \leq j \leq m$ such that $p \nmid c_j$ and
   \begin{equation*}
      p^i a^{(p^{n-i})}_{(c_1,\dots,c_m)} \equiv p^ia^{(p^{n-i-1})}_{(c_1/p,\dots,c_m/p)} \pmod{p^n}
   \end{equation*} 
   otherwise.
    Thus,  looking at the non zero coefficients of \eqref{expansionder}, we see that every $t_l$ appears a multiple of $p$ times among the $j_l$'s. Thus, if there exists some $t_l > r/p$, then $j_1+\dots+j_{p^{n-i}} \geq pt_l>r$ that is impossible.
  
  
By the argument above,  we can assume $p \mid c_l$ for every $l$. Write $c_l=pc_l'$. 
Thus, we have the following equality $\pmod{p^n}$:
\begin{multline*}
 p^i\sum_{\substack{t_1 <\dots<t_m \\ c_1t_1+\dots + c_mt_m=r \\ c_1+ \dots +c_m=p^{n-i} }}a^{(p^{n-i})}_{(c_1,\dots,c_m)}\partial_{n+1}^{[t_1]}(\tilde{f})^{c_1}\dots \partial_{n+1}^{[t_{m}]}(\tilde{f})^{c_m} \equiv_{p^n} \\
p^i\sum_{\substack{t_1<\dots<t_m \leq r/p \\ c_1't_1+\dots + c_m't_m=r/p \\ c_1'+ \dots +c_m'=p^{n-i-1} }}a^{(p^{n-i-1})}_{(c_1/p,\dots,c_m/p)}\partial_{n+1}^{[t_1]}(\tilde{f})^{pc_1'}\dots \partial_{n+1}^{[t_{m}]}(\tilde{f})^{pc_m'} 
\end{multline*}

that gives
  \begin{equation*}
   R\tilde{w}_{n} \partial_{n+1}^{[r]}(\ve^i([f]))=
   p^i\sum_{\substack{t_1<\dots<t_m \leq r/p \\ c_1't_1+\dots + c_m't_m=r/p\\ c_1'+ \dots +c_m'=p^{n-i-1} }}a^{(p^{n-i-1})}_{(c'_1,\dots,c'_m)}\partial_{n+1}^{[t_1]}(\tilde{f})^{pc_1'}\dots \partial_{n+1}^{[t_{m}]}(\tilde{f})^{pc_m'}.
     \end{equation*}
     In particular any term in the sum belongs to $\mathrm{Im}(F)$. 
    Together with the relation above and the \eqref{formulaw_nw_n-1}, it follows that
    \begin{multline*}
   \tilde{w}_{n-1} R \partial_{n+1}^{[r]}(\ve^i([f])) = F^{-1}R\tilde{w}_{n} \partial_{n+1}^{[r]}(\ve^i([f]))=\\
   \sum_{\substack{t_1<\dots<t_m \\ c_1't_1+\dots + c_m't_m=r/p \\ c_1'+ \dots +c_m'=p^{n-i-1} }}p^ia^{(p^{n-i-1})}_{(c'_1,\dots,c'_m)}\partial_{n}^{[t_1]}(\tilde{f})^{c_1'}\dots \partial_{n}^{[t_{m}]}(\tilde{f})^{c_m'} = \tilde{w}_{n-1}\partial_{n}^{[r/p]}(\ve^i([f])). 
 \end{multline*}
  This shows that 
   $\tilde{w}_{n-1} R \partial_{n+1}^{[r]}(\ve^i([f]))=\tilde{w}_{n-1}\partial_{n}^{[r/p]}(\ve^i([f])) \in A_{n+1}/{p^n}\simeq A_n$, thus by injectivity of $\tilde{w}_{n-1}$, we have $R \partial_{n+1}^{[r]}(\ve^i([f]))=\partial_{n}^{[r/p]}(\ve^i([f]))$.
   This shows the relation \eqref{derrestr}. 
   
    To prove the \eqref{derfrob}, 
   we can write analogously
   \begin{equation}
     \tilde{w}_n\partial_{n+1}^{[r]}(\F\ve^i([f]))= \partial_{n+1}^{[r]}(p^i(\tilde{f}^p)^{p^{n-i}})=
     p^i\sum_{\substack{t_1 <\dots< t_m \\ c_1t_1+\dots + c_mt_m=r \\ c_1+ \dots +c_m=p^{n-i} }}a^{(p^{n-i})}_{(c_1,\dots,c_m)}\partial_{n+1}^{[t_1]}(\tilde{f}^p)^{c_1}\dots \partial_{n+1}^{[t_{m}]}(\tilde{f}^p)^{c_m}.
   \end{equation}
   As before, 
\begin{equation}
     \tilde{w}_n\partial_{n+1}^{[r/p]}(\ve^i([f]))=\partial_{n+1}^{[r/p]}(p^i\tilde{f}^{p^{n-i}})=
   p^i\sum_{\substack{t_1<\dots<t_m \\ c_1t_1+\dots + c_mt_m=r/p \\ c_1+ \dots +c_m=p^{n-i} }}a^{(p^{n-i})}_{(c_1,\dots,c_m)}\partial_{n+1}^{[t_1]}(\tilde{f})^{c_1}\dots \partial_{n+1}^{[t_{m}]}(\tilde{f})^{c_m}.
   \end{equation} 
   Now by  \Cref{weakFpartial}, we have 
   \begin{equation}
   \partial_{n+1}^{[t_l]}(\tilde{f})^p\equiv \partial_{n+1}^{[pt_l]}(\tilde{f}^p) \pmod{p}.
   \end{equation}
     Thus, for $i \leq n$ we get
   \begin{align*}
    \tilde{w}_n\F\partial_{n+1}^{[r/p]}(\ve^i([f]))
    & = F\tilde{w}_n\partial_{n+1}^{[r/p]}(\ve^i([f]))
    \\& =p^i\sum_{\substack{t_1<\dots<t_m \\ c_1t_1+\dots + c_mt_m=r/p \\ c_1+ \dots +c_m=p^{n-i} }}a^{(p^{n-i})}_{(c_1,\dots,c_m)}\partial_{n+1}^{[t_1]}(\tilde{f})^{pc_1}\dots \partial_{n+1}^{[t_{m}]}(\tilde{f})^{pc_m} 
    \\&=p^i\sum_{\substack{t_1<\dots<t_m \\ c_1t_1+\dots + c_mt_m=r/p \\ c_1+ \dots +c_m=p^{n-i} }}a^{(p^{n-i})}_{(c_1,\dots,c_m)}\partial_{n+1}^{[pt_1]}(\tilde{f}^p)^{c_1}\dots \partial_{n+1}^{[pt_{m}]}(\tilde{f}^p)^{c_m}    
   \\ & = p^i\sum_{\substack{t_1<\dots<t_m \\ c_1t_1+\dots + c_mt_m=r \\ c_1+ \dots +c_m=p^{n-i} }} a^{(p^{n-i})}_{(c_1,\dots,c_m)}\partial_{n+1}^{[t_1]}(\tilde{f}^p)^{c_1}\dots \partial_{n+1}^{[t_{m}]}(\tilde{f}^p)^{c_m}
   \\ & = \tilde{w}_n\partial_{n+1}^{[r]}(\F\ve^i([f])),
   \end{align*}
   from which we deduce \eqref{derfrob}.
   
   To prove the \eqref{derversch}, we use that $R\F\ve=p$. Then,
   \begin{equation}
    R\F \circ (\ve \circ \partial_{n}^{[r]})= (\partial_{n}^{[r]} \circ R) \circ \F\ve = R \partial_{n+1}^{[pr]} \circ  \F \ve = R \F  \circ (\partial_{n+1}^{[r]} \circ \ve) 
   \end{equation}
 This means that 
   \begin{equation}
   \F \circ \ve \circ  \partial_{n}^{[r]} \equiv \F \circ \partial_{n+1}^{[r]} \circ \ve \pmod{\ve^n\w_{n+1}(A)}.
   \end{equation}
   By the injectivity of $\F$, follows that $\mathrm{Im}(\ve \circ  \partial_{n}^{[r]} -\partial_{n+1}^{[r]} \circ \ve) \subset \ve^n\w_{n+1}(A)$.
   Hence, since this holds for every $n$ and every $r$, also  $\mathrm{Im}(\ve \circ  \partial_{n+1}^{[pr]} -\partial_{n+2}^{[pr]} \circ \ve) \subset \ve^{n+1}\w_{n+2}(A)$, thus $$R_{n+2,n+1}\circ\ve \circ  \partial_{n+1}^{[pr]} - R_{n+2,n+1} \circ \partial_{n+2}^{[pr]} \circ \ve=0.$$ Since  $R_{n+2,n+1}\circ \ve=\ve \circ R$, by the relation \eqref{derrestr} follows that $\ve \circ  \partial_{n}^{[r]}\circ R -\partial_{n+1}^{[r]} \circ \ve \circ R=0$. This means that $\ve \circ  \partial_{n}^{[r]}=\partial_{n+1}^{[r]} \circ \ve \colon \wn(A)\rightarrow\w_{n+1}(A)$.\qedhere
   
    \end{proof}
\subsection{The ring of Witt differential operators}
   
  Let consider $\w_{n+1}(A)$ as a $\w(A)$-module, via the natural restriction $\w(A) \xrightarrow{R_{n+1}} \w_{n+1}(A)$.
  For any $m  \in \zz$, and any $a \in \w(A)$ we view $\F^{m-n}(a) \in \mathrm{End}_{\w(A)[1/p]}(\w_{}(A)[1/p])$ by left multiplication \footnote{For $m \in \mathbb{N}$, recall that $\Phi_A^{-m}=\frac{1}{p^m}\ve^m$.}. In particular, 
  \begin{lemma}
  Let $v_p(r) \leq n$ for some $r \geq 1$, and $1 \leq j \leq d$. The composition
  \begin{equation} \label{map_scalar_exten_witt_diff}
      \F^{v_p(r)-n}(a) \circ \partial_{j,n+1}^{[r]}=:\F^{v_p(r)-n}(a) \partial_{j,n+1}^{[r]} \colon \w_{n+1}(A) \rightarrow \ve^{n-v_p(r)}\w_{n+1}(A) \subset \w_{n+1}(A)
  \end{equation}
  is a well defined $\w_{n+1}(k)$-linear differential operator for any $a \in \w(A)$.
  \end{lemma}
  \begin{proof}
    Since for any $m \geq 0$  the equality  $x\ve^m(y)=\ve^m(\F^m(x)y)$ holds for any $x,y \in \w_{}(A)$, it follows that 
$$\F^{m-n}(a)(\ve^{n-m}\w(A)) \subset \ve^{n-m} \w(A),$$ where $\ve^{n-m}\w(A)=\w(A)$ if $m \geq n$. Thus, the restriction $R_{n+1}$ induces a $\w_{n+1}(A)$-linear map  $ \F^{m-n}(a) \in \mathrm{End}_{\w_{n+1}(A)}(\ve^{n-m}\w_{n+1}(A))$. Moreover, by the \eqref{diff_order_q_versch} it follows  that $\partial^{[r]}_{j,n+1 \mid \w_{n+1}(A)} \subset \ve^{n-v_p(r)}\w_{n+1}(A)$. This proves that the map \eqref{map_scalar_exten_witt_diff} is well defined, thus a $\w_{n+1}(k)$-linear differential operator, since $\partial_{j,n+1}^{[r]}$ is already.
\end{proof}
    
By convention, we extend the definition for $r=0$, by letting 
  \begin{equation*}
      \F^{v_p(0)-n}(a)\partial_{n+1}^{[0]}:= a\cdot Id_{\w_{n+1}(A)}.
  \end{equation*} 
  Furthermore, as seen in the proof of \Cref{uniquenesslift}, $\partial_{n+1}^{[\mathbf{r}]} \subset \ve^{n-v_p(\mathbf{r})}\w_{n+1}(A)$. Thus, every operator of the form $$ \F^{v_p(\mathbf{r})-n}(a)\partial_{n+1}^{[\mathbf{r}]}=\F^{v_p(\mathbf{r})-n}(a)\partial_{j_0,n+1}^{[r_{j_0}]}\prod_{j=1,j \neq j_0}^d \partial_{j,n+1}^{[r_j]},$$ with $r \geq 0$, and $1 \leq j_0 \leq d$ such that $v_p(\mathbf{r})=v_p(r_{j_0})$, defines an element of $\mathcal{D}(\w_{n+1}(A))$.
  
  \begin{lemma}\label{lemma_trivial_sum_diff}
      Any $\w_{n+1}(A)$-linear combination of the form $$\sum_{\substack{\mathbf{r} \in \mathbb{N}^d \\ v_p(\mathbf{r}) \leq n}}\F^{v_p(\mathbf{r})-n}
(a_{\mathbf{r}})\partial_{n+1}^{[\mathbf{r}]}+\sum_{\substack{\mathbf{r} \in \mathbb{N}^d \\ v_p(\mathbf{r}) > n}}a_{\mathbf{r}}\partial_{n+1}^{[\mathbf{r}]}$$ 
      is trivial if and only if $a_{\mathbf{r}} \in \ve^{v_p(\mathbf{r})+1}\w_{n+1
}(A)$ for $\mathbf{r} \in \mathbb{N}^d$ such that $v_p(\mathbf{r} ) \leq n$ and $a_{\mathbf{r}}=0$ otherwise.
      \end{lemma}
      \begin{proof}
      Indeed, as in the proof of  \Cref{uniquenesslift}, for any $\mathbf{r} \in \mathbb{N}^d$, there is an element $f_{\mathbf{r}}\in \w_{n+1}(A)$, such that 
      $$\partial^{[\mathbf{s}]}_{n+1}(f_{\mathbf{r}})=\left\{\begin{array}{ll}
           c_{\mathbf{r},\mathbf{s}}\ve^{n-v_p(\mathbf{r})}([\mathbf{t}]^{p^{-v_p(\mathbf{r})}(\mathbf{r}-\mathbf{s})}) & \text{if } n \geq v_p(\mathbf{s})>v_p(\mathbf{r}) \text{ or } v_p(\mathbf{r})< n < v_p(\mathbf{s})  \\
           e_{\mathbf{r},\mathbf{s}} \ve^{n-v_p(\mathbf{s})}([\mathbf{t}]^{p^{-v_p(\mathbf{s})}(\mathbf{r}-\mathbf{s})}) & \text{if } v_p(\mathbf{s}) \leq v_p(\mathbf{r}) \leq n
           \text{ or }   v_p(\mathbf{s})< n < v_p(\mathbf{r}) \\
      e_{\mathbf{r},\mathbf{s}} [\mathbf{t}]^{p^{-v_p(\mathbf{s})}(\mathbf{r}-\mathbf{s})}      & \text{if } n < v_p(\mathbf{s})\leq v_p(\mathbf{r}) \\
      c_{\mathbf{r},\mathbf{s}}[\mathbf{t}]^{p^{-v_p(\mathbf{r})}(\mathbf{r}-\mathbf{s})}      & \text{if } n < v_p(\mathbf{r})< v_p(\mathbf{s})
                 \end{array} \right.  $$
where $ c_{\mathbf{r},\mathbf{s}},e_{\mathbf{r},\mathbf{s}} \in \w_{n+1}(k)$  with $e_{\mathbf{r},\mathbf{r}}=1$. 
      It follows:
      \begin{multline*}
         \sum_{\substack{\mathbf{s} \\ v_p(\mathbf{s}) \leq n}}\F^{v_p(\mathbf{s})-n}(a_{\mathbf{s}})\partial_{n+1}^{[\mathbf{s}]}(f_{\mathbf{r}})+\sum_{\substack{\mathbf{s} \\ v_p(\mathbf{s}) > n}}a_{\mathbf{s}}\partial_{n+1}^{[\mathbf{r}]}(f_{\mathbf{r}})= \\
         =\left\{ \begin{array}{lc}
         \ve^{n-v_p(\mathbf{r})}(a_{\mathbf{r}})+ \text{linear combination of } \ve^i([\mathbf{t}]^l) \text{ s.t. } l \geq 1 & \text{if } v_p(\mathbf{r}) \leq n \\
         a_{\mathbf{r}}+ \text{linear combination of } \ve^i([\mathbf{t}]^l) \text{ s.t. } l \geq 1
         & \text{if } v_p(\mathbf{r}) > n.
         \end{array}\right.
      \end{multline*}
       Any trivial expression yields $\ve^{n-v_p(\mathbf{r})}(a_{\mathbf{r}})=0 \in \w_{n+1}(A)$ or $a_{\mathbf{r}}=0$.
  \end{proof}
    
 Let $X=\spec(A)$. If $A$ is smooth over $k$, then for any $\mathfrak{m} \in X$, there exists an open neighborhood $U$ of $\mathfrak{m}$  and local sections $z_1,\dots,z_m \in \Gamma(U,\ox)=A_f$ such that the scheme morphism 
    \begin{equation*}
      U \rightarrow \mathbb{A}_k^m=\spec(k[x_1,\dots,x_m]) 
    \end{equation*}
    induced by $A_f \ni z_i \mapsto x_i \in k[x_1,\dots,x_m] $, is étale.
    We say that the sections $(z_1,\dots,z_m)$ are \textit{local coordinates} of $X$ associated to a \textit{local chart} $U$.
    Notice that since composition of étale maps is étale, then any open $V \subset U$ is a local chart whenever it  $U$ is. In particular, for any smooth scheme $Y$ over $k$, the collection of those open affine subsets consisting of local charts forms a basis for the Zariski topology of $Y$. 
       If $A$ is smooth over $k$, then $A$ is locally étale over a polynomial algebra in the sense above. For any local chart $U=\spec(A_f)$, if $B$ denotes the corresponding polynomial algebra over $k$ generated by the local sections, one gets an identification of $k$-algebras $\mathcal{D}(A_{f})\simeq \mathcal{D}(B) \otimes_B A_{f}$ (cf. \cite[Theorem 3.2.5]{Will98}). By  \Cref{uniquenesslift}  we can lift any derivation of $\mathcal{D}(B)$ to $\mathcal{D}(\w_{n+1}(B))$. For any of such lifts, by  \Cref{lift_partial_H-S_étale}, since $\w_{n+1}(B) \rightarrow \w_{n+1}(A_f)$ is étale, corresponds a unique differential operator of $\mathcal{D}(\w_{n+1}(A_f))$. Hence,  we can lift any $k$-linear differential operator of $\mathcal{D}(A_{f})$ to a $\w_{n+1}(k)$-linear differential operator in $\mathcal{D}(\w_{n+1}(A_{f}))$. Since, $A \rightarrow A_f$, is also étale, the canonical map
       \begin{equation}
       \mathcal{D}(A)\rightarrow \mathcal{D}(A_{f})=\mathcal{D}(A)\otimes_A A_f
\end{equation}      
allows to lift a differential operator in  $\mathcal{D}(A)$ to a differential operator in $\mathcal{D}(\w_{n+1}(A_{f}))$ for any $f$ associated to a local chart in the covering of $X$. We are going to see that we can glue together suitable lifts in order to get a differential in $\mathcal{D}(\w_{n+1}(A))$ lifting the corresponding one in $\mathcal{D}(A)$.

If $(z_1,\dots,z_m)$ are local coordinates of $A$, we denote by $\partial_{i}^{[r]}:=\partial_{z_i}^{[r]}$ for any $r \geq 0$ and by $\partial_{i,n+1}^{[r]}$ its lift to $\w_{n+1}(A)$. Also, to simplify the notation, for any $a \in \w_{n+1}(A)$, we write for any $r \geq 0$,
\begin{equation}
a \bullet \partial_{i,n+1}^{[r]}:=\F^{v_p(r)-n}(a)\partial_{i,n+1}^{[r]}.
\end{equation}
      Notice that any sheaf for the Zariski topology of a smooth $k$-scheme $Y$ is determined by an open cover of local charts. Thus, we are lead to the following definition:
    \begin{defn}\label{def_witt_diff_smooth_algebras}
  For a smooth $k$-algebra $A$ of  dimension $m$, let $(z_1, \dots, z_m)$ be local coordinates associated to a local chart  $\spec(B)=U \subset \spec(A)$.
   For any $i$ such that $1 \leq i \leq m$, and $r \geq 0$, let  $\partial_{i}^{[r]}$ be the generators for the $B$-algebra $\mathcal{D}(B)$.  Let us consider the lifts $\partial_{i,n+1}^{[r]} \in \mathcal{D}(\w_{n+1}(B))$  obtained as in  \Cref{liftder}.
   
    We denote by $\mathcal{D}_{\w_{n+1}(A)}(U)$ the $\w_{n+1}(B)$-subalgebra of $\mathcal{D}(\w_{n+1}(B))$ generated by
   \begin{equation*}
    \begin{array}{lr}
    a \bullet \partial_{i,n+1}^{[r]},
   &  \text{  for any } a \in \w_{n+1}(A), \text{ } i=1,\dots,m, \text{ } r \geq 0.
   \end{array} 
   \end{equation*} Then, we define $\widetilde{\mathcal{D}_{\w_{n+1}(A)}}$ as the presheaf of $\w_{n+1}(k)$-algebras, given by  $$ \Gamma(U,\widetilde{\mathcal{D}_{\w_{n+1}(A)}}):=\mathcal{D}_{\w_{n+1}(A)}(U)=:\mathcal{D}_{\w_{n+1}(A)}(B)$$ for any local chart $U \subset \spec(A).$  
    \end{defn}
We will prove that the presheaf $\widetilde{\mathcal{D}_{\w_{n+1}(A)}}$ is indeed a sheaf in the following Lemma.
     \begin{rem}
   For $n=0$, the above definition coincides with  Grothendieck's sheaf of differential operators (defined for general commutative rings). For them, we can get rid of specifying such $U$. Thus, $ \widetilde{\mathcal{D}_{A}}=\widetilde{\mathcal{D}_{\w_1(A)}}=\mathcal{D}(A)$.   
   \end{rem}
 
  \begin{lemma}\label{lem_presh_of_wittdiff_is_sheaf}
  Let $X=\spec(A)$ for a smooth $k$-algebra $A$. Then,
  $\widetilde{\mathcal{D}_{\w_{n+1}(A)}}$ is a quasi-coherent $\w_{n+1}\ox$-mod.
  \end{lemma}
  \begin{proof}
  To prove the statement, it is sufficient (and necessary) to prove that for any local chart  $U=\spec(C) \subset X$ and any $f \in C$, $U_f=\spec(C_f)$  
  \begin{equation}\label{locDn}
  \mathcal{D}_{\w_{n+1}(A)}(U_f)=\mathcal{D}_{\w_{n+1}(A)}(U) \otimes_{\w_{}(C)} \w_{n+1}(C_f).
  \end{equation}
  Suppose the \eqref{locDn} holds. 
   Then take an open affine $U=\spec(A_g) \subset X$, with $g \in A$.
  An open cover of $U$ is given by a finite collection $(\spec(A_{gf_i}))_{i \in I}$ for $f_i \in A$ of spectra of localization at $gf_i \in A$ such that $(f_i)_{i \in I}$ generates $(1)=A_g$.  Set $B=A_g$. The sheaf condition requires to check that 
  \begin{equation}
  0 \rightarrow \mathcal{D}_{\w_{n+1}(A)}(U) \rightarrow \bigoplus_{i \in I}\mathcal{D}_{\w_{n+1}(A)}(U_{f_i}) \rightarrow \bigoplus_{i,j \in I} \mathcal{D}_{\w_{n+1}(A)}(U_{f_i}\times_U U_{f_j})
  \end{equation}
  is an exact sequence of $\w_{}(B)$-modules. Using that $B \rightarrow B_{f_i}$ are étale maps, we have 
  $\w_{n+1}(B_{f_i}\otimes B_{f_j})\simeq\w_{n+1}(B_{f_i}) \otimes_{\w_{n+1}(B)} \w_{n+1}(B_{f_j})\simeq\w_{n+1}(B)_{[f_if_j]}$. Together with the \eqref{locDn}, the above short exact sequence becomes,  
  \begin{equation}\label{sesDnlocsheaf}
  0 \rightarrow \mathcal{D}_{\w_{n+1}(A)}(B) \rightarrow \bigoplus_{i \in I}\mathcal{D}_{\w_{n+1}(A)}(B)_{[f_i]} \rightarrow \bigoplus_{i,j \in I} \mathcal{D}_{\w_{n+1}(A)}(B)_{[f_if_j]}
  \end{equation}
  It suffices to verify the exactness of \eqref{sesDnlocsheaf}  at the localization $B_{\mathfrak{m}}$ for any maximal ideal $\mathfrak{m} \subset B$. Thus, since $(f_i, \text{ } i \in I)=(1)=B$, we can further assume that there exists some $i \in I$ such that  $f_i \not \in \mathfrak{m}$, implying $f_i$ is invertible in $B_{\mathfrak{m}}$; we further assume that $f_1=1$. In this case, the \eqref{sesDnlocsheaf} is immediate. (compare with c.f. \cite[Tag 00EK, Lemma 10.24.1]{stacks-project}).
  
  Let us prove the relation \eqref{locDn}. It suffices to prove it for $U=X$ (i.e. $C=A$) assuming $X$ has  (local) coordinates $z_i$. Thus, let $f \in A$. 
 For any local coordinate $z_i$ of $A$, the natural map  $\mathcal{D}(\w_{n+1}(A)) \rightarrow \mathcal{D}(\w_{n+1}(A_f))$ sends $\partial_{{z_i}{},n+1}^{[r_i]} $ to $\partial_{\frac{z_i}{1},n+1}^{[r_i]}$ where $\partial_{\frac{z_i}{1},n+1}^{[r_i]}$ is the unique $r_i$-th differential operator lifting $\partial_{z_i,n+1}^{[r_i]} \in \mathcal{D}(\w_{n+1}(A))$ (since $\w_{n+1}(A) \rightarrow {\w_{n+1}(A_f)}=\w_{n+1}(A)_{[f]}$ is étale). Moreover, any  relation 
 $$\sum_{\substack{\mathbf{r} \\ v_p(\mathbf{r}) \leq n}} \Phi_{A_f}^{v_p(\mathbf{r})-n}(a_{\mathbf{r}})\partial_{\frac{\bm{z}}{1},n+1}^{[\mathbf{r}]}+\sum_{\substack{\mathbf{r} \\ v_p(\mathbf{r}) > n}} b_{\mathbf{r}}\partial_{\frac{\bm{z}}{1},n+1}^{[\mathbf{r}]}=0$$ 
 in $\mathcal{D}(\w_{n+1}(A_f))$ with $a_{\mathbf{r}},b_{\mathbf{r}} \in \w_{n+1}(A)$ yields, by  \Cref{lemma_trivial_sum_diff}, $a_{\mathbf{r}} \in \ve^{v_p(\mathbf{r})+1}(\w_{n+1}(A_f)) \cap \w_{n+1}(A)$, and $b_{\mathbf{r}}=0$. Thus, also $$\sum_{\mathbf{r}}\F^{v_p(\mathbf{r})-n}({a_{\mathbf{r}}})\partial_{\bm{z},n+1}^{[\mathbf{r}]}=0$$ in $\mathcal{D}(\w_{n+1}(A))$. Hence, it induces an injective map
   $$\mathcal{D}_{\w_{n+1}(A)}(A) \otimes_{\w_{n+1}(A)} \w_{n+1}(A_f) \hookrightarrow \mathcal{D}_{\w_{n+1}(A_f)}(A_f).$$ We notice, it is also surjective. Indeed, the isomorphism $\mathcal{D}(A) \otimes_A A_f \xrightarrow{\sim} \mathcal{D}(A_f)$ ensures that the collection $\{\partial_{\frac{z_i}{1},n+1}^{[r_i]} \mid r_i \geq 0\}$ is a set of generators over $\w_{n+1}(A_f)$ for $\mathcal{D}_{\w_{n+1}(A_f)}(A_f)$.       
  \end{proof}
  
   \begin{defn}\label{def_witt_diff_sheaf}
  For any smooth $k$-scheme $X$, and a covering  of affine local charts $\mathcal{U}:=\{U_i=\spec(A_i)\}_{i \in I}$ define $\mathcal{D}_{\w_{n+1}(X)}$ to be the unique quasi-coherent $\w_{n+1}\ox$-module such that $${\mathcal{D}_{\w_{n+1}(X)}}_{\mid U_i} \simeq \widetilde{\mathcal{D}_{\w_{n+1}(A_i)}}.$$
 \end{defn}
 \begin{rem}
 The  \Cref{def_witt_diff_sheaf} does not depend on the covering $\mathcal{U}$. Indeed, let $\{U_i=\spec(A_i)\}_{i  \in I}$ and $\{V_j=\spec(B_j)\}_{j \in J}$ be two coverings of $X$. Set $f_{ij}:\spec(C_{ij})=U_i \times_X V_j \rightarrow U_i$ and $g_{ij}:\spec(C_{ij})=U_i \times_X V_j \rightarrow V_j$ the respective open immersions, in particular they are étale maps.
 Let $\mathcal{F},\mathcal{G}$ be the unique sheaves of $\w_{n+1}\ox$-modules associated respectively to the covering $\{U_i\}_i$ and  $\{V_j\}_j$. Then, for any $i \in I, j \in J$
 $$(\mathcal{F}_{\mid U_i})_{\mid U_i \times_X V_j}= f_{ij}^{*}\widetilde{\mathcal{D}_{\w_{n+1}(A_i)}} \simeq \widetilde{\mathcal{D}_{\w_{n+1}(C_{ij})}} $$ and 
 $$(\mathcal{G}_{\mid V_j})_{\mid U_i \times_X V_j}= g_{ij}^{*}\widetilde{\mathcal{D}_{\w_{n+1}(B_j)}} \simeq \widetilde{\mathcal{D}_{\w_{n+1}(C_{ij})}}.$$
 Thus, there is an identity of sheaves 
 \begin{equation}
 \mathcal{F}_{\mid U_i \times_X V_j}=\mathcal{G}_{\mid U_i \times_X V_j}
 \end{equation}
 In particular they respect the cocycle condition for glueing sheaves with respect the covering $\{U_i \times_{X} V_j\}_{(i,j) \in I \times J}$ of $X$. It follows that $\mathcal{F} =\mathcal{G}$.
 \end{rem}
 \begin{rem}
 A definition of (completed) Witt differential operators appeared recently in the work of Dodd (c.f. \cite[Sec. 2]{dodd2024}) without truncation. Our definition for the truncated case is a subsheaf of the truncated version by Dodd. The main difference in \cite{dodd2024}, is the introduction of a canonical Hasse-Schmidt derivation, which is a specific lift over the Witt vectors of a Hasse-Schmidt derivation (i.e. operators satisfying \eqref{HS-formula}) in characteristic $p$, determined by an explicit formula (c.f. \cite[Corollary 2.6]{dodd2024}). Further, his sheaf of Witt differential operators is defined  intrinsically, avoiding local coordinates. However, a local description (c.f. \cite[Theorem 2.17]{dodd2024}) is in practice employed to give a presentation from which main properties follow.
 \end{rem}
 
\section{Hodge-Witt cohomology of Drinfeld's half space via local cohomology}\label{chap_drinf_hodge_witt}
\subsection{A spectral sequence for local cohomology}
Let $(X,\oo)$ be a ringed space, $\mathcal{F}$ be an $\oo$-module and $K^{\bullet}$ a bounded from below complex of $\oo$-modules. Let $\mathcal{I}^{\bullet}$ be an injective resolution of $\mathcal{F}$. Then, there is a first quadrant spectral sequence induced by the double complex $\mathrm{Hom}_{}(K^{\bullet},\mathcal{I}^{\bullet})$ \footnote{This set of homomorphism is taken in the category of complexes of sheaves}
\begin{equation}\label{generalss}
E_1^{r,s}=\ext^s_{}(K^{-r},\mathcal{F}) \implies \ext^{s+r}_{}(K^{\bullet},\mathcal{F}).
\end{equation}
Now let $(\mathcal{X},\oxx)$ be the $d$-dimensional Drinfeld upper half space over $k$, and $\mathcal{Y}$, its closed complement in $\mathbb{P}^d$. Set $\mathcal{O}:=\mathcal{O}_{\mathbb{P}_k^d}$, let $\mathcal{F}$ be an $\mathcal{O}$-module.
 Take an acyclic resolution 
 $ \zz_{\mathcal{Y}} \rightarrow J^{\bullet}$ of $\zz:=\zz_{\mathcal{Y}}$, the constant sheaf  over $\mathcal{Y}$ with value in $\zz$. Assume it is a finite resolution. Denote by $i: {\mathcal{Y}} \hookrightarrow \mathbb{P}^d$ the closed immersion, thus $i_{*}$ is exact. Then, if we take $K^{\bullet}=i_{*}J^{\bullet}$, it is a resolution of $i_{*}\zz$  and 
$$\ext^{r+s}_{}(K^{\bullet},\mathcal{F})=\ext^{r+s}_{}(i_{*}\zz,\mathcal{F})=\h^{s+r}_{\mathcal{Y}}(\mathbb{P}^d,\mathcal{F}).$$
The first equality follows by acyclicity of the complex $0 \rightarrow i_{*}\zz \rightarrow K^{\bullet} \rightarrow 0$ and the last equality follows from \cite[Proposition 2.3 bis. (21)]{sga2}. By assumptions, $J^{\bullet}$ is bounded and starts from degree 0. Then, the spectral sequence above rewrites as
\begin{equation}
E_1^{-r,s}=\ext^s_{}(K^{r},\mathcal{F}) \implies \h^{s-r}_{\mathcal{Y}}(\mathbb{P}^d,\mathcal{F}).
\end{equation}

\subsection{Orlik's acyclic resolution}
Here we recall an acyclic resolution of $\zz_{\mathcal{Y}}$ (cf. \cite[Section 2.1]{Orl08}).
For any $I \subset \Delta$, let $P_I \subset G$ be the associated parabolic subgroup. Let $\Delta \backslash I=\{\alpha_{i_0},\dots, \alpha_{i_r}\}$ with $i_0 < \dots <i_r$ and $\{e_0,\dots,e_d\}$ be the standard basis of $k^{d+1}$. Let $V_i=\sum_{s=0}^i ke_s$. Then, $Y_I:=\mathbb{P}(V_{i_0})$
is the closed $k$-subvariety of $\mathbb{P}^d$ stabilized by the action of $P_I$.
Notice that 
\begin{equation}
\mathcal{Y}= \bigcup_{I \subsetneq \Delta}\bigcup_{g \in G/{P_I}} g.Y_I , \quad Y_I \simeq \mathbb{P}^{i(I)}, \quad i(I):=\mathrm{min}\{j \mid \alpha_j \not \in I\}.
\end{equation}
Write $\Phi_{g,I}: g.Y_{I} \hookrightarrow Y$ for the  closed immersion given by the inclusion. Then, define the sheaves
\begin{equation}\label{orlikres}
 \quad J^r:=\bigoplus_{|I|=d-r-1}\bigoplus_{g \in G/{P_I}}(\Phi_{g,I})_{*}(\Phi_{g,I})^{-1}\zz, \quad r=0,\dots,d-1.
\end{equation}
For any $I \subset I' \subsetneq \Delta$ there are canonical inclusions $P_I \subset P_{I'}$ and closed immersions 
\begin{equation}
    \iota_{I,I'} \colon Y_{I} \hookrightarrow Y_{I'}.
\end{equation}
The projections 
\begin{equation}
    G/P_{I} \rightarrow G/P_{I'}, \quad gP_I \mapsto hP_{I'}
\end{equation}
induce also closed immersions
\begin{equation}
    \iota^{g,h}_{I,I'} \colon gY_I \rightarrow hY_{I'}
\end{equation}
such that $\Phi_{h,I'} \circ i_{I,I'}^{g,h}=\Phi_{g,I}$.  Furthermore, by functoriality, the map $i_{I,I'}^{g,h}$ induces a natural map of sheaves on $hY_{I'}$:
\begin{equation*}
    (\Phi_{h,I'})^{-1}\zz \rightarrow (i_{I,I'}^{g,h})_*(i_{I,I'}^{g,h})^{-1}(\Phi_{h,I'})^{-1}\zz
\end{equation*}
Then, applying the functor $(\Phi_{h,I})_*$ we get a map of sheaves on $\mathcal{Y}$: \begin{equation}
  p_{I,J}^{g,h} \colon  (\Phi_{h,I'})_{*} (\Phi_{h,I'})^{-1}\zz \rightarrow (\Phi_{g,I})_{*}(\Phi_{g,I})^{-1}\zz.
\end{equation}
Let 
\begin{equation}
    d_{I,I'}= \left\{\begin{array}{cc}
        (-1)^i\bigoplus_{(g,h) \in G/P_I \times G/P_{I'}} p_{I,J}^{g,h} & \text{ if } I'=I \sqcup \{\alpha_i\}, \\
       0  & \text{ otherwise} 
    \end{array}\right.
\end{equation}
where $p_{I,I'}^{g,h}$ is meant to be $0$ if $gP_I$ is not mapped to $hP_{I'}$. Then, the maps $d_{I,I'}$ induce
 morphisms $d^r : J^r\rightarrow J^{r+1}$, making $(J^{\bullet},d^{\bullet})$  a complex (cf. \cite[Section 2.1.1]{Ku16}). The following holds:
\begin{theorem}[{cf. \cite[Theorem 2.1.1]{Orl08}, \cite[Proposition 2.1.1.1]{Ku16}}]\label{theoacyclicity}
The complex of sheaves 
$0 \rightarrow \zz_{\mathcal{Y}} \rightarrow J^{\bullet} \rightarrow 0$  on ${\mathcal{Y}}$ is acyclic, i.e. it is an exact sequence in the category of sheaves.
\end{theorem} 
Observe that $(\Phi_{g,I})^{-1}\zz=\zz_{g.Y_I}$, that $i\circ \Phi_{g,I}$ is the closed immersion of $g.Y_I$ in $\mathbb{P}^d$, and $\ext^*_{}(-,\mathcal{F})$  commutes with direct sums. This in turn implies the following equality:
\begin{align*}
\ext^s_{}(i_{*}J^r,\mathcal{F}) & =  \bigoplus_{\substack{I \subset \Delta \\ |I|=d-r-1}}\bigoplus_{g \in G/{P_I}}\ext^s_{}((i \circ \Phi_{g,I})_{*}\zz_{g.Y_I}, \mathcal{F}) \\
 &=\bigoplus_{\substack{I \subset \Delta \\ |I|=d-r-1}}\bigoplus_{g \in G/{P_I}}\h^s_{g.Y_I}(\mathbb{P}^d,\mathcal{F}) \\
 &=\bigoplus_{\substack{I \subset \Delta \\ |I|=d-r-1}} \mathrm{Ind}^{G}_{P_I}\h^s_{Y_I}(\mathbb{P}^d,\mathcal{F}).
\end{align*}

\subsection{A spectral sequence for local cohomology of $\wn\oo_{\mathbb{P}^d}$-modules}
So at the end, the spectral sequence \eqref{generalss} above has the shape
\begin{equation}\label{ssext}
E_1^{-r,s}= \bigoplus_{|I|=d-r-1} \mathrm{Ind}^{G}_{P_I}\h^s_{Y_I}(\mathbb{P}^d,\mathcal{F}) \implies \h^{s-r}_{\mathcal{Y}}(\mathbb{P}^d,\mathcal{F}).
\end{equation}
Does this description hold in the context of Witt schemes? for which $\mathcal{F}$? Let $\mathcal{F}$ be a $\wn\oo$-module.
When $X$ is a $k$-scheme, the Witt scheme associated to $X$ is the ringed space $(|X|,\wn\ox)$. So there is a corresponding spectral sequence to \eqref{generalss}, in the category of sheaves on $\wn(X)$.
\begin{lemma}
Let $X$ be a $k$-scheme and $J^{\bullet}$ be a complex of sheaves on $\wn(X)$. Consider the natural closed immersion of schemes $\pi: X \hookrightarrow \wn(X)$. If $J^{	\bullet}$ is an acyclic complex of sheaves on $X$, then $\pi_{*}J^{\bullet}$ is an acyclic complex of sheaves on $\wn(X)$.
\end{lemma} 
\begin{proof}
We need to verify that $\pi_{*}J^{\bullet}$ viewed as a sequence of sheaves is exact. 
But $\pi$ is a closed immersion of schemes, thus $\pi_*$ is exact on the category of sheaves on $\wn(X)$, therefore $\pi_{*}J^{\bullet}$ is acyclic. \qedhere
\end{proof}
From  \Cref{theoacyclicity}, the following holds:
\begin{corol}
Let $J^{\bullet}$ be the complex \eqref{orlikres}. Then, $0 \rightarrow \pi_{*}\zz_{\mathcal{Y}} \rightarrow \pi_{*}J^{\bullet} \rightarrow 0$ is an acyclic complex 
of sheaves on $\wn(\mathcal{Y})$.
\end{corol}
Observe that $\pi$ is a closed immersion and a universal homeomorphism, since it is a nilpotent thickening. In particular, we have
$$\pi_{*}\zz_{\mathcal{Y}}=\zz_{\wn(\mathcal{Y})}. $$
This means that if $J^{\bullet}$ is an acyclic resolution of $\zz_{\mathcal{Y}}$ of sheaves on $\mathcal{Y}$, then $\pi_{*}J^{\bullet}$ is an acyclic resolution of $\zz_{\wn(\mathcal{Y})}$ of sheaves on  $\wn(\mathcal{Y})$. Moreover, when $i$ is a closed immersion, then $\wn(i)$ induces a closed immersion on the respective Witt schemes by \Cref{wn_preserve_immersions}.
Then, again we have identifications
\begin{equation}
\ext^{r+s}_{}(\wn(i)_{*}\pi_{*}J^{\bullet},\mathcal{F})=\ext^{r+s}_{}(\wn(i)_{*}\zz_{\wn({\mathcal{Y}})},\mathcal{F})=\h^{s+r}_{\wn(\mathcal{Y})}(\wn(\mathbb{P}^d),\mathcal{F})=\h^{s+r}_{\mathcal{Y}}(\mathbb{P}^d,\mathcal{F}),
 \end{equation}
 \begin{equation}
 \ext^{s}_{}(\wn(i)_{*}\pi_{*}J^{r},\mathcal{F})=\bigoplus_{|I|=d-r-1} \mathrm{Ind}^{G}_{P_I}\h^s_{Y_I}(\mathbb{P}^d,\mathcal{F}).
 \end{equation}
 In particular, when $\mathcal{F}=\wn\Omega^i_{\mathbb{P}^d}$, for any $i=0,\dots,d$, the spectral sequence \eqref{ssext} exists.
 
 Evaluating the spectral sequence \eqref{ssext}, we hope to compute $\h^*_{\mathcal{Y}}(\mathbb{P}^d,\mathcal{F})$. In turn this is related to $\h^0(\mathcal{X}, \mathcal{F})$ via the long exact sequence of the couple 
 \begin{equation*}
 \begin{tikzcd}
  {\mathcal{Y}} \arrow[r, closed]{}{} & \mathbb{P}^d & \arrow[l, open']{}{}  \mathcal{X}
\end{tikzcd}
\end{equation*}
 giving the following exact sequence 
 \begin{equation}
 0 \rightarrow \h^0(\mathbb{P}^d, \mathcal{F}) \rightarrow \h^0({\mathcal{X}}, \mathcal{F}) \rightarrow \tilde{\h}_{\mathcal{Y}}^1(\mathbb{P}^d, \mathcal{F})=\ker\big({\h}_{\mathcal{Y}}^1(\mathbb{P}^d, \mathcal{F})\rightarrow {\h}^1(\mathbb{P}^d, \mathcal{F})\big) \rightarrow 0,
 \end{equation}
 after noticing that ${\mathcal{X}}$ being affine and $\mathcal{F}$ quasi-coherent implies $\h^1({\mathcal{X}}, \mathcal{F})=0$. The local cohomology group $\tilde{\h}_{\mathcal{Y}}^1(\mathbb{P}^d, \mathcal{F})$ is related to the $E_1$ terms of the spectral sequence. To explain how, we need a geometric property of  Witt differentials. 
 \begin{lemma}\label{localcohomoofprojective}
Let $d \geq j$ be fixed non negative integers. Let $\mathcal{F}$ be one of the quasi-coherent $\wn\oo_{\mathbb{P}^d}$-modules $\wn\Omega^r_{\mathbb{P}^d}$ for any $r=0,\dots,d$, or $\wn\mathcal{L}$ associated to a line bundle $\mathcal{L}$ of $\mathbb{P}^d$. Then, 
\begin{enumerate}[label=\alph*)]
\item The local cohomology group sheaves $\mathcal{H}_{\mathbb{P}^j}^i(\mathcal{F})$ are trivial for any $i \not = d-j$. In particular,
\begin{equation}\label{triviallocal}
\h_{\mathbb{P}^j}^i(\mathbb{P}^d,\mathcal{F})=0, \quad \forall \text{ } i < d-j.
\end{equation}
\item $\h^i(\mathbb{P}^d \setminus \mathbb{P}^j,\mathcal{F})=0$ for any $i \geq d-j$.
\item $\h^i_{\mathbb{P}^j}(\mathbb{P}^d,\mathcal{F}) \simeq \h^i_{}(\mathbb{P}^d,\mathcal{F}),$ if $i>d-j$.
\end{enumerate}
 \end{lemma} 
 \begin{proof}
 \begin{enumerate}[label=\alph*)]
 \item The spectral sequence (cf. \cite[Theorem 2.6]{sga2})
 \begin{equation}
 E_2^{r,s}=\h^s(\mathbb{P}^d,\mathcal{H}_{\mathbb{P}^j}^r(\mathcal{F})) \implies \h_{\mathbb{P}^j}^{r+s}(\mathbb{P}^d,\mathcal{F}) 
 \end{equation}
 implies \eqref{triviallocal} from the triviality of $\mathcal{H}_{\mathbb{P}^j}^i(\mathcal{F})$ for $i \not= d-j.$ We do in detail the Witt differentials case, being the other one similar. Denote with $\mathcal{F}_n=\wn\Omega^r_{\mathbb{P}^d}$ for any $r=0,\dots,d$. Recall that from  \Cref{prop_exact_seq_gr}, we have a short  exact sequence,
 \begin{equation}\label{sesgr}
0 \rightarrow F_{X*}^{n+1} \frac{\Omega_{\mathbb{P}^d}^r}{B_n \Omega_{  \mathbb{P}^d}^r} \rightarrow \mathrm{gr}^n\mathcal{F}_{n+1}  \rightarrow F_{X*}^{n+1} \frac{\Omega_{\mathbb{P}^d}^{r-1}}{Z_n \Omega_{  \mathbb{P}^d}^{r-1}}\rightarrow 0.
\end{equation} 
where  $F_{X*}^{n+1} \frac{\Omega_{\mathbb{P}^d}^r}{B_n \Omega_{  \mathbb{P}^d}^r}$, $F_{X*}^{n+1} \frac{\Omega_{\mathbb{P}^d}^{r-1}}{Z_n \Omega_{  \mathbb{P}^d}^{r-1}}$ are locally free $\oo_{\mathbb{P}^d}$-modules of finite rank.  Since the result in a) holds in the case of coherent $\oo_{\mathbb{P}^d}$-modules (by arguments in \cite[Proposition 3.3 and Lemma 3.12]{sga2} ), it follows that $\mathcal{H}_{\mathbb{P}^j}^i(\mathrm{gr}^n\mathcal{F}_{n+1})=0$ if $i \not = d-j$,   by taking the long exact sequence associated to \eqref{sesgr}. The claim on $\mathcal{F}_n$  now follows by induction on $n$ for any $n \geq 1$,  by taking the associated long exact sequence  to
\begin{equation}\label{sesgrfil}
0 \rightarrow \mathrm{gr}^n\mathcal{F}_{n+1} \rightarrow \mathcal{F}_{n+1}\rightarrow \mathcal{F}_{n} \rightarrow 0.
\end{equation} 
  The case of Witt line bundles follows by analogy,  considering the short exact sequence \eqref{ses_witt_line_bundle}.
\item The corresponding result for coherent $\oo_{\mathbb{P}^d}$-modules $\mathcal{F}$ (i.e. when $n=1$) holds by computing the \v{C}ech cohomology for the covering $\mathcal{U}=\{D_+(z_r)\}_{j+1 \leq r \leq d}$ of $\mathbb{P}^d\backslash\mathbb{P}^j$. The resulting complex ${\check{C}}^{\bullet}(\mathcal{U},\mathcal{F})$ has degrees between $0$ and $d-j-1$. Therefore, for all $i \geq d-j$, the cohomology vanishes.
By considering the short exact sequence \eqref{sesgr} for  Hodge-Witt differentials, and \eqref{ses_witt_line_bundle} for Witt line bundles we see that b) follows by induction on $n$ (since for $n=1$ the vanishing holds).
\item By the long exact sequence associated to the couple $(\mathbb{P}^j,\mathbb{P}^d\backslash\mathbb{P}^j; \mathbb{P}^d)$, we have the exact sequence
\begin{equation}
    \h^{i-1}(\mathbb{P}^d\backslash\mathbb{P}^j,\mathcal{F}) \rightarrow\h^i_{\mathbb{P}^j}(\mathbb{P}^d,\mathcal{F}) \rightarrow \h^i(\mathbb{P}^d,\mathcal{F}) \rightarrow \h^i(\mathbb{P}^d\backslash\mathbb{P}^j,\mathcal{F}).
\end{equation}
Since $i-1 \geq d-j$, by part b) the outer terms of the above sequence are trivial, thus the map in the middle is an isomorphism.
 \end{enumerate}
 \end{proof}
 \begin{rem}
For a smooth $k$-scheme $X$ of dimension $d$, the canonical bundle $\omega_X=\ox(-d-1)$ and the $d$-th sheaf of differentials $\Omega_X^d$ agree. However, for $n>1$, in general $\wn\omega_X \neq \wn\Omega_{X}^d$.
\end{rem}

 For any $j=0,\dots,d$ the $E_1^{\bullet,j}$ terms of \eqref{ssext} have the property that $E_1^{-r,j}=0$ for any $r \geq j$: Indeed, if $I \subset \Delta=\{\alpha_0,\dots,\alpha_d\}$ is a subset of roots of $\G$, such that $|I|<d-j$, then also $i(I)<d-j$ and so $\h_{\mathbb{P}^{i(I)}}^j(\mathbb{P}^d,\mathcal{F})=0$ by the lemma above. 
 We wish to describe as explicitly as possible the $E_1$ page of the spectral sequence above. 
 \subsection{The generalized Steinberg modules over $\wn(k)$}\label{sec_generalized_steinb}
 For any  $I \subsetneq \Delta$, let us consider the $\zz$-module given by the following quotient:
 \begin{equation}
     {v}_{P_I}^G(\zz):=\mathrm{Ind}_{P_I}^G(\mathds{1}_{\zz})/(\sum_{\substack{G \supset Q \supsetneq P_I \\ Q \text{ parabolic sbgp}}} \mathrm{Ind}_{Q}^G(\mathds{1}_{\zz}))
 \end{equation}
 where $\mathds{1}_{\zz}$ denotes the ring $\zz$ as trivial $\zz[G]$-module. If $I=\emptyset$, $P_I=B$ and we denote ${v}_{B}^G(\zz)=:\mathrm{St}_G(\zz)$.
 The usual action of $G$ on the induction makes  any ${v}_{P_I}^G(\zz)$ a $\zz[G]$-module. Notice that ${v}_{P_I}^G(k)={v}_{P_I}^G(\zz) \otimes k$ gives us the generalized Steinberg representations of $G$. 
 This integral version already appears in \cite{ss91}, where in their setting $G$ is the group of points over a local field. In particular they consider the profinite topology on $G$ and its subgroups, consequently the induced representations they  study are smooth. Here there is no topology involved since we deal with a finite field, in particular $G$ is a finite group. Relating to their result, we can just consider our groups and subgroups equipped with the discrete topology, so that we can deduce the following properties.
 \begin{prop}[c.f. {\cite[Proposition 6.13]{ss91}}]\label{integralsteinberprop}
     For any $ I \subsetneq \Delta$, the integral generalized Steinberg modules ${v}_{P_I}^G(\zz)$ are finitely generated free $\zz$-modules. Moreover, for any $I$, there is a simplicial complex $\mathcal{T}_{\bullet}^{I}$, with the following properties: If $I= \emptyset$, $\mathcal{T}_{\bullet}^{I}$ is the combinatorial Tits building of $\mathrm{GL}_{d+1}(k)$; $\h^0(|\mathcal{T}_{\bullet}^{I}|, \zz)=\zz$, $\h^{d-1-|I|}(|\mathcal{T}_{\bullet}^{I}|, \zz)=v_{P_I}^G(\zz)$ and all other cohomology groups are trivial for any $I$.
    \end{prop}
    Before proceeding to the proof, we recall the following notations:
    For a Coxeter system $(W,S)$ and a subset $I \subset S$, the group $W_I \subset W$ is the subgroup of $W$ generated by the reflections associated to $I$.\\
    The set of reduced-$I$ elements of $W$, is the subset $W^I \subset W$ given by the representatives  $w \in W$ of the classes in the quotient $W/W_I$ such that $w$ has minimal length in the coset $wW_I$. (Every coset admits a unique reduced-$I$ element, cf. \cite[Lemma 3.2.1]{Digne_Michel_2020}).
    \begin{proof}
    The authors of the aforementioned Propositions  prove that such simplicial complex exists and it is acyclic. Then, by construction the simplicial integral cohomology has the desired properties. We recall the main point of the proof, specifying the stronger condition of using the discrete topology. Assume that $\Delta\backslash I=\{\alpha_{i_0},\dots, \alpha_{i_m}\}$. Let us consider the following simplicial sets: \\
    $Y_r^{I}:=$ simplicial set of $r+1$-tuples $(L_0,\dots,L_r)$ of lines in $k^{d+1}$ such that $\dim_{k} \sum_{i=1}^r L_i \leq j $ for some $j \in \{i_0+1,\dots, i_m+1\}$. \\
    $\mathcal{T}_r^{I}:=$ simplicial set of flags $(V_0 \subset \dots \subset V_r)$ of $k$-vector spaces in $k^{d+1}$ such that $\dim_{k}  V_i \in \{i_0+1,\dots, i_m+1\}$ for every $i=0,\dots,r$. \\
    $\mathcal{Z}_{rs}^{I}:=$ bisimplicial set of $(V_0 \subset \dots \subset V_r;L_0,\dots,L_s) \in \mathcal{T}_r^{I} \times Y_s^{I}  $  such that $\sum_{i=1}^s  L_i \subset V_0$. \\
        The face and degeneracy maps are given respectively by removing or doubling a vector space.
    We introduce also the following simplicial set: \\
     $\mathcal{NT}_r^{I}:=$ simplicial set of flags $(V_0 \subsetneq \dots \subsetneq V_r)$ of $k$-vector spaces in $k^{d+1}$ such that $\dim_{k}  V_i \in \{i_0+1,\dots, i_m+1\}$ for every $i=0,\dots,r$.\\ 
   $\mathcal{NT}_r^{I}$ is said to be the "normalization" of $\mathcal{T}_r^{I}$, where all the flags are assumed to not have repeated vector spaces. Notice that since $k$ is finite, the sets of vertices of the simplicial sets above are finite. 
    In particular,  for them the profinite topology coincides with the discrete topology.  The constant abelian sheaf $\underline{\zz}$ on any of this (discrete) simplicial sets assigns to any (finite) subset $U$, the corresponding abelian group $C(U,\zz)$ generated by all (set theoretical) maps $U \rightarrow \zz$.
     By a cosimplicial normalization theorem (see loc.cit. proof of Proposition 3.6), the natural inclusion
     $C(\mathcal{NT}^I_{\bullet}, \zz)\rightarrow C(\mathcal{T}^{I}_{\bullet},\zz)$ is a homotopic equivalence. Since for $r \geq d-|I|$, $\mathcal{NT}^I_r= \emptyset$, then $\h^{r}(|\mathcal{T}_{\bullet}^I|,\zz)=0$ for $r \geq d-|I|$.
      By loc. cit. Lemma 3.3 and Lemma 3.4,
    the natural maps $\mathcal{Z}_{\bullet, s}^{I} \xrightarrow{f_{\bullet}} Y_s^{I}$ and $ \mathcal{Z}_{r,\bullet}^{I} \xrightarrow{g_{\bullet}} \mathcal{T}_r^{I}$ induce respectively quasi-isomorphism of complexes $C(\mathcal{Z}_{\bullet, s}^{I}, \zz) \leftarrow C(Y_s^{I},\zz)$ and
    $C(\mathcal{Z}_{r,\bullet}^{I},\zz) \leftarrow C(\mathcal{T}_r^{I}, \zz)$ for any $r,s$\footnote{Here it is a sketched argument: the sheaves $\zz_{Y_{\bullet}^I}$ (resp. $\zz_{\mathcal{T}^I_{\bullet}}$) and $\zz_{\mathcal{Z}_{\bullet,\bullet}^I}$ are flasque and $f_{\bullet}$ (resp. $g_{\bullet}$) induces an acyclic resolution $\zz_{Y_{\bullet}^I} \rightarrow f_{\bullet,*}\zz_{\mathcal{Z}_{\bullet,\bullet}^I}$ of $\zz_{Y_{\bullet}^I}$ (resp. $\zz_{\mathcal{T}^I_{\bullet}}$): it is enough to check it on stalk, and then apply loc.cit. Lemma 3.3; Since flasque sheaves are global section-acyclic, the claim follows.}. In particular the cohomology of the total complex of $\mathcal{Z}_{\bullet,\bullet}^{I}$, computes the simplicial cohomology of $\mathcal{T}_{\bullet}^{I}$. Therefore, we have a second \footnote{I.e. whose filtration on the total complex is given by removing successive rows.} spectral sequence that read as
    \begin{equation}
        E_1^{r,s}:=h^s(C(\mathcal{Z}_{\bullet,r}^{I}, \zz)) \implies \h^{r+s}(|\mathcal{T}_{\bullet}^{I}|,\zz).
    \end{equation}
    Next step is proving that the following sequence induced by $f_{\bullet}$
    \begin{equation}
        0 \rightarrow C(Y_r^I,\zz) \rightarrow C(\mathcal{Z}^{I}_{0,r},\zz)\rightarrow \dots  \rightarrow C(\mathcal{Z}^{I}_{d-1-|I|-r,r},\zz)
    \end{equation}
    is exact.  With a bit of work that we here omit, it follows essentially by loc. cit. Lemma 3.3 again. Therefore, we deduce that $E_1^{r,s}=0$ if $r+s <d-|I|-1$ and $s>0$, $E_1^{r,0}=C(Y_r^I,\zz)$ for $r< d-1-|I|$, and $C(Y_{d-1-|I|}^{I},\zz) \subset E_1^{d-1-|I|,0} $. Furthermore, by loc. cit. Lemma 3.3 the complex
    \begin{equation}
        0 \rightarrow \zz \rightarrow C(Y_0^I,\zz) \rightarrow \dots \rightarrow C(Y^{I}_{d-1-|I|},\zz)
    \end{equation} is exact, therefore also $0 \rightarrow \zz \rightarrow E_1^{\bullet,0}$ for $\bullet \leq d-1-|I|$ it is so. It implies $E_2^{0,0}=\zz$ and $E_2^{r,s}=0$ for $0<r+s< d-1-|I|$, from which the vanishing result follows. To compute the highest cohomology group we consider the normalized simplicial complex (by homotopic invariance of simplicial cohomology). For any $J \subset I$, consider the flag 
    \begin{equation}
        \tau_J=(\sum_{i=0}^{j_0} k e_i \subsetneq \dots \subsetneq \sum_{i=0}^{j_r} k e_i )
    \end{equation}
    where ${e_i}$ for $i=0,\dots,d$ is the standard basis of $k^{d+1}$ and $\Delta \backslash J=\{\alpha_{j_0},\dots,\alpha_{j_r}\}$ with $j_0 < \dots < j_r$. Then the parabolic subgroup $P_J \subset G$ is the stabilizer of $\tau_J$. Moreover, the natural map 
    \begin{equation}
        \bigsqcup_{\substack{I \subset J \subset \Delta \\ |J|=d-1-r}} G/P_J \rightarrow \mathcal{NT}_r^{I}, \quad gP_j \mapsto g.\tau_J
     \end{equation}
     is a bijection. Hence, $$C(\mathcal{NT}_{r}^{I}, \zz)= \bigoplus_{\substack{I \subset J \subset \Delta \\ |J|=d-1-r}} \mathrm{Ind}_{P_J}^{G}(\mathds{1}_{\zz}).$$
     The cohomology group $\h^{d-1-|I|}(|\mathcal{T}_{\bullet}^{I}|,\zz)$ is equal to 
     $$\mathrm{coker}\Big(C(\mathcal{NT}_{d-2-|I|}^{I}, \zz) \rightarrow C(\mathcal{NT}_{d-1-|I|}^{I}, \zz)\Big)=\mathrm{coker}\Big(\bigoplus_{\alpha \in \Delta\backslash I}\ind_{P_{I \cup {\alpha}}}^{G}(\mathds{1}) \rightarrow \ind_{P_{I}}^{G}(\mathds{1})\Big)=v_{P_I}^{G}(\zz).$$
     To check the $\zz$-freeness of $v_{P_I}^G$, we prove that it has a finite descending filtration whose successive quotient are $\zz$-free. Let $W$ be the Weyl group of $G$ and let $W^I    \subset W$ be the subset of reduced-$I$ elements. Notice that for $I \subset J$, $W^J \subset W^I$. The Bruhat decomposition yields the equality
     \begin{equation}
         G/P_I=\bigsqcup_{w \in W^I} BwP_I/P_I
     \end{equation}
     and the surjections $C(w):=BwB/B \rightarrow BwP_I/P_I=:C_I(w)$ are actually bijections for any $w \in W^I$ (compare cf. \cite[Proposition 3.16 (ii)]{BorTit72} and \Cref{rem_iso_chamber_finite}). Fix an order on $W^I=\{w_1,\dots,w_m\}$ such that $a \leq b$ if and only if $l(w_a) \leq l(w_b)$. Denote by $F_I^r:=\{f \in \ind_{P_I}^G(\mathds{1}) : f(C_I(w_s))=0 \text{, } 1 \leq s \leq r\}$. Then, $\{0\}=F_I^m \subset F_I^{m-1} \subset \dots \subset F_I^0:=\ind_{P_I}^G(\mathds{1})$ has the property that for any $a < b$
     \begin{equation}
         F_I^a/F_I^b= C(\bigsqcup_{a <s \leq b} C_I(w_s),\zz)
     \end{equation}
     Indeed, any coset $f + F_I^b$ with $f \in F_I^a$ can be represented by an $f_{a,b}$ supported in $\bigsqcup_{a < s \leq b}  C_I(w_s)$: indeed, by taking $f_{a,b}$ be  the extension by $0$  of $f_{\mid \bigsqcup_{a< s \leq b}  C_I(w_s) }$, it follows that $f+F_I^b=f_{a,b}+(f-f_{a,b})+F_I^b=f_{a,b}+F_I^b$. Viceversa, the set $\{f_{a,b} \in C(\bigsqcup_{a < s \leq b} C_I(w_s),\zz) : f_{a,b}=0 \} \cap F_I^a\subset F_I^b$. The filtration $F_I^{\bullet}$
     induces a finite filtration $\Bar{F}_I^{\bullet}$ on the quotient $v_{P_I}^G$. We distinguish between the case where $w_r \in W^{I\cup \{\alpha\}}$ for some $\alpha \in \Delta\backslash I$ or not. In the first case we have a natural bijection 
     \begin{equation}
         C_{I}(w_r) \xrightarrow{\sim} C_{I\cup \{\alpha\}}(w_r)
     \end{equation}
     inducing an isomorphism
     \begin{equation}
         C(C_{I\cup \{\alpha\}}(w_r),\zz) \rightarrow C( C_{I}(w_r),\zz).
     \end{equation}
     If $f \in F_I^{a-1}$ then $f_{\mid C_I(w_r)} \in C( C_{I}(w_r), \zz) $ corresponds to a unique $\tilde{f} \in C(C_{I\cup \{\alpha\}}(w_r),\zz)$. Then, if $h \in \ind_{P_{I\cup\{\alpha\}}}^G(\mathds{1})$ denotes the extension by $0$ of $\tilde{f}$, by construction it follows that $f-h \in F_I^{r}$, thus we have the following equality of cosets 
     $$\bar{F}_I^{r-1} \ni f + \sum_{\alpha \in \Delta \backslash I}\ind_{P_{I\cup\{\alpha\}}}^G(\mathds{1})=(f-h)+\sum_{\alpha \in \Delta \backslash I}\ind_{P_{I\cup\{\alpha\}}}^G(\mathds{1}) \in \bar{F}_I^r.$$
     Therefore, $\bar{F}_I^{r-1}=\bar{F}_I^{r}$. Now, suppose  $w_r \not \in \bigcup_{\alpha \in \Delta\backslash I} W^{I\cup\{\alpha\}}$. We claim that $\bar{F}_I^{r-1}/\bar{F}_I^{r}={F}_I^{r-1}/{F}_I^{r}$. It is equivalent to check the relation
     $F_I^{r-1} \cap \sum_{\alpha \in \Delta \backslash I}\ind_{P_{I\cup\{\alpha\}}}^G(\mathds{1}) \subset F_I^{r}$.
     If $f \in F_I^{r-1}$ is such that $f=\sum_{\alpha \in \Delta \backslash I}f_{\alpha}$ with $f_{\alpha} \in \ind_{P_{I\cup\{\alpha\}}}^G(\mathds{1})$ then we can find a writing of $f=\sum_{\alpha \in \Delta\backslash I} g_{\alpha}$ with $g_{\alpha} \in F_I^{r-1} \cap \ind_{P_{I\cup\{\alpha\}}}^G(\mathds{1}) $. By induction we can suppose $f_{\alpha} \in F_I^{r-2}$. If at most one $\alpha$ exists such that $w_{r-1} \in W^{I\cup\{\alpha\}}$, the latter assertion can be proved with the same argument as in loc. cit. Prop. 4.4. Otherwise, if there exist different $\alpha \neq \beta$ such that $w_{r-1} \in W^{I\cup\{\alpha,\beta\}}$, as before we have a bijection 
     \begin{equation}
         C_{I \cup \{\alpha\}}(w_{r-1}) \xrightarrow{\sim}  C_{I \cup \{\alpha,\beta\}}(w_{r-1})
     \end{equation}
     and similarly to the argument above, we can find an $h \in F_I^{r-2} \cap \ind_{P_{I\cup\{\alpha,\beta\}}}^G(\mathds{1})$ such that $f_{\alpha}-h \in F_I^{r-1}$. Then we get a rewriting of $f$ replacing respectively $f_{\alpha}, f_{\beta}$ by  $g_{\alpha}:=f_{\alpha}-h$ and  $g_{\beta}:=f_{\beta}+h$ (for other $\gamma \in \Delta\backslash (I \cup\{\alpha,\beta\})$, let $g_{\gamma}:=f_{\gamma}$). Then the claim follows inductively. Lastly, notice that the condition on $w_{r}$ implies that for any $\alpha \in \Delta\backslash I$ there is some $s(\alpha)<r$ for which $C_{I \cup \{\alpha\}}(w_r)=C_{I \cup \{\alpha\}}(w_{s(\alpha)})$ : when $w_r \not \in W^{I \cup \{\alpha\}}$, then there exists some element in $(I \cup \{\alpha\}) \backslash I=\{\alpha\}$, i.e. $\alpha$ itself, such that $l(w_rs_{\alpha}) \leq l(w_r)$ where $s_{\alpha}$ is the simple reflection associated to $\alpha$. In particular the permutation $w_s$ (with $s \neq r$) of minimal length in $w_rW_{I \cup \{\alpha\}}$ must satisfy $l(w_s) \leq l(w_r)$, thus $s(\alpha)=s < r$.  Therefore, $F_I^{r-1} \cap \ind_{P_{I\cup\{\alpha\}}}^G(\mathds{1}) \subset F_I^{r}$. 
    This show that the successive quotients of $\bar{F}_I^{\bullet}$  are of the form $C(U,\zz)$  for some (finite) subset $U \subset G/P_I$, thus they are finitely generated free over $\zz$.   
    \end{proof}
    \begin{rem}\label{rem_iso_chamber_finite}
        If $F\colon \textbf{G}_{\bar{k}} \rightarrow \textbf{G}_{\bar{k}}$ denotes the standard geometric Frobenius on the ${\bar{k}}$-scheme $\G_{\bar{k}}$, then $\G_{\bar{k}}^F=G, \textbf{P}_{I,{\bar{k}}}^F=P_I, \textbf{B}_{\bar{k}}^F=B$. In the part of the proof related to verify $\zz$-freeness, we used the Bruhat decomposition associated to the $(B,N)$-pair of $\G_{\bar{k}}^F$, induced by that one on $\G_{\bar{k}}$. Then, the Weyl group we consider is $W=W_{\bar{k}}^F$, where $W_{\bar{k}}$ is the Weyl group for $\G_{\bar{k}}$.  For any $w \in W^{I}$ the bijection $C(w) \rightarrow C_I(w)$ then follows by the analogous statement for the schemes $\G_{\bar{k}}/\textbf{P}_{I,{\bar{k}}}$ (cf. \cite[Proposition 3.16 (ii)]{BorTit72}) and by the bijection $(\G_{\bar{k}}/\textbf{P}_{I,{\bar{k}}})^F=G/P_I$ (i.e comparing the two Bruhat decompositions, the cardinality of  $BwP_I/P_I$ and  $(\textbf{B}_{\bar{k}}w\textbf{P}_{I,{\bar{k}}}/\textbf{P}_{I,{\bar{k}}})^F$, we can deduce $|C_I(w)|=|(\textbf{B}_{\bar{k}}w\textbf{P}_{I,{\bar{k}}}/\textbf{P}_{I,{\bar{k}}})^F|$). 
    \end{rem}
   
    Hence by the proof of \Cref{integralsteinberprop}, it follows that the following augmented $G$-equivariant complex is acyclic:
    \begin{equation}\label{intgralsteinbcomplex}
        0 \rightarrow \zz \rightarrow \bigoplus_{\substack{I \subset J \subset \Delta \\ |J|=d-1}} \mathrm{Ind}_{P_J}^{G}(\mathds{1}_{\zz}) \rightarrow \dots \rightarrow  \mathrm{Ind}_{P_I}^{G}(\mathds{1}_{\zz}) \rightarrow {v}_{P_I}^G(\zz) \rightarrow 0.
    \end{equation}
    
   Set by definition ${\leftindex_{n} {v}}_{P_I}^G:={v}_{P_I}^G({\zz}) \otimes \wn(k)$. When $I=\emptyset$, and $G=\mathrm{GL}_j$, then ${\leftindex_{n} {v}}_{B}^{\mathrm{GL}_j}=:{\leftindex_{n} {\mathrm{St}}}_j$ is the Steinberg representation of $\mathrm{GL}_j$ over $\wn(k)$.
    \begin{corol}\label{corolacyclicsteinbwn}
        The following complex of $\wn(k)[G]$-modules is acyclic:
        \begin{equation}\label{wnsteinbcomplex}
        0 \rightarrow \wn(k) \rightarrow \bigoplus_{\substack{I \subset J \subset \Delta \\ |J|=d-1}} \mathrm{Ind}_{P_J}^{G}(\mathds{1}_{\wn(k)}) \rightarrow \dots \rightarrow  \mathrm{Ind}_{P_I}^{G}(\mathds{1}_{\wn(k)}) \rightarrow {\leftindex_n {v}}_{P_I}^G \rightarrow 0.
    \end{equation}
    \end{corol}
    \begin{proof}
  The complex \eqref{wnsteinbcomplex} is obtained by tensoring the complex \eqref{intgralsteinbcomplex} (namely $K^{\bullet}$) with $\wn(k)$.  The restriction of the functor $-\otimes_\zz \wn(k)$ on the full subcategory of projective $\zz$-modules is exact. Since $K^{\bullet}$ is an acyclic complex of free modules, it follows $K^{\bullet} \otimes_\zz \wn(k)$ is an acyclic complex of free $\wn(k)$-modules.
    \end{proof}
    \subsection{Computation of the $E_2$-page}
    In this section we will prove that similarly to the case of coherent $\mathcal{O}_{\mathbb{P}^d}$-modules cohomology \cite[cf. Theorem 2.1.2.1]{Ku16}, the computation of Hodge-Witt cohomology of the Drinfeld's upper half space over $k$, as well as the cohomology of Witt line bundles, can be described in terms of $\wn(k)[G]$-modules, given by certain local cohomology groups. This is reached by evaluating the $E_2$-page of \eqref{ssext}. When $\mathcal{L}=\ox(D)$ is a Witt line bundle for some Cartier Divisor $D$ on a $k$-scheme $X$, we denote $p\mathcal{L}:=\ox(pD)$.

    We recall a property of projective finitely generated modules over a ring $R$:
    \begin{lemma}
        Let $R$ be a commutative ring and $P,P'$ be $R$-modules. Assume that at least one between $P$ or $P'$ is finitely generated projective over $R$. Then,
        \begin{equation}
            \mathrm{Hom}_R(P \otimes_R P',R) \simeq \mathrm{Hom}_R(P,R ) \otimes_R \mathrm{Hom}_R(P',R)
        \end{equation}
        is a canonical isomorphism of $R$-modules
    \end{lemma}
    \begin{proof}
        Without loss of generality, we can assume $P$ be finitely generated projective. Then, $\mathrm{Hom}_R(P,-)$ is an exact endofunctor of $R$-mod and since $P$ is also finitely generated, there is a canonical bijection, functorial on $R$-modules $Q$ (cf. \cite[II, 4.2, Proposition 2 (i)]{Bou98}): 
        \begin{equation}
            \mathrm{Hom}_R(P,Q) \simeq \mathrm{Hom}_R(P,R) \otimes_R Q.
        \end{equation}
        
        Moreover, the functor $-\otimes_R P'$ is left adjoint to $\mathrm{Hom}_R(P', -)$. Therefore, we have the following natural identifications:
        \begin{equation*}
            \mathrm{Hom}_R(P \otimes_R P',R) \simeq \mathrm{Hom}_R(P,\mathrm{Hom}_R(P',R) ) \simeq \mathrm{Hom}_R(P,R) \otimes_R \mathrm{Hom}_R(P',R). \qedhere
        \end{equation*}
       \end{proof}
     We have the following: 
 \begin{prop} \label{prop_red_problem}
With the same notation of  \Cref{localcohomoofprojective}, the spectral sequence \eqref{ssext} degenerates at $E_2$. Moreover,
 \begin{equation}
 E_2^{0,j}=\h^j(\mathbb{P}^d,\mathcal{F}) \quad j \geq 2,
 \end{equation}
 and the terms $E_2^{-j+1,j}$ for $j \geq 1$ appear as an extension of certain $\wn(k)[G]$-modules:
 \begin{equation}\label{exte2term}
 0 \rightarrow E_{2,\sim}^{-j+1,j} \rightarrow E_{2}^{-j+1,j} \rightarrow E_{2,\mathrm{w.s.}}^{-j+1,j} \rightarrow 0,
 \end{equation}
 where, the following equality hold:
 \begin{equation}
 E_{2,\sim}^{-j+1,j}= \mathrm{Ind}^G_{P_{(d+1-j,j)}}(\tilde{\h}_{\mathbb{P}^{d-j}}^j(\mathbb{P}^d,\mathcal{F}) \otimes_{\wn(k)} {\leftindex_n {\mathrm{St}}}_j^{\vee})
 \end{equation}
 \begin{equation}
 E_{2,\mathrm{w.s.}}^{-j+1,j}=\h_{}^j(\mathbb{P}^d,\mathcal{F})\otimes_{\wn(k)} ({\leftindex_n {\upsilon}}^{G}_{P_{(d+1-j,1^{j})}})^{\vee},
 \end{equation}
 for any $1\leq j \leq d$, and finally
 \begin{equation}
 E_2^{0,1}=E_1^{0,1}=\mathrm{Ind}^G_{P_{(d,1)}}\h^1_{\mathbb{P}^{d-1}}(\mathbb{P}^d,\mathcal{F}).
 \end{equation}
 \end{prop}
 

 \begin{proof}
 We will prove that for any $j=1,\dots,d$, $E_1^{\bullet,j}$ defines an exact sequence of modules. Define $\Delta_j$ be the set of all subsets $I \subset \Delta$ such that $\alpha_0,\dots,\alpha_{d-j-1} \in I$, and $\alpha_{d-j} \not \in I $. When $j=d$, $\alpha_0 \not \in I$ is the only condition. By \Cref{localcohomoofprojective} a), if $i(I) < d-j$, then 
 \begin{equation}
    \h^j_{gY_I}(\mathbb{P}^d,\mathcal{F})=\h^j_{\mathbb{P}^{i(I)}}(\mathbb{P}^d,\mathcal{F})=0.
 \end{equation}
 It follows that we can write any $E_1^{\bullet,j}$ as
 \begin{equation}
     E_1^{\bullet,j}=\bigoplus_{\substack{I \subset \Delta\\ |I|=d-1+\bullet \\ i(I)=d-j}} \bigoplus_{g \in G/P_I}\h^j_{gY_I}(\mathbb{P}^d,\mathcal{F}) \oplus \bigoplus_{\substack{I \subset \Delta\\ |I|=d-1+\bullet \\ i(I)>d-j}} \bigoplus_{g \in G/P_I}\h^j_{gY_I}(\mathbb{P}^d,\mathcal{F}).
 \end{equation}
 The condition $i(I)=d-j$ is equivalent to $I \in \Delta_j$. Also, $i(I)>d-j$ is equivalent to $\alpha_0,\dots,\alpha_{d-j} \in I$.
 Then, for any $G$-equivariant quasi-coherent $\wn\oo_{\mathbb{P}^d}$-module $\mathcal{F}$, define
 \begin{equation}
 {}^nE^{\bullet,j}_{1,\sim}(\mathcal{F}):=\bigoplus_{\substack{ I \in \Delta_j\\|I|=d-1+\bullet}} \ind^{G}_{P_I}\tilde{\h}_{\mathbb{P}^{d-j}}^j(\mathbb{P}^d,\mathcal{F}),
 \end{equation}
 and 
 \begin{equation}
 {}^nE^{\bullet,j}_{1,\mathrm{w.s.}}(\mathcal{F}):=\bigoplus_{\substack{ I \subset \Delta\\|I|=d-1+\bullet \\ \alpha_0,\dots,\alpha_{d-j-1} \in I}} \ind^{G}_{P_I}\h^j(\mathbb{P}^d,\mathcal{F}) \quad (\alpha_0 \in I \text{ if } j=d.)
 \end{equation}
 be complexes where the differentials are induced by that one of $E_1^{\bullet,j}$. By construction, they satisfy a short exact sequence (of complexes) of $\wn(k)$-modules,
 \begin{equation}\label{sestildews}
 0 \rightarrow {}^nE_{1,\sim}^{\bullet,j} \rightarrow E_{1}^{\bullet,j} \rightarrow {}^nE_{1,\mathrm{w.s.}}^{\bullet,j} \rightarrow 0.
 \end{equation}
 Indeed, it is induced by the following short exact sequences of complexes:
 \begin{align*}
   &  0 \rightarrow \ind_{P_I}^{G}\tilde{\h}^j_{\mathbb{P}^{d-j}}(\mathbb{P}^d,\mathcal{F})\rightarrow \ind_{P_I}^{G}{\h}^j_{\mathbb{P}^{d-j}}(\mathbb{P}^d,\mathcal{F})  \rightarrow \ind_{P_I}^{G}{\h}^j_{}(\mathbb{P}^d,\mathcal{F})\rightarrow 0,   \text{ if } I \in \Delta_j\\
  &  0 \rightarrow 0 \rightarrow \ind_{P_I}^{G}{\h}^j_{}(\mathbb{P}^d,\mathcal{F})  \rightarrow \ind_{P_I}^{G}{\h}^j_{}(\mathbb{P}^d,\mathcal{F})\rightarrow 0,   \text{ if } \alpha_0,\dots,\alpha_{d-j} \in I \subset \Delta  \\
   & (0),  \text{ otherwise } 
 \end{align*}
 It suffices to prove exactness for ${}^nE_{1,\sim}^{\bullet,j}$ and ${}^nE_{1,\mathrm{w.s.}}^{\bullet,j}$. \footnote{cf. \cite[Proposition 2.2.4]{Orl08}, where the analogous property is fulfilled in the case of representations over a field ($n=1$ here) in characteristic $0$.}
 
 To start, we first claim that we have the following equality:
 \begin{equation}\label{tensorE_1sim}
     {}^nE_{1,\sim}^{\bullet,j}= \mathrm{Ind}_{P_{(d+1-j,j)}}^{G}\big( \bigoplus_{\substack{ I \subset \Delta_{GL_j}\\|I|=j-1+\bullet }}\ind^{\mathrm{GL}_j}_{P_I}(\mathds{1})\otimes \tilde{\h}_{\mathbb{P}^{d-j}}^j(\mathbb{P}^d,\mathcal{F}) \big)
 \end{equation}
 and 
  \begin{equation}\label{tensorE_1w.s.}
     {}^nE_{1,\text{w.s.}}^{\bullet,j}= \Big(\bigoplus_{\substack{ I \subset \Delta\\|I|=d-1+\bullet  \\ \alpha_0,\dots,\alpha_{d-j-1} \in I}}\ind^{G}_{P_I}(\mathds{1})\Big)\otimes {\h}_{}^j(\mathbb{P}^d,\mathcal{F}) 
 \end{equation}
 Moreover, if $ \mathcal{F}$ is any equivariant $\wn\oo_{\mathbb{P}^d}$-module, then the complexes above define exact sequence of $\wn(k)[G]$-modules. 
 
 Assume the equality \eqref{tensorE_1sim} and \eqref{tensorE_1w.s.} are verified. 
 Taking the tensor product of complexes with $ {\h}_{}^j(\mathbb{P}^d,\mathcal{F})$ and with $ \tilde{\h}_{\mathbb{P}^{d-j}}^j(\mathbb{P}^d,\mathcal{F})$ does not in general defines exact functors. However, by \Cref{corolacyclicsteinbwn},  the complex appearing inside $\mathrm{Ind}_{P_{(d+1-j,j)}}^{G}$ of \eqref{tensorE_1sim} defines an exact sequence of free $\wn(k)[P_{(d+1-j,j)}]$-modules. Further, the parabolic induction is an exact functor, therefore the sequence \eqref{tensorE_1sim} stays exact. In the relation \eqref{tensorE_1w.s.}, note that the index set of $I$ has the following property: for any such $I$ containing $\alpha_0,\dots,\alpha_{d-j-1}$, then $P_I \supset P_{(d+1-j,1^j)}$. Viceversa, any parabolic $Q$ such that $Q \supset P_{(d+1-j,1^j)} $ is of the form $P_I$ for some $I$ of such form. Hence, ${}^nE_{1,w.s}^{\bullet,j}$ coincides with the complex $C(\mathcal{NT}^{\{\alpha_0,\dots,\alpha_{d-j-1}\}}_{\bullet},\wn(k))$ tensor with ${\h}_{}^j(\mathbb{P}^d,\mathcal{F}) $. Therefore, exactness again follows from  \Cref{corolacyclicsteinbwn}. This proves the exactness property for $G$-equivariant (quasi-coherent) $\wn\oo_{\mathbb{P}^d}$-modules.
 
  The equality \eqref{tensorE_1w.s.} follows by definition, since ${\h}_{}^j(\mathbb{P}^d,\mathcal{F}) $ is a $G$-module. \\
Consider $I \in \Delta_j$. Then $P_{(d+1-j,1^j)}\subset P_I \subset P_{(d+1-j,j)}$ holds true and by transitivity of parabolic induction, we have $\ind_{P_I}^G=\ind_{P_{(d+1-j,j)}}^G\ind_{P_I}^{P_{(d+1-j,j)}}$ and
 $\ind_{P_{(d+1-j,1^j)}}^G=\ind_{P_{(d+1-j,j)}}^G\ind_{P_{(d+1-j,1^j)}}^{P_{(d+1-j,j)}}$. Moreover, the natural identifications of the quotients
 $$P_{(d+1-j,j)}/P_{(d+1-j,1^j)} \sim \mathrm{GL}_j/{(B \cap \mathrm{GL}_j)} $$ and 
$$P_{(d+1-j,j)}/P_{I} \sim \mathrm{GL}_j/{(P_I \cap \mathrm{GL}_j)}$$
induce for any $k[P_I]$-module $M$ (resp.  $k[P_{(d+1-j,1^j)}]$-module $M'$) isomorphism of representations  
\begin{equation} \label{identityparabind}
\ind_{P_I}^{P_{(d+1-j,j)}}(M)\simeq \ind_{P_I \cap \mathrm{GL}_j}^{\mathrm{GL}_j}(M), \qquad \ind_{P_{(d+1-j,1^j)}}^{P_{(d+1-j,j)}}(M') \simeq \ind_{B \cap \mathrm{GL}_j}^{\mathrm{GL}_j}(M').
\end{equation}
 Note that if $I=\{\alpha_0,\dots,\alpha_{d-j-1},\alpha_{i_0},\dots,\alpha_{i_{r+j-2}}\}$, for some $r$, then $P_I \cap \mathrm{GL}_j$ is the parabolic subgroup $Q_{\tilde{I}}$ of $GL_j \subset L_{(d+1-j,j)}$ associated to $\tilde{I}=\{\beta_{i_0-d+j+1},\dots,\beta_{i_{r+j-2}-d+j+1}\}$, where $\Delta_{GL_j}:=\{\beta_0,\dots,\beta_{j-1}\}$ is the usual system of simple roots for $GL_j$. Observing that $\h^j_{\mathbb{P}^{d-j}}(\mathbb{P}^d,\mathcal{F})$ is a $P_{(d+1-j,j)}$-module, then the relation \eqref{tensorE_1sim} readily follows. Moreover, the isomorphisms \eqref{identityparabind} do not depend on the base ring of the representations, in particular they hold true for any (equivariant) $\wn(k)$-module. 
   
   The degeneration is immediate. Indeed, $d_2^{-r,j}: E_2^{-r,j} \rightarrow E_2^{-r+2,j-1}$ is always the $0$ map: Since $E_1^{\bullet,j}$ defines  exact sequences for any $j=1,\dots,d$,  it follows that both terms $E_2^{-r,j}, E_2^{-r+2,j-1}$ are $0$ whenever $0 < r < j-1$; when $r=0$, or $r=j-1$ at least one of the terms is $0$.
 Then we get that 
 \begin{equation}
 \mathrm{gr}^{\bullet}\h^{s}_{Y}(\mathbb{P}^d,\mathcal{F})= \bigoplus_{j-r=s}E_2^{-r,j}.
 \end{equation}
 In particular, when $s > 1$, it follows
 \begin{equation}
 \mathrm{gr}^{\bullet}\h^{s}_{Y}(\mathbb{P}^d,\mathcal{F})= \bigoplus_{j \geq 1}E_2^{s-j,j}=E_2^{0,s}.
 \end{equation}
 This implies that 
 \begin{equation}
 E_2^{0,s}=\h^{s}_{Y}(\mathbb{P}^d,\mathcal{F}) \simeq \h^{s}(\mathbb{P}^d,\mathcal{F}) \quad (s > 1)
 \end{equation}
 where the last canonical isomorphism follows by the long exact sequence of the couple $(X,Y; \mathbb{P}^d)$.
 Notice that for $j \geq 1$, the map $d_1^{-j+1,j}$ induces a morphism of complexes relative to \eqref{sestildews} for $\bullet=-j+1, -j+2$. Applying the snake Lemma, we get the exact sequence \eqref{exte2term}, where 
 $$E^{-j+1,j}_{2,?}:=\ker\bigg({}^nE^{-j+1,j}_{1,?} \rightarrow {}^nE^{-j+2,j}_{1,?}\bigg), \quad ? \in \{\mathrm{w.s. },\sim\}.$$
  Also, notice that the collection of $\{P_I \cap \mathrm{GL}_j \mid I \in \Delta_j\}$ are the parabolic subgroups containing $B \cap \mathrm{GL}_j$, and those such that $|I|=d-j+1$, are the minimal ones.
So,
\begin{equation}\label{stein}
\mathrm{coker}\bigg(\bigoplus_{\substack{ I \in \Delta_j\\|I|=d+1-j}}\ind_{P_I \cap \mathrm{GL}_j}^{\mathrm{GL}_j}(\mathds{1}) \rightarrow \ind_{B \cap \mathrm{GL}_j}^{\mathrm{GL}_j}(\mathds{1})\bigg)={\leftindex_n {\mathrm{St}}}_j.
\end{equation}
 Since $\tilde{\h}_{\mathbb{P}^{d-j}}^j(\mathbb{P}^d,\mathcal{F})$ is a $P_{(d+1-j,j)}$-module, 
we see that
\begin{equation*}
    \begin{split}
 (E^{-j+1,j}_{2,\sim})^{\vee}
  & \simeq \mathrm{coker}\bigg(({}^nE^{-j+2,j}_{1,\sim})^{\vee} \rightarrow ({}^nE^{-j+1,j}_{1,\sim})^{\vee}\bigg)\\
&\simeq \left. \tilde{\h}_{\mathbb{P}^{d-j}}^j(\mathbb{P}^d,\mathcal{F})^{\vee} \otimes  \ind^{G}_{P_{(d-j+1,1^j)}} (\mathds{1}) \middle/ \Big(\tilde{\h}_{\mathbb{P}^{d-j}}^j(\mathbb{P}^d,\mathcal{F})^{\vee} \otimes \sum_{\substack{ I \in \Delta_j\\|I|=d-j+1}} \ind^{G}_{P_I} (\mathds{1})\Big)\right.  \\
 &  \simeq \mathrm{Ind}^G_{P_{(d+1-j,j)}}(\tilde{\h}_{\mathbb{P}^{d-j}}^j(\mathbb{P}^d,\mathcal{F})^{\vee} \otimes_{\wn(k)}{\leftindex_n {\mathrm{St}}}_j)
   \end{split}
\end{equation*}

implies
\begin{equation}
\mathrm{Ind}^G_{P_{(d+1-j,j)}}(\tilde{\h}_{\mathbb{P}^{d-j}}^j(\mathbb{P}^d,\mathcal{F}) \otimes_{\wn(k)} {{\leftindex_n {\mathrm{St}}}_j}^{\vee})=E^{-j+1,j}_{2,\sim}.
\end{equation}
In the last isomorphism, we use that the Steinberg module is a finitely generated free module over $\wn(k)$,  thus it is compatible with tensor product. For $E^{-j+1,j}_{2,\mathrm{w.s.}}$ a similar argument holds.  Notice that the cohomology  ${\h}_{}^j(\mathbb{P}^d,\mathcal{F})$ 
is a $G$-module (so also the dual is), and  then 
\begin{equation}
(E^{-j+1,j}_{2,\mathrm{w.s.}})^{\vee}=\mathrm{coker}\bigg( \bigoplus_{\substack{ I \subset \Delta\\|I|=d-j+1 \\ \alpha_0,\dots,\alpha_{d-j-1} \in I}} \ind^{G}_{P_I}(\mathds{1})    \rightarrow  \ind^{G}_{P_{(d-j+1,1^j)}}(\mathds{1})\bigg)\otimes ({\h}_{}^j(\mathbb{P}^d,\mathcal{F}))^{\vee}.
\end{equation}
The collection $\{P_I \mid I \subset \Delta, \alpha_0,\dots,\alpha_{d-j-1} \in I\}$ is the set of all parabolic subgroups of $G$ containing $P_{(d-j+1,1^j)}$. This means that 
\begin{equation}
\mathrm{coker}\bigg( \bigoplus_{\substack{ I \subset \Delta\\|I|=d-j+1 \\ \alpha_0,\dots,\alpha_{d-j-1} \in I}} \ind^{G}_{P_I}(\mathds{1})    \rightarrow  \ind^{G}_{P_{(d-j+1,1^j)}}(\mathds{1})\bigg)={\leftindex_n {\upsilon}}_{P_{(d-j+1,1^j)}}^G.
\end{equation}
Since ${\leftindex_n {\upsilon}}_{P_{(d-j+1,1^j)}}^G$ 
is finitely generated free over $\wn(k)$,  it implies
\begin{equation}
 E_{2,\mathrm{w.s.}}^{-j+1,j}=\h_{}^j(\mathbb{P}^d,\mathcal{F})\otimes_{\wn(k)} ({\leftindex_n {\upsilon}}^{G}_{P_{(d+1-j,1^{j})}})^{\vee}.
 \end{equation}
 \end{proof}
 
\section{Revisiting the crystalline Beilinson--Bernstein map}\label{sec_crys_bb}
In characteristic $0$, $\mathcal{D}_X$ is notably interesting for its relationship with Lie algebras representations. More precisely, If $K$ is a field of characteristic $0$ and $\textbf{G}$ is a reductive group over $K$, acting on a flag $K$-variety $X$, then there exists a map $$\phi^{\ox}: \mathcal{U}(\mathfrak{g}) \rightarrow \Gamma(X,\mathcal{D}_X)$$ obtained by "differentiating the action of $\textbf{G}$", that we call Beilinson-Bernstein (ab. BB) map, motivated by \cite{beiber81}. This notion comes precisely from the following. Let $Z$ be a representation of $\textbf{G}$, then we can associate with it a representation of the Lie algebra $\mathfrak{g}$: if $\zeta \in \mathfrak{g}$ and $f \in Z$, then 
\begin{equation}\label{eq:diff_act_char0}
    \zeta.f:=\frac{d}{d\epsilon}\left(\mathrm{exp}(\epsilon\zeta).f\right)_{\mid \epsilon=0}
\end{equation}
(see for example \cite[Sec. 4]{romanov21}). The operator $\frac{d}{d\epsilon}$ is a differential of order 1, so this action extends to the map $\phi^{\ox}$ above (comparing filtration of the enveloping algebra and filtration of the differential operators). In characteristic 0, the main feature  lies on the fact that $\phi^{\mathcal{O}_X}$ is surjective.  In a field $k$ of characteristic $p>0$, we can adapt this construction to produce a map $\dist(\textbf{G}_k) \rightarrow \Gamma(X,\mathcal{D}_X)$, but it is no longer surjective (cf. \cite[Theorem 3.11]{smith86}). Those maps are used, for example, in \cite{Orlik21} and \cite{smith86}. Then we can also describe a lift of this map, namely
$\w_n\langle \phi^{\ox}\rangle: \dist(\textbf{G}_k) \rightarrow \Gamma(X,\mathcal{D}_{\wn(X)})$, in order to investigate geometric representations over $\w_n(k)$.
\subsection{Description of the Beilinson--Bernstein map in positive characteristic }
For seek of completeness and since we are not able to find a reference where it is described, we will define a crystalline BB map in such way it agrees with the one used in \cite{Orlik21,smith86} similarly to the characteristic $0$ case (see \cite{BBL18,romanov21}). Let $k$ be a field of characteristic $p>0$. Let $\textbf{G}=\textbf{G}_{\mathbb{Z}}$ be an algebraic reductive group, and $X$ be a smooth $k$-variety equipped with an action of $\textbf{G}_k=\textbf{G} \otimes k$.

Let $\mathfrak{m}=\{f \in \Gamma(\textbf{G}_k,\mathcal{O}_{\textbf{G}_k}) \mid f(1)=0\}$ be the maximal ideal of $\mathcal{O}_{\textbf{G}_k,1}$.

\begin{defn}
    The Lie algebra $\mathfrak{g}_k$ is the tangent space at $1$ of $\textbf{G}_k$, i.e. $\mathfrak{g}_k=(\mathfrak{m}/{\mathfrak{m}^2})^{*}=\mathrm{Hom}_{k-lin.}(\mathfrak{m}/{\mathfrak{m}^2},k)$ as a $k$-module.
\end{defn} 
Let consider the first infinitesimal neighborhood of the identity $\spec(\mathcal{O}_{\textbf{G}_k}/{\mathfrak{m}^2})=:\textbf{G}_k^{(1)}$ . The $k$-algebra $\mathcal{O}_{\textbf{G}_k}/{\mathfrak{m}^2}$ is isomorphic to $k \oplus \epsilon \mathfrak{g}_k^*=:k[\G_k^{(1)}]$ where $\epsilon^2=0$.\footnote{Since we require that $\epsilon^2=0$, this is the same of giving the "extension by 0" multiplication structure on $k \oplus \mathfrak{g}_k^*$, i.e. $(c,\zeta)\cdot (c',\eta):=(cc',c\eta+c'\zeta)$.} 
\begin{lemma}[{cf. \cite[Proposition 3.4]{milneLAG}}]
      There are natural bijections 
\begin{equation}ker(\G_k(k[\epsilon]/(\epsilon^2))\xrightarrow{\epsilon \mapsto 0} \G_k(k)) \leftrightarrow \mathrm{Der}_k(k[\G_k^{(1)}],k) \leftrightarrow (\mathfrak{m}/{\mathfrak{m}^2})^{*}, \quad 1+\epsilon\zeta \mapsto \zeta
\end{equation}
Moreover,  the group structure of 
$\ker(\G_k(k[\epsilon]/(\epsilon^2))\xrightarrow{\epsilon \mapsto 0} \G_k(k))$ corresponds to the additive structure of $(\mathfrak{m}/{\mathfrak{m}^2})^{*}$, while the $k$-linear structure  is given by $\lambda.(1+\epsilon\zeta):=1+\epsilon\lambda\zeta$ , for any $\lambda \in k$.
\end{lemma}
The Lie algebra structure on $(\mathfrak{m}/{\mathfrak{m}^2})^{*}$ is given in the following way: 
 
 Let $\zeta, \eta \in \mathrm{Hom}_{k-lin.}(\mathfrak{m}/{\mathfrak{m}^2},k)=\mathrm{Der}_k(k[\G_k^{(1)}],k)$, and let $\Delta$ be the comultiplication of the algebraic group $\G_k^{(1)}$. Then, we define 
 \begin{equation}
 \zeta.\eta \colon k[\G_k^{(1)}] \xrightarrow{\Delta} k[\G_k^{(1)}] \otimes_k k[\G_k^{(1)}] \xrightarrow{\zeta \otimes\eta} k \otimes k \simeq k
 \end{equation}
 where the last map is the multiplication in $k$. Set $[\zeta,\eta]:=\zeta.\eta-\eta.\zeta \in (\mathfrak{m}/{\mathfrak{m}^2})^{*} $.
 \begin{lemma}
     The $k$-module $(\mathfrak{m}/{\mathfrak{m}^2})^{*}$ endowed with the operation $[-,-]$ is a Lie algebra. 
 \end{lemma}
 \begin{proof}
     This is cf. \cite[Part I, 7.7]{Jantzen03}.
 \end{proof}
   It follows that $1+\epsilon[-,-]$ makes $\ker(\textbf{G}_k(k[\epsilon]/(\epsilon^2))\xrightarrow{\epsilon \mapsto 0} \textbf{G}_k(k))$ a Lie algebra too.
 
We consider the restriction of the action of $\textbf{G}_k$ on $X$  to the action $\sigma \colon \G_k^{(1)} \times_k X \rightarrow X$.
Let $\mm$ be a quasi-coherent $\G$-linearizable $\oo_X$-module. The restriction on $\G_k^{(1)}$ induces a linearization
\begin{equation}
   \phi \colon \sigma^*\mm \rightarrow \mathrm{pr}_1^*\mm. 
\end{equation}
\begin{prop}\label{lie_alg_map_induced_by_lineariz}
    The isomorphism $\phi$ induces a Lie algebra homomorphism
    \begin{equation}
       \rho: \mathfrak{g}_k \rightarrow \mathrm{End}_k(\mm),
    \end{equation}
    such that, for any open affine $U \subset X$, $\zeta \in \mathfrak{g}_k$, $m \in \Gamma(U,\mm)$, $s \in \oo_X(U)$, we have
    \begin{equation}\label{leibn_rel}
        \rho(\zeta)(sm)=s\rho(\zeta)(m) + \rho(\zeta)(s)m.
    \end{equation}
\end{prop}
\begin{proof}
For any open affine $U \subset X$, we have that
\begin{equation}
    \Gamma(\G_k^{(1)} \times U,p_1^*\mm)=k[\G_k^{(1)}] \otimes_k \Gamma(U,\mm)
\end{equation}
Now, consider the isomorphisms of schemes $ \varepsilon_i \colon \G_k^{(1)} \times X \rightarrow \G_k^{(1)} \times X $ for $i \in \{1,2\}$, given on points by $\varepsilon_1(g,x)=(g,gx)$ and $\varepsilon_2(g,x)=(g,g^{-1}x)$ for $g \in \G_k^{(1)}, \text{ } x \in X$. Clearly, $\varepsilon_1=\varepsilon_2^{-1}$ and $\mathrm{pr}_1 \circ \varepsilon_1=\sigma$. Therefore, we have the following equalities:
\begin{align*}
    \Gamma(\G_k^{(1)} \times U,\sigma^*\mm)&=\Gamma(\G_k^{(1)} \times U,\varepsilon_1^*\mathrm{pr}_1^*\mm)\\
    &=\Gamma(\G_k^{(1)} \times U,(\varepsilon_{2})_*\mathrm{pr}_1^*\mm) \\
    &=\Gamma(\varepsilon_2^{-1}(\G_k^{(1)} \times U),\sigma^*\mm) \\
    &=\Gamma(\varepsilon_1(\G_k^{(1)} \times U),\sigma^*\mm) \\
    &=\Gamma(\G_k^{(1)} \times U,\mathrm{pr}_1^*\mm)
\end{align*}
where the last equality is satisfied since as a topological space $\G_k^{(1)}$ is a singleton and $g.U=U$ for $g \in \G_k^{(1)}$, since $\G_k^{(1)}$ is a neighbourhood of the identity of $\G$. Therefore, $\phi$ induces an isomorphism
\begin{equation}
    k[\G_k^{(1)}] \otimes_k \Gamma(U,\mm) \xrightarrow{\phi} k[\G_k^{(1)}] \otimes_k \Gamma(U,\mm).
\end{equation}
The composition with the natural morphism $\sigma^{\#} \colon  \Gamma(U,\mm) \rightarrow k[\G_k^{(1)}] \otimes_k \Gamma(U,\mm)$ induces a  map \footnote{The \eqref{eq:counit},\eqref{eq:comult} follow by the definition of linearization:  in particular, $(1 \times Id_X)^*\phi=Id_{\mm} \colon \mm \rightarrow \mm$ implies \eqref{eq:counit} and $pr_{12}^*\phi \circ (1 \times \sigma)^*\phi=(m \times id_X)^*\phi$ implies \eqref{eq:comult}. }
\begin{equation}
    \Gamma(U,\mm) \xrightarrow{\phi \circ \sigma^{\#}} k[\G_k^{(1)}] \otimes_k \Gamma(U,\mm)
\end{equation}
such that 
  \begin{align}\label{eq:counit}
  &(1 \otimes Id) \circ \phi \circ \sigma^{\#}=Id, \\
  \label{eq:comult}
  &(\phi \circ \sigma^{\#})\otimes Id_{k[\G_k^{(1)}]} \circ (\phi \circ \sigma^{\#})=Id_{\Gamma(U,\mm)} \otimes \Delta \circ (\phi \circ \sigma^{\#})
  \end{align} where $1 \colon k[\G^{1}_k] \rightarrow k$ is the $k$-rational point corresponding to the neutral element of $\G_k^{(1)}$ (In particular $\phi \circ \sigma^{\#}$ is a comodule map over the Hopf algebra $k[\G_k^{(1)}]$). For any $\zeta \in \mathfrak{g}_k$, we define $\rho$ given by the following map: 
\begin{equation}
    \rho(\zeta) \colon \Gamma(U,\mm) \xrightarrow{\phi \circ \sigma^{\#}} k[\G_k^{(1)}] \otimes_k \Gamma(U,\mm) \xrightarrow{\zeta \otimes Id} \Gamma(U,\mm).
\end{equation}
Let $s \in \ox(U)$ and $m \in \Gamma(U,\mm)$, then we have
\begin{align*}
    \rho(\zeta)(sm)&=(\zeta \otimes Id ) \circ (\phi \circ \sigma^{\#}(sm)) \\
    &=(\zeta \otimes Id)(\sigma^{\#}(s)\cdot (\phi \circ \sigma^{\#})(m))\\
    &=(\zeta \otimes Id)(\sigma^{\#}(s))\cdot (1 \otimes Id) \circ (\phi \circ \sigma^{\#})(m)+ (1 \otimes Id)(\sigma^{\#}(s))\cdot (\zeta \otimes Id) \circ (\phi \circ \sigma^{\#})(m)\\
    &=\rho(\zeta)(s)\cdot m+s \cdot \rho(\zeta)(m).
\end{align*}
Since $\phi \circ \sigma^{\#}$ is a comodule map, it follows that $\rho$ is a Lie algebra map by \cite[Part I, 7.11 (2)]{Jantzen03}.
\end{proof}
In particular, for $\mm=\ox$ the \eqref{lie_alg_map_induced_by_lineariz} is the Lie algebra homomorphism $$\mathfrak{g}_k \rightarrow \mathrm{Der}_k(\Gamma(U,\ox)), \quad \zeta \mapsto (f \mapsto  (Id \otimes \zeta) \circ \sigma^{\#} (f)), \text{ for any } \zeta \in \mathfrak{g}_k, f \in \Gamma(U,\ox) $$
that we call \textit{differentiated contragradient action} (cf. \cite[I.2.7, I.2.8, I.7.11 (1)]{Jantzen03}).
\begin{rem}
In particular, note that the Lie algebra of an algebraic group could act on a scheme, even if the algebraic group does not. For example, in the case of $\mathcal{X} \subset \mathbb{P}^d_k $, the projective space is equipped with an action of $ GL_{d+1,k}$, thus $\h^0(\mathcal{X},\oo_{\mathbb{P}^d_k})$ inherits an action of $\mathfrak{gl}_{d+1,k}$ by the construction above. However, $GL_{d+1,k}$  does not act on $\mathcal{X}$ (not to be confused with the finite group of $k$-rational points $GL_{d+1}(k)$).
\end{rem}
 There is an alternative description of the Lie algebra map above, better adapted for later computations:
\begin{defn}\label{def_cano_diff}
Let $\sigma$ be an action of $\textbf{G}_k$ on a smooth $k$-variety $X$. The \textit{ canonical differential action} of $\sigma$ is the map $\partial\sigma:\mathfrak{g}_k \rightarrow \mathrm{End}_k(\Gamma(U,\ox))$ given by 
\begin{equation}
    \partial\sigma(\zeta)(f)={\frac{d}{d\epsilon}\Big(\sigma(1+\epsilon\zeta)(f)\Big)}_{\mid \epsilon=0}
\end{equation}
for $f \in \Gamma(U,\ox)$ and $\zeta \in \mathfrak{g}_k$.
\end{defn}
\begin{rem}
The map appearing \Cref{def_cano_diff} reminds the shape of the differential action in characteristic 0 \eqref{eq:diff_act_char0}. Here, we formally truncate the expansion of the exponential function at the first order. That is the reason for calling it "canonical differential action".  
\end{rem}
\begin{lemma}
For any $\zeta \in \mathfrak{g}_k$, the map $\partial\sigma(\zeta)$ is a $k$-linear derivation. 
\end{lemma}
\begin{proof}
Let $\zeta,\zeta' \in \mathfrak{g}_k$ and $f,f' \in \Gamma(U,\ox)$. Since $\sigma(1+\epsilon \zeta)$ is a $k$-algebra homomorphism, then 
\begin{multline}
\partial\sigma(\zeta)(ff')={\frac{d}{d\epsilon}\Big(\sigma(1+\epsilon\zeta)(ff')\Big)}_{\mid \epsilon=0}=
{\frac{d}{d\epsilon}\Big(\sigma(1+\epsilon\zeta)(f)\sigma(1+\epsilon\zeta)(f')\Big)}_{\mid \epsilon=0}\\
=\sigma(1)(f){\frac{d}{d\epsilon}\Big(\sigma(1+\epsilon\zeta)(f')\Big)}_{\mid \epsilon=0}
+ \sigma(1)(f'){\frac{d}{d\epsilon}\Big(\sigma(1+\epsilon\zeta)(f)\Big)}_{\mid \epsilon=0}\\
=f\partial\sigma(f')+f'\partial\sigma(f).
\end{multline}
This shows that $\partial\sigma(\mathfrak{g}_k) \subset \mathrm{Der}_k(\Gamma(U,\ox))$.
\end{proof}
\begin{lemma}
    The maps \eqref{lie_alg_map_induced_by_lineariz} and $\partial\sigma$ agree when $\mm=\ox$. In particular, $\partial\sigma$ is a Lie algebra homomorphism.
\end{lemma}
\begin{proof}
Set $\oo_U:=\Gamma(U,\ox)$ and $\oo_{U_{\epsilon}}:=\oo_{U}[\epsilon]/(\epsilon^2)$. Let $\zeta \in \mathfrak{g}_k$. We view $\zeta_k \colon k[\textbf{G}_k^{(1)}] \rightarrow k$ as a derivation. Equivalently, it is a $k$-algebra homomorphism $1+ \epsilon \zeta \colon k[\textbf{G}_k^{1}] \rightarrow k[\epsilon]/{(\epsilon^2)}$. Then, the statement follows by the commutative diagram
\begin{equation*}
 \begin{tikzcd}
\mathcal{O}_{U_{\epsilon}} \arrow[r, "\sigma^{\#}_{\epsilon}"] & {\mathcal{O}_{U_{\epsilon}} \otimes k[\textbf{G}_k^{(1)}]} \arrow[r, "Id\otimes(1+\epsilon\zeta)"] & {\mathcal{O}_{U_{\epsilon}} \otimes k[\epsilon]/{(\epsilon^2)}=\mathcal{O}_U \otimes k[\epsilon]/{(\epsilon^2)} \otimes k[\epsilon]/{(\epsilon^2)}} \arrow[d, "Id \otimes m"] \\
\mathcal{O}_U \arrow[u, hook] \arrow[r, "\sigma^{\#}"]         & {\mathcal{O}_U \otimes k[\textbf{G}^{(1)}_k]} \arrow[r] \arrow[rd, "Id \otimes\zeta"]              & {\mathcal{O}_U \otimes k[\epsilon]/{(\epsilon^2)}} \arrow[d, "\frac{d}{d\epsilon}_{\mid \epsilon=0}"]                                                                         \\
                                                               &                                                                                   & \mathcal{O}_U                                                                                                                                                                
\end{tikzcd}
\end{equation*}
Here, $m \colon k[\epsilon]/{(\epsilon^2)} \otimes  k[\epsilon]/{(\epsilon^2)} \rightarrow
k[\epsilon]/{(\epsilon^2)} $ denotes the multiplication,
${\frac{d}{d\epsilon}}_{\mid \epsilon=0}: \mathcal{O}_U[\epsilon]/{(\epsilon^2)} \rightarrow \mathcal{O}_U$ is the map sending $a+b\epsilon \mapsto b$ , $\mathcal{O}_{U_{\epsilon}}:=\mathcal{O}_U \otimes k[\epsilon]/(\epsilon^2)$ and $\sigma_{\epsilon}^{\#}:=\sigma^{\#}\otimes k[\epsilon]/(\epsilon^2)$.  Then, by definition  we can read
$$\sigma(1+\epsilon\zeta) \colon  \mathcal{O}_U \rightarrow \mathcal{O}_U \otimes k[\textbf{G}_k^{(1)}] \rightarrow{} \mathcal{O}_U [\epsilon]/{(\epsilon^2)},$$ 
following the vertical left, top horizontal, and top vertical right arrows.  In the bottom, we have the contragradient action of $\textbf{G}_k^{(1)}$.  Since $1+ \epsilon\zeta \colon k[\textbf{G}_k^{(1)}] \rightarrow k[\epsilon]/{(\epsilon^2)}$  corresponds to the derivation $\zeta : k[\textbf{G}_k^{(1)}] \rightarrow k$, compatibly with $k[\epsilon]/{(\epsilon^2)} \xrightarrow{{\frac{d}{d\epsilon}}_{\mid \epsilon=0}} k $,  then the diagram above commutes.
\end{proof}

\begin{rem}
From now on, when we talk about "differentiating an action" we mean to consider the canonical differential action. 
\end{rem}

 Thus $$\partial\sigma: \mathfrak{g}_k \rightarrow \Gamma(U,\mathrm{Fil}_1\mathcal{D}_X) \subset \Gamma(U,\mathcal{D}_X) $$  is a morphism of Lie algebras and extends to a morphism of  associative algebras  $\mathcal{U}(\mathfrak{g}_k) \to \Gamma(X,\mathcal{D}_X) $ sending $\mathcal{U}_m(\mathfrak{g}_k)$ to $\Gamma(U,\mathrm{Fil}_m\mathcal{D}_X)$. Further, there is another natural morphism of filtered algebras $\mathcal{U}(\mathfrak{g}_k) \to \dist(\textbf{G}_k)$ induced (by the universal property of enveloping algebras) by the natural inclusion of $\mathfrak{g}_k \subset \dist(\textbf{G}_k)$. Here,  the filtration on $\dist(\G_k)$ is given by distributions of order $n$, i.e. $\fil^n\dist(\G_k)=\dist_n(\G_k)$.
 
 Choose an ordered basis of $\mathfrak{g}_{\zz}$, namely 
 $$\mathcal{B}_{\textbf{G}}=\{L_{\alpha},H_{\beta},Y_{-\alpha},h_1,\dots,h_m \mid \alpha \in \Phi^+ , \beta \in \Delta \}.$$ It is formed in the following way:
 
 For any $\alpha \in \Phi^+$, $L_{\alpha} \in (\mathfrak{g}_{\zz})_{\alpha}, Y_{-\alpha} \in  (\mathfrak{g}_{\zz})_{-\alpha} $ are the generators of the corresponding one dimensional weight spaces associated to the roots $\pm\alpha$ and $H_{\alpha}=[L_{\alpha},Y_{-\alpha}]$.  The subset $\mathcal{B}_{\textbf{T}}=\{H_{\beta},h_1,\dots,h_m \mid \beta \in \Delta\}$ form a basis for $\mathrm{Lie}(T_{\zz}) \subset \mathfrak{g}_{\zz}$. Also, it contains a semisimple part, given by the subset $\mathcal{B}^{ss}_{\textbf{T}} \subset \mathcal{B}_{\textbf{T}}$ formed by the elements $H_{\beta}$ with $\beta \in \Delta$. 
 \begin{lemma}[{cf. \cite[Part II,1,1.12]{Jantzen03} and \cite[Ch. 2, Corollary to Lemma 5]{Ste16}}]\label{lemma_distrib_of_G}
      The distribution algebra $\dist(\G)$ is the $\zz$-subalgebra of  $\mathcal{U}(\mathfrak{g}_{\zz})\otimes \mathbb{Q}$ generated by $L_{\alpha}^{[a_{\alpha}]}:=\frac{1}{a_{\alpha}!}L_{\alpha}^{a_{\alpha}}$,  $ Y_{-\alpha}^{[b_{\alpha}]}:=\frac{1}{b_{\alpha}!}Y_{-\alpha}^{b_{\alpha}}$, and $\binom{h_j}{c_j}$ for ${\alpha} \in \Delta$ and $j=1,\dots,m$ such that $h_j \in \mathcal{B}_{\textbf{T}} \backslash \mathcal{B}^{ss}_{\textbf{T}}, a_{\alpha}, b_{\alpha}, c_j \in \mathbb{N}.$
 \end{lemma} 
Thus, by \Cref{lemma_distrib_of_G} we get the following:
\begin{prop}\label{prop_bb_map_constr}
    Assume that $X$ has a smooth lift over $\zz$, i.e. there exists a smooth $\zz$-scheme $X_{\zz}$ such that $X \simeq X_{\zz} \times k$, and the action $\sigma$ lifts to an action $\sigma_{\zz}$ of $\G$ on $X_{\zz}$. Then, for any open $U \subset X$, there is a unique well defined morphism of $k$-algebras 
 \begin{equation}\label{bbmap}
    \phi^{\ox}: \dist(\G_k)=\dist(\G) \otimes k \rightarrow \Gamma(U,\mathcal{D}_X),
\end{equation}
sending
\begin{align*}
\frac{1}{a_{\alpha}!}L_{\alpha}^{a_{\alpha}} \otimes 1 &\mapsto \frac{1}{a_{\alpha}!}\partial\sigma(L_{\alpha} \otimes 1)^{a_{\alpha}}, & \text{for all } a_{\alpha} \in \mathbb{N}\\
\frac{1}{b_{\alpha}!}Y_{-\alpha}^{b_{\alpha}} \otimes 1 &\mapsto \frac{1}{b_{\alpha}!}\partial\sigma(Y_{-\alpha} \otimes 1)^{b_{\alpha}}, & \text{ for all  } {\alpha}  \in \Delta, b_{\alpha} \in \mathbb{N} \\
\binom{h_j}{c_j} \otimes 1 &\mapsto \binom{\partial\sigma(h_j \otimes 1)}{c_j},  & \text{ for all  } h_j \in \mathcal{B}_{\textbf{T}} \setminus \mathcal{B}^{ss}_{\textbf{T}}, c_j \in \mathbb{N}.
\end{align*}
\end{prop}
\begin{proof}
Let $\dim(X)=d$. Let $\tilde{X}:=X_{\zz} \times \mathbb{Q}$ and $\sigma_{\mathbb{Q}}=\sigma_{\zz} \otimes \mathbb{Q} $. Similarly,  $\tilde{U}=U_{\zz} \times \mathbb{Q}$ for some $U_{\zz} \subset X_{\zz}$ open such that $U_{\zz} \times k \simeq U$.
We have that $\dist(\G) \subset \dist(\G) \otimes \mathbb{Q}=\mathcal{U}(\mathfrak{g}_{\mathbb{Q}})$ is a free $\zz$-module.
Moreover, 
    $\Gamma(U,\mathcal{D}_X) \simeq S_d(k)=k[z_1,\dots,z_d]\langle \partial_{z_1}^{[s_1]},\dots \partial_{z_d}^{[s_d]}\rangle=S_d(\zz) \otimes k$. 

    We notice that $S_d(\mathbb{Q})\simeq \mathcal{U}_{\oo_{\tilde{X}}(U)}(\mathcal{T}_{\tilde{X}}(\tilde{U}))$.\footnote{The notation $\mathcal{U}_{\mathcal{A}}(\mathcal{L})$ means the universal enveloping algebra of a Lie-Rinehart algebra $\mathcal{L}$ over a commutative $k$-algebra $\mathcal{A}$ (cf. \cite[Section 2]{Rinehart}).}
    
     By the universal property of enveloping Lie algebras, the differential of the action $\sigma_{\mathbb{Q}}$:
     \begin{equation}
         \partial\sigma_{\mathbb{Q}} \colon \mathfrak{g}_{\mathbb{Q}} \rightarrow \mathcal{T}_{\tilde{X}}(\tilde{U}),
     \end{equation}
     is a Lie algebra homomorphism that lifts uniquely to a map of associative $\mathbb{Q}$-algebras
     \begin{equation}
         \phi \colon \mathcal{U}(\mathfrak{g}_{\mathbb{Q}}) \rightarrow \mathcal{U}_{\oo_{\tilde{X}}(U)}(\mathcal{T}_{\tilde{X}}(\tilde{U})).
     \end{equation}
      Moreover, for any $\alpha \in \Delta$ and any $m \in \mathbb{N}$,
      $m!L_{\alpha}^{[m]}=L_{\alpha}^m$,  $m!Y_{-\alpha}^{[m]}=Y_{-\alpha}^m$, $m!\binom{h_i}{m}=h_i(h_i-1)\cdots (h_i-m+1)$ hold in $\dist(\G)$. Applying $\phi$, we get respectively the equality $m!\phi(L_{\alpha}^{[m]})=\partial\sigma_{\zz}(L_{\alpha})^m$,  $m!\phi(Y_{-\alpha}^{[m]})=\partial\sigma_{\zz}(Y_{-\alpha})^m$, $m!\phi(\binom{h_i}{m})=\partial\sigma_{\zz}(h_i)\partial\sigma_{\zz}(h_i-1)\cdots \partial\sigma_{\zz}(h_i-m+1)$ therefore the algebra map $\phi$ induces an associative $\zz$-algebra map 
      \footnote{The corresponding elements $\phi(L_{\alpha}^{[m]})=\frac{1}{m!}\partial\sigma(L_{\alpha})^m$ and $\phi(Y_{-\alpha}^{[m]})=\frac{1}{m!}\partial\sigma(Y_{-\alpha})^m$ are well defined, since for any $\delta \in \mathcal{T}_{X_{\zz}}(U_{\zz})$, and $m \in \mathbb{N}$, $\delta^{[m]}$ is a $\zz$-linear combination of elements of the form $\prod_{i=1}^d \partial_{z_i}^{[r_i]}$ with $r_i \geq 0$. Also $\phi(\binom{h_i}{m_i})=\binom{\partial\sigma(h_i)}{m_i}$ is well defined: Indeed, since $h_i \in \mathrm{Lie}(T)$, then $\partial\sigma(h_i)$ acts diagonally  on $\Gamma(U_{\zz}, \oo_{X_{\zz}})$ (c.f. \cite[II, 1.19]{Jantzen03}) by looking at the weight space decomposition.}
\begin{equation}
    \phi_{\zz} \colon \dist(\G) \rightarrow \Gamma(U_{\zz}, \mathcal{D}_{X_{\zz}}).
\end{equation}
 Hence, we get the searched map after tensor by $k$.
\end{proof}
\begin{defn}
We call the morphism \eqref{bbmap} \textit{the Beilinson-Bernstein map} (BB ) w.r.t. $\ox$.
\end{defn}
We have by definition a commutative diagram
\begin{equation*}
\begin{tikzcd}
\mathcal{U}(\mathfrak{g}_k) \arrow[rr, "\partial\sigma"] \arrow[rd] &                                     & {\Gamma(X,\mathcal{D}_X)} \\
                                                                    & \dist(\G_k) \arrow[ru, "\phi^{\ox}"] &                
\end{tikzcd}
\end{equation*}
of associative $k$-algebras.
 Hence, the map \eqref{bbmap} is a morphism of filtered associative algebras, so for any distribution $\zeta$ of order $n$, and any distribution $\eta$ of order $m$, 
$$ \phi^{\ox}([\zeta,\eta])=[ \phi^{\ox}(\zeta), \phi^{\ox}(\eta)] \in \Gamma(X,\mathrm{Fil}^{n+m-1}\mathcal{D}_X).$$

\begin{lemma}
If $X$ is a smooth scheme over $k$, the natural inclusion (and homeomorphism) $i \colon X \hookrightarrow \wn(X)$, induces a surjective morphism of sheaves:
\begin{equation}
i^{*} \colon \mathcal{D}_{\wn(X)} \twoheadrightarrow i_{*}\mathcal{D}_{X}
\end{equation}
\end{lemma}
\begin{proof}
By the  \Cref{uniquenesslift}, it follows that the induced map on the stalks of any point of $X$ is surjective, thus the map of sheaves is surjective as well.
\end{proof}
We can define a section for this map (in the category of sheaves of sets).
\begin{prop}\label{prop_teich_lift_diff}
 There exists a map of sheaves of sets
   \begin{equation}
[.] \colon i_{*}\mathcal{D}_X \rightarrow \mathcal{D}_{\wn(X)}
\end{equation} 
such that $i^* \circ [.]=id$. 
\end{prop}
\begin{proof}
Let $U \subset X$ be an open subset.  We define the map 
\begin{equation}
[.]_U \colon \Gamma(U,\mathcal{D}_{X}) \rightarrow \Gamma(U,\mathcal{D}_{\wn(X)})
\end{equation}
such that if $\partial \in \Gamma(U,\mathcal{D}_{X}) $ is one of the $\partial_{\bm{z}}^{[\mathbf{r}]}$, for  local coordinates $\bm{z}$, then 
$[\partial]_U:=\partial_{n,\bm{z}}^{\mathbf{r}}$ as defined in  \Cref{liftder}. Otherwise, there is a unique way to write any $\partial \in \Gamma(U,\mathcal{D}_{X}) $ as 
\begin{equation}
\partial=\sum_{\mathbf{r}} b_{\mathbf{r},U}\partial_{\bm{z}}^{[\mathbf{r}]},\quad  b_{\mathbf{r},U} \in \ox(U).
\end{equation}
In that case, let 
\begin{equation}
[\partial]_U:=\sum_{\mathbf{r}} [b_{\mathbf{r},U}]\partial_{\bm{z},n}^{[\mathbf{r}]}
\end{equation}
where $[.]$ is the Teichm\"uller map $\ox(U) \rightarrow \wn\ox(U)$.
We have the following commutative diagram, for any $V \subset U $:
\begin{equation}
\begin{tikzcd}
\Gamma(U,\mathcal{D}_X) \arrow[d, "\rho_{UV}"] \arrow[r, "{[.]_U}"] & \Gamma(U,\mathcal{D}_{\wn(X)}) \arrow[d, "\rho^{(n)}_{UV}"] \\
\Gamma(V,\mathcal{D}_X) \arrow[r, "{[.]_V}"]                        & \Gamma(V,\mathcal{D}_{\wn(X)})                             
\end{tikzcd}
\end{equation}
Indeed, it is clear for the differential operators of the form $\partial_{\bm{z}}^{[\mathbf{r}]}$,  since they uniquely determine lifts $\partial_{\bm{z},n}^{[\mathbf{r}]}$. In the general case, it follows since the restriction of $[.]_U$ to $\wn\ox$ is the Teichm\"uller map, for which such a diagram is commutative. 
In particular, the maps $[.]_U$ define a map of sheaves with the desired property.
\end{proof}
\begin{defn}
For any smooth $k$-scheme $X$, we define the map of sets
\begin{equation}\label{teich_lift_diff}
    \w_n\langle \phi^{\ox}\rangle: \dist(\textbf{G}_k) \xrightarrow{\phi^{\ox}} \Gamma(X,\mathcal{D}_X)\xrightarrow{[.]}\Gamma(X,\mathcal{D}_{\wn(X)}).
\end{equation}
\end{defn}
\subsection{Examples}\label{sec_ex_bb}
Here, we want to describe the BB map for some reductive group $\textbf{G}$, acting on the flag variety $X=\textbf{G}/\textbf{B}$, where $\textbf{B}$ is a fixed Borel subgroup. We suppose that all groups are defined over $\mathbb{Z}$. The index $(-)_k$ means we are tensoring with $k$, as done before.
\begin{ex}[$\G=\mathrm{SL}_2$]
For $\textbf{G}= \mathrm{SL}_2$, we have $X=\mathbb{P}^1$. The Lie algebra $\mathfrak{g}=\mathfrak{sl}_2$ is generated by the matrices 
\begin{equation}
    L=\begin{pmatrix} 0 & 1 \\
                      0 & 0\end{pmatrix}, \quad 
    H=\begin{pmatrix} 1 & 0 \\
                      0 & -1\end{pmatrix}, \quad
    Y= \begin{pmatrix} 0 & 0 \\
                       1 & 0\end{pmatrix}.            
\end{equation}
The Kostant's $\mathbb{Z}$-form of $\mathfrak{g}=\mathfrak{sl}_2$ is the $\mathbb{Z}$-algebra $U_{\mathbb{Z}}$ generated by
\begin{equation*}L^{[a]}:=\frac{L^a}{a!},\quad \binom{H}{b}:=\frac{H(H-1)\dotsc (H-b+1)}{b!}, \quad Y^{[c]}:=\frac{Y^c}{c!},
\end{equation*} where $a,b,c \in \mathbb{Z}_{\geq 0}$. 
\begin{lemma}\label{lemma_global_D_P^1}
    The $\mathbb{Z}$-algebra
       generated by the following sections of $\mathcal{D}_{\mathbb{P}^1}$ (not a priori globally defined)
    $$ z^r \Big(\frac{\partial}{\partial z}\Big)^{[s]}, \qquad \text{for all } r, s \in \mathbb{N} \text{ such that }  0 \leq r \leq 2s,$$
    is a $\zz$-subalgebra of $\Gamma(\mathbb{P}^1,\mathcal{D}_{\mathbb{P}^1})$.
\end{lemma}
\begin{proof}
    If $[x_0 : x_1]$ are the coordinates of a point in $X$, the standard covering of $X$ is given by\footnote{$X_f$ for a regular function $f$ denotes the standard open $X \setminus Z(f)$.} $(U_0=X_{x_0}\simeq \mathbb{A}^1,U_{\infty}=X_{x_1}\simeq \mathbb{A}^1)$. Denote with $z$ the local coordinate of $U_0$ and $w$ for the local coordinate of $U_{\infty}$. On the intersection we have $w=\frac{1}{z}$. Then we have
$$\Gamma(U_0,\mathcal{D}_{\mathbb{P}^1})=\mathbb{Z}[z]\Big\langle \Big(\frac{\partial}{\partial z}\Big)^{[s]}, s \geq 0\Big\rangle,$$ i.e. the crystalline  Weyl algebra $S_1(\zz)$ is generated by $z$ and the differential operator $\frac{\partial}{\partial z}$ ( cf.  \Cref{smith2.7}).
Similarly,
$$\Gamma(U_{\infty},\mathcal{D}_{\mathbb{P}^1})=\mathbb{Z}[w]\Big\langle \Big(\frac{\partial}{\partial w}\Big)^{[s]}, s \geq 0\Big\rangle.$$ On the intersection $U_0 \cap U_{\infty}$, we have that $$\frac{\partial}{\partial w}=-z^2\frac{\partial}{\partial z}.$$ For any $m \in \zz$, we have that $z^r \Big(\frac{\partial}{\partial z}\Big)^{[s]}(z^m)=\binom{m}{s}z^{r-s+m} \neq 0$ if and only if  $|m| \geq s$. The equality $\Gamma(\mathbb{P}^1,\mathcal{D}_{\mathbb{P}^1})=\Gamma(U_0,\mathcal{D}_{\mathbb{P}^1})     \cap \Gamma(U_{\infty},\mathcal{D}_{\mathbb{P}^1})$ where the intersection is taken in $\mathcal{D}(\zz[z,1/z])$, means that $\Gamma(\mathbb{P}^1,\mathcal{D}_{\mathbb{P}^1})$  is generated by the operators of the form $z^r \Big(\frac{\partial}{\partial z}\Big)^{[s]}$ sending $\zz[z]$ to $\zz[z]$ and $\zz[w]$ to $\zz[w]$. Therefore, $z^r \Big(\frac{\partial}{\partial z}\Big)^{[s]}$ is a global section if and only if $r$ and $s$ are such that
\begin{equation}
   \left\{ \begin{array}{lr}
        r-s+m \geq 0 & \forall m \geq s  \\
         r-s+m \leq 0 & \forall m \leq -s 
    \end{array}\right.
\end{equation}
that is equivalent to $0 \leq r \leq 2s$.
\end{proof}

To understand the BB map, we have to look at the action of $\textbf{G}$ on $\ox$: we choose the natural action given by 
\begin{equation*}\label{action_SL_2onX}
    \sigma(g)(f)([x_0:x_1]):=f(g^{-1} \cdot [x_0:x_1]) \qquad g \in \mathrm{SL}_2, f \in \ox
\end{equation*}
\begin{prop}(cf. \cite[p. 175]{smith86}) 
For $\mathbf{G}=\mathrm{SL}_2$ and $X=\mathbf{G}/\mathbf{B} \simeq \mathbb{P}^1$ the BB map is the map of filtered $k$-algebras given on generators by:
\begin{equation}
\begin{array}{cc}
        \phi^{\ox}: \dist(\mathbf{G}_k)\simeq U_{\mathbb{Z}} \otimes k \longrightarrow \Gamma(\mathbb{P}_k^1,\mathcal{D}_{\mathbb{P}_k^1}) \\
        L   \otimes 1 \longmapsto z^2\frac{\partial}{\partial z} \\
        Y \otimes 1 \longmapsto -\frac{\partial}{\partial z}. 
    \end{array}
\end{equation}
\end{prop}
\begin{proof}
We need to compute the canonical differential action of $\sigma$ for the elements $L,Y \in \mathfrak{sl}_2$. We can make such computation on the chart $U_0$. For any $Q(z) \in \oo_{\mathbb{P}^1}(U_0)$, we get 
 \begin{align*}\partial\sigma(L)(Q(z)) &=
 \frac{d}{d\epsilon}_{\mid \epsilon=0}\sigma{(1+\epsilon L_{})}(Q(z)) \\
 &=\frac{d}{d\epsilon}_{\mid \epsilon=0}Q\left(\frac{z}{1-\epsilon z}\right) \\
 &=\frac{z^2}{1-\epsilon z}_{\mid \epsilon=0}\frac{\partial Q}{\partial z}\left(\frac{z}{1-\epsilon z}\right)_{\mid \epsilon=0} \\
 &=z^2\frac{\partial Q}{\partial z}(z)
 \end{align*} and
\begin{align*}
\partial\sigma(Y)(Q(z)) &=\frac{d}{d\epsilon}_{\mid \epsilon=0}\sigma{(1+\epsilon Y_{})}(Q(z)) \\&=\frac{d}{d\epsilon}_{\mid \epsilon=0}Q({z-\epsilon}) \\ &=-\frac{\partial Q}{\partial z}({z -\epsilon})_{\mid \epsilon=0}\\ &=-\frac{\partial Q}{\partial z}(z). \qedhere
\end{align*}

\end{proof}
\end{ex}
\begin{ex}[$\G=\mathrm{SL}_{d+1,\zz}$]
The example of $\mathrm{SL}_2$ reflects the situation in general for the group $\textbf{G}=\mathrm{SL}_{d+1,\zz}$. Let $\Delta$ be a basis of a root system $\Phi=\Phi^+ \cup \Phi^-$ of $\textbf{G}$. Consider $\mathfrak{t}=\mathrm{Lie}(T)$, where $T$ is a fixed maximal torus, $\mathfrak{g}_{\alpha}= \mathrm{Lie}(\mathbb{G}_a)$, where $\mathbb{G}_a$ is the one dimensional subgroup associated to the root $\alpha \in \Phi^+$ (resp. $\Phi^-$) of the unipotent group $U$ (resp. $U^{-}$) . Then we denote by 
\begin{equation*}
    L_{\alpha} \in \mathfrak{g}_{\alpha}, \text{ } \alpha \in \Phi^+,\quad 
    Y_{\alpha} \in  \mathfrak{g}_{\alpha}, \text{ } \alpha \in\Phi^{-},\quad
    H_{\alpha}=[L_{\alpha},Y_{-\alpha}] \in \mathfrak{t}, \text{ } \alpha \in \Phi^{+}
  \end{equation*}
their generators.
Then, $\mathfrak{g}$ is the $\zz$-module spanned by $L_{\alpha},Y_{-\alpha}, H_{\alpha}$ with $\alpha \in \Phi^+$.

 For any $\alpha \in \Phi^{+}$, the triple $\{L_{\alpha},Y_{-\alpha},H_{\alpha}\}$ generates a copy of $\mathfrak{sl}_{2,\mathbb{Z}}$ in $\mathfrak{g}$.
In the case of $\mathrm{SL}_{d+1}$ we have $d(d+1)$ roots, $d(d+1)/2$ of them are the positive roots, each one corresponding to an injective morphisms of Lie algebras  $s_{ij} \colon \mathfrak{sl}_2 \hookrightarrow \mathfrak{sl}_{d+1}$ for any $\alpha_{ij} \in \Phi^+$.  The variety $X=G/B \simeq \mathbb{P}^d$ is covered by the standard open cover of affine schemes: denote with $[x_0,\dots,x_d] \in \mathbb{P}^d$ point coordinates, then 
\begin{equation*}
    (V_0=X_{x_0}\simeq \mathbb{A}^d,\dots,V_{d}=X_{x_d}\simeq \mathbb{A}^d)
\end{equation*}
is the standard cover of $X$. For each $ 0 \leq i\not=j \leq d$ the intersection $V_i \cap V_j$  admits local coordinates $$z_i=(z_{0i}:=\frac{x_0}{x_i},\dots,\hat{\frac{x_i}{x_i}},\dots,z_{di}:= \frac{x_d}{x_i})$$ and 
$$z_j=(z_{0j}=\frac{x_0}{x_j},\dots,\hat{\frac{x_j}{x_j}},\dots, z_{dj}=\frac{x_d}{x_j})$$
such that 
\begin{equation} \label{glu}
z_{ij}=\frac{1}{z_{ji}}.
\end{equation}
So we can see $z_{ij}$ and $z_{ji}$ respectively as local coordinates of a copy of $\mathbb{A}^1  \subset V_j$ and $\mathbb{A}^1 \subset V_i$. After glueing those affine lines along \eqref{glu} we get a copy of $\mathbb{P}^1$. Now for every $\alpha_{ij} \in \Phi^{+}$, we get an $\mathfrak{sl}_2$-triple together with a BB map of $\mathbb{P}^1$. Then consider the closed immersion $$l \colon \mathbb{P}^1 \longrightarrow{} \mathbb{P}^d, \quad [x_0:x_1] \longmapsto [0:\dots:x_0:\dots:x_1:\dots: 0] $$ where, in  the coordinate of $\mathbb{P}^d$, $x_0, x_1$ are respectively in the $i$-th and $j$-th position. We have a natural morphism of $\mathcal{O}_{\mathbb{P}^d}$-mod, $l_*\mathcal{T}_{\mathbb{P}^1} \rightarrow l_*l^*\mathcal{T}_{\mathbb{P}^d}$. Moreover, let $\mathcal{I} \subset \oo_{\mathbb{P}^d}$ the sheaf ideal cutting out $\mathbb{P}^1$ via the closed immersion $l$. Then, there is a natural  $\mathcal{O}_{\mathbb{P}^d}$-module isomorphism $\mathcal{T}_{\mathbb{P}^d}/\mathcal{I}\mathcal{T}_{\mathbb{P}^d} \sim l_*l^*\mathcal{T}_{\mathbb{P}^d}$. Thus we get a natural map $$\Gamma(\mathbb{P}^1,\mathcal{T}_{\mathbb{P}^1})=\Gamma(\mathbb{P}^d,l_*\mathcal{T}_{\mathbb{P}^1}) \rightarrow \Gamma(\mathbb{P}^d,\mathcal{T}_{\mathbb{P}^d}/\mathcal{I}\mathcal{T}_{\mathbb{P}^d})$$ 
mapping (locally) $\frac{\partial}{\partial{z_{}}} \mapsto \frac{\partial}{\partial{z_{ij}}} $, where $z$ is a local coordinate for $\mathbb{A}^1 \subset \mathbb{P}^1$. Thus, we have a commutative diagram 
\begin{equation*}
\begin{tikzcd}
\mathfrak{sl}_2 \arrow[r, "\partial\sigma"] \arrow[d, "s_{ij}", hook] & {\Gamma(\mathbb{P}^1,\mathcal{T}_{\mathbb{P}^1})} \arrow[d] \\
\mathfrak{sl}_{d+1} \arrow[r, "\partial\sigma"]                       & {\Gamma(\mathbb{P}^d,\mathcal{T}_{\mathbb{P}^d})}          
\end{tikzcd}
\end{equation*}
of Lie algebra morphisms, extending by construction to the analogous one for $\phi^{\mathcal{O}_X}$.
\end{ex}

\begin{ex}[$\G=\mathrm{GL}_{d+1,k}$, $X=\mathbb{P}_k^d$]\label{action_of_G_on_Pd}
In the case $\textbf{G}=\mathrm{GL}_{d+1,k}$, we consider the flag variety $\G_k/\mathbf{P}_{(1,d),k}$ given by taking the quotient over the parabolic subgroup associated to the partition $(1,d)$. 
Let \begin{equation}
    \sigma \colon \textbf{G}_k \times \mathbb{
P}_k^d \rightarrow \mathbb{P}_k^d
\end{equation} be the action given on closed points by, 
$$(g, [z_0:\dots:z_d]) \mapsto [z_0:\dots:z_d]g^{-1} \text{ for any } g \in \textbf{G}_k.$$ 
Thus, the induced contragradient action is given by $g.f(z)=f(zg)$ for $f \in \mathcal{O}_{\mathbb{P}_k^d}$.
The Lie algebra $\mathfrak{g}_k$ acts on $\Gamma(U,\oo_{\mathbb{P}^d})$, for any open $U \subset \mathbb{P}^d$, via the canonical differential action $\partial\sigma$.

For any $\alpha \in \Phi^+$ we compute the image of $L_{\alpha} \in \mathfrak{g}_{\alpha}$ and $Y_{-\alpha} \in \mathfrak{g}_{-\alpha}$ under $\partial\sigma$.
Let $f \in \oo_{\mathbb{P}^d}(U)$. Write $\alpha=\alpha_{ij}$ for some $0 \leq i<j \leq d$. For any $0 \leq l \leq d$, put $V_l:=D_+(z_l)$.

It suffices to compute $\partial\sigma(L_{\alpha})(f)$, resp. $\partial\sigma(Y_{-\alpha})(f)$ in the local chart $V_i$, resp. $V_j$. Therefore, we have
\begin{align}
    \partial\sigma(L_{\alpha})(f)_{\mid V_i}&= \frac{d}{d\epsilon}_{\mid \epsilon=0}\sigma(1+L_{\alpha_{ij}}\epsilon)(f) \\
    &=\frac{d}{d\epsilon}_{\mid \epsilon=0}f(z_{0i},\dots,\epsilon+z_{ji},\dots,z_{di}) \\
    &=\frac{\partial}{\partial z_{ji}}f,
\end{align}
where we wrote $f$ as a function in the local coordinates $(z_{li})_{l \neq i}$ of $V_i$. Similarly we obtain
\begin{equation}
    \partial\sigma(Y_{-\alpha})(f)_{\mid V_j}=\frac{\partial}{\partial z_{ij}}f.
\end{equation}
On the intersection $V_i \cap V_j$ we have the relation $z_{ij}=z_{ji}^{-1}$, thus follows that $\frac{\partial}{\partial z_{ji}}f=-z_{ij}^2\frac{\partial}{\partial z_{ij}}f$ on $V_i \cap V_j$. 
Following this computation, we can consider $L_{\alpha},Y_{-\alpha}$ as derivations after taking their images under $\phi^{\mathcal{O}_{\mathcal{X}}}$,
i.e. we let
\begin{equation}\label{def_y_as_derivative}
\begin{array}{rcl}
 -z_{ij}^2\frac{\partial}{\partial z_{ij}}=\phi^{\mathcal{O}_{\mathcal{X}}}(L_{\alpha_{ij}})=: y_{\alpha_{ij}} =: y_{ij} & \mbox{if} & \alpha_{ij} \in \Phi^+ \\
\frac{\partial}{\partial z_{ij}}=\phi^{\mathcal{O}_{\mathcal{X}}}(Y_{-\alpha_{ij}})=: y_{-\alpha_{ij}}=: y_{ji} & \mbox{if} & \alpha_{ij} \in \Phi^+
\end{array}
\end{equation}
  Then the description of $\Gamma(\mathbb{P}_k^d,\mathcal{D}_{\mathbb{P}_k^d}) $ can be given locally as crystalline Weyl algebra generated by $z_{ij}$ and differentials $y_{\alpha_{ij}}^{[m]}$. 
More precisely, 
\begin{lemma}
For any $0\leq j\leq d$,
\begin{equation*} 
   \Gamma( V_j,\mathcal{D}_{\mathbb{P}_k^d})\simeq k[z_{ij} \mid i \neq j ]\big\langle y_{\mathrm{sgn}(i-j)\alpha_{ij}}^{[m]} \mid m \geq 0,\text{ } \alpha_{ij} \in \Phi\big\rangle.
\end{equation*}

\end{lemma}
\begin{proof}
 The open subvariety $V_j=D_+(z_j) \subset \mathbb{P}_k^d $ is a local chart of $\mathbb{P}_k^d$  with local coordinates $z_i=z_{ij}$, for $0\leq i \neq j \leq d$, isomorphic to $\mathbb{A}_k^d$. Also, $y_{\mathrm{sgn}(i-j)\alpha_{ij}}=\frac{\partial}{\partial z_{ij}}$ if $i<j$ or $-z_{ij}^2\frac{\partial}{\partial z_{ij}}$ if $i>j$. Then, the result follows by   \Cref{smith2.7}. 
\end{proof}
\end{ex}
\section{Applications}\label{sec_main_appl}
In this section let $k$ be a finite field of characteristic $p$ and $\mathbb{P}^d:=\mathbb{P}_k^d.$ Here, we are going to apply the content of  \Cref{secwittdiffop} to lift the operators above to differential operators over the Witt scheme $\wn(\mathbb{P}^d)$.  The action of $\G_k=:\G$ on $\mathbb{P}^d$, as described by $\sigma$ in  \Cref{action_of_G_on_Pd}, is considered. 
If $g \in G$ is a $k$-rational point and $U \subset \mathbb{P}^d$ is an open  such that $\sigma_g(U)=U$, then the ${G}$-linearization \eqref{isofasci} induces an isomorphism \footnote{As abuse of notation, here we use $\sigma_g$ to denote the action on the structure sheaf $\oo_{\mathbb{P}^d}$.}
\begin{equation}
\mathcal{D}_{\mathbb{P}^d}(U) \rightarrow \mathcal{D}_{\mathbb{P}^d}(U), \quad \eta \mapsto {}^{\sigma_g}\eta: (f \mapsto (\sigma_{g^{-1}}\circ \eta \circ \sigma_g)(f)).
\end{equation}  

For any $\alpha \in \Phi^+$, $s \in \mathbb{N}$,  let $\w_n(y_{-\alpha_{}}^{[s]})$ (resp. $\w_n(y_{\alpha_{}}^{[s]})$) be the  image of  $Y_{-\alpha_{}}^{[s]} \in \dist(\G)$ (resp. $L_{\alpha}^{[s]}$) under the map  \eqref{teich_lift_diff}. Those are elements of $\Gamma(\mathbb{P}^d,\mathcal{D}_{\wn(\mathbb{P}^d)})$. 

For any $\alpha \in \Phi$, write $\alpha=\alpha_{ij}$ with $i \neq j$ ($0 \leq i,j \leq d$). We introduce $$z_{\alpha}=z_{ij}:=\frac{z_i}{z_j} \in k(\mathbb{P}^d)=\mathrm{Frac}(\oo_{\mathbb{P}^d}(V_{0,\dots,d})), $$
where $V_{0,\dots,d}=\cap_{l=0}^d V_l$,  $V_l:=D_+(z_l) \subset \mathbb{P}^d$ for any $0 \leq l \leq d$ and $k(\mathbb{P}^d)$ is the function field of $\mathbb{P}^d$.

The sheaves $\wn\mathcal{O}_{V_j}$ have a natural (left) $\mathcal{D}_{\wn(V_j)}$-module structure, given by evaluating a differential operator to functions, thus the canonical map $ \mathcal{D}_{\wn(\mathbb{P}_k^d)} \rightarrow \mathcal{D}_{\wn(V_j)}$ induces a (left) $\mathcal{D}_{\wn(\mathbb{P}_k^d)}$-mod structure.

 Let $\mathbf{P}=\mathbf{P}_j:=\mathbf{P}_{(j+1,d-j)}$ be the maximal parabolic subgroup of $\G$ associated to the partition $d+1=(j+1)+(d-j)$ with Levi decomposition $\mathbf{P}_j=\mathbf{L}_j\mathbf{U}_j$.

\begin{lemma}
For any $\beta \in \Phi^+$, the differential operators $z_{\beta}^{p-1}y_{-\beta}^{[p]} \in  \mathcal{D}(\mathbb{A}^d)$ are elements of $\mathcal{D}(\mathbb{P}^d)$.
\end{lemma}
\begin{proof}
We notice that if $\beta=\alpha_{il}$ (with $i<l$), we have  $ z_{\beta}^{p-1}y_{-\beta}^{[p]} \in \mathcal{D}(D_+(z_l)) \simeq \mathcal{D}(\mathbb{A}^d)$ and under the identification $k[\frac{z_0}{z_l},\dots,\frac{z_d}{z_l}]\simeq k[\frac{z_i}{z_l}]\otimes_k k[\frac{z_r}{z_l}\mid r \neq i,l]$ we have $z_{\beta}^{p-1}y_{-\beta}^{[p]}(k[\frac{z_0}{z_l},\dots,\frac{z_d}{z_l}])=z_{\beta}^{p-1}y_{-\beta}^{[p]}(k[\frac{z_i}{z_l}])\otimes 1$. Thus we can assume $d=1$. Hence, let $z_{\beta}=z$ and $y_{-\beta}=\partial_z$. Then, it follows by
\Cref{lemma_global_D_P^1}.
\end{proof}

  We can get a generalization of \cite[Proposition 2.1.5.3]{Ku16} and \cite[Proposition 3.11]{Orlik21} in the case $\mathcal{E}=\mathcal{O}_{\mathbb{P}_k^d}$, by the following. Let $\mathcal{D}_n:=\Gamma(\mathbb{P}^d,\mathcal{D}_{\wn(\mathbb{P}^d)})$.
\begin{prop}\label{mainfinitegeneration}
Assume that $\mathrm{char}(k)\neq 2$. Then, the $P_j$-module $\tilde{\h}^{d-j}_{\mathbb{P}_k^j}(\mathbb{P}_k^d,\wn(\mathcal{O}_{\mathbb{P}_k^d}))$ admits a submodule $N_{n,j}$ that is a finitely generated $P_j$-module over $\wn(k)$ 
and a $\wn(k)$-linear epimorphism of $\mathcal{D}_n$-modules
\begin{equation}\label{epiwn}
    \rho_{n,j}: \mathcal{D}_n \otimes_{\wn(k)}N_{n,j} \twoheadrightarrow \tilde{\h}^{d-j}_{\mathbb{P}_k^j}(\mathbb{P}_k^d,\wn(\mathcal{O}_{\mathbb{P}_k^d})).
\end{equation}
\end{prop}
Before proving the proposition we need some technical preparation. Let $S \subset \zz^{d+1}$ be a finite set. Let $\wn(k)[\underline{T}]=\wn(k)[T_0,\dots,T_d]$ and $\underline{T}:=(T_0,\dots,T_d)$.
Consider the polynomial algebra  $\wn(k)[\underline{T}^{\underline{m}} \mid \underline{m} \in S]$ as $\wn(k)$-module. Let
$\wn(k)[\underline{T}^{\underline{m}} \mid \underline{m} \in S ]^{}_{r}$
be the $\wn(k)$-submodule of $\wn(k)[\underline{T}^{\underline{m}} \mid \underline{m} \in S ]$ consisting of homogeneous polynomials of degree $r$. Denote by 
 $$ \wn(k)[\underline{T}^{\underline{m}} \mid \underline{m} \in S]^{\ge 1}_{\leq p^l}:=\bigcup_{r=0}^{l}\wn(k)[\underline{T}^{\underline{m}} \mid \underline{m} \in S ]_{p^r}.$$ 
 The latter is a free $\wn(k)$-module with a $\wn(k)$-basis consisting of elements of the form 
 $$\prod_{\underline{m} \in S}   (\underline{T}^{\underline{m}})^{i_{\underline{m}}}$$
  with $\sum_{\underline{m} \in S}i_{\underline{m}} = p^r$, $0 \leq r \leq l$
  and $i_{\underline{m}}\geq 0$ for any $\underline{m} \in S$.
  \begin{defn}
      Let  $k[\underline{z}]^{\leq}_{n}(S;1) \subset \wn(k[\underline{z}])$ be the $\wn(k)$-submodule
      generated by the elements of the form $\ve^{l}(\prod_{\underline{m} \in S}   (\underline{T}^{\underline{m}})^{i_{\underline{m}}})$ where $T_i:=[z_i]$ for $0\leq i \leq d$ and
      $\prod_{\underline{m} \in S}   (\underline{T}^{\underline{m}})^{i_{\underline{m}}} \in  \wn(k)[\underline{T}^{\underline{m}} \mid \underline{m} \in S]^{\ge 1}_{\leq p^l} $
  \end{defn}
  
   \begin{rem}
       Since $S$ is finite, $k[\underline{z}]^{\leq}_{n}(S;1)$ is a finitely generated $\wn(k)$-module.
   \end{rem} 
   
  \begin{rem}
We may choose a bijection between the sets $S$ and $\{1,\dots, |S|\}$. Thus, for any $\underline{m} \in S$ corresponds a unique natural number $s=s(\underline{m})$. Under this bijection, the variable $\underline{T}^{\underline{m}}$ corresponds to a variable $X_s$. Then, $k[\underline{z}]_n^{\leq}(S;1)$ is generated over $\wn(k)$ by $\ve^l(f)$ where $f$ runs over homogeneous monomials $f \in \w_{n-l}(k)[X_1,\dots,X_{|S|}]$ of degree $p^r$ for some $0 \leq r \leq l$.
 \end{rem}
  \begin{defn}
      Let $A$ be a unitary commutative $k$-algebra. 
      Let $a_1,\dots,a_r \in A$ be pairwise distinct. For any $i \in \mathbb{N}$, define $A_n(\{a_1,\dots,a_r\};i)$ be the $\wn(k)$-submodule of $\wn(A)$ generated by $\ve^l(\prod_{j=1}^r[a_j]^{m_j})$ for $0 \leq l < n$ with $\sum_{j=1}^r m_j=p^li$.
  \end{defn}
  \begin{rem}
      Let $A=k[\underline{z}]$ as above and $S \subset \zz^{d+1}$ be a finite set. Then,  for each $\underline{m} \in S$ corresponds a distinct monomial $\underline{z}^{\underline{m}} \in A$. We see by construction that \begin{equation}
          k[\underline{z}]_n(\{\underline{z}^{\underline{m}} \mid \underline{m} \in S\};1) \subset k[z]_n^{\leq}(S;1)
      \end{equation} 
      is a $\wn(k)$-submodule. The symbol "$\leq$" underlines the condition on the degree of an element of the form $\ve^l(f)$.  
  \end{rem}


\begin{lemma}\label{lemmaA_n(a,b;i)}
Let $A$ be a unitary commutative $k$-algebra. Let $a,b \in A$ two distinct elements and $[-] \colon A \rightarrow \wn(A)$ be the Teichm\"uller map.
 Then, $[a+b]^i \in A_n(\{a,b\};i).$
\end{lemma}   
\begin{proof}
If $q(a,b)=\sum_{j=0}^{r}c_ja^jb^{r-j} \in k[a,b]$ is a homogeneous polynomial of degree $r\geq 1$,  there exist polynomials $q_1(a,b),\dots, q_{n-1}(a,b) \in k[a,b] $ such that
$$[q(a,b)]=\left(\sum_{j=0}^{r}[c_j][a]^j[b]^{r-j}\right)+\sum_{i=1}^{n-1}\ve^i([q_{i}(a,b)]).$$ Moreover, for any $i=1,\dots,n-1$,  $q_i(a,b)$ is homogeneous of degree $p^ir$. Indeed, the existence of some polynomials is clear;  we need to prove that are of that kind. Applying the $i$-th ghost map $w_i$  to both sides of the latter equality, we get the relations
$$\left(\sum_{j=0}^{r}c_j^{p^i}a^{p^ij}b^{p^ir-p^ij}\right)+pq_1(a,b)^{p^{i-1}}+p^2q_2(a,b)^{p^{i-2}} + \dots + p^{i}q_{i}(a,b)= q(a,b)^{p^i}$$
for each $i=1,\dots,n-1$. By induction on $i$, it shows that $q_{i}(a,b)$ is a homogeneous polynomial of degree $p^ir$.
Now consider $q(a,b)=(a+b)^i$. By the reasoning above we can write
$$ [a+b]^i=t_0(a,b)+\sum_{l=1}^{n-1}\ve^l([q_{l}(a,b)])$$
with $\mathrm{deg}(q_l(a,b))=p^li$,  $t_0(a,b) \in A_n(\{a,b\};i) $ and $[q_l(a,b)] \in \w_{n-l}(k[a,b])$. Then, we can apply the same argument to each $q_l(a,b)$:

 We have an equality
$$ [q_l(a,b)]= t_{l,0}(a,b)+\sum_{l'=1}^{n-l-1}\ve^{l'}([q_{l,l'}(a,b)])$$ 
where $\mathrm{deg}(q_{l,l'}(a,b))=p^{l+l'}i$, $t_l$ is a homogeneous polynomial of degree $p^li$ in the variables $[a]$ and $[b]$, and $[q_{l,l'}(a,b)] \in \w_{n-l-l'}(k[a,b])$.

 Thus we can write 
$$ \ve^l([q_l(a,b)])= t_l(a,b)+\sum_{l'=l+1}^{n-1}\ve^{l'}([q_{l,(l'-l)}(a,b)])$$ 
where $\mathrm{deg}(q_{l,(l'-l)}(a,b))=p^{l'}i$, $t_l=\ve^l(t_{l,0}) \in A_n(\{a,b\};i)$, and $[q_{l,(l'-l)}(a,b)] \in \w_{n-l'}(k[a,b])$ for $n > l' \geq l+1$. Hence, we can iterate the same argument to every polynomial $q_{l,l'}(a,b)$, until we get an expression of $[a+b]^i$ in $A_n(\{a,b\};i)$.
 \end{proof}
 \begin{corol}\label{corolA_n(S;i)}
 Let $A$ be a unitary commutative $k$-algebra with $\mathrm{char}(k) \neq 2$. Let $a_1,\dots,a_r \in A$ pairwise distinct.
 Then, $[\sum_{j=1}^r a_j]^i \in A_n(\{a_1,\dots,a_r\};i).$
 \end{corol}
 \begin{proof}
  We proceed by induction for $r \geq 2$, using the previous lemma as the base case.  Then, assume $r>2$. Up to permutations of the indices of $a_i$, we can assume $b_r=a_1+\dots + a_{r-1} \neq a_r$ (otherwise $2a_j=\sum_{i=1}^r a_i$ for any $j$, contradicting the hypothesis of having distinct $a_j$'s). By Lemma  
   \ref{lemmaA_n(a,b;i)}, it follows that
 $ [a_r+b_r]^i=[\sum_{j=1}^r a_j]^i \in A_n(\{b_r,a_r\}; i)$.  Thus $[a_r+b_r]^i$ is a linear combination of elements of the form $\ve^l([a_1+\dots+a_{r-1}]^{j}[a_r]^{p^li-j})$ for $0 \leq l < n$.  By inductive hypothesis $[a_1+\dots+a_{r-1}]^{j} \in A_{n-l}(\{a_1,\dots,a_{r-1}\};j)$, thus it is a linear combination of elements of the form $\ve^s(\prod_{u=1}^{r-1} [a_u]^{m_u})$ with  $\sum m_u=p^sj$ and  $0 \leq s < n-l$. Since $\ve^l(\ve^s(\prod_{u=1}^{r-1} [a_u]^{m_u})[a_r]^{p^li-j})=\ve^{l+s}([a_r]^{p^{l+s}i-p^sj}\prod_{u=1}^{r-1} [a_u]^{m_u}) \in A_n(\{a_1,\dots,a_r\}; i)$, it follows also that $[a_1+ \dots + a_r]^i \in A_n(\{a_1,\dots, a_r\}; i)$.
 \end{proof}
 
\begin{lemma}\label{multvs_i}
Let $A$ be a $k$-algebra, and $m,r$ be non negative integers. Let $a_1,\dots,a_m \in A$. Let $\underline{d}_1,\dots,\underline{d}_r \in \mathbb{N}^{m}$ and denote by $[\underline{a}]^{\underline{d}_i}:=[a_1]^{(\underline{d}_i)_1}\dots [a_m]^{(\underline{d}_i)_m} \in \wn(A)$ for any $1 \leq i \leq r$. Let $s_1,\dots,s_r \in \mathbb{N}$. Assume $s_r=\mathrm{max}_i \{s_i\}$. Then
\begin{equation}
\prod_{i=1}^r \ve^{s_i}([\underline{a}]^{\underline{d}_i})=p^{s_1+\dots+s_{r-1}}\ve^{s_r}\big([\underline{a}]^{\sum_{i=1}^r p^{s_r-s_i}\underline{d}_i}\big)
\end{equation} 
where $p^{s_r-s_i}\underline{d}_i \in \mathbb{N}^m$ is defined by $(p^{s_r-s_i}\underline{d}_i)_j:=p^{s_r-s_i}(\underline{d}_i)_j$ for any $1 \leq j \leq m$.
\end{lemma} 
\begin{proof}
Without loss of generality, we can assume $s_1 \leq  \dots \leq s_r$. We can proceed by induction on $r \geq 1$. The base case is trivially true. Then, we have
\begin{align*}
\prod_{i=1}^r \ve^{s_i}([\underline{a}]^{\underline{d}_i}) 
&=p^{s_1+\dots+s_{r-2}}\ve^{s_{r-1}}([\underline{a}]^{\sum_{i=1}^{r-1} p^{s_{r-1}-s_i}\underline{d}_i})\ve^{s_r}([\underline{a}]^{\underline{d}_r})  \\&=
p^{s_1+\dots+s_{r-2}}\ve^{s_{r-1}}\big([\underline{a}]^{\sum_{i=1}^{r-1} p^{s_{r-1}-s_i}\underline{d}_i}F^{s_{r-1}}\ve^{s_r}([\underline{a}]^{\underline{d}_{r}})\big) \\&=
p^{s_1+\dots+s_{r-1}}\ve^{s_{r-1}}\big([\underline{a}]^{\sum_{i=1}^{r-1} p^{s_{r-1}-s_i}\underline{d}_i}\ve^{s_r-s_{r-1}}([\underline{a}]^{\underline{d}_{r}})\big)\\&=
p^{s_1+\dots+s_{r-1}}\ve^{s_{r-1}}\big(\ve^{s_r-s_{r-1}}([\underline{a}]^{d_r+p^{s_r-s_{r-1}}\sum_{i=1}^{r-1} p^{s_{r-1}-s_i}\underline{d}_i})\big)
\\&=p^{s_1+\dots+s_{r-1}}\ve^{s_r}\big([\underline{a}]^{\sum_{i=1}^r p^{s_r-s_i}\underline{d}_i}\big).
\end{align*}
\end{proof}
 For any $\underline{d}=(d_1,\dots,d_m) \in \mathbb{N}^m$ let the \textit{degree of} $\underline{d}$ be  $\mathrm{deg}(\underline{d}):=\sum_{i=1}^m d_i \in \mathbb{N}$.
 
 The following is immediate:
  \begin{corol}\label{corolmultvs_i}
 Let $m \in \mathbb{N}$, $A=k[\underline{z}]$ for $\underline{z}=(z_1,\dots,z_{m})$, $a_i=z_i$ with $[z_i]=T_i$, $i=1,\dots,m$. Finally let $s_i$, $\underline{d}_i,r$ as in  \Cref{multvs_i}. Then, $\prod_{i=1}^r \ve^{s_i}(\underline{T}^{\underline{d}_i})=\ve^{s}(f)$ for a homogeneous monomial $f \in \w_{n-s}(k)[\underline{T}]$ of degree  $\mathrm{deg}(f)=\sum_{i=1}^rp^{s-s_i}\mathrm{deg}(\underline{d}_i)$,  where $s=\mathrm{max}_i\{s_i\}$. $\square$
 \end{corol}

\begin{proof}[Proof of  \Cref{mainfinitegeneration}]
Rephrasing the statement, it is equivalent to say that the cohomology group
$\tilde{\h}^{d-j}_{\mathbb{P}_k^j}(\mathbb{P}_k^d,\wn(\mathcal{O}_{\mathbb{P}_k^d}))$ is generated as $\mathcal{D}_n^{}$-module by a $P_j$-module $N_{n,j}$ finitely generated over $\wn(k)$. We define the following sets
\begin{equation*}
\begin{array}{lc}
     I:=\{(m_0,\dots,m_d) \in \zz^{d+1} \mid m_0,\dots,m_j \geq 0, m_{j+1},\dots,m_d<0, \sum_{i}m_i=0\}, \\
     I_j:=\{(m_i)_i \in I \mid m_i=-1 \quad \forall i \geq j+1\}.
\end{array}
\end{equation*} 
Notice that the set $I$ is infinite and closed under the sum (component-wise), while $I_j$ is a finite set. Firstly, we need to prove that in the category of $\wn(k)$-modules the following holds:

    \begin{equation}\label{redloccohomcomputation}
     \tilde{\h}^{d-j}_{\mathbb{P}_k^j}(\mathbb{P}_k^d,\wn(\mathcal{O}_{\mathbb{P}_k^d}))= \left\{\begin{array}{lcr}
     \wn(k[z_1,\dots,z_d])/{\wn(k)} & \mbox{if} & d-j=1, \\
     \wn\Big(\bigoplus_{(m_0,\dots,m_d) \in I}k \cdot z_0^{m_0}\dotsc z_j^{m_j}z_{j+1}^{m_{j+1}}\dotsc z_d^{m_d}\Big) & \mbox{if} & d-j \geq 2, \\
     0 & \mbox{if} & d=j,
     \end{array}\right.
     \end{equation}
      Notice that  $\bigoplus_{(m_0,\dots,m_d) \in I}k \cdot z_0^{m_0}\dotsc z_j^{m_j}z_{j+1}^{m_{j+1}}\dotsc z_d^{m_d}$ is not a unitary ring. Thus we set by definition, in the second line of  \eqref{redloccohomcomputation},
     $$\wn\Big(\bigoplus_{(m_0,\dots,m_d) \in I} k \cdot z_0^{m_0}\dotsc z_j^{m_j}z_{j+1}^{m_{j+1}}\dotsc z_d^{m_d}\Big):= \text{ $\wn(k)$-module generated by  } \ve^l([\underline{z}^{\underline{m}}]) \text{,  }$$
    where  $\underline{z}^{\underline{m}}:=\prod_{i=0}^d z_i^{m_i} $ and $0 \leq l < n$. 
     To check the \eqref{redloccohomcomputation}, we use the  the long exact sequence
     \begin{multline*}
     {\h}^{d-j-1}_{}(\mathbb{P}_k^d,\wn(\mathcal{O}_{\mathbb{P}_k^d}))\rightarrow
     {\h}^{d-j-1}_{}(\mathbb{P}_k^d\backslash \mathbb{P}_k^j,\wn(\mathcal{O}_{\mathbb{P}_k^d}))\rightarrow
     {\h}^{d-j}_{\mathbb{P}_k^j}(\mathbb{P}_k^d,\wn(\mathcal{O}_{\mathbb{P}_k^d}))\rightarrow\\ \rightarrow
     {\h}^{d-j}_{}(\mathbb{P}_k^d,\wn(\mathcal{O}_{\mathbb{P}_k^d}))\rightarrow
     {\h}^{d-j}_{}(\mathbb{P}_k^d\backslash \mathbb{P}_k^j,\wn(\mathcal{O}_{\mathbb{P}_k^d}))=0.
     \end{multline*}
    If $d-j \geq 2$, then ${\h}^{d-j-1}_{}(\mathbb{P}_k^d,\wn(\mathcal{O}_{\mathbb{P}_k^d}))=0$, thus $$\tilde{\h}^{d-j}_{\mathbb{P}_k^j}(\mathbb{P}_k^d,\wn(\mathcal{O}_{\mathbb{P}_k^d}))\simeq {\h}^{d-j-1}_{}(\mathbb{P}_k^d\backslash \mathbb{P}_k^j,\wn(\mathcal{O}_{\mathbb{P}_k^d})).$$  If $d-j=1$, then ${\h}^{d-j-1}_{}(\mathbb{P}_k^d,\wn(\mathcal{O}_{\mathbb{P}_k^d}))=\wn(k)$. Hence, we have the short exact sequence 
    \begin{equation*}
    0 \rightarrow \wn(k) \rightarrow \h^0(\mathbb{P}_k^d\backslash \mathbb{P}_k^{d-1},\wn(\mathcal{O}_{\mathbb{P}_k^d})) \rightarrow \tilde{\h}^{1}_{\mathbb{P}_k^{d-1}}(\mathbb{P}_k^d,\wn(\mathcal{O}_{\mathbb{P}_k^d})) \rightarrow 0.
    \end{equation*}
     Since $\mathbb{P}_k^d\backslash \mathbb{P}_k^{d-1}\simeq \mathbb{A}_k^d$, we get the first line of \eqref{redloccohomcomputation}. The case $d=j$ is trivial. Now consider the open cover of $$\mathbb{P}_k^d\backslash \mathbb{P}_k^{j}= \bigcup_{i=j+1}^{d}D_{+}(z_i).$$
    
  The \v{C}ech complex associated to this covering has its highest components of degree $d-j-2$ and $d-j-1$, given by
  \begin{equation}
  C^{d-j-2}=\bigoplus_{i_1=j+1}^d \Gamma\left(\bigcap_{\substack{i \neq i_1 \\ j+1 \leq i \leq d}}D_{+}(z_i),\wn(\mathcal{O}_{\mathbb{P}_k^d}) \right) \xrightarrow{\partial^{d-j-2}} C^{d-j-1}= \Gamma\left(\bigcap_{i=j+1}^d D_{+}(z_i),\wn(\mathcal{O}_{\mathbb{P}_k^d})\right)
  \end{equation}
  Then, ${\h}^{d-j-1}_{}(\mathbb{P}_k^d\backslash \mathbb{P}_k^j,\wn(\mathcal{O}_{\mathbb{P}_k^d}))=\mathrm{coker}(\partial^{d-j-2})$.
  We have that 
  $$\Gamma\left(\bigcap_{\substack{i \neq i_1 \\ j+1 \leq i \leq d}}D_{+}(z_i),\wn(\mathcal{O}_{\mathbb{P}_k^d}) \right)=\wn((k[z_0,\dots,z_d]_{z_{j+1}\dots \hat{z_{i_1}}\dots z_d})^0)$$
  where the superscript $0$ means that we are considering the degree $0$ elements in the localization. By the natural inclusion $$\wn((k[z_0,\dots,z_d]_{z_{j+1}\dots \hat{z_{i_1}}\dots z_d})^0) \subset \wn(k[z_0,\dots,z_d])_{[z_{j+1}\dots \hat{z_{i_1}}\dots z_d]},$$ we see that $\wn((k[z_0,\dots,z_d]_{z_{j+1}\dots \hat{z_{i_1}}\dots z_d})^0)$ is generated as $\wn(k)$-module by the elements of the form
  \begin{equation}
  \ve^l([z_0^{m_0}\dots z_j^{m_j}z_{j+1}^{m_{j+1}+p^lm'_{j+1}}\dots z_{d}^{m_{d}+p^lm'_{d}}]) \text{ such that } \left\{\begin{array}{lc}
  0 \leq l \leq n-1, \\
  \sum_{i=0}^d m_i+p^l\sum_{i=j+1}^d m'_{i}=0, \\
  m_{i_1}+p^lm'_{i_1} \geq 0, \\
  m'_{j+1}, \dots, m'_d \leq 0, \\
  m_0,\dots,m_d \geq 0.
  \end{array}\right.
  \end{equation}
  Analogously, 
  \begin{equation*}
  \Gamma\left(\bigcap_{i=j+1}^d D_{+}(z_i),\wn(\mathcal{O}_{\mathbb{P}_k^d})\right)=\wn((k[z_0,\dots,z_d]_{z_{j+1}\dots z_d})^0) \subset \wn(k[z_0,\dots,z_d])_{[z_{j+1}\dots z_d]}
  \end{equation*} is generated by 
  \begin{equation}
  \ve^l([z_0^{m_0}\dots z_j^{m_j}z_{j+1}^{m_{j+1}+p^lm'_{j+1}}\dots z_{d}^{m_{d}+p^lm'_{d}}]) \text{ such that } \left\{\begin{array}{lc}
  0 \leq l \leq n-1, \\
  \sum_{i=0}^d m_i+p^l\sum_{i=j+1}^d m'_{i}=0, \\
  m'_{j+1},\dots,m'_d \leq 0, \\
  m_0,\dots,m_d \geq 0
  \end{array}\right.
  \end{equation}
  as $\wn(k)$-module. 
  It follows that $\mathrm{coker}(\partial^{d-j-2})$ is generated as $\wn(k)$-module by the elements  
  \begin{equation}
  \ve^l([z_0^{m_0}\dots z_j^{m_j}z_{j+1}^{m_{j+1}+p^lm'_{j+1}}\dots z_{d}^{m_{d}+p^lm'_{d}}]) \text{ such that } \left\{\begin{array}{lc}
  0 \leq l \leq n-1, \\
  \sum_{i=0}^d m_i+p^l\sum_{i=j+1}^d m'_{i}=0,\\
  m'_{j+1}, \dots, m'_d \leq 0, \\
  m_i < p^l|m'_{i}|,j+1 \leq i \leq d, \\
  m_0,\dots,m_d \geq 0
  \end{array}\right.
  \end{equation}
Finally, we see that these elements generate (since for any fixed $l$, $m_i-p^l|m'_i|$ take all the negative integer values) $\wn\Big(\bigoplus_{(m_0,\dots,m_d) \in I} k \cdot z_0^{m_0}\dotsc z_j^{m_j}z_{j+1}^{m_{j+1}}\dotsc z_d^{m_d}\Big)$. 

Set $Y=\bigcap_{i=j+1}^d D_{+}(z_i)$ 
and
$Y_{i_1}=\bigcap_{\substack{i \neq i_1 \\ j+1 \leq i \leq d}}D_{+}(z_i)$ for any $j+1 \leq i_1 \leq d$.

We regard the $\wn(k)$-modules $C^{d-j-2}, C^{d-j-1}$ respectively with the natural structure of 
$\mathcal{D}$-module, induced respectively by 
$ \Gamma(\mathbb{P}_k^d,\mathcal{D}_{\wn(\mathbb{P}_k^d)})\rightarrow \Gamma(Y_{i_1},\mathcal{D}_{\wn(\mathbb{P}_k^d)})$ and $ \Gamma(\mathbb{P}_k^d,\mathcal{D}_{\wn(\mathbb{P}_k^d)})\rightarrow \Gamma(Y,\mathcal{D}_{\wn(\mathbb{P}_k^d)})$. Thus also the local cohomology $ \tilde{\h}^{d-j}_{\mathbb{P}_k^j}(\mathbb{P}_k^d,\wn(\mathcal{O}_{\mathbb{P}_k^d}))$ inherits a $\Gamma(\mathbb{P}^d,\mathcal{D}_{\wn(\mathbb{P}_k^d)})$-module structure.\\
From now on, let  $ 0 \leq a \leq j$ and $j+1 \leq b \leq d$.\\
Write a Teichm\"uller representative $T_{ab}$ of $z_{ab}$ as $\frac{T_a}{T_b}:=T_{ab}=:T_aT_b^{-1}$. If $\underline{m}:=(m_0,\dots,m_d)$, an element $$\underline{T}^{\underline{m}}:=T_0^{m_0}\dotsc T_j^{m_j}T_{j+1}^{m_{j+1}}\dotsc T_d^{m_d} \in \wn\Big(\bigoplus_{(m_0,\dots,m_d) \in I} k \cdot z_0^{m_0}\dotsc z_j^{m_j}z_{j+1}^{m_{j+1}}\dotsc z_d^{m_d}\Big)$$ is well-defined whenever $(m_0,\dots,m_d) \in I$. With this notation, the action of $\mathcal{D}_n^{}$ on an element of $\wn\Big(\bigoplus_{(m_0,\dots,m_d) \in I} k \cdot z_0^{m_0}\dotsc z_j^{m_j}z_{j+1}^{m_{j+1}}\dotsc z_d^{m_d}\Big)$ induces an action on the set $I$.  As matter of notation, write just $m$ to denote an element of $I_j$.  
 Let $N_{n,j}$  be the finitely generated $\wn(k)$-module $k[\underline{z}]^{\leq}_{n}(I_j;1)$.  For $n=1$, $N_{1,j}=\bigoplus_{(m_0,\dots,m_j,-1,\dots,-1) \in I_j} k \cdot z_0^{m_0}\dotsc z_j^{m_j}z_{j+1}^{-1}\dotsc z_d^{-1}$  is a $P_j$-module.
 
 If $n >1$, any $x \in N_{n,j}$ can be written as
 $x= \sum_{l=0}^{n-1} a_{l,\underline{i}}\ve^{l}(\prod_{m \in I_j}   (\underline{T}^{{m}})^{i_{{m}}})$
 where $\underline{i}=(i_m)_{m \in I_j}$ and $a_{l,\underline{i}} \in \wn(k)$. 
  If $g \in P_j$, then \footnote{The operation $g.$ denotes the action of $G$} 
  \begin{equation}\label{pjmodwn}
    g.\ve^{l}\left(\prod_{m \in I_j}   (\underline{T}^{m})^{i_{m}}\right)=\ve^{l}\left(\prod_{m \in I_j}   [g.\underline{z}^{m}]^{i_{m}}\right)=\ve^l\left(\prod_{m \in I_j} \left[\sum_{m' \in I_j}b_{m',i_{m} }\underline{z}^{m'}\right]^{i_{m}}\right)
  \end{equation}
  where $b_{m,i_m} \in k$ and $\sum_{m \in I_j}i_{m}=p^r$ for some $0\leq r \leq l$. For any $m \in I_j$ the corresponding $\left[\sum_{m' \in I_j}b_{m',i_{m} }\underline{z}^{m'}\right]^{i_{m}}$ lies in $k[\underline{z}]_{n-l}(I_j;i_{m})$ by  \Cref{corolA_n(S;i)}, thus it is a linear combination of elements of the form 
  $$\ve^{s_{m}}\left(\prod_{{m'} \in I_j}   (\underline{T}^{m'})^{r(m',s_{m} ,i_{m})}\right) $$
   where $r(m',s_{m} ,i_{m}) \in \mathbb{N}$ are such that  $\sum_{m' \in I_j}r(m',s_{m} ,i_{m})=p^{s_{m}}i_{m}$.
   Choose a bijection of $I_j$ with $\{1,2\dots,|I_j|\}$ and denote by $s_1,s_2,\dots,s_{|I_j|}$ (resp. $i_1,\dots,i_{|I_j|}$) the corresponding elements $s_{m}$ (resp. $i_{m}$ ) via the chosen bijection.
   The expression in the right hand side of \eqref{pjmodwn} is a linear combination of elements of the form 
   \begin{equation}\label{big_prod_v^l}
   \ve^l\left( \prod_{u=1}^{|I_j|}\ve^{s_{u}}\left(\prod_{{m'} \in I_j}   (\underline{T}^{m'})^{r(m',s_{u} ,i_{u})}\right) \right).
\end{equation}     
Let $s=\mathrm{max}_u s_u$ and $r(s_{u}, i_{u}) \in \mathbb{N}^{|I_j|} $ be given by letting  $r(s_{u} ,i_{u})_v:=r(v,s_{u} ,i_{u})$. By \Cref{corolmultvs_i} applied to $A=k[\underline{z}^{m} \mid m \in I_j]$ and $a_i=\underline{z}^{m}$ for $i=1,\dots,|I_j|$, where $m$ is the unique element of $I_j$ corresponding to $i$ under the chosen bijection, the quantity \eqref{big_prod_v^l} is equal to $\ve^{l+s}(f)$ for a homogeneous monomial $f \in \w_{n-l-s}(k)[\underline{T}^{m} \mid m \in I_j]$ of degree
 $$\mathrm{deg}(f)=\sum_{u=1}^{|I_j|} p^{s-s_u}\mathrm{deg}(r(s_{u} ,i_{u}))=\sum_{u=1}^{|I_j|} p^{s-s_u}\left(\sum_{v=1}^{|I_j|}r(v,s_{u} ,i_{u})\right)= \sum_{u=1}^{|I_j|}p^{s-s_u}p^{s_u}i_u=p^{s+r}.$$
 Since $s+r \leq s+l$, it shows that $\ve^{l+s}(f) \in k[\underline{z}]_n^{\leq}(I_j;1)$,
     thus  $g.\ve^{l}\left(\prod_{m \in I_j}   (\underline{T}^{m})^{i_{m}}\right)\in N_{n,j}.$ 
      It follows that $N_{n,j}$ is a $P_j$-module for any $n \geq 1$. 
Furthermore, since operators in $\mathcal{D}_n$ are compatible with Verschiebung maps (by the \eqref{derversch}), we need to prove the proposition for $n=1$. We achieve the result by applying the following procedure.

We say that a monomial $\underline{z}^{\underline{m}}=\prod_{m_i \in \underline{m} }z_i^{m_i} $ contains $z_i^{m_i}$ if it appears in the product with $m_i \neq 0$. Also, every monomial is parametrized by a scalar and the vector $\underline{m} \in I$, that we refer as the \textit{associated vector} to $\underline{z}^{\underline{m}}$. We recall that the differential operators $y_{il}^{[r]}$ act in such a way:
$$y^{[r]}_{il} : z_i^{m_i}z_l^{m_l}\prod_{\substack{ \underline{m} \in I \\ s\neq i,l}}z_s^{m_s} \longmapsto \binom{m_l}{r}z_i^{m_i+r}z_l^{m_l-r}\prod_{\substack{ \underline{m} \in I \\ s\neq i,l}}z_s^{m_s} $$
for any $i \neq l$, where for $m<0$, $\binom{m}{r}:=\frac{m(m-1)\dots (m-r+1)}{r!}$. 
\begin{enumerate}
    \item[Step 1:] Apply $y_{ab},\dots,y_{ab}^{p-1} \in \mathcal{D}^{}_{1}$ to (the monomials with associated vector) $m \in I_j$: starting by the element $z_0^{m_0}\dots z_j^{m_j}z_{j+1}^{-1}\dots z_{d}^{-1}$, one gets all monomials with associated vectors $(m_0,\dots,m_j,m_{j+1},\dots,m_d)$ such that  $m_{j+1},\dots,m_d \in \{-1,\dots,-p\}$ ; Note that $(-1)\dots (-m) \not \equiv 0 \pmod p$ if $m < p$.
    \item[Step 2:] Apply $y_{ab}^{[p]} \in \mathcal{D}^{}_{1}$ to $m \in I_j$:\\ We have $y_{ab}^{[p]}( z_a^{m_a} z_{b}^{-1}\prod_{\substack{ \underline{m} \in I \\ i \neq a,b}}z_i^{m_i})= (-1)^p z_a^{m_a+p} z_{b}^{-1-p}\prod_{\substack{ \underline{m} \in I \\ i \neq a,b}}z_i^{m_i}$ with $m_a+p \geq p$; 
    we wish to produce those 
    $$ z_a^{m_a} z_{b}^{-1-p}\prod_{\substack{ \underline{m} \in I \\ i \neq a,b}}z_i^{m_i}$$
     for which  $m_a < p$.  If $j=0$, the latter condition is empty since in this case $m_0 \geq p+1$ (since $a=0$, in this case $m_0 \geq 0$ and $m_i <0$ for any $i \neq 0$, with $m_0+(-p-1)\sum_{i \neq0,b}m_i=0$). Thus, we can assume $j>0$. For $x \leq j$ and $x \neq a$, apply $T_{ax}^{p-1}y_{xa}^{[p]} \in \mathcal{D}_1^{}$:
          \begin{multline*}T_{ax}^{p-1}y_{xa}^{[p]}(z_a^{m_a+p} z_{b}^{-1-p}\prod_{\substack{ \underline{m} \in I \\ i \neq a,b}}z_i^{m_i})=z_a^{p-1}z_x^{1-p}\binom{m_a+p}{p}z_x^{m_x+p}z_a^{m_a}z_b^{-1-p}\prod_{\substack{ \underline{m} \in I \\ i \neq a,b,x}}z_i^{m_i} \\
     =\binom{m_a+p}{p}z_a^{m_a+p-1}z_b^{-1-p}z_x^{m_x+1}\prod_{\substack{ \underline{m} \in I \\ i \neq a,b,x}}z_i^{m_i}.
     \end{multline*}
     Since we can assume $m_a<p$, we have $p \nmid \binom{m_a+p}{p}$.

     More generally, we have by induction that 
     
    $$(T_{ax}^{p-1}y_{xa}^{[p]})\big( z_a^{m_a+p-s} z_{b}^{-1-p}\prod_{\substack{ \underline{m} \in I \\ i \neq a,b}}z_i^{m_i}\big)= u \cdot z_a^{m_a+p-s-1} z_{b}^{-1-p}z_x^{m_x+1}\prod_{\substack{ \underline{m} \in I \\ i \neq a,b,x}}z_i^{m_i} $$ for $1 \leq s \leq m_a$ and $u= \binom{m_a+p-s}{p} \in \mathbb{Z}_{(p)}^{\times}$; For $p \geq s > m_a$, we have that $m_a+p-s < p$, thus we can apply $$y_{xa}(z_a^{m_a+p-s-1}z_b^{-1-p}z_x^{m_x+1}\prod_{\substack{ \underline{m} \in I \\ i \neq a,b,x}}z_i^{m_i})=(m_a+p-s)z_a^{m_a+p-s-1}z_b^{-1-p}\prod_{\substack{ \underline{m} \in I \\ i \neq a,b}}z_i^{m_i}$$ having non zero coefficient.
    
    \item[Step 3:] Restart from applying Step 1 to the associated vectors of the form  $(m_0,\dots,m_j,-(p+1),\dots,-(p+1)) \in I$ in place of $m \in I_j$. 
\end{enumerate}
It remains to show that this algorithm is well-defined and generates all elements in $I$. In Step 1 we produce all associated vectors in $I$ such that $|m_{j+1}|,\dots,|m_d| \leq p$.  Elements obtained in Step 1 and Step 2 form the subset of $I$ where $|m_{j+1}|,\dots,|m_d| \leq p+1$. In this way, we see that at any iteration $r$, Step 3 is well-defined, since $p+1,2p+1,\dots rp+1$ are not $0$ modulo $p$, and it generates all associated vectors in $I$ such that $|m_{j+1}|,\dots,|m_d| \leq rp+1$. Thus the union for each $r>0$ gives $I$.
\end{proof}
\begin{rem}
We recall that $\wn$ is not defined in the category of modules. In particular, we cannot take any "image" of $\wn$ in order to define a certain $N_{n,j}$ as the natural "lift" of some $k$-module. The definition of $k[\underline{z}]_n^{\leq}(I_j;1)$ try to solve this problem in the sense that for $n=1$ it coincides with the $P_j$-module $N_j$ defined in (\cite[Section 3]{Orlik21}).
\end{rem}



  



\appendix
\section{Witt vectors}\label{appendix_witt_vectors}
In the following appendix we recall definition and properties of Witt vectors. Every result mentioned here can be found in \cite[Ch. 0,1]{Illusie79},\cite[Appendix A]{Langer2004DERC} or more generally in \cite{Hes15}. All rings and algebras considered are commutative and unitary.\\
 Let $A$ be a ring. The ring of $p$-typical Witt vectors $\wn(A)$ of length $n \geq 1$  is the ring object defined by the following property:
 \begin{prop}[cf. {\cite[Proposition 1.2]{Hes15}}] For any  ring $A$, there exists a unique ring $\wn(A)$ functorial in the category of rings such that:
 \begin{itemize}
 \item[a)] The underlying set is $\wn(A)= A \times A \times \cdots \times A$, where the cartesian product is taken $n$-times;
 \item[b)] For any $i=0,\dots,n-1$, the \textit{$i$-th ghost map} $w_i$  defined by
 \begin{equation}
 w_i \colon \wn(A) \rightarrow A,\quad (a_1,\dots,a_n) \mapsto \sum_{j=0}^i p^ja_{j+1}^{p^{i-j}}
\end{equation}
is a ring homomorphism.
\end{itemize}  
\end{prop}  
  The  ring map $$ w=(w_0,\ldots,w_{n-1}): \w_n(A) \rightarrow A^n$$ where on the target we consider the product ring structure, is called the \textit{ghost map}. \\
   There are unique ring homomorphisms, functorial in $A$, called Frobenius maps, $$F: \w_{n}(A) \rightarrow \w_{n-1}(A),$$ such that $w \circ F=F^{w} \circ w$ where $F^{w}: A^{n} \rightarrow A^{n-1}$ is the shift map
$$
F^{w}\left( a_{1}, a_{2}, \ldots, a_n\right)=\left(a_{2}, a_{3}, \ldots, a_n\right).
$$
Notice that $w_{n-1}=F^{n-1} \colon \wn(A) \rightarrow A$. \\
 For any $n$, there are surjective ring homomorphisms (restrictions) $R: \w_{n+1}(A) \rightarrow \w_{n}(A)$ defined by
$$R\left( a_{1}, \ldots, a_{n+1}\right)=\left( a_{1}, \ldots, a_{n}\right), \quad (a_1,\dots,a_{n+1}) \in \w_{n+1}(A) .$$
Notice that $w_0=R^{n-1} \colon \wn(A) \rightarrow A$. Moreover, $RF=FR$. \\
 Also, there are natural additive maps, functorial in $A$, called Verschiebung maps $\ve \colon \w_{n-1}(A) \rightarrow \w_{n}(A)$ defined by
$$
\ve \left( a_{1}, a_{2}, \ldots, a_n\right)=\left(0, a_{1}, a_{2}, \ldots, a_n\right), \quad (a_1,\dots,a_n) \in \wn(A).
$$
We have that $RV=VR$. \\
Finally one has the multiplicative map, called the Teichmüller map, $[\text{ }]=[\text{ }]_n :A \rightarrow \w_{n}(A)$,
 where $[a]=(a, 0,0, \ldots,0) \in \wn(A)$ for any $a \in A$. The Teichm\"uller map is compatible with $R$.
 
These maps are related by the following identities:
\begin{enumerate}
  \item $F(\ve(a))=p a, \quad a \in \wn(A).$
  \item $a \ve\left(a^{\prime}\right)=\ve\left(F(a) a^{\prime}\right), \quad a,a' \in \wn(A).$ In particular, $\ve\wn(A) \subset \wn(A)$ is an ideal.
  \item $F([a])=\left[a^{p}\right], \quad a \in A.$
  \item If $pA=0$, then $\ve(F(a))=pa, \quad a \in \wn(A).$
\end{enumerate}
Furthermore, for any $n,r$, we have the following exact sequence:
\begin{equation}\label{wittexact}
    0 \rightarrow \mathrm{W}_{r}(A) \xrightarrow{\ve^n}\mathrm{W}_{n+r}(A)\xrightarrow{R^{r}} \w_n(A) \rightarrow 0.
\end{equation}
Therefore, the ring of \textit{Witt vectors} $\w(A):=\varprojlim_{n} (\cdots\rightarrow\w_{n+1}(A) \xrightarrow{R} \wn(A) \rightarrow \cdots)$ is separated and complete for the $\ve$-filtration defined by the ideals $\ve^n\w(A)$ for any $n \geq 0$. \\
The definition of the maps $\ve, F, [\text{ }]$ extends compatibly on the projective limit $\w(A)$. The Formulas 1-4 above hold also for $\w(A)$. Moreover, there is an exact sequence:
\begin{equation}
 0 \rightarrow \mathrm{W}_{}(A) \xrightarrow{\ve^n}\mathrm{W}_{}(A)\xrightarrow{R} \w_n(A) \rightarrow 0,
\end{equation}
where here $R$ is the canonical restriction defined by the projective limit. \\
Assume that $A$ is a $k$-algebra, where $k$ is a perfect field of characteristic $p$.\\ Then, $\wn(A)$ (resp. $\w(A)$, let us say for $n=\infty$ by definition) is canonically a $\wn(k)$-algebra (resp. a $\w(k)$-algebra ). Let $\phi_k \colon k \rightarrow k$, $\phi_A \colon A \rightarrow A$  be the Frobenius morphism $x \mapsto x^p$, and set $\Phi=\wn(\phi_k)$, $\Phi_A=\wn(\phi_A)$. Then, if $n < \infty$,  the relation  $F=R \circ \Phi_A$ holds. If $n= \infty$, then $F=\Phi_A$.  Furthermore, by the Formula 2), the Verschiebung map $\ve$ is a $\Phi^{-1}$-semilinear map over $\wn(k)$ (resp. over $\w(k)$ if $n=\infty$). \\
The following properties hold:
\begin{prop}
 Let $R$ be a ring such that $p \in R$ is nilpotent. Set $n \geq 1$ be a an integer. Then,
\begin{itemize}
\item[a)] Let $A$ be a ring. If $S \subset A$ is a multiplicative subset, the image of $S$ in $\wn(A)$ under $[\text{ }]$ is a multiplicative subset $[S]$ and there is an identification between the localizations $\wn(S^{-1}A)=[S]^{-1}\wn(A)$.
\item[b)] If $ f \colon A \rightarrow B$ is an étale morphism of $R$-algebras, then $\wn(f)$ is an étale morphism of $\wn(R)$-algebras.
\item[c)] If $A,B$ are $R$-algebras such that $A$ is étale over $R$, then the canonical map $\wn(B) \rightarrow \wn(A \otimes_R B)$ is étale and induces an isomorphism $\wn(A) \otimes_{\wn(R)} \wn(B) \xrightarrow{ \sim} \wn(A \otimes_R B)$.
\end{itemize}
\end{prop}
\begin{proof}
The assertion a) follows by arguments in \cite[Ch. 0, Sec. 1.5]{Illusie79}; b) is \cite[Proposition A.8]{Langer2004DERC} and c) is \cite[Corollary A.12]{Langer2004DERC}. 
\end{proof}
If $(X,\ox)$ is a $k$-scheme, then the presheaf $$ \wn\ox \colon U \mapsto \wn(\ox(U)),$$ where $U \subset X$ is open, is a sheaf of $\wn(k)$-algebras and the locally ringed space $\wn(X):=(|X|,\wn\ox)$ is a $\wn(k)$-scheme. The maps $F,\ve,R,[\text{ }],w$ sheafify on $\wn\ox$ (cf. \cite[Sec. 1.5]{Illusie79}).
\begin{prop}[cf. {\cite[Proposition 1.5.6]{Illusie79}}]\label{wn_preserve_immersions}
 If $X \rightarrow Y$ is an open (resp. closed) immersion of $k$-schemes, then $\wn(f)$ is an open (resp. closed) immersion of $\wn(k)$-schemes.
 \end{prop}

\section{Lifting properties}
In this appendix we collect some results about lifting properties of differential operators over smooth algebras. All rings and algebras are commutative with unit.
\begin{defn}[{\cite[Ch. 0, Definition 19.3.1]{EGA4}}]
 A ring map $f \colon A \rightarrow B$ is  formally smooth (resp.formally  étale), or equivalently  $B$ is a formally smooth $A$-algebra, (resp. formally étale $A$-algebra) if for any commutative diagram of ring maps 
\begin{equation}
    \begin{tikzcd}
B \arrow[r] \arrow[rd, dashed] & C/I         \\
A \arrow[r, "g"] \arrow[u, "f"]                       & C \arrow[u]
\end{tikzcd}
\end{equation}
where $I$ is a nilpotent ideal, there exists (resp. exists and it is unique) a ring map $\tilde{g} \colon B \rightarrow C$ making the diagram above commutative.\\
We say that $f$ is smooth (resp. étale) if $f$ is formally smooth (resp. formally étale) and of finite presentation.
\end{defn}

\begin{prop}[{\cite[Theorem 18.1.2]{EGA4}} Topological invariance of étale site]\label{top_inv_étale}
i) Let $S \xrightarrow{\pi} T$ a surjective ring map such that $\ker{(\pi)}$ is nilpotent. Then, for any étale $T$-algebra $T'$ there exists a unique étale $S$-algebra $S'$ and a ring map $S' \xrightarrow{\pi'} T'$ such that the following 
\begin{equation}
    \begin{tikzcd}
S' \arrow[r, "\pi'", dotted] & T'          \\
S \arrow[r, "\pi"] \arrow[u] & T \arrow[u]
\end{tikzcd}
\end{equation}
is a pushout diagram in the category of commutative rings.\\
ii) (Second formulation) For $S \xrightarrow{\pi} T$ as in i), the base change functor
\begin{equation}
    \{\text{ \'Etale } S-\text{algebras}\} \rightarrow \{\text{ \'Etale } T-\text{algebras}\}, \text{ } S' \mapsto S' \otimes_S T=:T'
\end{equation}
is an equivalence between the categories of étale $S$-algebras and étale $T$-algebras.
\end{prop}
Let $R$ be a base ring.
\begin{lemma}\label{lift_H-S_étale}
Let $A$ be an $R$-algebra and $D \colon A \rightarrow A[\epsilon]/(\epsilon^r)$ be a $R$-algebra homomorphism such that 
$  A \xrightarrow{D} A[\epsilon]/(\epsilon^r) \xrightarrow{\epsilon \mapsto 0} A$ is the identity of $A$. Then, for any $R$-algebra $B$ that is an étale $A$-algebra, there exists a unique $R$-algebra homomorphism 
\begin{equation}
    \tilde{D} \colon B \rightarrow B[\epsilon]/(\epsilon^r)
\end{equation} 
such that  $B \xrightarrow{\tilde{D}} B[\epsilon]/(\epsilon^r) \xrightarrow{\epsilon \mapsto 0} B$ is the identity of $B$
and compatible with $D$. In particular, $\tilde{D}=D \otimes_A B.$
\end{lemma}
\begin{proof}
We have that $B[\epsilon]/(\epsilon^r)$ is an $A$-algebra via the pushout square
\begin{equation}
    \begin{tikzcd}
B \arrow[r, "D \otimes_A B"]                                     & {B[\epsilon]/(\epsilon^r)}           \\
A \arrow[r, "D"] \arrow[u, "f"] \arrow[ru, "g"] & {A[\epsilon]/(\epsilon^r)} \arrow[u]
\end{tikzcd}
\end{equation}
 By definition of $D$ we have a commutative diagram
\begin{equation}
    \begin{tikzcd}
B \arrow[r, "id"]               & B                                                           \\
A \arrow[r, "g"] \arrow[u, "f"] & {B[\epsilon]/(\epsilon^r)} \arrow[u, "\epsilon \mapsto 0"']
\end{tikzcd}
\end{equation}
Thus, since $f$ is formally étale there exists a unique ring map $\tilde{D} \colon B \rightarrow B[\epsilon]/(\epsilon^r)$ such that the diagram above commutes.
\end{proof}
\begin{prop}\label{partial_H-S_to_ring}
    Let $A$ be a torsionfree $R$-algebra. For any integer $r \geq 0$, let $\partial^{[r]} \colon A \rightarrow A$ be additive maps, such that $\partial:=\partial^{[1]}$ is a $R$-linear derivation and for any integers $r,s \geq 0$ the following relations hold:
    \begin{equation}\label{rel_partial_H-S}
        \partial^{[0]}=id_A, \text{ } \partial^{[r]} \circ \partial^{[s]}=\binom{r+s}{r}\partial^{[r+s]}.
    \end{equation}
    Then, any $\partial^{[r]}$ only depends on $\partial$ and the map 
    \begin{equation}
        D^{(r)} \colon A \rightarrow A[\epsilon]/(\epsilon^{r+1}), \text{ } x \mapsto \sum_{i=0}^r \epsilon^i\partial^{[i]}(x) 
    \end{equation}
    is a $R$-algebra homomorphism such that $A \xrightarrow{D^{(r)}} A[\epsilon]/(\epsilon^{r+1}) \xrightarrow{\epsilon \mapsto 0} A$ agrees with $id_A$.
\end{prop}
\begin{proof}
    Any $\partial^{[r]}$ only depends on $\partial$: Indeed, by the relation $\partial^{[r-1]} \circ \partial=r\partial^{[r]}$, this follows by induction on $r$, being trivially satisfied for $r=1$ .\\
    For any $x,y \in A$, we claim that for every $r \geq 1$,
    \begin{equation}\label{HS-formula}
        \partial^{[r]}(xy)=\sum_{i=0}^r \partial^{[i]}(x)\partial^{[r-i]}(y).
    \end{equation}
    For $r=1$, the latter is just the Leibniz rule for $\partial$. We proceed by induction on $r$. We have that 
    \begin{align*}r\partial^{[r]}(xy)=\partial^{[r-1]}(\partial(xy)) \\ &= \partial^{[r-1]}(x\partial(y)+y\partial(x))\\
    & =\partial^{[r-1]}(x\partial(y))+\partial^{[r-1]}(y\partial(x))\\
    & =\sum_{i=0}^{r-1} \partial^{[i]}(x)\partial^{[r-i-1]}(\partial(y))+\partial^{[i]}(y)\partial^{[r-i-1]}(\partial(x))\\
    &= \sum_{i=0}^{r-1} (r-i)\partial^{[i]}(x)\partial^{[r-i]}(y)+(r-i)\partial^{[i]}(y)\partial^{[r-i]}(x)\\
    &= r\sum_{i=0}^r \partial^{[i]}(x)\partial^{[r-i]}(y).
    \end{align*}
    By assumptions on $A$, the latter equality implies the  \eqref{HS-formula}.
    Then, we have that $D^{(r)}(1)=1$ and the following relation for $D^{(r)}$ follows:
    \begin{align*}
        D^{(r)}(xy)  &=\sum_{i=0}^r\epsilon^i\partial^{[i]}(xy)
                   \\&= \sum_{i=0}^r\epsilon^i\sum_{s=0}^i\partial^{[s]}(x)\partial^{[i-s]}(y)=\sum_{i=0}^r\sum_{s=0}^i\epsilon^s\partial^{[s]}(x)\epsilon^{i-s}\partial^{[i-s]}(y)
                   \\&= D^{(r)}(x)D^{(r)}(y). \qedhere
    \end{align*}
\end{proof}
\begin{rem}
By induction on $s$ applied to the product $x_1\cdots x_s$, with $x_1,  \dots,x_s \in A$,  the  \eqref{HS-formula} generalizes to
\begin{equation}\label{general_H-S_formula}
    \partial^{[r]}(x_1\cdots x_s)=\sum_{i_1+\cdots i_s=r} \partial^{[i_1]}(x_1)\cdots \partial^{[i_s]}(x_s).
\end{equation}
\end{rem}
\begin{corol}\label{lift_partial_H-S_étale} 
Let $A$ be a  $R$-algebra and $\partial_A^{[r]} \colon A \rightarrow A$ be a collection of $R$-linear maps satisfying \eqref{general_H-S_formula} and \eqref{rel_partial_H-S}. Then, for any $R$-algebra $B$ that is an étale $A$-algebra, there is a unique collection of additive maps $\partial_B^{[r]} \colon B \rightarrow B$ such that $\partial_B^{[1]}$ is a $R$-linear derivation, the relations \eqref{general_H-S_formula} hold and they are compatible with $\partial_A^{[r]}$.
\end{corol}
\begin{proof}
    By  \Cref{lift_H-S_étale} the corresponding $R$-algebra homomorphism $D_A^{(r)}$ defined as in  \Cref{partial_H-S_to_ring} lifts to a unique $R$-algebra homomorphism $D_B^{(r)}$. Then $\partial^{[r]}_B$ is determined by writing $D^{(r)}_B(x)=x+\epsilon\partial_B(x)+\cdots+\epsilon^r \partial_B^{[r]}(x)$. The relations \eqref{general_H-S_formula} then hold. 
\end{proof}
\begin{corol}\label{lift_derivation_over_smooth}
 Let $T$ be an $R$-algebra.
\begin{itemize}
\item[i)] Let $\tilde{T}, \tilde{R}$ be  torsionfree rings equipped with  surjective ring maps $\tilde{T} \rightarrow T, \tilde{R} \rightarrow R $. Assume, the diagram
\begin{equation}
\begin{tikzcd}
\tilde{T} \arrow[r]  & T \\
\tilde{R}  \arrow[u] \arrow[r] &  R \arrow[u] 
\end{tikzcd}
\end{equation}
is a pushout square. Then, any collection of maps $\partial_{\tilde{T}}^{[r]}\colon \tilde{T} \rightarrow \tilde{T}$ as in  \Cref{partial_H-S_to_ring} induces a  collection of maps
$\partial_{{T}}^{[r]}\colon {T} \rightarrow {T}$, compatible with  $\partial_{\tilde{T}}^{[r]}$. In particular, $\partial_{{T}}^{[r]}$ satifies the relation \eqref{general_H-S_formula} for any $r \geq 0$;
\item[ii)]Let $S \xrightarrow{f} T$ be a surjective $R$-algebra homomorphism with nilpotent kernel. Assume that $S$ is a smooth $R$-algebra. Any collection of maps
$\partial_{{T}}^{[r]}\colon {T} \rightarrow {T}$ satisfiyng the relation \eqref{general_H-S_formula} lifts to some collection of maps $\partial_S^{[r]} \colon S \rightarrow S$ compatible with $\partial_T^{[r]}$.
\end{itemize} 
\end{corol}
\begin{proof}
\textit{i)}: By  \Cref{partial_H-S_to_ring}, let $D_{\tilde{T}}^{(r)}$ denote the corresponding $\tilde{R}$-algebra homomorphism to $\partial_{\tilde{T}}^{[r]}$. The base change of $D_{\tilde{T}}^{(r)}$ along $\tilde{T} \rightarrow T$, induces a $R$-algebra homomorphism $D_{{T}}^{(r)} \colon T \rightarrow T[\epsilon]/(\epsilon^{r+1})$. Writing $D_{T}^{(r)}(x)=x+\epsilon\partial_T(x)+\cdots+\epsilon^r \partial_T^{[r]}(x)$ for any $x \in T$, we get a collection of maps $\partial_T^{[r]}$ with the desired properties. \\
\textit{ii)}: The assumptions ensure that the map $D_T^{(r)}:=\sum_{i=0}^r\epsilon^i\partial_{T}^{[r]} \colon T \rightarrow T[\epsilon]/(\epsilon^{r+1})$ is an $R$-algebra morphism. Since $S$ is smooth, it is in particular formally smooth. Moreover, the natural map $$S[\epsilon]/(\epsilon^{r+1}) \xrightarrow{f \otimes 1} T[\epsilon]/(\epsilon^{r+1})$$ induced by the tensor product of $f$ with $R[\epsilon]/(\epsilon^{r+1})$, is surjective with nilpotent kernel, thus by formally smoothness the following solid square
\begin{equation}
    \begin{tikzcd}
S \arrow[r, "D_T^{(r)} \circ f"] \arrow[rd, "\exists ", dashed] & T[\epsilon]/(\epsilon^{r+1})         \\
R \arrow[r] \arrow[u]                       & S[\epsilon]/(\epsilon^{r+1}) \arrow[u, "f \otimes 1" ']
\end{tikzcd}
\end{equation}
admits a $R$-algebra homomorphism $D_S^{(r)} \colon S \rightarrow S[\epsilon]/(\epsilon^{r+1})$. Now as in \textit{i)}, the $\partial_S^{[r]} \colon S \rightarrow S$ are defined by expanding the writing of $D_S^{(r)}(x)$ for any $x \in S$ and they are by construction compatible with $\partial_T^{[r]}$. 
\end{proof}
 
 \begin{corol}\label{cor_smooth_lift_der_over_wn}
 Let $T'$ be a smooth $k$-algebra. Then the following holds:
 \begin{itemize}
     \item[i)] There exist a collection of $k$-linear maps $\partial_{T'}^{[r]} \colon T' \rightarrow T'$
     satisfying \eqref{general_H-S_formula}
     \item[ii)] There exists a projective system of smooth $\wn(k)$-algebras $S_n'$ such that $S_1'=T'$ and for any $r \geq 0$, a projective system of $\wn(k)$-linear maps $\partial_{S'_n}^{[r]}\colon S_n'\rightarrow S_n'$ satisfying \eqref{general_H-S_formula}, lifting $\partial_{T'}^{[r]}$. 
 \end{itemize}
 \end{corol}
\begin{proof}
We verify the statement for $T'=T=k[\underline{x}]$ being a polynomial algebra over $k$. \\
Let $R=k$, $\tilde{R}=\w(k)$. The polynomial algebra  $\tilde{T}:=\w(k)[\underline{x}]$ over $\w(k)$ is a torsionfree lift of $T$. The assumptions of  \Cref{lift_derivation_over_smooth} i) are then satisfied. By \Cref{partial_H-S_to_ring}, any $\w(k)$-linear derivation $\partial_{\tilde{T}}$ of $\tilde{T}$ determines uniquely a collection of $\w(k)$-linear maps $\partial_{\tilde{T}}^{[r]}$, thus by  \Cref{lift_derivation_over_smooth} i), they induce $k$-linear maps $\partial_{T}^{[r]}$  satisfying the relation \eqref{general_H-S_formula}. Let $S_n=\wn(k)[\underline{x}]$ and consider $T$ as $\wn(k)$-algebra via the restriction $\wn(k) \rightarrow k$. Then, the map $S_n \rightarrow T$ satisfies the assumptions of  \Cref{lift_derivation_over_smooth} ii), thus  the collection $\partial_T^{[r]}$ lift to $\wn(k)$-linear maps $\partial_{S_n}^{[r]}$ satisfying the relation \eqref{general_H-S_formula}. Also, the $\w_{n+1}(k)$-algebra map $S_{n+1} \rightarrow S_n$ satisfies the condition of  \Cref{lift_derivation_over_smooth} ii), therefore for any $r \geq 0$ we can find a projective system of $\wn(k)$-linear maps $\{\partial_{S_n}^{[r]}\}_n$ with respect to the surjections $S_{n+1} \rightarrow S_n$, satisfying the \eqref{general_H-S_formula}. 

A generic smooth $T'$ is étale over a polynomial algebra $T$, and such $\partial_{T'}^{[r]}$ are induced uniquely by $\partial_{T}^{[r]}$ by \Cref{lift_partial_H-S_étale} (and do not depend on the choice of local coordinates).

Now, assume that $T'$ is  étale over some polynomial algebra $T=k[\underline{x}]$ over $k$. Then, there exists a projective system of smooth lifts $S'_n$ over $\wn(k)$: Indeed, the natural ring map $S_n=\wn(k)[\underline{x}] \rightarrow \w_{n-1}(k)[\underline{x}]=S_{n-1}$ satisfies the condition of Proposition \ref{top_inv_étale}, therefore applying the same proposition successively for any $n >1$, we get a projective system of étale $\wn(k)[\underline{x}]$-algebras $S'_n$. In particular, $S'_n$ is smooth over $\wn(k)$. A general $T'$ smooth over $k$ is locally étale over a polynomial algebra $T$. Hence, there exist locally such étale $T$-algebras $S'_n$. Since $\spec(T')$ is smooth affine, the obstruction class to glue the $S'_n$'s is $0$ (cf. \cite[Theorem 6.3]{SGA1}). Therefore, by glueing we get a smooth scheme over $\wn(k)$ lifting $T$.\footnote{Note that this lift is affine, since any scheme admitting a surjective integral map from an affine scheme is affine (c.f. \cite[Tag 05YU, Lemma 32.11.2]{stacks-project}).}

 Thus, by applying  \Cref{lift_derivation_over_smooth} ii) and induction on $n$ to the ring maps $S'_{n+1} \rightarrow S'_n$,  we get a projective system of $\wn(k)$-linear lifts $\{\partial_{S'_n}^{[r]}\}_n$ of $\partial_{T'}^{[r]}$ for any $r \geq 0$.
\end{proof}

 


 \nocite{*}


\printbibliography
\printindex
\end{document}